\documentclass[11pt]{article}
\usepackage[margin=1in]{geometry}
\usepackage{lmodern}
\usepackage{amsmath}
\usepackage{color}
\usepackage{mathtools}
\mathtoolsset{showonlyrefs}
\usepackage{amsfonts,amstext, amssymb,amsthm, amsmath, rotating, latexsym,wasysym}
\usepackage{graphicx,bm,bbm,bigints}% ,
\usepackage{adjustbox}

\usepackage[colorinlistoftodos]{todonotes}

\usepackage{footnote}
\makesavenoteenv{tabular}
\makesavenoteenv{table}
\usepackage{nameref,xr}
\usepackage[caption=false]{subfig}
\usepackage[subfigure]{tocloft}

\usepackage{multirow}
\usepackage{booktabs}
\usepackage{epstopdf}
\usepackage{bigstrut}
\usepackage{enumitem}
\DeclareRobustCommand{\rchi}{{\mathpalette\irchi\relax}}
\newcommand{\irchi}[2]{\raisebox{\depth}{$#1\chi$}}
\theoremstyle{plain}
\newtheorem{thm}{Theorem}[section]
\newtheorem{lemma}[thm]{Lemma}

\newtheorem{prop}[thm]{Proposition}
\newtheorem{remark}[thm]{Remark}
\newtheorem{example}[thm]{Example}

\makeatletter
\def\namedlabel#1#2{\begingroup
    #2%
    \def\@currentlabel{#2}%
    \phantomsection\label{#1}\endgroup
}
\makeatother

\makeatletter

\newcommand{\Rome}[1]{\expandafter\@slowromancap\romannumeral #1@}
\makeatother

\newcommand{\be}{\begin{equation}}
\newcommand{\ee}{\end{equation}}

\DeclareMathOperator*{\argmin}{arg\,min}

\newcommand{\R}{\mathbb{R}}

\newcommand{\D}{\mathcal{D}}
\newcommand{\p}{\mathbb{P}}

\newcommand{\w}{\mathcal{W}}
\newcommand{\g}{\mathbb{G}}
\newcommand{\f}{\mathcal{F}}
\newcommand{\E}{\mathbb{E}}
\newcommand{\M}{\mathcal{M}}

\newcommand{\h}{\mathcal{H}}
\newcommand{\C}{\mathfrak{C}}

\newcommand{\A}{\mathfrak{A}}

\newcommand{\Var}{\mathrm{Var}}
\definecolor{lgray}{gray}{0.70}
\usepackage{tikz}

\providecommand{\norm}[1]{\left\lVert#1\right\rVert}

\usepackage[linesnumbered,ruled,vlined]{algorithm2e}

\usepackage[numbers]{natbib}
\usepackage{blindtext}
\usepackage{nicefrac}

\usepackage{libertine}

\usepackage[font=small,labelfont=bf]{caption}

\usepackage[colorlinks,citecolor=blue,urlcolor=blue,linkcolor=red,naturalnames=true]{hyperref}
\usepackage{authblk}

\newcommand{\blind}{1}

\begin{document}
%% ADD before submitting
\numberwithin{equation}{section}
\def\spacingset#1{\renewcommand{\baselinestretch}%
{#1}\small\normalsize} \spacingset{1}
\spacingset{1.5} % DON'T change the spacing!

\if1\blind
{
\title{\bf Semiparametric Efficiency in Convexity Constrained Single Index Model}
\author[1]{Arun K. Kuchibhotla\thanks{Email: {\tt             arunku@cmu.edu}.}}
 \author[2]{ Rohit K. Patra\thanks{Email: {\tt rohitpatra@ufl.edu}.}} 
 \author[3]{Bodhisattva Sen\thanks{Email:{\tt bodhi@stat.columbia.edu} Supported by NSF Grants DMS-17-12822 and AST-16-14743.}}
\affil[1]{Carnegie Mellon University }
 \affil[2]{University of Florida}
 \affil[3]{Columbia University}
    \date{}
    \maketitle
  } \fi
  
  \if0\blind
  {    \title{Semiparametric Efficiency in Convexity Constrained Single Index Model}
  \date{}
        \maketitle
  } \fi

%  \author{Arun K. Kuchibhotla, Rohit K. Patra, and Bodhisattva Sen}
% \maketitle

\begin{abstract}
We consider estimation and inference in a single index regression model with an unknown convex link function.  We introduce a convex and Lipschitz constrained least squares estimator (CLSE) for both  the parametric and the nonparametric components given independent and identically distributed observations. We prove the consistency and find the rates of convergence of the CLSE when the errors are assumed to have only $q \ge 2$ moments and are allowed to depend on the covariates. When $q\ge 5$, we establish $n^{-1/2}$-rate of convergence and asymptotic normality of  the estimator of the parametric component.  Moreover,  the CLSE is proved to be semiparametrically efficient if the  errors happen to be homoscedastic. \if1\blind{We develop and implement a numerically stable and computationally fast algorithm to compute our proposed estimator in the R package~\texttt{simest}}\fi\if0\blind{We develop and implement a numerically stable and computationally fast algorithm to compute our proposed estimator in an R package}\fi. We illustrate our methodology through extensive simulations and data analysis. Finally,  our proof of efficiency is geometric and provides a general framework that can be used to prove efficiency of estimators in a wide variety of semiparametric models even when they do not satisfy the efficient score equation directly. 
 % Our proposed algorithm works even when $n$ is modest and $d$ is large (e.g., $n = 500$, and $d = 100$). 
 \end{abstract}
\noindent%
\section{Introduction}\label{sec:intro}
Suppose  we have $n$ i.i.d.~observations $\{(X_i,Y_i)\in   \rchi\times \R, 1\le i\le n\}$
from the following single index regression model:
\begin{equation}\label{eq:simsl}
Y = m_0(\theta_0^{\top}X) + \epsilon,
\end{equation}
where $X \in   \rchi \subset \R^d$ ($d \ge 1$) is the predictor,  $Y \in \R$ is the response variable, and $\epsilon$ satisfies $\mathbb{E}(\epsilon|X) = 0$ and $\mathbb{E}(\epsilon^2|X) < \infty$ almost everywhere (a.e.) $P_X$, the distribution of $X$. We assume that the real-valued link function  $m_0$ and $\theta_0 \in \R^d$ are the unknown parameters of interest. 

Single index models are ubiquitous in regression because they provide convenient dimension reduction and interpretability. The single index model circumvents the curse of dimensionality encountered in estimating the fully nonparametric regression function $\E(Y|X = \cdot)$ by assuming that the link function depends on $X$ only through a one dimensional projection, i.e., $\theta_0^\top X$; see e.g.,~\cite{Powelletal89}. Moreover, the coefficient vector $\theta_0$ provides interpretability~\citep{liracine07} and the one-dimensional nonparametric link function $m_0$ offers some flexibility in modeling. The above model has received a lot of attention in statistics in the last few decades; see e.g.,~\cite{Powelletal89,LiDuan89,ICHI93,HardleEtAl93,Hristacheetal01,DelecroixEtal06,cuietal11,Patra16} and the references therein. The above papers propose estimators for the single index model under the assumption that $m_0$ is smooth (i.e., two or three times differentiable). 

 However, quite often in the context of a real application, qualitative assumptions on $m_0$ may be available. 
% {\color{red}The implicit assumption here that the conditional expectation of $Y$ given $X$ is only a function of $\theta_0^{\top}X$, although reduces the flexibility of full non-parametric model, compensates by a lot in quality of estimation and interpretation. Here the quality of estimation is in terms of the rate at which the regression function can be estimated. It is well-known from the results in \cite{Gyorfi02} that without the single index restriction, the best rate of convergence is $n^{-2/(4+d)}$ but with this restriction we are able to obtain $n^{-2/5}$. The single index model assumption also allows for recognizing relative importance of predictor variables in predicting $Y$.}
% In this paper, we assume further that $m_0$ is a convex Lipschitz function. This assumption is motivated by the fact that in a wide range of applications in various fields the regression function is known to be convex or concave.
 For example, in microeconomics, production and utility functions are often assumed to be concave and nondecreasing; concavity indicates decreasing marginal returns/utility~\citep{Varian84,Matzkin91,liracine07}. In finance, the relationship between call option prices and strike price are often known to be convex and decreasing~\citep{AD03}; in stochastic control, value functions are often assumed to be convex~\citep{K11}.  The following two real-data examples further illustrate that convexity/concavity constraints arise naturally in many applications.

\begin{figure}[ht!]
\centering
\includegraphics[width=.9\textwidth]{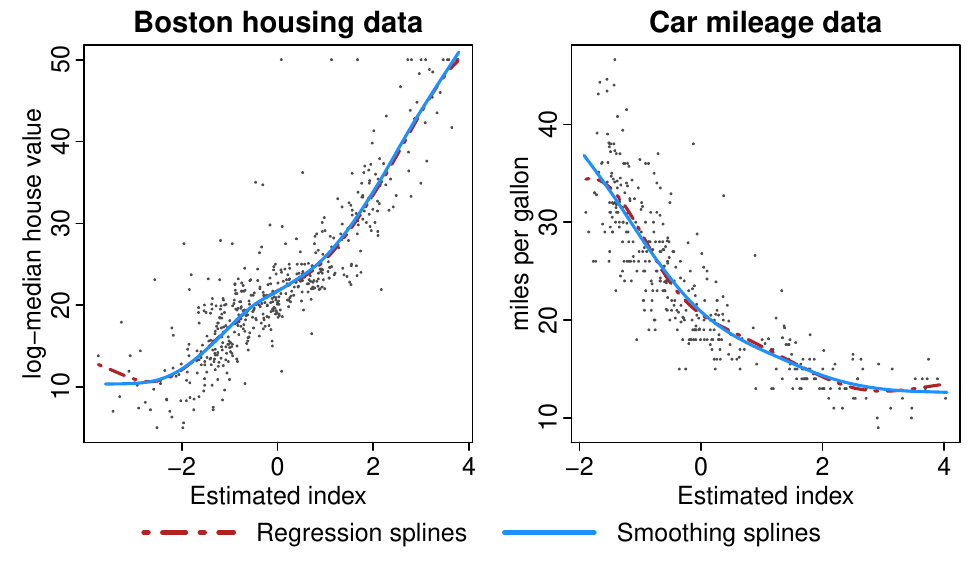}
%% scale=.8
  \caption[]{Scatter plots of $\{(X_i^\top\hat{\theta}, Y_i)\}_{i=1}^n,$ where $\hat{\theta}$ is the estimator of $\theta_0$ proposed in~\cite{MR2529970}. The plot is overlaid with the smoothing and regression splines based  function estimators of $m_0$ proposed in~\cite{Patra16} and~\cite{MR2529970}, respectively. Left panel: Boston housing data (see Section~\ref{sec:boston}); right panel: the car mileage data (see Section~\ref{sec:car}).}
  \label{fig:real_data_plot_prelim}
\end{figure}

\begin{example}[Boston housing data]\label{ex:boston}
\citet{harrison1978hedonic} studied the effect of different covariates on real estate price in the greater Boston area.  The response variable $Y$ was the log-median value of homes in each of the $506$ census tracts in the Boston standard metropolitan area. A single index model is  appropriate for this dataset; see e.g.,~\cite{gu2015oracally,MR2529970,MR2589322,MR2787613}. The above papers considered the following covariates in their analysis: average number of rooms per dwelling, full-value property-tax rate per $10000$ USD, pupil-teacher ratio by town school district, and proportion of population that is of ``lower (economic) status'' in percentage points. In the left panel of Figure~\ref{fig:real_data_plot_prelim}, we provide the scatter plot of  $\{(Y_i, \hat{\theta}^\top X_i)\}_{i=1}^{506}$, where $\hat{\theta}$ is the estimate of $\theta_0$ obtained in~\cite{MR2529970}. We also plot estimates of $m_0$ obtained from~\cite{Patra16} and~\cite{MR2529970}. The plot suggests a convex and nondecreasing relationship between the log-median home prices and the index, but the fitted link functions satisfy these shape constraints only approximately.

\end{example}

\begin{example}[Car mileage data]\label{ex:car}
  \citet{cars_1983} consider a dataset containing mileages of different cars. The data contains mileages of $392$ cars as well as the following covariates:  displacement,   weight,  acceleration,   and horsepower.  \citet{MR2957294} and \cite{Patra16} have fit a partial linear model and a  single index model, respectively. In the right panel of Figure~\ref{fig:real_data_plot_prelim}, we plot the estimators  proposed in~\cite{Patra16} and~\cite{MR2529970}. Both of these works consider estimation in the single index model under only smoothness assumptions. The ``law of diminishing returns'' suggests $m_0$ should be convex and nonincreasing. However, as observed in Figure~\ref{fig:real_data_plot_prelim}, the estimators based only on smoothness assumptions satisfy this shape constraint only approximately.

\end{example}
% subsection example_2_car_mileage_data (end)
In both of the examples, the smoothing based estimators do not incorporate the known shape of the nonparametric function. Thus the estimators are not guaranteed to be convex (or monotone) in finite samples. Moreover, the choice of the tuning parameter in smoothness based estimators is tricky as different values for the tuning parameter lead to very different shapes. This unpredictable behavior makes the smoothness  based  estimators of $m_0$ less interpretable, and motivates the study of a convexity constrained single index model. {We discuss these two datasets and our analysis in more detail in Sections~\ref{sec:boston} and~\ref{sec:car}.}

In this paper, we propose constrained least squares estimators for $m_0$ and $\theta_0$ that is guaranteed to satisfy the inherent convexity constraint in the link function everywhere. The proposed methodology is appealing for two main reasons: (1) the estimator is interpretable and takes advantage of naturally occurring qualitative constraints; and (2) unlike smoothness based estimators, the proposed estimator is highly robust to the choice of the tuning parameter without sacrificing efficiency. 

In the following, we conduct a systematic study of the computation, consistency, and rates of convergence of the estimators, under mild assumptions on the covariate and error distributions. We further prove that the estimator for the finite-dimensional parameter $\theta_0$ is asymptotically normal. Moreover, this estimator is shown to be semiparametrically efficient if the errors happen to be homoscedastic, i.e., when $\mathbb{E}(\epsilon^2|X) \equiv \sigma^2$ a.e.~for some constant $\sigma^2$. It should be noted that in the examples above the link function is also known to be monotone. To keep things simple, we focus on  only convexity constrained single index model. However, \textit{all} our results continue to hold under the additional monotonicity assumption, i.e., our conclusions hold for convex/concave and nondecreasing/nonincreasing $m_0$. {More generally, our results continue to hold under \textit{any} additional shape constraints; see Remarks~\ref{rem:Monotone},~\ref{rem:Efficency_additional}, and~\ref{rem:Add_mono} and Section~\ref{sec:real_data_analysis} in the paper for more details.}
 % 

% estimation and inference in shape-restricted single index models has not received much attention.

One of the main contributions of this paper is our novel geometric proof of the semiparametric efficiency of the constrained least squares estimator. Note that proving semiparametric efficiency of constrained (and/or penalized) least squares estimators often requires a delicate use of the structure of the estimator of the nonparametric component (say $\hat m$) to construct {\it least favorable paths}; see e.g.~\cite{VANC},~\cite[Chapter 9.3]{VdV02}, and~\cite{MR1394975} (also see Example~\ref{ex:coxmodel}). In contrast, our approach is based on the following simple observation. For a traditional smoothness based estimator $\hat{m}$, the path $t \mapsto \hat{m} + t a$ will belong to the (function) parameter space for {\it any} smooth ``perturbation'' $a$ (for small enough $t \in (-1,1)$). However this is no longer true when the underlying parameter space is constrained. But,  observe that the projection of $\hat{m} + t a$ onto the constrained function space certainly yields a ``valid" path. Our proof technique is based on differentiability properties of the path $t\mapsto \Pi (\hat{m} + ta)$, where $\Pi$ denotes the $L_2$-projection onto the (constrained) function space. This general  principle is applicable to other shape constrained semiparametric models, because differentiability of the projection operator is well-studied  in the context of constrained optimization algorithms; see~Section~\ref{sub:_semiparametric_eff_shape} below for a more detailed discussion. {Also see Example~\ref{ex:coxmodel}, where we discuss the applicability of our technique in (re)proving the semiparametric efficiency of the nonparametric maximum likelihood estimator in the Cox proportional hazard model under current status censoring~\cite{MR1394975}.}
 %Due to the assumed shape constraint on $m_0$, the parametric submodels for the link function are nonlinear and the nuisance tangent space is intractable.}
To be more specific, we study the following Lipschitz constrained convex least squares estimator (CLSE): 
\begin{equation}\label{eq:CLSE}
(\check{m}_{L},\check{\theta}_{L}) \coloneqq \argmin_{(m,\theta) \in \M_{L} \times\Theta} Q_n(m, \theta),
\end{equation}
where
\begin{equation}\label{eq:CLSE-2}
Q_n(m, \theta)\coloneqq \frac{1}{n} \sum_{i=1}^n \{Y_i - m(\theta^\top X_i)\}^2
\end{equation}
and  $\M_L$ denotes the class of all $L$-Lipschitz real-valued convex functions on $\R$ and
\be\label{smooth:eq:theta}
\Theta \coloneqq \{\eta=(\eta_1, \ldots, \eta_d) \in\mathbb{R}^d: |\eta|=1\mbox{ and }\eta_{1}\ge 0 \} \subset S^{d-1}.
\ee 
Here $| \cdot |$ denotes the usual Euclidean norm, and  $S^{d-1}$ is the Euclidean unit sphere in $\mathbb{R}^d$. The norm-1 and the positivity constraints are necessary for identifiability of the model\footnote{\label{foo:ident} Without any sign or scale constraint on $\Theta$ no $(m_0, \theta_0)$ will be identifiable. To see this, fix any $(m_0, \theta_0)$ and define $m_1(t) \coloneqq m_0(-2t)$ and $\theta_1 = -\theta_0/2$, then  $m_0(\theta_0^\top \cdot)\equiv m_1(\theta_1^\top \cdot)$; see~\cite{carroletal97}, \cite{cuietal11}, and~\cite{MR2369025} for identifiability of the model~\eqref{eq:simsl}. Also see Section~\ref{sec:ident} for further discussion.}. 

The Lipschitz constraint in~\eqref{eq:CLSE} is not restrictive as all convex functions are Lipschitz in the interior of their domains. Furthermore in shape-constrained single index models,  the Lipschitz constraint is known to  lead to computational advantages~\citep{kalai2009isotron,kakade2011efficient,lim2014convergence,ganti2015learning,mazumder2019computational}. 
% It should be noted that  the additional Lipschitz assumption is not generally required when studying the nonparametric convex regression model (see e.g.,~\cite{MR2850215}) but this restriction eases the theoretical analysis in this semiparametric model.
{ Additionally on the theoretical side, the Lipschitzness assumption allows us to control the behavior of  the estimator near the boundary of its domain. This control is crucial for establishing semiparametric efficiency. To the best of our knowledge, this is the first work proving semiparametric efficiency for an estimator in a \textit{bundled parameter} problem (where the parametric and nonparametric components are intertwined; see \cite{huang1997interval}) where the nonparametric estimate is shape constrained and non-smooth.} % The asymptotic normality of $\check{\theta}_{L,n}$ can be readily used to  construct confidence intervals for ${\theta}_0$.} 
%As any convex function is Lipschitz in the interior of its domain, Our estimator $(\check{m}_{L},\check{\theta}_{{L}})$ can be thought of as a shape-constrained nonparametric LSE for the parameters $m_0$ and $\theta_0$ in model~\eqref{eq:simsl}. 
Note that the convexity constraint in~\eqref{eq:CLSE} leads to a convex piecewise affine estimator $\check{m}_{L}$ for the link function $m_0$; see Section~\ref{sec:CLSE} for a detailed discussion. 

% ~\citet{groeneboom2016current} propose a $\sqrt{n}$-consistent and asymptotically normal but \textit{inefficient} estimator of the index vector in the current status model based on the (non-smooth) maximum likelihood estimate (MLE) of the nonparametric component under just monotonicity constraint.  They also propose two other estimators of the index vector based on kernel smoothed versions of the MLE for the nonparametric component.  Although these estimators do not achieve the efficiency bound their asymptotic variances can be made arbitrarily close to the efficient variance.

% \subsection{Summary of our contributions} % (fold)
% \label{sub:contributions_of_the_paper}

% In this paper, we study the CLSE~\eqref{eq:CLSE}. 
Our theoretical and methodological study can be split in two broad categories. In Section~\ref{sec:CLSE}, we find the rate of convergence of the CLSE as defined in~\eqref{eq:CLSE}, whereas in Section~\ref{sec:SemiInf} we establish the asymptotic normality and semiparametric efficiency of $\check{\theta}_{L}$. Suppose that $m_0$ is $L_0$-Lipschitz, i.e., $m_0 \in \M_{L_0}$. If the tuning parameter $L$ is chosen such that $L \ge L_0$, then under mild distributional assumptions on $X$ and $\epsilon$, we show that $\check{m}_{L}$ and $\check{m}_{L}(\check{\theta}_{L}^\top \cdot)$ are minimax rate optimal for estimating $m_0$  and $m_0(\theta_0^\top \cdot)$, respectively; see Theorems~\ref{thm:rate_m_theta_CLSE} and~\ref{thm:ratestCLSE}. We also allow for the tuning parameter $L$ to depend on the data and show that the rate of convergence of $\check{m}_{L}(\check{\theta}_{L} \cdot)$ is uniform in $ L\in [L_0, nL_0]$, up to a $\sqrt{\log\log n}$ multiplicative factor; see Theorem~\ref{thm:UniformLRate}. This result justifies the usage of a data-dependent choice of $L$, such as cross-validation. Additionally, in Theorem~\ref{thm:rate_derivCLSE}, we find the rate of convergence of $\check{m}'_{L}.$ In Section~\ref{sec:SemiInf}, we establish that if $L\ge L_0$, then  $\check{\theta}_{L}$  is $\sqrt{n}$-consistent and $n^{1/2}(\check{\theta}_L - \theta_0)$ is asymptotically normal with mean $0$ and finite variance; see Theorem~\ref{thm:Main_rate_CLSE}. The asymptotic normality of $\check{\theta}_{L}$ can be readily used to  construct confidence intervals for ${\theta}_0$. Further, we show that  if the errors happen to be homoscedastic, then $\check{\theta}_{L}$ is semiparametrically efficient. 

Our contributions on the computational side are two fold. In Section~\ref{sec:compute} of the supplementary file, we propose an alternating descent algorithm for estimation in the single index model~\eqref{eq:simsl}. Our descent algorithm works as follows: when $\theta$ is fixed, the $m$ update is obtained by solving a quadratic program with linear constraints, and when $m$ is fixed, we update $\theta$ by taking a small step on the Stiefel manifold $\Theta$ with a guarantee of descent. \if1\blind  {We implement the proposed algorithm in the R package \texttt{simest}}\fi
  \if0\blind  {We implement the proposed algorithm in the R package {\texttt{***}}}\fi. Through extensive simulations (see Section~\ref{sec:Simul_Cvx} and Section~\ref{app:add_simul} of the supplementary file) we show that the finite sample performance of our estimators is robust to the choice of the tuning parameter $L$.  Thus we think the practitioner can choose $L$ to be very large without sacrificing any finite sample performance. Even though the minimization problem is non-convex, we illustrate that the proposed algorithm (when used with multiple random starting points) performs well in a variety of simulation scenarios when compared to existing methods.

\subsection{Semiparametric efficiency and shape constraints} % (fold)
\label{sub:_semiparametric_eff_shape}
Although estimation in single index models under smoothness assumptions is well-studied (see e.g.,~\cite{Powelletal89,LiDuan89,ICHI93,HardleEtAl93,Hristacheetal01,DelecroixEtal06,MR2529970,cuietal11} and the references therein), estimation and efficiency in shape-restricted single index models have not received much attention. The earliest reference on this topic we could find was the work of \citet{VANC}, where the authors considered a penalized likelihood approach in the current status regression model (which is similar to the single index model) with a monotone link function.~\citet{CHSA} consider maximum likelihood estimation in a generalized additive index model (a more general model than \eqref{eq:simsl}) and only prove consistency of the proposed estimators.  In~\citet{2016arXiv161006026B}, the authors study model~\eqref{eq:simsl} under  monotonicity constraint and prove $n^{1/3}$-consistency of the LSE of $\theta_0$; however they do not obtain the limiting distribution of the estimator of $\theta_0.$ \citet{2017arXiv171205593B} propose a tuning parameter-free $\sqrt{n}$-consistent (but not semiparametrically efficient) estimator for the index parameter in the monotone  single index model. 

% subsection _se_ (end)

In this paper, we show that $\check{\theta}_L$ is semiparametrically efficient under homoscedastic errors. Our proof of the semiparametric efficiency is novel and can be applied to other semiparametric models when the estimator does not readily satisfy the efficient score equation. In fact, we provide a new and general technique for establishing semiparametric efficiency of an estimator when the nuisance tangent set is not the space of all square integrable functions. The basic idea is as follows. Suppose $\ell_{\theta_0, m_0}(y, x)$ represents the semiparametrically efficient influence function, meaning that the ``best'' estimator $\tilde{\theta}$ of $\theta_0$ satisfies the following asymptotic linear expansion:
\begin{equation}\label{eq:eff_intro}
\eta^\top (\tilde{\theta} - \theta_0) = \frac{1}{n}\sum_{i=1}^n \eta^\top \ell_{\theta_0, m_0}(Y_i, X_i) + o_p(n^{-1/2}),
\end{equation}
for every $\eta\in \R^d$. A crucial step in establishing that $\check{\theta}_L$ satisfies~\eqref{eq:eff_intro} is to show for any $\eta\in\mathbb{R}^d$,
\[n^{-1}\sum_{i=1}^n  \eta^{\top}\ell_{\check{\theta}_L, \check{m}_L}(Y_i, X_i) = o_p(n^{-1/2}),\]
 i.e., $\check{\theta}_L$ is an \emph{approximate zero} of the efficient score equation~\cite[Theorem 6.20]{VdV02}. Because $(\check{m}_L,\check{\theta}_L)$ minimizes $(m,\theta) \mapsto Q_n(m,\theta)$ over $\mathcal{M}_L\times\Theta$, the traditional way to prove the approximate zero property is to use the fact that $\partial Q_n(\check{m}_L + t a, \check{\theta}_L + t\eta)/\partial t|_{t = 0} =0$ for all perturbation ``directions'' $(a, \eta)$ and find an $a$ such that the derivative of $t\mapsto Q_n(\check{m}_L + ta, \check{\theta}_L + t\eta)$ at $t=0$  is $n^{-1}\sum_{i=1}^n \eta^\top \ell_{\check{\theta}_L, \check{m}_L}(Y_i, X_i)$; see e.g.,~\cite{NeweyStroker93}. {In fact, using this method one can often show that the estimator satisfies the efficient score equation \textit{exactly}.} If $\check{m}_L+ta$ is a valid path (i.e., $\check{m}_L+ t a \in \M_L$ for all $t$ in some neighborhood of zero) for an arbitrary but ``smooth'' $a$ then it is relatively straightforward to establish the approximate zero property~\citep{NeweyStroker93}.\footnote{ As $\theta\in \Theta$ is restricted to have norm $1$, ${\theta} + t \eta$ does not belong to the parametric space for $t\neq 0$ and $\eta^\top {\theta} \ne 0$. However, this can be easily remedied by considering another path that is differentiable and has the same ``direction''; we define such a path in~\eqref{eq:path_para}.} 
However, this approach does not work when the nonparametric function $m_0$ is constrained. This is because under constraints,  $\check{m}_L + t a$ might not be a valid path for arbitrary but smooth $a$.
% The novelty of our proposed approach lies in using the fact that  $\check{m}_L + t a$ might not be a valid path but $\Pi_{\mathcal{M}_L}(\check{m}_L + ta)$ is always a valid path, here $\Pi_{\mathcal{M}_L}(f)$ for any $f\in L_2$ denotes the projection of $f$ into $\mathcal{M}_L$.
 The novelty of our proposed approach lies in observing that in contrast to $t\mapsto\check{m}_L + ta$, $t \mapsto\Pi_{\mathcal{M}_L}(\check{m}_L + ta)$ is always a valid path for every smooth $a$; here $\Pi_{\mathcal{M}_L}(f)$ is the $L_2$-projection of $f$ onto $\mathcal{M}_L$. Thus if $t\mapsto \Pi_{\mathcal{M}_L}(\check{m}_L + ta)$ is differentiable, then $\partial Q_n(\Pi_{\mathcal{M}_L}(\check{m}_L + ta), \check{\theta}_L + t\eta)/\partial t|_{t = 0} = 0$ for any perturbation $(a,\eta)$. Then establishing that $\check{\theta}_L$ is an approximate zero boils down to finding an $a$ such that 
 \[
 \frac{\partial }{\partial t}Q_n(\Pi_{\mathcal{M}_L}(\check{m}_L + ta), \check{\theta}_L + t\eta)\Big|_{t = 0}= n^{-1}\sum_{i=1}^n \eta^\top  \ell_{\check{\theta}_L, \check{m}_L}(Y_i, X_i) +o_p(n^{-1/2}).
 \]
  Differentiability of projection operators is well-studied; e.g., see~\cite{dharanipragada1996quadratically,fitzpatrick1982differentiability,mccormick1972gradient,shapiro1994existence,sokolowski1992shape} for sufficient conditions for a general projection operator to be differentiable. The generality and the usefulness of our technique can be understood from the fact that no specific structure of $\check{m}_L$ or $\mathcal{M}_L$ is used in the previous discussion; we elaborate on this in Section~\ref{sec:SemiCLSE}. On the other hand, existing methods (see e.g.,~\cite{VANC}) require delicate (and not generalizable) use of the structure of the nonparametric estimator to create valid paths around the nonparametric function; see e.g.,~\cite{VANC} for semiparametric efficiency in  current status regression, and~\cite[Chapter 9.3]{VdV02} and~\cite{MR1394975} for efficiency in the Cox proportional hazard model with current status data; see Example~\ref{ex:coxmodel}.

% Our method does not require the user to find the exact characterization of the nuisance tangent space but only a large enough subset of the nuisance tangent space (so that $\ell_{\check{m},\check{\theta}}$ is close to this subset, as above). This can be very useful because in many  cases when the estimator is non-smooth and/or the estimator satisfies shape constraints. Because in such cases the nuisance tangent space (at the estimator)  can be quite complicated and hard to characterize (see cox prop hazard model; will look after class).

% The proof of Theorem~\ref{thm:Main_rate_CLSE} is complicated as the nuisance tangent space at $\check{m}_{L}$ is not easy to quantify.  This behavior is due to the fact that $\check{m}_{L}$ (a piecewise affine function) lies on the boundary of the set of all convex functions, i.e., the nuisance tangent space at $\check{m}_{L}$ is not $L_2(\Lambda)$ {\clr What is $\Lambda$?}; see Section~\ref{sec:SemiCLSE} for more details. {\clr Anything lost if we remove this$\rightarrow$ Our analysis is further complicated by the fact that $m_0$ and $\theta_0$ are \textit{bundled} (where the parametric and nonparametric components are intertwined; see~\cite{huang1997interval})}.  

% \todo[inline]{In summary...}

% subsection contributions_of_the_paper (end)

\subsection{Organization of the exposition} % (fold)
\label{sub:description_of_exposition}

% subsection description_of_exposition (end)
Our exposition is organized as follows: in Section \ref{sec:Estim}, we introduce some notation and formally define the CLSE. In Section \ref{sec:CLSE}, we state our assumptions, prove consistency, and give rates of convergence for the  CLSE. 
% In Section~\ref{sec:AsymRegFunEstimate} we analyze $\check{m}_{L_n,n}(\check{\theta}_{L_n,n}^\top \cdot)$, the estimator for the regression function $m_0(\theta_0^\top \cdot)$,  and in Section~\ref{sec:SepPara} we analyze the CLSE for $m_0$ and $\theta_0$ separately. 
In Section~\ref{sec:SemiInf}, we detail our new method to prove semiparametric efficiency of the CLSE.  We  use this to prove $\sqrt{n}$-consistency, asymptotic normality,  and efficiency (when the errors happen to be homoscedastic) of the CLSE of $\theta_0$. We discuss an algorithm to compute the proposed estimator in Section~\ref{sec:compute}. In Section~\ref{sec:Simul_Cvx}, we provide an extensive simulation study and compare the finite sample performance of the proposed estimator with existing methods in the literature. In Section~\ref{sec:real_data_analysis}, we analyze the Boston housing data \cite{harrison1978hedonic} and the car mileage data \cite{cars_1983} introduced in Examples~\ref{ex:boston} and~\ref{ex:car} in more details. In both of the cases,  we show that the natural shape constraint leads to stable and interpretable estimates.  Section~\ref{sec:discussion} provides a brief summary of the paper and discusses some open problems. 

{Section numbers in the supplementary file are prefixed with ``S.". Section~\ref{app:sketchCLSE} of the supplementary file provides some insights into the proof of Theorem~\ref{thm:Main_rate_CLSE},  one of our main results. 
Section~\ref{app:add_simul} provides further simulation studies. Section~\ref{sec:proof:Estims}  provides additional discussion on the identifiability of the parameters. Sections~\ref{app:proof:existanceCLSE}--\ref{sec:proof_semi} contain the proofs of our results.  Section~\ref{sec:proof_of_eq:app_score_equation}  completes our novel proof of semiparametric efficiency sketched in Section~\ref{sec:SemiCLSE}.}

\section{Notation and Estimation}\label{sec:Estim}
\subsection{Preliminaries} \label{sec:prelim}
% \subsection{Notation}

In what follows, we assume that we have i.i.d.~data $\{(X_i,Y_i)\}_{i=1}^{n}$ from~\eqref{eq:simsl}. We start with some notation. Let  $\rchi \subset \mathbb{R}^d$ denote the support of $X$ and define
\begin{equation}\label{eq:Defn_Q_D}
D \coloneqq\text{conv} \{\theta^{\top}x :\,  x\in\rchi, \theta\in \Theta\},\quad D_{\theta}\coloneqq \{ \theta^\top x: x \in \rchi\}, \quad\text{and}\quad D_{0} \coloneqq D_{\theta_0},
\end{equation}
 where $\text{conv}(A)$ denotes the convex hull of the set $A$. Let $\M_L$ denote the class of real-valued convex functions on $D$ that are uniformly Lipschitz with Lipschitz bound $L.$ For any $m\in \M_L$, let $m^\prime$ denote the nondecreasing right derivative of the real-valued convex function $m$. Because $m$ is a uniformly Lipschitz function with Lipschitz constant $L$, without loss of generality, we can assume that $|m'(t)|\leq L$, for all $t\in D.$ We use $\p$ to denote the probability of an event and $\E$ for the expectation of a random quantity. For any $\theta\in\Theta$, let $P_{\theta^\top X}$ denote the distribution of $\theta^\top X$.  For $g : \rchi \to \mathbb{R}$, define $  \|g\|^2 \coloneqq \int g^2(x) dP_X(x).$ Let $P_{\epsilon, X}$ denote the joint distribution of $(\epsilon, X)$ and let $P_{\theta,m}$ denote the joint distribution of $(Y,X)$ when $Y = m(\theta^\top X) +\epsilon,$ where $\epsilon$ is defined in~\eqref{eq:simsl}. In particular,  $P_{\theta_0, m_0}$ denotes the joint distribution of $(Y,X)$ when  $X\sim P_X$ and $(Y,X)$ satisfies \eqref{eq:simsl}.  For any set $I \subseteq \mathbb{R}^p$ ($p\ge 1$) and any function $g: I \to \mathbb{R}$, we define $\|g\|_{\infty} \coloneqq \sup_{u \in I} |g(u)|$ and $\|g\|_{I_1} \coloneqq \sup_{u \in I_1} |g(u)|,$ for $I_1 \subseteq I.$  The notation $a\lesssim b$ is used to express that $a\le C b$ for some  constant $C>0$.  For any function $f:\rchi\rightarrow\R^r, r\ge1$, let $\{f_i\}_{1\le i \le r}$ denote each of the components of $f$, i.e.,  $f(x)= (f_1(x), \ldots, f_r(x))$ and  $f_i: \rchi \to \R$. We define $\| f\|_{2, P_{\theta_0, m_0}}\coloneqq \sqrt{ \sum_{i=1}^r \|f_i\|^2}$  and $\| f\|_{2, \infty}\coloneqq \sqrt{ \sum_{i=1}^r \|f_i\|^2_\infty}.$ For any  function $g:D\rightarrow\R$ and $\theta\in\Theta$, we define $(g\circ\theta)(x) \coloneqq g(\theta^{\top}x),$ for all $ x \in \rchi.$ We use the following (standard) empirical process theory notation. For any function $f: \R\times \rchi\rightarrow \R$, $\theta \in \Theta,$ and $m :\R \rightarrow \R$, we define \[P_{\theta, m} f \coloneqq \int f(y,x) dP_{\theta, m}(y,x).\] Note that   $P_{\theta, m} f$ can be a random variable when $\theta$ or $m$ or both are random. Moreover, for any function $f: \R\times \rchi\to \R$, we define $\p_nf \coloneqq n^{-1}\sum_{i=1}^n f(Y_i, X_i)$ and $\g_n f\coloneqq {\sqrt{n}} (\p_n -P_{\theta_0,m_0}) f.
$
\subsection{Identifiability}\label{sec:ident} We now discuss the identifiability of $m_0\circ\theta_0$ and $(m_0, \theta_0)$. Letting $Q(m,\theta) \coloneqq \E[Y -m(\theta^\top X)]^2,$ observe that $(m_0, \theta_0)$ minimizes $Q(\cdot, \cdot).$ In fact we can show in Section~\ref{proof:lem:Ident_cvxSim}, that
\begin{equation}\label{eq:true_mimina}
\inf_{\{(m, \theta):\; m\circ\theta \in L_2(P_X) \text{ and } \|m\circ\theta - m_0\circ\theta_0\| > \delta\}} \big[Q(m,\theta) - Q(m_0,\theta_0)\big] > \delta^2, \quad \text{for any} \quad \delta > 0.
\end{equation}
This implies that $m_0\circ\theta_0$ is always identifiable and further, one can hope to consistently estimate $m_0\circ\theta_0$ by minimizing  the sample version of $Q(m, \theta)$; see~\eqref{eq:CLSE}.

Note that the identification of $m_0\circ\theta_0$ does not guarantee that both $m_0$ and $\theta_0$ are separately identifiable. Hence, in what follows, when dealing with the properties of separated parameters, we will directly assume:
\begin{enumerate}[label=\bfseries (A\arabic*)]
\setcounter{enumi}{-1}
 \item  The parameters $m_0\in \M_{L_0}$ and $\theta_0\in \Theta$ are separately identifiable, i.e., $m\circ \theta = m_0\circ\theta_0$ for some $(m, \theta)\in \M_{L_0}\times \Theta$ implies that $m = m_0$ and $\theta= \theta_0$.\label{a0}
 \end{enumerate} 
 \citet{ICHI93} has found general sufficient conditions on the distribution of $X$ under which~\ref{a0} holds; these sufficient conditions allow for some components of $X$ to be discrete, also see~\citet[Pages 12--17]{Horowitz98} and~\citet[Proposition 8.1]{liracine07}. When $X$ has a density with respect to Lebesgue measure,~\citet[Theorem 1]{LinKul07} find a simple sufficient condition for~\ref{a0}. We  discuss and compare these two sufficient conditions in Section~\ref{sec:IdentDisc} of the supplementary file.

%{\color{red} \cite{ICHI93} gives conditions for identifiability of $\theta$ and $m$.  Should we just use the same conditions?}

%\section{Two estimators} \label{sec:Estim}
%In Lemma~\ref{lem:Ident_cvxSim}, we show that $(m_0, \theta_0)$ minimizes the population version of the square error loss. In this section, we use $Q_n(m,\theta),$ the empirical version of the square error loss, to formally define the CLSE and the PLSE.
% using the notation developed in Section \ref{sec:prelim}.
\section{Convex and Lipschitz constrained LSE} \label{sec:CLSE}
Recall that CLSE is defined as the minimizer of $(m, \theta) \mapsto Q_n(m, \theta)$ over $\M_L\times \Theta$. Because $Q_n(m, \theta)$ depends   only on the values of the function at  $\{{\theta} ^\top X_i\}_{i=1}^n$,  it is immediately clear that  the minimizer $\check{m}_{{L}}$ is  unique only at $\{\check{\theta}_{{L}} ^\top X_i\}_{i=1}^n$. Since $\check{m}_{{L}}$ is restricted to be convex, we interpolate the function linearly between $\check{\theta}_{{L}} ^\top X_i$'s  and extrapolate the function linearly outside the data points.\footnote{\label{foo:extrapolate}Linear interpolation/extrapolation does not violate the convexity or the $L$-Lipschitz property} Thus $\check{m}$ is piecewise affine. In Section~\ref{app:proof:existanceCLSE} of the supplementary file,  we prove the existence of the minimizer in \eqref{eq:CLSE}. The optimization problem~\eqref{eq:CLSE} might not have a unique minimizer and the results that follow hold true for any global minimizer. 
 % In Sections \ref{sec:AsymRegFunEstimate} and  \ref{sec:SemiCLSE} we study the  $(\check{m}_{{L}},\check{\theta}_{{L}})$ is a consistent estimator of $(m_0, \theta_0)$ and study its asymptotic properties, respectively.
\begin{remark}\label{rem:compute}
 For every fixed $\theta$,  $m (\in  \M_L) \mapsto Q_n(m,\theta)$ has a unique minimizer. The minimization over the class of uniformly Lipschitz functions is a quadratic program with linear constraints and can be computed easily; see Section~\ref{sec:CLSE_comp}.
\end{remark}

\subsection{Asymptotic analysis of the regression function estimate}\label{sec:AsymRegFunEstimate}
In this section, we study the asymptotic behavior of $\check{m}_{{L}}\circ\check{\theta}_{{L}}$. We will now list the assumptions under which we study the rates of convergence of the CLSE for the regression function.
\begin{enumerate}[label=\bfseries (A\arabic*)]
\setcounter{enumi}{0}
\item \label{aa1_new} The unknown convex link function $m_0$ is bounded by some constant $M_0$ $ (\ge 1)$ on $D$ and is uniformly Lipschitz with Lipschitz constant $L_0$. 
\item \label{aa1} The support of $X$, $\rchi$, is a  subset of $\mathbb{R}^d$ and  $\sup_{x\in \rchi} |x| \le T,$ for some finite $T\in \R$.
\item \label{aa2}The error $\epsilon$ in model~\eqref{eq:simsl} has finite  $q$th moment, i.e., $K_q\coloneqq \big[\E(|\epsilon|^q)\big]^{1/q}  < \infty$ where $q\ge 2$. Moreover, $\mathbb{E}(\epsilon|X) = 0,$ $P_X$ a.e.~and $\sigma^2(x)\coloneqq \E(\epsilon^2|X=x) \le \sigma^2 < \infty$ for all $x\in \rchi.$ 

\end{enumerate}

The above assumptions deserve comments. \ref{aa1} implies that the support of the covariates is bounded. 
% However, using arguments similar to Remark 4.2 (and Section 9.2) of~\cite{2016arXiv161006026B} one can relax assumption~\ref{aa1} and allow for distribution $X$ to be sub-Gaussian. {\clr Write some arguments?} 
In assumption~\ref{aa2}, we allow $\epsilon$ to be heteroscedastic and $\epsilon$ can depend on $X$. Our assumption on $\epsilon$ is more general than those considered in  the shape constrained literature, most works assume that all moments of $\epsilon$ are finite and ``well-behaved'', see e.g., \cite{2017arXiv171205593B},~\cite{Hristacheetal01}, and~\cite{Xiaetal02}.

Theorem \ref{thm:rate_m_theta_CLSE} (proved in Section~\ref{proof:rate_m_theta_CLSE}) below  provides an upper bound on the rate of convergence of $\check{m}_L\circ\check{\theta}_L$ to $m_0\circ\theta_0$ under the $L_2(P_X)$ norm. The following result is a finite sample result and shows the explicit dependence of the rate of convergence on $L = L_n, d,$ and $q$.

\begin{thm}\label{thm:rate_m_theta_CLSE}
Assume \ref{aa1_new}--\ref{aa2}. Let $\{L_n\}_{n\ge 1}$ be a fixed sequence such that $L_n\ge L_0$ for all $n$ and let 
\begin{equation}\label{eq:r_n_def}
r_n\coloneqq\min \left\{  \frac{n^{2/5}}{d^{2/5}L_n},  \frac{n^{1/2-1/2q}}{L_n^{(3q+ 1)/(4q)}} \right\}.
\end{equation}
 Then for every $n\ge 1$ and $u\ge 1$, there exists a constant $\mathfrak{C}\ge 0$ depending only on $\sigma, M_0, L_0, T,$ and  $K_q$, and constant $C$ depending only on $K_q, \sigma$, and $q,$ such that
  \[ \sup_{\theta_0, m_0,\epsilon, X} \p\left( r_n \|\check{m}_{L_n}\circ\check{\theta}_{L_n} - m_0\circ\theta_0\| \ge u \mathfrak{C}\right) \le \frac{C}{u^{q}}+ \frac{\sigma^2}{n},\]
where the supremum is taken over all $\theta_0\in \Theta$ and all joint distributions of $(\epsilon,X)$ and parameters $m_0$ for which assumptions~\ref{aa1_new}--\ref{aa2} are satisfied with constants $\sigma, M_0, L_0, T,$ and  $K_q$. In particular if $q\ge 5$, $d = O(1),$ and $L_n = O(1)$ as $n\to\infty$, then $\|\check{m}_{L_n}\circ\check{\theta}_{L_n} - m_0\circ\theta_0\| = O_p({n^{-2/5}}).$
\end{thm}
Note that~\eqref{eq:r_n_def} allows for the dimension $d$ to grow with $n$ and $\theta_0$ to change with $n$. For example if $L_n\equiv L$ for some fixed $L\ge L_0$, then we have that $\|\check{m}_{L_n}\circ\check{\theta}_{L_n} - m_0\circ\theta_0\| =o_p(1)$  if $d = o(n^{1 - 1/q})$. In  the rest of the paper, we assume that $d$ is fixed.  In Proposition~\ref{lem:MinimaxLowerbound} in Section~\ref{sec:minimax_lower_bound}, we find the minimax lower bound for the single index model ~\eqref{eq:simsl}, and show that $\check{m}_{L}\circ\check{\theta}_L$ is minimax rate optimal when $q\ge 5$.

  The next result shows that the rates  in Theorem~\ref{thm:rate_m_theta_CLSE} are in fact uniform (up to a $\sqrt{\log\log n}$ factor) in $ L \in [L_0, n L_0]$.
  % \footnote{The dependence on $L$ in the second term of the rate can be improved via a slightly more careful calculation.} 
This uniform-in-$L$ result is important for the study of the estimator with a data-driven choice of $L$ such as cross-validation or Lepski's method~\cite{lepski1997optimal}. Theorem~\ref{thm:rate_m_theta_CLSE} alone cannot provide such a rate guarantee because it requires $L$ to be non-stochastic.

\begin{thm}\label{thm:UniformLRate}
Under the assumptions of Theorem~\ref{thm:rate_m_theta_CLSE}, the CLSE satisfies 
\begin{equation}\label{eq:UniformLRate}
\sup_{L_0 \le L \le nL_0}\,\min\left\{\frac{n^{2/5}}{L}, \frac{n^{1/2 - 1/(2q)}}{\sqrt{L}}\right\}\|\check{m}_L\circ\check{\theta}_L - m_0\circ\theta_0\| = O_p\left(\sqrt{\log\log n}\right).
\end{equation}
\end{thm}
\begin{remark}[Diverging $L$]\label{rem:DependenceonL}   The dependence on $L$ in Theorems~\ref{thm:rate_m_theta_CLSE} and~\ref{thm:UniformLRate} suggest that the estimator may not be consistent if $L \equiv L_n$ diverges too quickly with the sample size. The simulation in Section~\ref{sub:robustness_of_choice_of_} suggests that the estimation error has negligible dependence on $L$ and that the dependence on $L$ in Theorems~\ref{thm:rate_m_theta_CLSE} and~\ref{thm:UniformLRate} might be sub-optimal. We believe this discrepancy is due to the lack of available technical tools to prove uniform boundedness of the estimator $\check{m}_{n,L}$ in terms of $L$. At present, we are only able to prove that with high probability, $\|\check{m}_{n,L}\|_{\infty} \le LT + M_0 + 1$ for all $L \ge L_0$; see Lemma~\ref{lem:Upsilion_ep}. If one can prove $\|\check{m}_{n,L}\|_{\infty} \le C$ for all $L \ge L_0$, with high probability, for a constant $C $ independent of $L$, then our proofs can be modified to remove the dependence on $L$ in Theorems~\ref{thm:rate_m_theta_CLSE} and~\ref{thm:UniformLRate}. 
\end{remark}

% \begin{remark}\label{rem:Tailbound}
% The results above are possibly the first (in the shape constrained literature) deriving the rates for the (constrained) LSE  when the errors have only finite moments and can depend on the covariates. See \cite{han2017sharp,2018arXiv180502542H} for a discussion for results when errors are independent of the covariates. Also, see Section~\ref{sec:MajorGeneral} for more details. 
% \end{remark}

\subsection[Asymptotic analysis of the separated parameters]{Asymptotic analysis of~$\check{m}$ and $\check{\theta}$}\label{sec:SepPara}
In this section we establish the consistency and find rates of convergence of $\check{m}_{L_n}$ and $\check{\theta}_{L_n}$ separately. In Theorem~\ref{thm:rate_m_theta_CLSE} we proved that $\check{m}_{L_n}\circ\check{\theta}_{L_n}$ converges in the $L_2(P_{\theta_0, m_0})$ norm but that does not guarantee that $\check{m}_{L_n}$ converges to $m_0$ in the $\|\cdot\|_{D_0}$ norm. A typical approach for proving consistency of $\check{m}_{L_n}$ is to prove that $\{\check{m}_{L_n}\}$ is precompact in the $\|\cdot\|_{D_0}$ norm ($D_0$ is defined in~\eqref{eq:Defn_Q_D}); see e.g.,~\cite{2016arXiv161006026B,VANC}. The Arzel\`{a}-Ascoli theorem establishes that the necessary and sufficient condition for compactness (with respect to the uniform norm) of an arbitrary class of continuous functions on a bounded domain is that the function class be uniformly bounded and equicontinuous. However, if $L_n$ is allowed to grow to infinity, then it is not clear whether the sequence of functions $\{\check{m}_{L_n}\}$ is equicontinuous.  Thus to study the asymptotic properties of $\check{m}_{L_n}$ and $\check{\theta}_{L_n},$ we assume that $L_n\equiv L \ge L_0$, is a fixed constant.  For the rest of paper, we will use $\check{m}$ and $\check\theta$ to denote $\check{m}_L$ (or $\check{m}_{{L_n}}$) and $\check\theta_L$ (or $\check{\theta}_{{L_n}}$), respectively. The next theorem (proved in Section~\ref{proof:thm:uucons}) establishes consistency of $\check{m}$ and $\check{\theta}$ separately. Recall that $m_0^\prime$ denotes the nondecreasing right derivative of the  convex function $m_0$. 
 \begin{thm}\label{thm:uucons} Suppose the assumptions of Theorem~\ref{thm:rate_m_theta_CLSE} and~\ref{a0} hold. Then, for any fixed $L\ge L_0$ and any compact subset $C$ in the interior of $D_0$, we have
\[
|\check{\theta} - \theta_0| = o_p(1), \qquad  \|\check{m} - m_0\|_{D_0} = o_p(1), \quad \text{and} \quad \|\check{m}' - m'_0\|_{C} = o_p(1). \]
%\todo[inline]{Proof has gap}
% Under assumptions \ref{c1}--\ref{c4}, the constrained LSE satisfies $\check{\theta}_n\overset{P}{\to}\theta_0$ and $\|\check{m}_n\circ\theta_0 - m_0\circ\theta_0\| = O_p(n^{-2/5}L_n^{4/5}).$
\end{thm}

Fix an orthonormal basis  $\{e_1,\ldots, e_d\}$ of $\R^d$ such that $e_1 =\theta_0.$ Define $H_{\theta_0}\coloneqq [e_2,\ldots,e_d]\in\mathbb{R}^{d\times(d-1)}$. We will use the following two additional assumptions to establish upper bounds on the rate of convergence of $\check{m}$ and $\check{\theta}.$

% \todo[inline]{Do we need uniform consistency on $D_0$ or just compact subsets would do? \\Ans: We need uniform consistency. Otherwise we would not be able to establish rates of convergence for $\check{m}'$.}

\begin{enumerate}[label=\bfseries (A\arabic*)]
\setcounter{enumi}{3}
% \item  is a nonsingular matrix. \label{aa4}
\item $ H_{\theta_0}^\top\E\big[\Var(X|\theta_0^\top X) \{m_0'(\theta_0^{\top}X)\}^2 \big]H_{\theta_0}$ is a positive definite matrix.
% The conditional distribution of $X$ given $\theta_0^{\top}X$ is nondegenerate.
% {\clr Where is this used, and how is this connected to the semi parametric part. {\cln Used to Theorem~\ref{thm:ratestCLSE}, see~\eqref{eq:G1SquareBound}, not the semiparametric part}}
 \label{aa5_new}
% \item $\Var(X)$ is a positive definite matrix. \label{aa5}

\item  \label{aa6} The density of $\theta_0^\top X$ with respect to the Lebesgue measure is bounded above by $\overline{C}_d < \infty$. 
% \item There exist positive constants $\overline{C}_d, \underline{C_d},$ and $r$,  such that for every $\theta\in B_{\theta_0}(r),$ the density of $\theta^\top X$ with respect to the Lebesgue measure  on $D_\theta$ is bounded below by $\underline{C_d}$ and bounded above by $\overline{C}_d$.

\end{enumerate}
Assumption~\ref{aa5_new}, is used to find the rate of convergence for $\check\theta$ and $\check{m}$ separately and is widely used in all works studying root-$n$ consistent estimation of $\theta_0$ in the single index model, see e.g.,~\cite{Powelletal89,ICHI93,Patra16,2017arXiv171205593B}; also see~Remark~\ref{rem:SingularInformation}. \ref{aa6} is mild, and is satisfied if $X=(X_1,\ldots, X_d)$ has a continuous covariate $X_k$ such that: (1) $X_k$ has a bounded density; and (2) $\theta_{0,k}>0$. Compare assumption~\ref{aa6} with~\citep{ICHI93,cuietal11,2017arXiv171205593B,MR2529970,wang2015spline} where it is assumed that $\theta^\top X$ has a density bounded away from zero for all $\theta$ in a neighborhood of $\theta_0$.
Assumption~\ref{aa6} is used to find rates of convergence of the derivative of the estimators of $m_0$. In Theorem~\ref{thm:ratestCLSE}, we only use the fact that $\theta_0^\top X$ is absolutely continuous with respect to Lebesgue measure. The following result (proved in Section~\ref{sec:proof_ratestPLSE}) establishes upper bounds on the rate of convergence of $\check{\theta}$ and $\check{m}$ respectively.
 \begin{thm}\label{thm:ratestCLSE}
If assumptions~\ref{a0}--\ref{aa6} hold, $q\ge 5$, and $L\ge L_0$, then we have 
\[
|\check{\theta} - \theta_0| = O_p(n^{-2/5}) \quad \text{and} \quad  {\int (\check{m}(t)- m_0 (t))^2 dP_{\theta_0^{\top}X}(t) dt = O_p(n^{-4/5})}.\]
% \begin{equation}\label{eq:IntegralFtheta_011}
% \int_{D_0} \left(\check{m}_n(t) - m_0(t)\right)^2f_{\theta_0^{\top}(X)}(t) dt = O_p(n^{-4/5}).
% \end{equation}}
\end{thm}
\begin{remark}\label{rem:SingularInformation}
Note that, under homoscedastic errors in~\eqref{eq:simsl}, the efficient information for $\theta_0$ is a scalar multiple of $H_{\theta_0}^\top\E\big[\Var(X|\theta_0^\top X) \{m_0'(\theta_0^{\top}X)\}^2 \big]H_{\theta_0}=: \mathcal{I}_0$; see Section~\ref{sec:eff_score}. If $\mathcal{I}_0$ is  not positive definite, then there is zero information for $\theta_0$ along some directions. In that case, we can show that $| \mathcal{I}_0^{1/2} (\check{\theta} - \theta_0)| = O_p(n^{-2/5})$; see~\eqref{eq:step_th6} in the supplementary file. 
\end{remark}
A simple modification of the proof of Proposition~\ref{lem:MinimaxLowerbound} will prove that $\check{m}$ is also minimax rate optimal. {Under additional smoothness assumptions on $m_0$, in the following theorem (proved in Section~\ref{proof:DerivCLSE}) we show that $\check{m}'$,  the right derivative of $\check{m},$ converges to $m'_0$ in both the $L_2$  and the supremum norms.

\begin{thm}\label{thm:rate_derivCLSE} Suppose assumptions  of Theorem~\ref{thm:ratestCLSE} hold and $m_0'$ is $\nicefrac{1}{2}$-H\"{o}lder continuous on $D_0$, then 
\begin{equation} \label{eq:deriv_theta_0}
{\|\check{m}^\prime \circ \theta_0- m'_0\circ\theta_0\| = O_p\big(n^{-2/15}\big)}\quad \text{and} \quad  \|\check{m}^\prime \circ \check\theta- m'_0\circ\check\theta\|= O_p\big(n^{-2/15}\big).
\end{equation}
Further, if $m_0$ is twice continuously differentiable and assumption~\ref{bb2'} (in Section~\ref{sec:SemiInf}), then for any compact subset $C$ in the interior of $D_0$, we have % {\clr such that $ f_{\theta_0^{\top}(X)}$ is bounded away from $0$ on $C$, } we have that
\begin{equation}\label{eq:sup_rate}
\sup_{t\in C} |\check{m} (t)- m_0 (t)| = O_p(n^{-8/(25 + 5\beta)})%n^{-2/5} n^{(2 + 2\beta)/(25 + 5\beta)}
 \quad \text{and} \quad \sup_{t\in C} |\check{m}^\prime (t)- m'_0 (t)| = O_p(n^{-4/(25 + 5\beta)}).%n^{-1/5} n^{(1+\beta)/(25 + 5\beta)}
\end{equation}
\end{thm}
% The fact that $\check{m}^\prime$ is a step function complicates the proof of the above result (given in Section~\ref{proof:DerivCLSE} of the supplementary file). 
\begin{remark}\label{rem:TwiceCont}
As in~\eqref{eq:deriv_theta_0},~\eqref{eq:sup_rate} can also be proved under $\gamma$-H\"{o}lder continuity of $m_0'$, but in this case the rate of convergence depends on $\gamma$ explicitly. Assumption~\ref{bb2'} allows for the density of $\theta_0^{\top}X$ to be zero at some points in its support; see Section~\ref{sec:SemiInf} for a detailed discussion. Further if the density of $\theta_0^{\top}X$ is bounded away from zero, then $\beta$ can be taken to be $0$. 
\end{remark}
\begin{remark}\label{rem:lowerQ}
The condition $q\ge5$ in Theorems~\ref{thm:ratestCLSE} and~\ref{thm:rate_derivCLSE} can be relaxed at the expense of slower rates of convergence. { In fact, by following the arguments in the proofs, we can show, with $p_n := \max\{n^{-2/5}, n^{-1/2 + 1/(2q)}\}$ for any $q\ge2$, that $|\check{\theta} - \theta_0| = O_p(p_n)$, and 
\[
\|\check{m}\circ \theta_0- m_0\circ \theta_0 \|  = O_p(p_n), \quad \|\check{m}^\prime \circ \theta_0- m'_0 \circ \theta_0\| = O_p(p_n^{1/3})\quad \text{and} \quad  \|\check{m}^\prime \circ \check\theta- m'_0\circ\check\theta\|= O_p(p_n^{1/3}).
\]}
% where
%  \[
%   p_n \coloneqq \min\left\{{n^{2/5}},{n^{1/2 - 1/(2q)}}\right\}.
% \]
\end{remark}
% }% In fact, the obtained rate need not be optimal, but is sufficient for our purposes (in deriving the efficiency of $\check\theta$; see Section~\ref{sec:SemiCLSE}).
\begin{remark}[Additional shape constraints on the link function]\label{rem:Monotone} {It might often be the case that in addition to convexity, the practitioner is interested in {imposing additional shape constraints (such as monotonicity, unimodality, or $k$-monotonicity~\cite{MR3881209}) on $m_0$.} For example, in the datasets considered in Examples~\ref{ex:boston} and~\ref{ex:car}, the link function is {plausibly} both convex \textit{and} monotone; see~\cite{CHSA} for further motivation {on} additional shape constraints. The conclusions (and proofs) of {Theorems~\ref{thm:rate_m_theta_CLSE} and~\ref{thm:UniformLRate}--\ref{thm:rate_derivCLSE}} also hold for the  CLSE under additional constraints {on the link function}. An intuitive explanation is that the parameter space $\M_L$ is only reduced by imposing {additional constraints} on the link function and this can only give better rates (if not the same). In case of an additional monotonicity constraint on $m_0$, one can modify the proof of Proposition~\ref{lem:MinimaxLowerbound} to show that {the} rate obtained in Theorem~\ref{thm:rate_m_theta_CLSE} is in fact minimax optimal for {the} the CLSE  (under further monotonicity constraint). }

\end{remark}

\section{Semiparametric inference for the CLSE} \label{sec:SemiInf}
The main result in this section shows that $\check{\theta}$ is $\sqrt{n}$-consistent and asymptotically normal; see Theorem~\ref{thm:Main_rate_CLSE}. Moreover,~$\check{\theta}$ is shown to be semiparametrically efficient for $\theta_0$ if the errors happen to be homoscedastic. The asymptotic analysis of $\check \theta$  is involved  as $\check{m}$ is a piecewise affine function and hence not differentiable everywhere.
%As stated in Theorem \ref{thm:existanceCLSE}, the estimate $\check{m}$ of the link function $m_0$ is a piecewise affine function and hence not differentiable everywhere. This complicates the asymptotic analysis of $\check{\theta}$, relative to $\hat{\theta}$.

Before deriving the limit law of $\check{\theta}$, we introduce some  notations and assumptions.
 Let $p_{\epsilon,X}$ denote the joint density (with respect to some dominating measure on $\R \times \rchi$) of $(\epsilon, X)$. Let $p_{\epsilon|X} (\cdot,x)$ and $p_X(\cdot)$ denote the corresponding conditional probability density of $\epsilon$ given $X = x$ and the marginal density of $X$, respectively. In the following we list additional assumptions used in Theorem~\ref{thm:Main_rate_CLSE}. Recall $D$ and $D_0$ from~\eqref{eq:Defn_Q_D} and let  $\Lambda$ denote the Lebesgue measure.
% \begin{enumerate}[label=\bfseries (B\arabic*)]
% \setcounter{enumi}{0}
%  \item For every $\theta\in B_{\theta}(r, \theta_0)$, $D_\theta$ is strict subset of $D^{(r)}$. For the rest of the paper we redefine $D\coloneqq D^{(r)}$.\label{bb1}
%  \end{enumerate}
 \begin{enumerate}[label=\bfseries (B\arabic*)]
\setcounter{enumi}{0}
 \item $m_0\in \M_{L_0}$ and  $m_0$ is $(1+ \gamma)$-H\"{o}lder continuous on $D_0$ for some $\gamma>0$. Furthermore, $m_0$ is strongly convex on $D$, i.e., there exists a $\kappa_0>0$ such that  $m_0(t)- \kappa_0 t^2$ is convex.   \label{bb1}

 \item  There exists $\beta\ge 0$ and $\underline{C}_d>0$ such that $\p(\theta_0^\top X \in I) \ge \underline{C}_d\,\Lambda(I)^{1 + \beta},$  for all intervals $I \subset D_0$. \label{bb2'}
  
% \item Assume that there exist  $r>0$ such that for all $\theta \in S^{d-1}\cap B_{\theta_0}(r)$ we have
% $$D_\theta \subsetneq D.$$\label{b1}
\end{enumerate}
For every $\theta\in\Theta$, define  $h_\theta(u) \coloneqq \E[X |\theta^\top X=u]$.
%  \begin{align}
% h_\theta(u) &\coloneqq \E[X |\theta^\top X=u].\label{eq:h_beta}
% \end{align}
 \begin{enumerate}[label=\bfseries (B\arabic*)]
\setcounter{enumi}{2}
{\item 
The function $u\mapsto h_{\theta_0}(u)$ is $1/2$-H\"{o}lder continuous and for a constant $\bar{M}>0$,\label{bb2}
\begin{equation}\label{eq:L_2Lip}
\E\Big(|h_{\theta}(\theta_0^{\top}X) - h_{\theta_0}(\theta_0^{\top}X)|^2\Big) \le \bar{M}|\theta - \theta_0|\quad\mbox{for all}\quad \theta\in\Theta.
\end{equation}}
\item The  density $p_{\epsilon|X} (e,x)$ is differentiable with respect to $e$ for all $x\in \rchi$.
\label{bb3}
\end{enumerate}

Assumptions \ref{bb1}--\ref{bb3} deserve comments. \ref{bb1} is much weaker than the  standard assumptions used in semiparametric inference in  single index models~\cite[Theorem 3.2]{VANC}. Assumption~\ref{bb2'} is an improvement compared to the assumptions in the existing literature. Assumption~\ref{bb2'} pertains to the distribution of $\theta_0^\top X$ and is inspired by~\cite[assumption~(D)]{MR2369025}. In contrast, most existing works require the density of $\theta_0^\top X$ to be bounded away from zero (i.e., $\beta=0$); see e.g.,~\cite[Assumption 5.3(II)]{ICHI93},~\cite[Assumption (d)]{cuietal11},~\cite[Lemma F.3]{2017arXiv171205593B},~\cite[Assumption~A2]{MR2529970},~\cite[Assumption (A2)]{wang2015spline}. Our assumption is significantly weaker because it allows the density of $\theta_0^\top X$ to be zero at some points in its support. For example, when $X\sim\text{Uniform}[0,1]^d$, the density of $\theta_0^\top X$ might not be bounded away from zero~\citep[Figure 1]{MR2369025}, but~\ref{bb2'} holds with $\beta=1.$
  Assumption~\ref{bb2} can be favorably compared to those in \cite[Theorem 3.2]{VANC},~\cite[Assumption (A5)]{groeneboom2016current},~\cite[Assumption (A5)]{2017arXiv171205593B}, and~\cite[Assumption G2 (ii)]{song2014semiparametric}.  We use the smoothness assumption~\ref{bb2} when establishing semiparametric efficiency of $\check\theta$. The Lipschitzness assumption~\eqref{eq:L_2Lip} can be verified by using the techniques of~\cite{alonso1998lp}, when $u\mapsto h_{\theta} (u)$ is $1/2$-H\"{o}lder continuous for all $\theta$ in a neighborhood of $\theta_0$ and the H\"{o}lder constants  are uniformly bounded in $\theta$.
% \subsection{Disc} % (fold)
% \label{sub:disc}

% subsection disc (end)
{In general, establishing semiparametric efficiency of an estimator proceeds in two steps. Let $\hat{\xi}$ and $\hat\gamma$ denote the estimators of a parametric component $\xi_0$ and a nuisance component $\gamma_0$ in a general semiparametric model. In a broad sense, the proof of semiparametric efficiency} of $\hat{\xi}$ involves two main steps: (i) finding the efficient score of the model at the truth (call it ${\ell}_{\xi_0, \gamma_0}$); and (ii) proving that $(\hat{\xi},\hat\gamma)$ satisfies $\p_n {\ell}_{\hat\xi, \hat\gamma} = o_p(n^{-1/2})$; see~\cite[pages 436-437]{VdV02} for a detailed discussion. In the Sections~\ref{sec:eff_score} and~\ref{sec:SemiCLSE}, we discuss steps (i) and (ii) in our context, respectively.
%{\clr \textbf{Things to fix}
% \begin{itemize}
%\item change $D$ to ${D}_r$ and $D_\beta$ to ${D}_\beta$.
%\item Fix the order of $(t,\beta)$ in the subscript.
%\item Check what assumptions you need
%\item change $t$ to $s$ in the earlier part of section 4.
%\item  Horowitz 1992 page 509 after assumption 3. similar compactness for the form of $\beta$  arguments are used
%\item Give justifiction for reparamterization
%\end{itemize}
%}

\subsection{Efficient score}\label{sec:eff_score}
In this subsection we  calculate the efficient score for the model:
\begin{equation}\label{eq:Score_model}
Y=m(\theta^\top X)+\epsilon,
\end{equation}
where $m, X,$ and $\epsilon$ satisfy assumptions~\ref{bb1}--\ref{bb3}.  First observe that the parameter space $\Theta$ is a closed subset of $\R^d$ and the interior of $\Theta$  in $\R^d$ is the empty set. Thus to compute the score for model~\eqref{eq:Score_model}, we construct a path on the sphere. We use $\R^{d-1}$ to parametrize the paths for model~\eqref{eq:Score_model} on $\Theta$ when $\theta_{0,1} >0$. For each $\eta\in\mathbb{R}^{d-1}, $ $s \in \R$, and $|s| \le |\eta|^{-1}$,  define the following path , with ``direction'' $\eta$, through $\theta$ (which lies on the unit sphere)
\be\label{eq:path_para}
\zeta_s(\theta,\eta)\coloneqq\sqrt{1-s^2|\eta|^2}\, \theta + s H_\theta \eta,
\ee
 where for every $\theta\in \Theta$, $H_\theta \in \R^{d\times(d-1)}$ is such that for every $\eta \in \R^{d-1}$,  $| H_\theta \eta| =|\eta|$ and $H_\theta\eta$ is orthogonal to $\theta$. Furthermore, we need $\theta \mapsto H_\theta$ to satisfy some smoothness properties; see Lemma 1 of~\cite{Patra16} for such a construction. {Note that, if $\theta_{0,1}=0$, then for any $s$ in a neighborhood of zero, there exists an $\eta\in \R^{d-1}$ such that  $\zeta_s(\theta_0,\eta) \notin \Theta$. Thus, if $\theta_{0,1}=0$, then $\theta_0$ lies on the ``boundary'' of $\Theta$ and the existing semiparametric theory breaks down. Therefore, for the rest of the paper, we assume that $\theta_{0,1}$ is strictly positive.}

%  \begin{enumerate}[label=(H\arabic*)]
% \item $\xi\mapsto H_\theta\xi$ are bijections from $\mathbb{R}^{d-1}$ to the hyperplanes $\{x\in\mathbb{R}^d: \theta^\top  x=0\}$. \label{H1}
% \item The columns of $H_\theta$ form an orthonormal  basis for $\{x\in\mathbb{R}^d: \theta^\top  x=0\}$.\label{H2}
% \item   $\|H_\theta - H_{\theta_0}\|_2\le|\theta- \theta_0 |.$ Here for any matrix $A\in \R^{d_1\times d_2}$, $\|A\|_2\coloneqq \sup_{\{ b\in \R^{d_2}:\, |b|=1\}} |Ab|$.\label{H3}
% \item For all  distinct $\eta,\beta \in \Theta\setminus\theta_0$, such that $|\eta-\theta_0|\le1/2$ and $|\beta-\theta_0|\le1/2,$\label{H4}
% \begin{equation}\label{eq:H_lip_gen}
% \|H_{\eta}-H_{\beta}\|_2 \le 8(1+8/\sqrt{15})\frac{ |\eta-\beta| }{ |\eta-\theta_0|+|\beta-\theta_0|}.
% \end{equation}
% \end{enumerate}%
% See Lemma 1 of~\cite{Patra16} for a construction of a class of matrices satisfying the above properties.
 The log-likelihood  of model~\eqref{eq:Score_model} is $l_{\theta, m}(y,x)= \log[ p_{\epsilon|X} (y-m(\theta^\top x),x) p_X(x)].$ For  any $\eta \in S^{d-2}$, consider the path defined as $s \mapsto \zeta_s(\theta,\eta)$. Note that by the definition of $H_\theta$, $s \mapsto \zeta_s(\theta,\eta)$ is  a valid path in $\Theta$ through $\theta$; i.e., $\zeta_0(\theta,\eta)=\theta$ and $\zeta_s(\theta,\eta) \in \Theta$ for every $s$ in some neighborhood of $0$. Thus the score for the parametric submodel is
\begin{equation}\label{eq:paraScore}
\left. \frac{\partial l_{\zeta_s(\theta,\eta), m} (y,x) }{\partial s}\right\vert_{s=0}= \eta^\top S_{\theta,m}(y,x),
\end{equation}
where
\begin{equation}\label{eq:S_def_118}
S_{\theta,m}(y,x):=   -\frac{p'_{\epsilon|X} \big(y-m(\theta^\top x),x\big)}{p_{\epsilon|X} \big(y-m(\theta^\top x),x\big)}m^\prime(\theta^\top x) H_\theta ^\top x.
\end{equation}

% \begin{remark}\label{rem:qmd}
% % {\clr Before computing the score for the above submodel, note that $m \in \M_L$ need not be differentiable everywhere.  showing that the underlying class of distributions is differentiable in quadratic mean requires some careful analysis; in Remark~\ref{rem:qmd} (in Section~\ref{sec:qmd} of the supplementary file) we show this for the model with Gaussian errors. With this in mind,  the score for the above submodel is }
% \end{remark}
The next step in computing the efficient score for model~\eqref{eq:Score_model} at $(m,\theta)$ is to compute the nuisance tangent space of the model (here the nuisance parameters are $p_{\epsilon|X}, p_X$, and $m$). To do this define a parametric submodel for  the unknown nonparametric components:
\begin{align} \label{eq:nonp_path}
\begin{split}
m_{s,a}(t)&=m(t) - s a(t), \quad p_{\epsilon|X;s, b}(e, x) = p_{\epsilon|X}(e, x) (1 + s b(e, x)), \quad  p_{X; s,q}(x) =p_X(x)(1+s q(x)),
\end{split}
\end{align} where $s\in \R$,   $b: \R \times \rchi\to\R$ is a bounded function such that $\E(b(\epsilon, X) |X)=0$ and $\E(\epsilon b(\epsilon, X) |X)=0$,  $q: \rchi \to \R$ is a bounded function such that $\E(q(X))=0,$ and $a\in\D_m$, with 
\begin{align}\label{eq:d_def_m}
\begin{split}
\D_m := \big\{f \in L_2(\Lambda): f'(\cdot)\;&\text{exists and }m_{s,f}(\cdot) \in \M_{L}\; \text{ for all } s \in B_0(\delta)\;\text{for some }\delta >0\big\}.
\end{split}
\end{align}
% \todo{Why $\M_{L_0}$}
  Note that when $m$ satisfies~\ref{bb1} then $\D_m$ reduces to $ 
 \D_m=  \big\{f \in L_2(\Lambda): f'(\cdot)\;\text{exists}\}.$
    % Since $m_0$ satisfies assumption~\ref{bb1}, we can find  $\delta$ (small enough) such that for all $s\in B_{0}(\delta)$, $m_{s,a}\in \M_L$ and  $p_{\epsilon|X;s, b}$ and $p_{X; s,q}$ are valid densities.
    Thus
% \begin{equation}\label{eq:NonBdry}
     $\overline{\mathrm{lin}}\,\D_m = L_2(\Lambda).$
% \end{equation}
 Theorem~4.1 of~\cite{NeweyStroker93} (also see~\citet[Proposition 1]{MaZhusemi13})  shows that when  the parametric score is $\eta^\top S_{\theta,m} (\cdot, \cdot)$ and the nuisance tangent space corresponding to $m$ is $L_2(\Lambda)$, then the efficient score for model~\eqref{eq:Score_model} is 
\begin{equation}\label{eq:eff_score_smooth_def}
\frac{1}{\sigma^2(x)}(y-m(\theta^\top x)) m^\prime(\theta^\top x)  H_\theta^\top \left\lbrace x -  \frac{\E(\sigma^{-2}(X)X|\theta^\top X=\theta^\top x)}{\E(\sigma^{-2}(X)|\theta^\top X=\theta^\top x)}\right\rbrace.
\end{equation}
Note that the efficient score depends on $p_{\epsilon|X}$ and $p_X$ only through $\sigma^2(\cdot)$. However if the errors happen to be homoscedastic (i.e., $\sigma^2(\cdot)\equiv\sigma^2$) then the \textit{efficient score} is $\ell_{\theta,m} (x,y)/\sigma^2$, where 
\begin{equation} \label{eq:EffScoreCLSE}
 \ell_{\theta,m} (x,y) := (y-m(\theta^\top x)) m^\prime(\theta^\top x) H_\theta^\top [  x -h_{\theta}(\theta^\top x) ].
 \end{equation}
 As $\sigma^2(\cdot)$ is unknown we restrict ourselves to efficient estimation under homoscedastic error; see Remark~\ref{rem:Efficency_general} for a brief discussion.  
% m_0^\prime(\theta^\top x) m^\prime(\theta^\top x)
 
 % $\psi_{\check\theta,\check m}$ is ``well-behaved''  and satisfies

\subsection{Efficiency of the CLSE} \label{sec:SemiCLSE}
The $\sqrt{n}$-consistency, asymptotic normality, and  efficiency (when the errors are homoscedastic) of $\check{\theta}$ will be established \textit{if} we could show that 
 \begin{equation}\label{eq:surrogate}
\sqrt{n}\,\p_n \ell_{\check\theta,\check m} = o_{p}(1)
\end{equation}
   and the class of functions $\ell_{\theta,m}$ indexed by  $(\theta, m)$ in a ``neighborhood'' of $(\theta_0,m_0)$ satisfies some technical conditions; see e.g., \citet[Chapter 6.5]{VdV02}.  As discussed in Section~\ref{sub:_semiparametric_eff_shape}, because $(\check{m},\check{\theta})$ minimizes $(m,\theta) \mapsto Q_n(m,\theta)$ over $\mathcal{M}_L\times\Theta$, the traditional way to prove~\eqref{eq:surrogate} is to use the fact that $\partial Q_n(\check{m}_{s,a}, \zeta_s(\theta,\eta))/\partial s|_{s = 0} =0$ for any $(a, \eta)$ such that $s\mapsto (\check{m}_{s,a}, \zeta_s(\theta,\eta))$ is a valid path (i.e., $a\in \overline{\mathrm{lin}}\,\D_{\check{m}}$). One then finds $(a, \eta) \in \D_{\check{m}}\times \R^{d-1}$ such that the derivative of $s\mapsto Q_n(\check{m}_{s,a}, \zeta_s(\theta,\eta))$ at $s=0$  is approximately $n^{-1}\sum_{i=1}^n \eta^\top \ell_{\check{\theta}, \check{m}}(Y_i, X_i)$; such an $(a, \eta)$ is called the  (approximate) \textit{least favorable submodel}; see~\citet[Section 9.2]{VdV02}. In Section~\ref{sec:eff_score}, we saw that if $m$ is strongly convex then $\overline{\mathrm{lin}}\,\D_m = L_2(\Lambda)$. However $\check{m}$ is piecewise affine and we can only show that $\overline{\mathrm{lin}}\,\D_{\check{m}} \subset L_2(\Lambda)$. Thus $s \mapsto \check{m}_{s,a}$ is valid path only if $ a\in \D_{\check{m}}$; see~\cite{VANC} for another example where $\overline{\mathrm{lin}}\,\D_{\check{m}} \neq L_2(\Lambda)$. In such cases it is hard to find the least favorable submodel as often the step to compute the least favorable model involves computing projection onto  $\overline{\mathrm{lin}}\,\D_{\check{m}}$; see e.g.,~\cite{Newey90}. Thus when $\overline{\mathrm{lin}}\,\D_{\check{m}}$ is not $L_2(\Lambda)$ (or a very simple subspace of $L_2(\Lambda)$), the standard linear path arguments fail to find the least favorable submodel.  To overcome this,~\cite{VANC} use a very complicated and non-linear path; see Section 6.2 of~\cite{VANC}; also see~\cite{Patra16}.

 Our proposed technique crucially relies on the observation that $s \mapsto \Pi_{\M_L}(\check{m}_{s,a})$ is a valid path for every $a\in L_2(\Lambda)$.  Thus if $s\mapsto \Pi_{\mathcal{M}_L}(\check{m}_{s,a})$ is differentiable, then establishing that $\check{\theta}$ is an approximate zero boils down to finding an $a\in L_2(\Lambda)$ such that 
 \begin{equation}\label{eq:approx_pi_score}
 \frac{\partial }{\partial s}Q_n(\Pi_{\mathcal{M}_L}(\check{m}_{s,a}), \zeta_s(\theta,\eta))\Big|_{s = 0}= n^{-1}\sum_{i=1}^n \eta^\top \ell_{\check{\theta}, \check{m}}(Y_i, X_i) +o_p(n^{-1/2}).
 \end{equation}
 for every $\eta \in \R^{d-1}. $  In Section~\ref{sec:proof_of_eq:app_score_equation}, we show $s \mapsto \Pi_{\M_L}(\check{m}_{s, a})$ is differentiable if $a \in \mathcal{X}_{\check{m}}$, where 
 \begin{equation}\label{eq:X_m_def}
 \mathcal{X}_{\check{m}} :=\big\{a\in L_2(\Lambda): a \text{ is a piecewise affine continuous function with kinks at }\{\check{t}_i\}_{i=1}^\mathfrak{p}\big\},
 \end{equation}
 and $\{\check{t}_i\}_{i=1}^\mathfrak{p}$ are the set of  kinks of $\check{m}$. For  a piecewise affine function, a kink is a point where the slope changes. Furthermore, in Theorem~\ref{thm:projection}, we find an $a\in \mathcal{X}_{\check{m}}$ that satisfies~\eqref{eq:approx_pi_score}. The advantage of the technique proposed here is that the construction of approximate least favorable submodel is  analytic and does not rely on the ability of the user to ``guess'' the least favorable submodel; see e.g.,~\cite[Section 9.2-9.3]{VdV02} and~\cite{VANC}. The above discussion and~\cite[Theorem 6.20]{VdV02} lead to our main result (Theorem~\ref{thm:Main_rate_CLSE}) of this section.  Recall $S_{\theta_0,m_0}$ and $\ell_{\theta,m}$  defined in~\eqref{eq:paraScore} and~\eqref{eq:EffScoreCLSE}, respectively.

\begin{thm} \label{thm:Main_rate_CLSE}
 Assume \ref{a0}--\ref{aa6} and \ref{bb1}--\ref{bb3} hold. Let $\theta_{0,1}>0$, $q\ge 5$, and $L\ge L_0$.  
% Define the function
% \be \label{eq:EffScoreCLSE}
% {\ell}_{\theta,m}(y,x) :=\big(y-m(\theta^\top x)\big) m^\prime(\theta^\top x)  H_\theta^\top \left\lbrace x - h_\theta(\theta^\top x)\right\rbrace.
% \ee
If $\gamma > 1/2+\beta/8$ and $V_{\theta_0,m_0}:=P_{\theta_0,m_0}({\ell}_{\theta_0,m_0} S^\top_{\theta_0,m_0})$ is a nonsingular matrix in $\R^{(d-1) \times (d-1)}$, then
\begin{equation}\label{eq:globalEffCLSE}
\sqrt{n} (\check{\theta}- \theta_0)\stackrel{d}{\rightarrow}   N(0,H_{\theta_0} V_{\theta_0,m_0}^{-1}  {I}_{\theta_0,m_0}(H_{\theta_0} V_{\theta_0,m_0}^{-1})^\top), 
\end{equation}
where  ${I}_{\theta_0,m_0} := P_{\theta_0,m_0} ({\ell}_{\theta_0,m_0}{\ell}^\top_{\theta_0,m_0})$. Further, if $\sigma^2(\cdot)\equiv\sigma^2$, then $V_{\theta_0,m_0}= {I}_{\theta_0,m_0}$ and
\begin{equation}\label{eq:localeffestimCLSE}
\sqrt{n} (\check{\theta}-\theta_0) \stackrel{d}{\rightarrow} N(0, \sigma^4 H_{\theta_0}{I}^{-1}_{\theta_0,m_0} H_{\theta_0}^\top).
\end{equation}
\end{thm}
\begin{remark}\label{rem:gammaBeta}
{If $m_0$ is twice continuously differentiable then $\gamma=1$. Hence, $\gamma > 1/2+\beta/8$ is equivalent to assuming $\beta \in [0,4)$. Note that $\beta>0$ allows for covariate distributions for which the density of $\theta_0^\top X$ can go to zero. In Theorem~\ref{thm:Main_rate_CLSE}, to keep notations in the proof simple, we assume that $q\ge 5$. However, by using Remark~\ref{rem:lowerQ}, this condition can be weakened to $q\ge4$. In Section~\ref{sec:Degenracy}, we show that the limiting variances in Theorem~\ref{thm:Main_rate_CLSE} are unique and do not depend on the particular choice of $\theta\mapsto H_\theta$.}
\end{remark}
\paragraph{Sketch of the proof.} The proof follows along the lines of Theorem 6.20 of~\cite{VdV02}. The main novelty in the proof is a new mechanism to verify that the estimator satisfies the score equation~\eqref{eq:surrogate}.  However to simplify the algebra involved,\footnote{All the proofs will go through with $\ell_{\theta,m}$ instead of $\psi_{\theta, m}$. However, usage of $\ell_{\theta, m}$ will require more remainder terms to be controlled and thus will lead to more tedious proofs.} we will work with 
\begin{equation} \label{eq:App_score}
 \psi_{\theta,m} (x,y) := (y-m(\theta^\top x))m^\prime(\theta^\top x) H_\theta^\top [  x -h_{\theta_0}(\theta^\top x)],
 \end{equation}
 a slight modification of $\ell_{\theta, m}$. The only difference between $\ell_{\theta, m}$ and $\psi_{\theta, m}$ is the last term ($h_{\theta}(\theta^\top X))$. 
 In Section~\ref{app:sketchCLSE} of the supplementary file we show that     \begin{equation}\label{eq:App_score_equation}
   \sqrt{n}\, \p_n \psi_{\check{\theta},\check{m}} =o_p(1),
   \end{equation}
implies 
   \begin{equation}\label{eq:RAL_1}
        \sqrt{n}V_{\theta_0,m_0}  H_{\theta_0}^\top(\check{\theta}- \theta_0)  ={}\g_n \psi_{\theta_0,m_0} + o_p(1+\sqrt{n} |\check{\theta}-\theta_0|).
       \end{equation}    The conclusion of the proof follows by  observing that 
$\psi_{\theta_0,m_0} = {\ell}_{\theta_0,m_0}$. 
We will now give a brief sketch of the proof of~\eqref{eq:App_score_equation}. Define for every $(m, \theta)$,  $\eta\in \R^{d-1}$, $a: D\to \R$, and $t\in \R$, 
\[\zeta_t(\theta,\eta):=\sqrt{1-t^2|\eta|^2}\, \theta + t H_\theta \eta \qquad\text{and}\qquad  \xi_t(u;a, {m}) := \Pi_{\M_L}({m} - t a)(u).
\]
Observe that $(\check{m}, \check{\theta})$ is the minimizer of $(m,\theta) \mapsto Q_n(m,\theta)$ and  $t \mapsto (\zeta_t(\check\theta,\eta), \xi_t(u;a, \check{m}))$ is a valid path in $\M_L\times \Theta$ through $( \check \theta, \check{m})$. Thus  $t=0$ is the minimizer of  $t\mapsto Q_n(\zeta_t(\check\theta,\eta), \xi_t(\cdot;a, \check{m}))$ for every  $\eta\in \R^{d-1}$ and $a: D\to \R$.
Hence if $t\mapsto Q_n(\zeta_t(\check\theta,\eta), \xi_t(\cdot;a, \check{m}))$ is differentiable then 
\begin{equation}\label{eq:derivative_Q_n_equals_zero} \frac{\partial}{\partial t} Q_n(\zeta_t(\check\theta,\eta), \xi_t(\cdot;a, \check{m}))\Big|_{t=0}=0.\end{equation}
Furthermore, if functions $a_1, a_2, \ldots, a_K$ (for some $K\ge1$) are such that $t\mapsto Q_n(\zeta_t(\check\theta,\eta), \xi_t(\cdot;a_j, \check{m}))$ is differentiable for all $1\le j\le K$, then
\[
\sum_{j=1}^{K} \alpha_j\frac{\partial}{\partial t} Q_n(\zeta_t(\check\theta,\eta), \xi_t(\cdot;a_j, \check{m}))\Big|_{t=0}=0,
\]
for any $\alpha_1,\ldots,\alpha_K\in\mathbb{R}$. Note that the proof of~\eqref{eq:App_score_equation} will be complete, if we can show that for every $\eta\in S^{d-2}$, there exist a $K\ge1$ and functions $a_j: D\to \R, 1\le j\le K$ such that $t\mapsto \Pi_{\M_L}(\check{m} - t a_j)(u)$ is differentiable and 
\begin{equation}\label{eq:psi_approx_phi}
\eta^\top   \p_n \psi_{\check{\theta}, \check{m}} = \sum_{j=1}^K \alpha_j\frac{\partial}{\partial t} Q_n(\zeta_t(\check\theta,\eta), \xi_t(\cdot;a_j, \check{m}))\Big|_{t=0} + o_p(n^{-1/2}).
\end{equation}
This means that it is enough to consider the approximation of $\eta^{\top}\mathbb{P}_n\psi_{\check{\theta},\check{m}}$ by the linear closure of $\{\partial Q_n(\zeta_t(\check\theta,\eta), \xi_t(\cdot;a, \check{m}))/\partial t|_{t=0}: t\mapsto Q_n(\zeta_t(\check\theta,\eta), \xi_t(\cdot;a, \check{m}))\mbox{ is differentiable at } t=0\}$. {Instead of fully characterizing the linear closure set, we find a large enough subset that suffices for our purpose using the following steps.}
\begin{enumerate}
\item We find a set of perturbations $a$ such that $t\mapsto \xi_t(\cdot;a, {m})$ is differentiable. Recall $\mathcal{X}_{\check{m}}$ defined in~\eqref{eq:X_m_def}. In Lemma~\ref{lem:projection_is_identity} (stated and proved in the supplementary file), we show that $\mathcal{X}_{\check{m}} \subseteq \{a:D\to\mathbb{R}\,|\,t\mapsto\xi_t(\cdot; a, \check{m}) \text{ is differentiable at } t=0\}.$
\item For every such $a\in\mathcal{X}_{\check{m}}$, in Lemma~\ref{lem:deriv}, we show that
\[
-\frac{1}{2}\frac{\partial}{\partial t}Q_n(\zeta_t(\check\theta,\eta), \xi_t(\cdot;a, \check{m}))\Big|_{t = 0} ~=~ \mathbb{P}_n\left[\big(y-\check{m}(\check\theta^\top x)\big)\Big\{ \eta^\top  \check{m}'(\check\theta ^\top x) H_{\check\theta}^\top x - a(\check\theta^\top x)\Big\}\right].
\]
\end{enumerate}
 % For every $a\in\mathcal{X}_{\check{m}}$, the following provides an expression for the derivative on the left hand side of~\eqref{eq:derivative_Q_n_equals_zero}.  
Thus to prove~\eqref{eq:psi_approx_phi}, it is enough to show that
\[
\inf_{a\in\overline{\mathrm{lin}}(\mathcal{X}_{\check{m}})}\left|\eta^{\top}\mathbb{P}_n\psi_{\check{\theta},\check{m}} - \mathbb{P}_n\left[(y-\check{m}(\check\theta^\top x))\{ \eta^\top  \check{m}'(\check\theta ^\top x) H_{\check\theta}^\top x - a(\check\theta^\top x)\}\right]\right| = o_p(n^{-1/2}),
\] 
where $\psi_{\theta, m}$ is defined in~\eqref{eq:App_score}.
In more general constraint spaces, one might need to use the generality of $\overline{\mathrm{lin}}(\mathcal{X}_{\check{m}})$ but in our case, it suffices to work with $\mathcal{X}_{\check{m}}$; see Theorem~\ref{thm:projection}. \qed

 % A full proof is provided in Section~\ref{sec:proof_of_eq:app_score_equation} of the supplementary file. 

% \end{proof}

\begin{remark}[Efficiency under heteroscedasticity]\label{rem:Efficency_general} It is important to note that~\eqref{eq:eff_score_smooth_def}, the efficient score, depends on $\sigma^2(\cdot)$.  Without additional assumptions, estimators of $\sigma^2(\cdot)$  will have poor finite sample  performance (especially if $d$ is large) which in turn will lead to  poor finite sample performance of the weighted LSE; see \citet[pages 93-95]{Tsiatis06}. 
\end{remark}

\begin{remark}[Efficiency under additional shape constraints]\label{rem:Efficency_additional}
As discussed in Remark~\ref{rem:Monotone}, it might be the case that the practitioner is interested in {imposing} additional shape constraints such as monotonicity, unimodality, or $k$-monotonicity (in addition to convexity). If $m_0$ satisfies these constraints {in a strict sense} (i.e., $m_0$ is strictly monotone {or} $k$-monotone) then the discussion in~Section~\ref{sec:eff_score} implies that {the efficient score (at the truth) is still~\eqref{eq:eff_score_smooth_def} even under the additional shape constraints}. This is true, because $\overline{\mathrm{lin}}\,\D_{m_0} = L_2(\Lambda)$ even under these additional shape constraints on link functions, as $m_0$ does not lie {on} the ``boundary'' of the parameter space. In fact, under these additional constraints, the proof of Theorem~\ref{thm:Main_rate_CLSE} can be used with minor modifications to show that CLSE of $\theta_0$ satisfies~\eqref{eq:globalEffCLSE}.
\end{remark}

To further illustrate the usefulness of our new approach we discuss the proof of semiparametric efficiency in the Cox proportional hazards model under  current status censoring~\cite{MR1394975,VdV02}. 

\begin{example}[Cox proportional hazards model with current status data]\label{ex:coxmodel}
Suppose that we observe a random sample of size $n$ from the distribution of $X = (C, \Delta, Z)$, where $\Delta = 1\{T \le C\}$, such that the survival time $T$ and the  observation time $C$ are independent given $Z \in \R^d$, and that $T$ follows a Cox proportional hazards model with parameter $\theta_0$ and cumulative hazard function $\Lambda_0$; e.g., see~\cite[Section 2]{MR1394975} for a more detailed discussion of this model.~\citet{MR1394975} shows that $\hat\Lambda$, the nonparametric maximum likelihood estimator (NPMLE) of $\Lambda_0$, is a right-continuous step function with possible discontinuities only at $C_1,\ldots, C_n$ (the observed  censoring/inspection times).~\citet{MR1394975} also proves that $\hat{\theta}$ (the NPMLE for $\theta_0$) is an efficient estimator for $\theta_0$. However just as in the single index model, the proof of efficiency is complicated due to the fact that $s\mapsto \hat{\Lambda}+ s h$ will not necessarily be a valid hazard function for every smooth $h(\cdot)$.\footnote{$\hat{\Lambda}+ s h$ is not guaranteed to be monotone as $\hat{\Lambda}$  is a nondecreasing piecewise constant function and not strictly increasing.}  To establish~\eqref{eq:surrogate} for the above model,~\citet[pages 563-564]{MR1394975} ``guesses'' an approximately least favorable path (also see~\cite[pages 439-441]{VdV02}). However, using the arguments above we can easily see that $ s \mapsto \Pi(\hat{\Lambda}+ s h)$ is differentiable if $h$ is a piecewise constant function with possible discontinuities only at the points of discontinuities of $\hat{\Lambda}$. Then using the property that $\|\hat{\Lambda}-\Lambda_0\| =o_p(n^{-1/3}),$ one can establish a result similar to~\eqref{eq:approx_pi_score}. A similar strategy can be used to establish efficiency in the current status regression model in~\citet{VANC}.
 \end{example}

 \subsection{Construction of confidence sets and validating the asymptotics}\label{sec:cfd_set}
% \clr{Theorem~\ref{thm:Main_rate_PLSE} proves asymptotic normality of the estimator $\hat{\theta}$ of $\theta_0$. In order to do inference with our estimator, one needs an estimate of the asymptotic variance. Under homoscedasticity, a simple plug-in approach works.
Theorem~\ref{thm:Main_rate_CLSE} shows that when the errors happen to be homoscedastic the CLSE of $\theta_0$ is $\sqrt{n}$-consistent and asymptotically normal with covariance matrix:
\begin{equation}\label{eq:finite_var}
\Sigma^0:= \sigma^4 H_{\theta_0} P_{\theta_0,m_0}[ {\ell}_{\theta_0,m_0}(Y,X) {\ell}^\top_{\theta_0,m_0}(Y,X)]^{-1} H_{\theta_0}^\top,
\end{equation}
where $\ell_{\theta_0, m_0}$ is defined in~\eqref{eq:EffScoreCLSE}. This result can be used to construct confidence sets for $\theta_0.$ However since $\Sigma^0$ is unknown, we propose using the following plug-in estimator of $\Sigma^0$:
\begin{equation}\label{eq:V_hat}
\check{\Sigma} := \check{\sigma}^4 H_{\check{\theta}} \big[\p_n\big( {\ell}_{\check{\theta},\check{m}}(Y,X) {\ell}^\top_{\check{\theta},\check{m}}(Y,X)\big)\big]^{-1} H_{\check{\theta}}^\top,
\end{equation}
where $\check{\sigma}^2: =  \sum_{i=1}^n [Y_i-\check{m}(\check{\theta}^\top X_i)]^2/n$. Note that Theorems~\ref{thm:ratestCLSE} and~\ref{thm:rate_derivCLSE} imply consistency of $\check\Sigma$.
% To be precise, it can be easily shown (using Therorems~\ref{thm:rate_m_theta_CLSE}--\ref{thm:rate_derivCLSE}) that one can consistently estimate $\Sigma^0$ (the asymptotic variance in \eqref{eq:localeffestimPLSE}) by the following plug-in estimator,
%  is a consistent estimator of $\sigma^2$.

 For example one can construct the following $1-2\alpha$ confidence interval for $\theta_{0,i}$:
\begin{equation}\label{eq:Conf_int}
\bigg[\max \left\{-1, \check{\theta}_i - \frac{z_{\alpha}}{\sqrt{n}} \left(\check{\Sigma}_{i,i}\right)^{1/2}\right\} ,\; \min\left\{1, \check{\theta}_i  + \frac{z_{\alpha}}{\sqrt{n}} \left(\check{\Sigma}_{i,i}\right)^{1/2} \right\}\bigg],
\end{equation}
where $z_{\alpha}$ denotes the upper $\alpha$th-quantile of the standard normal distribution. The truncation guarantees that confidence interval is a subset of the parameter set. 

{We now give an illustrative simulation example.  We generate $n$ i.i.d.~observations from the model: $Y=(\theta_0^\top X)^2 + N(0, .3^2),$ where $ X \sim \text{Uniform} [-1, 1]^3$ and $\theta_0 = (1,1,1)/\sqrt{3},$ for $n$ increasing from $50$ to $1000.$  For the above model, $\Sigma^0_{1,1}$ is $0.22$.\footnote{To compute the limiting variance in~\eqref{eq:finite_var}, we used a Monte Carlo approximation of $P_{\theta_0,m_0}[ {\ell}_{\theta_0,m_0}(Y,X) {\ell}^\top_{\theta_0,m_0}(Y,X)]$ with sample size $2\times 10^5$ and true $(m_0, \theta_0, P_X)$. The limiting covariance matrix $\Sigma^0 = 0.33 {I}_3 - 0.11 J_3$, where $I_3$ is the $3\times 3$ identity matrix and $J_3$ is the $3\times 3$ matrix of all ones.} In the left panel of Figure~\ref{fig:QQplot}, we present the Q-Q plot of $\sqrt{n}[\Sigma^0_{1,1}]^{-1/2}(\check{\theta}_{1}-\theta_{0,1})$  based on 800 replications; on the $x$-axis we have the quantiles of the standard normal distribution. The Q-Q plot validates the asymptotic normality and shows that the sample variance of the CLSE converges to the limiting variance found in Theorem~\ref{thm:Main_rate_CLSE}.  In the right panel of Figure~\ref{fig:QQplot}, we present empirical coverages (from $800$ replications) of $95\%$ confidence intervals based on  the CLSE constructed via~\eqref{eq:Conf_int}.
% We use this to validating the asymptotic properties derived in~Theorem~\ref{eq:globalEffCLSE}. We use~\eqref{eq:Conf_int} to construct a confidence interval and  for a simulated dataset as the sample size increases
}

\begin{figure}[!h]
  \begin{minipage}{0.50\linewidth}
    \centering
  \includegraphics[width=\textwidth]{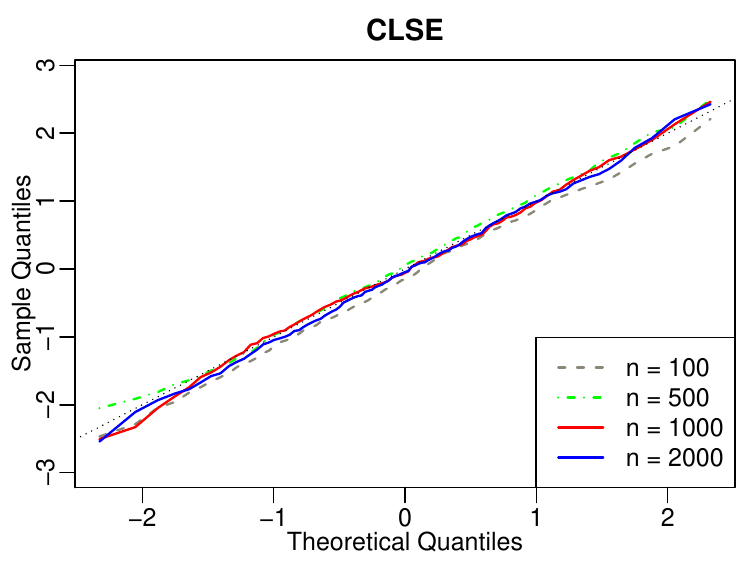}%
      \par

  \end{minipage}%
  \begin{minipage}{0.50\linewidth}
    \centering
\begin{tabular}{rcc}
  \toprule
  % \multicolumn{2}{c}{$n$}
  % \cmidrule(r){2-3} \cmidrule(rl){4}\\
   \multirow{2}{*}{$n$} &  \multicolumn{2}{c}{CLSE}\\
   \cmidrule(lr){2-3} 
 & Coverage & Avg Length\\
  \hline
  50  & 0.92 & 0.30 \\
  100  & 0.91 & 0.18 \\
  200  & 0.92 & 0.13 \\
  500  & 0.94 & 0.08 \\
  1000  & 0.93 & 0.06 \\
  % 2000  & 0.92 & 0.04 \\
   \hline
   
\end{tabular}

    \par
\end{minipage}
\caption{Summary of $\check{\theta}$ (over  800 replications) based on $n$ i.i.d.~observations from the model~\ref{sec:cfd_set}. Left panel:  Q-Q plots for $\sqrt{n}\left[\Sigma^0_{1,1}\right]^{-1/2}(\check{\theta}_{1}-\theta_{0,1})$  for  $n\in \{100, 500,1000, 2000\}$. The dotted black line corresponds to the $y=x$ line; right panel: estimated coverage probabilities and average lengths of nominal $95\%$ confidence intervals for the first coordinate of $\theta_0$. }%}
\label{fig:QQplot}
\end{figure}

% \todo[inline]{Remark about small $\hat{\lambda}$ choice doesn't matter and choice of $L$}

\section{Simulation study}\label{sec:Simul_Cvx}
{In Section~\ref{sec:compute} of the supplementary file, we develop an alternating minimization algorithm to compute the CLSE~\eqref{eq:CLSE}. 
  In this section we illustrate the finite sample performance of the CLSE using the implementation in the \texttt{R} package \if1\blind{\texttt{simest}.} \fi  \if0\blind{\texttt{***}}\fi}. We also compare its performance with other existing estimators, namely, the \texttt{EFM} estimator (the estimating function method; see \cite{cuietal11}), the \texttt{EDR} estimator (effective dimension reduction; see~\citet{Hristacheetal01}), and the estimator proposed in~\cite{Patra16} with the tuning parameter chosen by generalized cross-validation (\cite{Patra16}; we denote this estimator by \texttt{Smooth}).
% We use $(\hat{m}_{SS}, \hat{\theta}_{SS})$ to denote the estimator proposed in Section~\ref{smooth:sec:prelim} when the tuning parameter is chosen by generalized cross-validation; \cite{wahba90}. $(\hat{m}, \hat{\theta})$, $(\check{m}, \check{\theta})$, and $(m^\dagger, \theta^\dagger)$ denote the PLSE, the CLSE, and the LSE estimators, respectively.
We use \texttt{CvxLip} to denote the CLSE.

\subsection{Another convex constrained estimator}
Alongside these existing estimators, we also numerically study  another natural estimator under the convexity shape constraint --- the convex LSE --- denoted by \texttt{CvxLSE} below. This estimator is obtained by minimizing the sum of squared errors  subject to only the convexity constraint. Formally, the \texttt{CvxLSE} is
\begin{equation} \label{eq:ultimate_est}
(m^\dagger_n, \theta^\dagger_n):= \argmin_{(m, \theta)\in\, \mathcal{C}\times\Theta} Q_n(m, \theta).
\end{equation}
{The computation of \texttt{CvxLSE} is discussed in Remark~\ref{rem:CvxLSE} and is implemented in the \texttt{R} package \if1\blind{\texttt{simest}.} \fi  \if0\blind{\texttt{***}}\fi. }However, theoretical analysis of this estimator is difficult because of various reasons; see Section~\ref{sec:discussion_on_the_theoretical_analysis_of_the_texttt_CvxLSE} of the supplementary file for a brief discussion. %Although we present some evidence that our estimators $\check{\theta}_n$ and $\hat{\theta}_n$ are robust to the choice of their respective tuning parameters $L$ and $\lambda$, an application of this methods still requires the user to choose them. In this regard, $\theta^{\dagger}_n$ could be more bankable if its theory checks out since there is no tuning parameter to pick.
In our simulation studies we observe that the performance of  \texttt{CvxLSE} is very similar to that of \texttt{CvxLip}.
%, and \texttt{CvxLSE} to denote the convex LSE estimator proposed in~\eqref{eq:ultimate_est}

In what follows, we will use $(\tilde{m},\tilde{\theta})$ to denote a generic estimator that will help us describe the quantities in the plots; e.g., we use $ \|\tilde{m}\circ \tilde{\theta}-m_0\circ\theta_0\|_n =[\frac{1}{n}\sum_{i=1}^n (\tilde{m}( \tilde{\theta}^\top x_i)-m_0(\theta_0^\top x_i))^2]^{1/2}$ to denote the in-sample root mean squared estimation error of $(\tilde{m},\tilde{\theta})$, for all the estimators considered. From the simulation study it is easy to conclude that the proposed estimators have superior finite sample performance in most sampling scenarios considered.

\begin{figure}[!ht]
% \captionsetup[subfigure]{labelformat=empty}
\centering
\includegraphics[width=.75\textwidth]{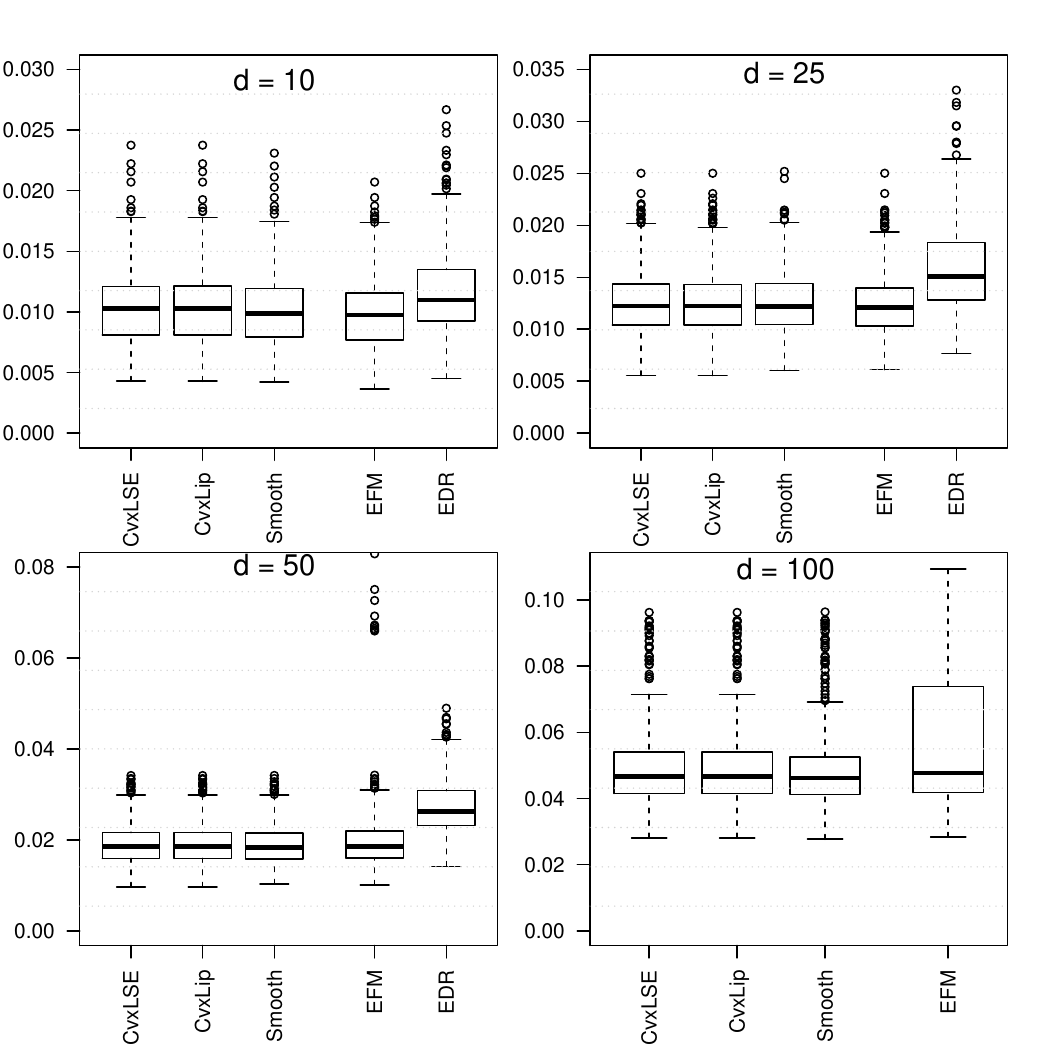}\\ [-2ex]%%
\caption{Boxplots of $\sum_{i=1}^d |\tilde{\theta}_i-\theta_{0,i}|/d$ (over 500 replications) based on $100$ observations from the simulation setting in Section~\ref{sub:increasing_dimension} for dimensions $10,$ $25,$ $50,$ and $100,$ shown in the top-left, the top-right, the bottom-left, and the bottom-right panels, respectively. The bottom-right panel doesn't include EDR as the R-package \texttt{EDR} does not allow for $d=100.$}
\label{fig:Ex3Cuietal_homo}
%% label for entire figure
\end{figure}

\subsection{Increasing dimension} % (fold)
\label{sub:increasing_dimension}
To illustrate the behavior/performance of the estimators  as $d$ grows, we consider the following single index model $Y=  (\theta_0^\top X)^2 + t_6, \text{ where } \theta_0= (2, 1, \mathbf{0}_{d-2})^\top/\sqrt{5} \text{ and } X\in \R^d \sim\text{Uniform}[-1,5]^d,$
where $t_6$ denotes the Student's $t$-distribution with $6$ degrees of freedom.
In each replication we observe $n= 100$ i.i.d.~observations from the model. It is easy to see that the performance of all the estimators worsen as the dimension increases from $10$ to $100$ and \texttt{EDR} has the worst overall performance; see Figure~\ref{fig:Ex3Cuietal_homo}. However when $d=100$, the convex constrained estimators have significantly better performance. This simulation scenario is similar  to the one considered in Example 3 of Section 3.2 in~\cite{cuietal11}.
% subsection increasing_dimension (end)
% the following works.
% \vspace{-.1in}

% subsection verifying_the_asymptotics_ (end)
% subsection subsection_name (end)

\subsection{Choice of \texorpdfstring{$L$}{Lg}} % (fold)
\label{sub:robustness_of_choice_of_}
In this subsection, we consider a simple simulation experiment to demonstrate that the finite sample performance of the CLSE is robust to  the choice of tuning parameter. We generate an i.i.d.~sample (of size $n= 500$) from the following model:
\begin{equation}\label{eq:robust}
Y=(\theta_0^\top X)^2 + N(0, .1^2), \quad \text{where}\; X \sim \text{Uniform} [-1, 1]^4\; \text{and} \; \theta_0 = (1,1,1,1)^\top/2.
\end{equation}
Observe that, we have   $-2 \le \theta^\top X \le 2$ and  $L_0:= \sup_{t \in [-2,2]} m_0'(t)= 4$ as $m_0(t)= t^2.$ To understand the effect of $L$ on the performance of the CLSE,  we  show the box plot of $\sum_{i=1}^4|\check{\theta}_i-\theta_{0,i}|/4$  as $L$ varies from $3\, (< L_0)$ to $10$ in Figure~\ref{fig:robust}. Figure~\ref{fig:robust} also includes the \texttt{CvxLSE} which corresponds to $L=\infty$.  The plot clearly show that the performance of \texttt{CvxLip}  is not significantly affected by the particular choice of the tuning parameter. The observed robustness in the behavior of the estimators can be attributed to the stability endowed  by the convexity constraint.
\begin{figure}[!ht]
% \captionsetup[subfigure]{labelformat=empty}
\centering
\includegraphics[width=.5\textwidth]{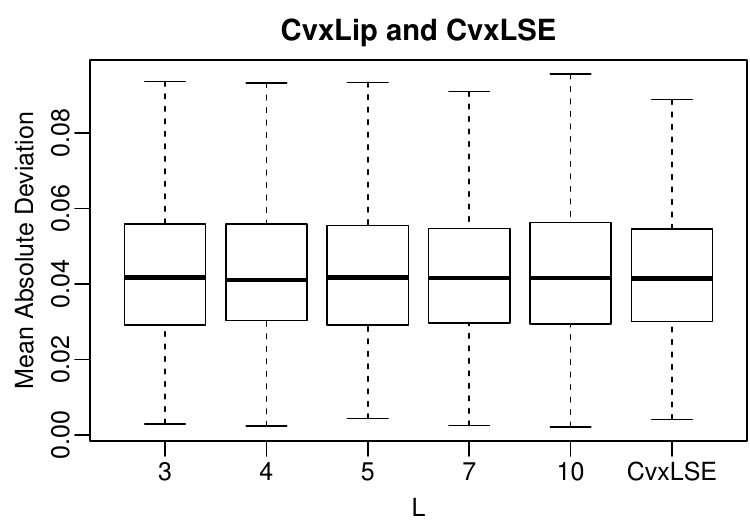}\\ %%
% \caption[Boxplots of estimates to study their robustness]{Box plots of $\frac{1}{4}\sum_{i=1}^4|\tilde{\theta}_i-\theta_{0,i}|$ (over $1000$ replications) for the model~\eqref{eq:robust} ($d=4$ and $n=500$) as the tuning parameter varies. Left panel: \texttt{CvxPen} when $\lambda_n= \exp(-T)\times n^{1/5}$ for $T={\{0,0.7, 1, 2, 5,7\}}$; right panel: \texttt{CvxLip} for $L= \{3, 4, 5,7, 10\}$ and \texttt{CvxLSE}. }
\caption[Boxplots of estimates to study their robustness]{Box plots of $\frac{1}{4}\sum_{i=1}^4|\tilde{\theta}_i-\theta_{0,i}|$ (over $1000$ replications) for the model~\eqref{eq:robust} ($d=4$ and $n=500$) \texttt{CvxLip} for $L= \{3, 4, 5,7, 10\}$ and \texttt{CvxLSE} (i.e., $L=\infty$). }
% \clr{Both plots have $1000$ replications}.}
\label{fig:robust}
%% label for entire figure
\end{figure}
\section{Real data analysis} % (fold)
\label{sec:real_data_analysis}

In this following we analyze the two real datasets discussed in Examples~\ref{ex:boston} and~\ref{ex:car}.
% and apply the developed methodology for prediction and estimation.

\subsection{Boston housing data}\label{sec:boston}

% \citet{harrison1978hedonic} studied the effect of different covariates on real estate price in the greater Boston area.  The response variable $Y$ was the log-median value of homes in each of the $506$ census tracts in the Boston standard metropolitan area. A single index model is  appropriate for this dataset; see e.g.,~\cite{gu2015oracally,MR2529970,MR2589322,MR2787613}. The above papers considered the following covariates in their analysis: average number of rooms per dwelling, full-value property-tax rate per $10000$ USD, pupil-teacher ratio by town school district, and proportion of population that is of ``lower (economic) status'' in percentage points. In the left panel of Figure~\ref{fig:real_data_plot_prelim}, we provide a scatter plot of  $\{(Y_i, \hat{\theta}^\top X_i)\}_{i=1}^{506}$, where $\hat{\theta}$ is an estimate of $\theta_0$ {\clr obtained in~\cite{Patra16}}. We also plot estimates of $m_0$ {\clr obtained from~\cite{Patra16} and~\cite{MR2529970}}. The plot suggests a convex and nondecreasing relationship between the log-median home prices and the index, but the fitted link functions satisfy these shape constraints only approximately. 

We briefly recall the discussion in Example~\ref{ex:boston}. The Boston housing dataset was collected by \cite{harrison1978hedonic} to study the effect of different covariates on the real estate price in the greater Boston area. The dependent variable $Y$ is the  log-median value of homes in each of the $506$ census tracts in the Boston standard metropolitan area. \citet{harrison1978hedonic} observed $13$ covariates and fit a linear model after taking $\log$ transformation for $3$ covariates and power transformations for three other covariates; also see \cite{MR2589322} for a discussion of this dataset.

\citet{breiman1985estimating} did further analysis to deal with multi-collinearity of the covariates and selected four variables using a penalized stepwise method. The chosen covariates were:  average number of rooms per dwelling (RM), full-value property-tax rate per $10,000$ USD (TAX), pupil-teacher ratio by town school district (PT), and proportion of population that is of ``lower (economic) status'' in percentage points (LS). Following \cite{MR2529970} and \cite{MR2787613}, we take logarithms of LS and TAX to reduce sparse areas in the dataset.  Furthermore, we have scaled and centered each of the covariates to have mean $0$ and variance $1.$  \citet{MR2529970} fit a nonparametric additive regression model to the selected variables and obtained an $R^2$ (the coefficient of determination) of $0.64$. \citet{MR2589322} fit a single index model to this data using the set of covariates suggested in \cite{MR1624402}. In~\cite{gu2015oracally}, the authors create $95\%$ uniform confidence band for the link function and reject the null hypothesis that the link function is linear. Both in~\cite{gu2015oracally} and~\cite{MR2589322}, the fitted link function is approximately nondecreasing and  convex; see Figure~2 of \cite{MR2589322} and Figure~5 of~\cite{gu2015oracally}.  This motivates us to fit a \textit{nondecreasing} and convex single index model to the Boston housing dataset. 
% Moreover, the shape constraints add interpretability to the estimators of both $\theta_0$ and $m_0$.}
 In particular, we consider the following estimator:
\begin{equation}\label{eq:CLSE_nondec}
(\hat{m}_{L},\hat{\theta}_{L}) \coloneqq \argmin_{\substack{\theta \in  \Theta \\  m \in \M_L \cap \,\mathcal{N} } }\; \sum_{i=1}^n (Y_i - m(\theta^\top X_i))^2,
\end{equation}
where $\mathcal{N}$ is the set of real-valued nondecreasing functions on $D$.
Following the discussions in~Remarks~\ref{rem:Monotone} and~\ref{rem:Efficency_additional}, we observe that the results in this paper also hold for $(\hat{m}_L,\hat{\theta}_L)$. The computation of the CLSE under the additional monotonicity constraint is discussed in~Remark~\ref{rem:Add_mono} and implemented in the accompanying R package. 

We summarize our results in Table~\ref{tab:real_dat}. We call $(\hat{m}_L,\hat{\theta}_L)$, the \texttt{MonotoneCLSE}. In Figure~\ref{fig:real_data_plot}, we plot the scatter plot of $\{(\hat{\theta}^\top_L X_i, Y_i)\}_{i=1}^{506}$  overlaid with the plot of $\hat{m}_L(\cdot)$ and the regression splines based estimator of~\cite{MR2529970}.  For \texttt{MonotoneCLSE} and \texttt{CvxLip}, we chose $L= 30$ (an arbitrary but large number). We also observe that the $R^2$ for the monotonicity and convexity constrained (\texttt{MonotoneCLSE}) and just convexity constrained single index models (\texttt{CvxLip} and \texttt{CvxLSE}), when using all the available covariates, is approximately $0.80$. {To further understand the predictive properties of the estimators under different smoothness and shape constraints, in Table~\ref{tab:real_dat} we report the $5$-fold cross-validation error averaged over 100 random partitions.  The large cross-validation error for the \texttt{CvxLSE} is due to over-fitting of $m_n^\dagger$ at the boundary of its support; see~Figure~\ref{fig:SimplModel} for an illustration of this boundary effect.}

 %Inclusion of all the extra variables leads to only a minor increase in $R^2$ at the cost of interpretability; \cite{MR2589322} also reached to a similar conclusion.}
% \todo[inline]{({\color{red} Cross-validated prediction error??}).}
% \begin{itemize}
%   \item  \cite{breiman1985estimating} did further analysis to reduce the number of ``impottant variable to 4'' see Wang, Yang 09
%   \item Wang and and Yang take log for  \texttt{LS} \texttt{TAX} so  as not to have gaps. They fit a additive regression model and get a correlation of .80. We get a correlation of .87
%   \item Need to understand \cite{MR2589322}
%   \item see \verb|Boston_MaY_final_analysis.R| file in Smooth single index model folder.
% \end{itemize}

\subsection{Car mileage data}\label{sec:car}
First, we briefly recall the discussion in Example~\ref{ex:car}. We consider the car mileage dataset of~\citet{cars_1983} for a second application for the convex single index model. We model the mileage ($Y$) of $392$ cars using the covariates ($X$):  displacement (Ds),   weight (W),  acceleration (A),   and horsepower (H).  \citet{MR2957294} fit a partial linear model to this this dataset, while \cite{Patra16} fit a single index model (without any shape constraint). The ``law of diminishing returns'' suggests $m_0$ should be convex and nonincreasing. However, the estimators based only on smoothness assumptions satisfy these shape constraints only approximately.  In the right panel of Figure~\ref{fig:real_data_plot}, we fit a convex and nonincreasing single index model. 
 % Our analysis reinforces the notion that cars with higher acceleration and/or weight generally have lower mileage.

We have scaled and centered each of covariates to have mean $0$ and variance $1$ for our analysis, just as in Section~\ref{sec:boston}. We performed a test of significance for $\theta_0$ using the plug-in variance estimate in Section~\ref{sec:cfd_set}. The covariates  A, Ds, and  H were found to be significant and each of them had $p$-value less than $10^{-5}$. In the right panel of  Figure~\ref{fig:real_data_plot}, we have the scatter plot of $ \{(\hat{\theta}_L^\top X_i, Y_i)\}_{i=1}^{392}$ overlaid with the plot of $\hat{m}_L(\cdot)$ and regression splines based estimator obtained in~\cite{MR2529970}; here $\hat{\theta}_L$ is defined as in~\eqref{eq:CLSE_nondec} but $\mathcal{N}$ now denotes the class of real-valued \textit{nonincreasing} functions on $D$. Table~\ref{tab:real_dat} lists different estimators for $\theta_0$ and their respective $R^2$ and cross-validation errors.
  \begin{table}[h]
\caption[Estimates of $\theta_0$  and generalized $R^2$ for the datasets in Sections~\ref{sec:boston} and \ref{sec:car}. {\texttt{EFM} and \texttt{EDR} do not provide a function estimator and hence we do not show an $R^2$ value.}]{Estimates of $\theta_0$  and generalized $R^2$ for the datasets in Sections~\ref{sec:boston} and \ref{sec:car}. {\texttt{EFM} and \texttt{EDR} do not provide a function estimator and hence we do not show an $R^2$ value. { CV-error denotes out of  5-fold cross validation averaged over 100 random partitions.} }}\label{tab:real_dat}
\centering
\begin{adjustbox}{max width=\textwidth}

\begin{tabular}{l*{13}{c}}
 \toprule
\multirow{2}{*}{Method} &\multicolumn{6}{c}{Boston Data}& &\multicolumn{6}{c}{Car mileage data}  \bigstrut \\
\cmidrule(rl){2-7} \cmidrule(rl){9-13}
  &RM& $\log(\text{TAX})$ & PT &$\log(\text{LS})$& $R^2$&  CV-error & & Ds & W  &A&H& $R^2$& CV-error \\

 \midrule
\texttt{LM}\footnote{\texttt{LM} denotes the linear regression model.}  &  2.34 & $-0.37$ & $-1.55$ & $-5.11$ &  0.73 & 20.75& & $-0.63$ & $-4.49$ & $-0.06$ & $-1.68$ &  0.71 & 18.61\\
\texttt{Smooth} &  0.44 & $-0.18$ & $-0.27$ & $-0.83$ &  0.77 & 17.80& &0.42 &  0.18 &  0.11 &  0.88 &  0.76 &15.29\\
%\texttt{CvxPen} &  0.48 & $-0.19$ & $-0.25$ & $-0.82$ &  0.77 &  &0.45 &  0.15 &  0.13 &  0.87 &  0.76 \\
\texttt{MonotoneCLSE} &  0.49 & $-0.21$ & $-0.25$ & $-0.81$ &  0.80 & 17.93& &0.44 &  0.17 &  0.13 &  0.87 &  0.76 & 15.34\\
\texttt{CvxLip} &  0.48 & $-0.23$ & $-0.26$ & $-0.80$ &  0.80 & 17.93&  &0.44 &  0.18 &  0.12 &  0.87 &  0.76& 15.22 \\
% [1]  0.4796955 -0.2341851 -0.2631070 -0.8036319
\texttt{CvxLSE} &  0.43 & $-0.20$ & $-0.28$ & $-0.84$ &  0.80 & 21.44&  &0.39 &  0.14 &  0.12 &  0.90 &  0.77 &16.38\\
\texttt{EFM} &  0.48 & $-0.19$ & $-0.21$ & $-0.83$ &  --- & --- & &0.44 &  0.18 &  0.13 &  0.87 &--- &  --- \\
\texttt{EDR} &  0.44 & $-0.14$ & $-0.18$ & $-0.87$ &  --- &--- &  &0.33 &  0.11 &  0.15 &  0.93 & --- & --- \\

\bottomrule

\end{tabular}
\end{adjustbox}

\end{table}

   % we display the estimates of $\theta_0$  based on the methods considered in the paper. The MAVE, the EFM estimator, and the PLSE give similar estimates while the EDR gives a different estimate of the index parameter.
\begin{figure}[h!]
\centering
\includegraphics[width=1\textwidth]{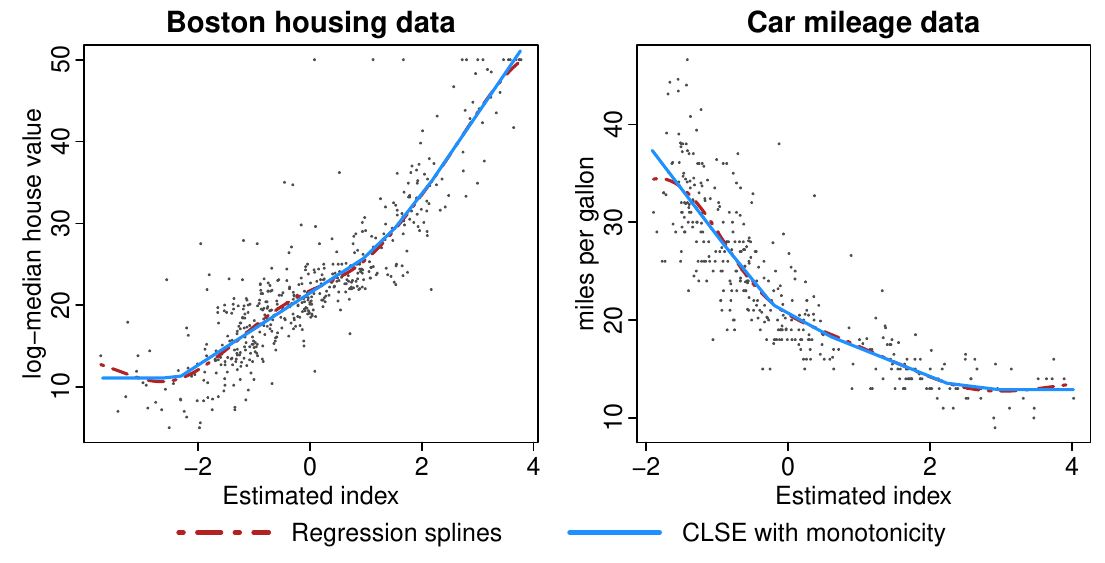}
%% scale=.8
  \caption[]{Scatter plots of $\{(X_i^\top\check{\theta}, Y_i)\}_{i=1}^n$ overlaid with the plots of function estimates proposed in~\cite{MR2529970} (red, dot-dashed) and monotonicity constrained CLSE proposed in this paper (blue, solid) for the two real datasets considered. Left panel: Boston housing data (Section~\ref{sec:boston}), nondecreasing CLSE; right panel: the car mileage data (Section~\ref{sec:car}), nonincreasing CLSE.}
  \label{fig:real_data_plot}
\end{figure}

\section{Discussion} \label{sec:discussion}
  In this paper we have proposed and studied a Lipschitz constrained LSE in the convex single index model. Our estimator  of the regression function is minimax rate optimal (Proposition~\ref{lem:MinimaxLowerbound}) and  the estimator of the index parameter is semiparametrically efficient when the errors happen to be homoscedastic (Theorem~\ref{thm:Main_rate_CLSE}). This work represents the first in the literature of semiparametric efficiency of the LSE when the nonparametric function estimator is non-smooth and parameters are bundled. Our proof of semiparametric efficiency is geometric and provides a general framework that can be used to prove efficiency of estimators in a wide variety of semiparametric models even when the estimators do not satisfy the efficient score equation directly; see sketch of proof of Theorem~\ref{thm:Main_rate_CLSE} and Example~\ref{ex:coxmodel} in Section~\ref{sec:SemiCLSE}. 

  % {\clr 
  % I don't think we should stress too much about why CvxLse is hard. I guess it would be enough to say ``An interesting future direction would be the study of CvxLSE introduced in the simulation section. This requires a control of the function estimator near the boundary of its domain and generalizations of our theorems 3.2, 3.5 and 3.7.'' Something from below is
  % ``if one can prove results similar to Theorems~\ref{thm:rate_m_theta_CLSE}--\ref{thm:rate_derivCLSE} for the convex LSE, then the techniques used in Section~\ref{sec:SemiInf} can be readily applied to prove asymptotic normality of $\theta^\dagger$. These challenges make the study of the convex LSE a very interesting problem for future research''}
  % % {\clg The new result allows the errors to have only finite moments and to depend on the covariates. }

{Theorem~\ref{thm:rate_m_theta_CLSE} proves the worst case rate of convergence for the CLSE.  It is well-known in convex regression that if the true regression function is piecewise linear, then the LSE converges at a much faster (near parametric) rate~\cite{MR3881209}. This behavior is called the \emph{adaptation} property of the LSE.
It is natural to wonder if such a property also holds for $\check{m}\circ\check\theta$. In Section~\ref{sub:investigation_of_the_adaptation_of_the_clse} of the supplementary file, we investigate the behavior of  $\check{m}\circ\check\theta$ and $\check{\theta}$ (as sample size increases) when $m_0$ is piecewise linear. The simulation suggests that $\check{m}\circ\check{\theta}$ converges at a near parametric rate when $m_0$ is piecewise linear. However a formal proof of this is beyond the scope of this paper as it requires different techniques. Furthermore, the asymptotic behavior of $\check \theta$ in this setting is an open problem. 
% the effect of the adaptive behavior of the CLSE on the asymptotic limit of $\check{\theta}$ is an open problem.

   % but it will require substantially different techniques. Furthermore, the effect of the adaptation behavior of the CLSE on the asymptotic limit of $\check{\theta}$ is an open problem, even in  the much simpler partial linear regression model. In Section~\ref{sub:investigation_of_the_adaptation_of_the_clse} of the supplementary file, we investigate the behavior of  $\check{m}\circ\check\theta$ and $\check{\theta}$ as sample size increases. However, given the focus on semiparametric efficiency in this paper, we refrain from further studying the adaptation behavior here. 
   }

% \begin{supplement}
%  \sname{Supplementary material for}
%  \stitle{Efficient Estimation in Convex Single Index Models}
%  \slink[doi]{10.1214/00-AOASXXXXSUPP}
%  \sdatatype{.pdf}
%  \sdescription{Due to space constraints, some finite sample examples and  the proofs of the results in the paper are relegated to the supplementary file.}
% \end{supplement}
% \noeqref{eq:deriv_theta_0}
% % \noeqref{eq:deriv_theta_hat}
% % \noeqref{eq:mprime_sup_theta}\noeqref{eq:BdryCond}
% \noeqref{eq:sup_rate}
% \noeqref{eq:true_mimina}
% \noeqref{eq:True_inf}
% \noeqref{eq:RAL_1}
% \noeqref{eq:rn_Major}
% \noeqref{eq:RAL_1}

\bibliographystyle{chicago}
{\bibliography{SigNoise}}
\newpage
\setcounter{section}{0}
\setcounter{equation}{0}
\setcounter{figure}{0}
\renewcommand{\thesection}{S.\arabic{section}}
\renewcommand{\theequation}{E.\arabic{equation}}
\renewcommand{\thefigure}{S.\arabic{figure}}
\cftsetindents{section}{1em}{2.5em}
\cftsetindents{subsection}{1.5em}{3em}

  \begin{center}
  \Large {\bf Supplement to ``Semiparametric Efficiency in Convexity Constrained Single Index Model''}
  \end{center}
       
\begin{abstract}
Section~\ref{sec:compute} proposes an alternating minimization algorithm to compute the estimators proposed in the paper. Section~\ref{app:sketchCLSE}  provides some insights into the proof of Theorem~\ref{thm:Main_rate_CLSE}. Section~\ref{sec:Degenracy} shows that the asymptotic variance in~Theorem~\ref{thm:Main_rate_CLSE} is the Moore-Penrose inverse of the efficient information matrix.
Section~\ref{app:add_simul}  provides further simulation studies. Section~\ref{sec:proof:Estims}  provides additional discussion on our identifiability assumptions. Section~\ref{sec:minimax_lower_bound} finds the minimax lower bound for the model~\eqref{eq:simsl} under~\ref{aa1_new}--\ref{aa2} and shows that the CLSE is minimax rate optimal when $q\ge 5.$ Section~\ref{sec:maximal_inequalities_for_heavy_tailed} provides new maximal inequalities that allow for unbounded 
errors. 
% Section~\ref{app:proof:existanceCLSE}  proves Proposition~\ref{thm:existanceCLSE}. 
 These maximal inequalities are used in Section~\ref{app:AsymCLSE_Proof} to allow for heavy-tailed and heteroscedastic errors. These results are also of independent interest. Sections~\ref{app:proof:existanceCLSE}--\ref{sec:proof_semi} contain the proofs omitted from the main text.  Section~\ref{app:AsymCLSE_Proof}  proves the results in Section~\ref{sec:CLSE}.  Section~\ref{sec:proof_of_eq:app_score_equation} completes the proof of the approximate zero property in~\eqref{eq:App_score_equation}. Sections~\ref{app:stepsCLSE} and \ref{sec:proof_semi} complete the proofs of the steps in Section~\ref{app:sketchCLSE}. Section~\ref{rem:bin} provides a comment regarding the computation of the function estimate in the CLSE when there are ties.
 \end{abstract}
{
  \hypersetup{linkcolor=blue}
  \input{JASA_supp4.toc}
}
\newpage
\section{Alternating minimization algorithm}\label{sec:compute}
% \todo{Remove PLSE connections}
In this section we describe an algorithm for computing the estimator defined in \eqref{eq:CLSE}.  As mentioned in Remark~\ref{rem:compute}, the minimization of the desired loss function for a fixed $\theta$ is a convex optimization problem; see Section~\ref{sec:CLSE_comp} below for more details. With the above observation in mind, we propose the following general alternating minimization algorithm to compute the proposed estimator. \if1\blind  {The algorithms discussed here are implemented in our R package \texttt{simest}~\cite{simest}.}\fi \if0\blind  {The algorithms discussed here are implemented in our R package \texttt{***}.}\fi 

We first introduce some notation. Let $(m, \theta)\mapsto \C(m, \theta)$ denote a nonnegative criterion function, e.g., $\C(m, \theta) = Q_n(m, \theta)$. And suppose, we are interested in finding the minimizer of $\C(m, \theta)$ over $(m, \theta) \in \A\times \Theta$, e.g., in our case $\A$ is $\M_L$. For every $\theta \in \Theta$, let us define
 \begin{equation} \label{eq:profile}
  m_{\theta, \A}:= \argmin_{m \in \A} \C(m,\theta).
 \end{equation}
 Here, we have  assumed that for every $\theta \in \Theta$, $m \mapsto \C(m, \theta)$ has a unique minimizer in $\A$ and $m_{\theta, \A}$ exists. The general alternating scheme is described in Algorithm~\ref{algo:b}.
 \begin{algorithm}
%\DontPrintSemicolon % Some LaTeX compilers require you to use \dontprintsemicolon    instead
\KwIn{Initialize $\theta$ at $\theta^{(0)}$.}
\KwOut{$(m^*, \theta^*) := \argmin_{(m,\theta)\in\A\times\Theta}\C(m, \theta)$.}
At iteration $k \ge 0$, compute $m^{(k)} := m_{\theta^{(k)},\A} =  \argmin_{m \in \A} \C(m,\theta^{(k)})$.\\
Find a point $\theta^{(k+1)}\in\Theta$ such that
\[
\C(m^{(k)}, \theta^{(k+1)}) \le \C(m^{(k)}, \theta^{(k)}).
\]
In particular, one can take $\theta^{(k+1)}$ as a minimizer of $\theta\mapsto\C(m^{(k)}, \theta)$.\\
Repeat steps 1 and 2 until convergence.
\caption{Alternating minimization algorithm}
\label{algo:b}
\end{algorithm}
%\begin{enumerate}
%\item Start with a initial estimate of $\theta$, say, $\theta^{(0)}$.
%\item At iteration $k$, compute $m^{(k)}:= m_{\theta^{(k)}, \A}.$
%\item Find a point $\theta^{(k+1)} \in \Theta$ such that
%\begin{equation}\label{eq:descent}
% \C(m^{(k)},\theta^{(k+1)}) \le  \C(m^{(k)},\theta^{(k)}).
%\end{equation}
%In particular, one can take $\theta^{(k+1)}$ as a minimizer of $\C(m^{(k)},\theta)$.
%\item Repeat steps 2 and 3 until convergence.
%\end{enumerate}

Note that, our assumptions on $\C$ does not imply that  $\theta \mapsto \C(m_{\theta, \A}, \theta)$  is a convex function. In fact in our case the ``profiled" criterion function $\theta \mapsto \C(m_{\theta, \A}, \theta)$ is not convex. Thus the algorithm discussed above is not guaranteed to converge to a global minimizer. However, the algorithm guarantees that the criterion value is nonincreasing over iterations, i.e., $\C(m^{(k+1)},\theta^{(k+1)})\le \C(m^{(k)},\theta^{(k)})$ for all $k\ge 0.$ { To lessen the chance of getting stuck at a local minima,  we use multiple random starts for $\theta^{(0)}$ in Algorithm~\ref{algo:b}. Further, following the idea of~\cite{dumbgen2013stochastic}, we use other existing $\sqrt{n}$-consistent estimators of $\theta_0$ as warm starts; see Section~\ref{sec:Simul_Cvx} for examples of such estimators.} In the following section,  we discuss an algorithm to compute $m_{\theta,\mathfrak{A}}$, when $\C(m, \theta)=Q_n(m, \theta)$ and $\mathfrak{A} =\M_L$.% while in Section~\ref{sec:PLSE_comp} we discuss the computation of $m_{\theta, \mathcal{R}}$ when $\C(m, \theta)=\mathcal{L}_n(m, \theta; \lambda).$
\subsection{Strategy for estimating the link function} \label{sec:CLSE_comp}
In this subsection, we describe an algorithm to compute  $m_{\theta, \M_L}$ as defined in~\eqref{eq:profile}. We use the following notation. Fix an arbitrary $\theta\in\Theta$. Let $(t_1, t_2, \cdots, t_n)$ represent the vector $(\theta ^\top x_1, \cdots, \theta ^\top x_n)$ with sorted entries so that $t_1 \le t_2 \le \cdots \le t_n$. Without loss of generality, let $y := (y_1, y_2, \ldots, y_n)$ represent the vector of responses corresponding to the sorted $t_i$.

% \subsubsection{Lipschitz constrained least squares (CLSE)} \label{sec:CLSE_comp}
When $\C(m, \theta)=Q_n(m, \theta)$, we consider the problem of minimizing $\sum_{i=1}^n \{y_i - m(t_i)\}^2$ over $m\in\M_L$.  Note that the loss depends  only on the values of the function at the $t_i$'s and the minimizer is only unique at the data points.  Hence, in the following we identify $m:= (m(t_1), \ldots, m(t_n)):= (m_1, \ldots, m_n)$ and interpolate/extrapolate the function linearly between and outside the data points; see footnote~\ref{foo:extrapolate}. Consider the general problem of minimizing
\[
(y - m)Q(y - m) = |Q^{1/2}(y - m)|^2,
\]
for some positive definite matrix $Q$. In most cases $Q$ is the $n\times n$ identity matrix; see Section~\ref{rem:bin} of the supplementary file for other possible scenarios. Here $Q^{1/2}$ denotes the square root of the matrix $Q$ which can be obtained by Cholesky factorization.

The Lipschitz constraint along with convexity (i.e., $m\in\mathcal{M}_L$) reduces to imposing the following linear constraints:
\begin{equation}
-L \le \frac{m_2 - m_1}{t_2 - t_1}  \le \frac{m_3 - m_2}{t_3 - t_2} \le \cdots \le \frac{m_{n} - m_{n-1}}{t_n - t_{n-1}} \le L.\footnote{In Section~\ref{rem:bin} of the supplementary file, we discuss a solution for scenarios with ties.}\label{lipcons}
\end{equation}
In particular, the minimization problem at hand can be represented as
\begin{equation}
\mbox{minimize }|Q^{1/2}(m - y)|^2 \qquad \mbox{ subject to } \qquad Am \ge b,\label{p1}
\end{equation}
for $A$ and $b$ written so as to represent~\eqref{lipcons}. It is clear that the entries of  $A$ involve $1/(t_{i+1} - t_i), 1\le i\le n-1$. If the minimum difference is close to zero, then the minimization problem~\eqref{p1} is ill-conditioned and can lead to numerical inaccuracies. For this reason, in the implementation we have added a pre-binning step in our implementation; see Section~\ref{rem:bin} of the supplementary for details.
\begin{remark}[Additional monotonicity assumption]\label{rem:Add_mono}
{Note that if $m$ is additionally monotonically nondecreasing, then
\[
m_1 \le m_2 \le \cdots \le m_n \quad\Leftrightarrow\quad A'm \ge \textbf{0}_{n-1},
\]
where $\textbf{0}_{n-1}$ is the zero vector of dimension $n-1$, $A'\in\mathbb{R}^{(n-1)\times n}$ with $A'_{i,i} = -1, A'_{i,i+1} = 1$ and all other entries of $A'$ are zero. Thus, the problem of estimating convex Lipschitz function that is additionally monotonically nondecreasing can also be reduced to problem~\eqref{p1} with another matrix $A$ and vector $b$.}

\end{remark}

In the following we reduce the optimization problem~\eqref{p1} to a nonnegative least squares problem, which can then be solved efficiently using the \texttt{nnls} package in R. Define $z := Q^{1/2}(m - y)$, so that $m = Q^{-1/2}z + y$. Using this, we have $Am\ge b$ if and only if $AQ^{-1/2}z \ge b - Ay.$ Thus, \eqref{p1} is equivalent to
\begin{equation}\label{p2}
\mbox{minimize }|z|^2\mbox{ subject to }Gz\ge h,
\end{equation}
where $G := AQ^{-1/2}$ and $h := b - Ay.$ An equivalent formulation is
\begin{equation}\label{p3}
\mbox{minimize }|Eu - \ell|,\mbox{ over }u\succeq0,\mbox{ where } E := \begin{bmatrix}G^{\top}\\h^{\top}\end{bmatrix}\mbox{ and }\ell := [0,\ldots,0,1]^{\top}\in\mathbb{R}^{n+1}.
\end{equation}
Here $\succeq$ represents coordinate-wise inequality. A proof of this equivalence can be found in \citet[page 165]{LAW}; see \cite{CHEN} for an algorithm to solve \eqref{p3}.
 % can be found in \cite{LAW} and \cite{CHEN}

If $\hat{u}$ denotes the solution of \eqref{p3}  then the solution of \eqref{p2} is given as follows. Define $r := E\hat{u} - \ell$. Then $\hat{z}$, the minimizer of \eqref{p2}, is given by $\hat{z} := (-r_1/r_{n+1},\ldots,  -r_n/r_{n+1})^\top$\footnote{Note that~\eqref{p2} is a Least Distance Programming (LDP) problem and \citet[page 167]{LAW} prove that $r_{n+1}$ cannot be zero in an LDP with a feasible constraint set.}. Hence the solution to \eqref{p1} is given by $\hat{y} = Q^{-1/2}\hat{z} + y$.

\begin{remark}\label{rem:CvxLSE} Recall, the \texttt{CvxLSE} defined in~\eqref{eq:ultimate_est}. The \texttt{CvxLSE} can be computed via~Algorithm~\ref{algo:b} with $\A = \mathcal{C}$. To compute $m^{(k)}$ in Step 1 of~Algorithm~\ref{algo:b}, we can use strategy developed in~Section~\ref{sec:CLSE_comp} with ~\eqref{lipcons} replaced by
 % the algorithm developed in~Section~\ref{sec:CLSE_comp}The convexity constraint (i.e., $m\in\mathcal{C}$) can be represented by
 the following set of $n-2$ linear constraints:
\begin{equation}
\frac{m_2 - m_1}{t_2 - t_1} \le \frac{m_3 - m_2}{t_3 - t_2} \le \cdots \le \frac{m_{n} - m_{n-1}}{t_n - t_{n-1}}.\label{convcons}
\end{equation}
% where we use the notation of Section~\ref{sec:compute}.
 Similar to the CLSE, this reduces the computation of $m$ (for a given $\theta$) to solving a quadratic program with linear inequalities; see Section~\ref{sec:CLSE_comp}. The algorithm for computing $\theta^{(k+1)}$ developed below works for both \texttt{CvxLip} and $\texttt{CvxLSE}$.
\end{remark}
\subsection{Algorithm for computing \texorpdfstring{$\theta^{(k+1)}$}{Lg}} \label{sec:ThetaGradDesc}
In this subsection we describe an algorithm to find the minimizer $\theta^{(k+1)}$ of $\C(m^{(k)},\theta)$ over $\theta\in\Theta$. Recall that $\Theta$ is defined to be the ``positive" half of the unit sphere, a $d-1$ dimensional manifold in $\R^d$. Treating this problem as minimization over a manifold, one can apply a gradient descent algorithm by moving along a geodesic; see e.g.,~\citet[Section 3.3]{SAM}. But it is computationally expensive to move along a geodesic and so we follow the approach of \cite{WEN} wherein we move along a retraction with the guarantee of descent. To explain the approach of \cite{WEN}, let us denote the objective function by $f(\theta)$, i.e., in our case $f(\theta)=\C(m^{(k)},\theta)$. Let $\alpha \in \Theta$ be an initial guess for $\theta^{(k+1)}$ and define 
\[
  g := \nabla f(\alpha) \in \R^d \quad\mbox{and}\quad A := g\alpha^{\top} - \alpha g^{\top},
\]
where $\nabla$ denotes the gradient operator.
% In the following we use $g$ to denote both the function and the
Next we choose the path $\tau\mapsto \theta(\tau),$ where
\[
\theta(\tau) := \left(I + \frac{\tau}{2}A\right)^{-1}\left(I - \frac{\tau}{2}A\right)\alpha = \frac{1 + \frac{\tau^2}{4}[(\alpha^{\top}g)^2 - |g|^2] + \tau \alpha^{\top}g}{1 - \frac{\tau^2(\alpha^{\top}g)^2}{4} + \frac{\tau^2|g|^2}{4}}\alpha - \frac{\tau}{1 - \frac{\tau^2(\alpha^{\top}g)^2}{4} + \frac{\tau^2|g|^2}{4}}g,
\]
for $\tau\in\mathbb{R}$, and find a choice of $\tau$ such that $f(\theta(\tau))$ is as much smaller than $f(\alpha)$ as possible; see step 2 of Algorithm~\ref{algo:b}. It is easy to verify that
$$
\frac{\partial f(\theta(\tau))}{\partial \tau}\bigg|_{\tau = 0} \le 0;$$
see Lemma 3 of \cite{WEN}. This implies that $\tau \mapsto f(\theta(\tau))$ is a nonincreasing function in a neighborhood of $0$. %Given a value of $g$ and $\alpha$, $\theta(\tau)$ has the following closed form expression:
%\begin{equation}\label{eq:theta_tau_def}
%\theta(\tau) = \frac{1 + \frac{\tau^2}{4}[(\alpha^{\top}g)^2 - |g|^2] + \tau \alpha^{\top}g}{1 - \frac{\tau^2(\alpha^{\top}g)^2}{4} + \frac{\tau^2|g|^2}{4}}\alpha - \frac{\tau}{1 - \frac{\tau^2(\alpha^{\top}g)^2}{4} + \frac{\tau^2|g|^2}{4}}g;
%\end{equation}
%see Lemma 4 of \cite{WEN}.
Recall that for every $\eta\in\Theta$, $\eta_1$ (the first coordinate of $\eta$) is nonnegative. For $\theta(\tau)$ to lie in $\Theta$,  $\tau$ has to satisfy the following inequality
\begin{equation}\label{quad:ineq}
\frac{\tau^2}{4}[(\alpha^{\top}g)^2 - |g|^2] + \tau \left(\alpha^{\top}g - \frac{g_1}{\alpha_1}\right) + 1 \ge 0,
\end{equation}
where $g_1$ and $\alpha_1$ represent the first coordinates of the vectors $g$ and $\alpha$, respectively. This implies that a valid choice of $\tau$ must lie between the zeros of the quadratic expression on the left hand side of \eqref{quad:ineq}, given by
\[
2\frac{\left(\alpha^{\top}g - g_1/\alpha_1\right)\pm\sqrt{\left(\alpha^{\top}g - g_1/\alpha_1\right)^2 + |g|^2 - (\alpha^{\top}g)^2}}{|g|^2-(\alpha^{\top}g)^2}.
\]
Note that this interval always contains zero.  Now we can perform a simple line search for $\tau\mapsto f(\theta(\tau))$, where $\tau$ is in the above mentioned interval, to find $\theta^{(k+1)}$. \if1\blind  {We implement this step in the R package \texttt{simest}.} \fi  \if0\blind{We implement this step in the R package \texttt{***}.} \fi

\section{Main components in the proof of Theorem~\ref{thm:Main_rate_CLSE}} % (fold)
\label{app:sketchCLSE}
  In this section prove that~\eqref{eq:App_score_equation} implies~\eqref{eq:RAL_1}.

   \begin{enumerate}[label=\bfseries Step \arabic*]

  \item  In Theorem \ref{LIP:thm:nobiasCLSE} we show that $\psi_{\check{\theta},\check{m}}$ is approximately unbiased in the sense of \cite{VdV02}, i.e.,
  \begin{equation} \label{eq:nobias_main_th}
   \sqrt{n} P_{\check{\theta}, m_0} \psi_{\check{\theta},\check{m}} =o_p(1).
  \end{equation}Similar conditions have appeared before in proofs of asymptotic normality of maximum likelihood estimators (e.g., see \cite{MR1394975}) and the construction of efficient one-step estimators (see~\cite{MR913573}). The above condition essentially ensures that  ${\psi}_{\theta_0, \check{m}}$ is a good ``approximation'' to ${\psi}_{\theta_0, m_0}$; see~Section 3 of \cite{MurphyVaart00} for further discussion.\label{item:step3}
  % \todo[inline]{We need to fix the proof of Step 2 and 3 for because the definition of $\psi$ changed. It is NOT
  % $\psi_{\theta,m} (x,y) := (y-m(\theta^\top x)) H_\theta^\top [ m^\prime(\theta^\top x) x - { m_0^\prime(\theta^\top x)} h_{\theta_0}(\theta^\top x)].$

  % }
   \item We prove \label{item:step4}
  \begin{equation} \label{eq:Emp_proc_zero}
   \g_n ( \psi_{\check{\theta},\check{m}}- \psi_{\theta_0,m_0}) =o_p(1)
  \end{equation}
  in Theorem \ref{LIP:thm:ConsistencyofG_n}. Furthermore, as $\psi_{\theta_0,m_0}= {\ell}_{\theta_0,m_0},$ we have
  $P_{\theta_0,m_0} [\psi_{\theta_0, m_0}]=0.$
  Thus, by \eqref{eq:App_score_equation} and \eqref{eq:nobias_main_th}, we have  that \eqref{eq:Emp_proc_zero} is equivalent to
  \begin{equation}\label{eq:Final_eqpart1}
   \sqrt{n} (P_{\check{\theta}, m_0} - P_{\theta_0, m_0}) \psi_{\check{\theta},\check{m}} =\g_n {\ell}_{\theta_0,m_0} + o_p(1).
  \end{equation}
  \item To complete the proof, it is now enough to show that \label{item:step5}
  \begin{equation} \label{eq:ParaScore_approx}
  \sqrt{n} (P_{\check{\theta}, m_0} - P_{\theta_0, m_0}) \psi_{\check{\theta},\check{m}} = \sqrt{n}V_{\theta_0,m_0}  H_{\theta_0}^\top (\check{\theta}- \theta_0) + o_p(\sqrt{n} |\check{\theta} -\theta_0|).
  \end{equation}
  A proof of \eqref{eq:ParaScore_approx} can be found in the proof of Theorem 6.20 in \cite{VdV02}; also see~\citet[Section 10.4]{Patra16}. Lemma~\ref{thm:nobiasCLSE_part2} in  Section~\ref{app:nobias_part2CLSE} of the supplementary file
    proves that $(\check\theta, \check{m})$ satisfy the required conditions of Theorem 6.20 in \cite{VdV02}.
   \end{enumerate}
   \noindent  Observe that  \eqref{eq:Final_eqpart1} and \eqref{eq:ParaScore_approx} imply
   \begin{align} \label{eq:final_thm}
    \begin{split}
  \sqrt{n}V_{\theta_0,m_0}  H_{\theta_0}^\top(\check{\theta}- \theta_0)  ={}&\g_n \ell_{\theta_0,m_0} + o_p(1+\sqrt{n} |\check{\theta}-\theta_0|),\\
   \Rightarrow  \sqrt{n} H_{\theta_0}^\top (\check{\theta}- \theta_0) ={}& V_{\theta_0,m_0}^{-1}\g_n \ell_{\theta_0,m_0} + o_p(1) \stackrel{d}{\rightarrow} V_{\theta_0,m_0}^{-1} N(0, {I}_{\theta_0,m_0}).
   \end{split}
   \end{align}
  The proof of the theorem will be complete, if we can show that  $$\sqrt{n}(\check{\theta}- \theta_0)=H_{\theta_0} \sqrt{n} H_{\theta_0}^\top(\check{\theta}- \theta_0) +o_p(1),$$ the proof of which can be found in Step 4 of Theorem 5 in~\cite{Patra16}. 
\begin{figure}[!ht]
% \captionsetup[subfigure]{labelformat=empty}
\centering
\includegraphics[width=.85\textwidth]{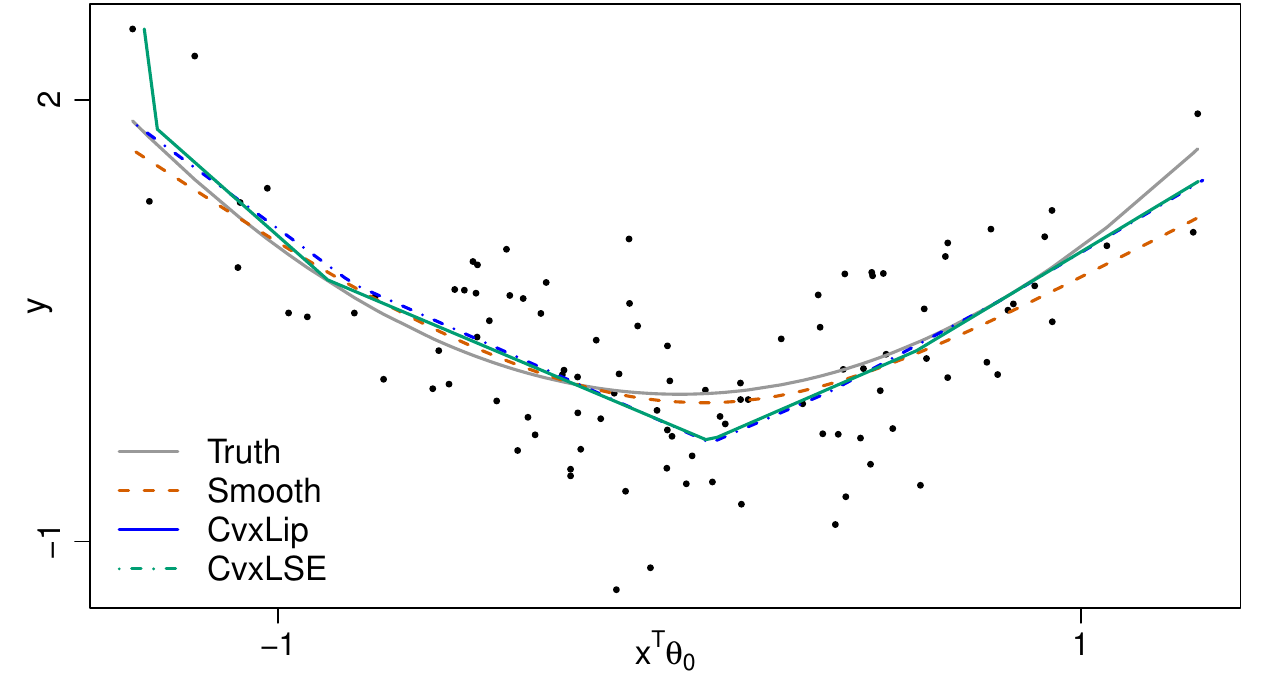}\\ %%
\caption[An illustrative example]{Function estimates for  the model $Y=(\theta_0^\top X)^2+ N(0,1)$, where $\theta_0=(1,1,1,1,1)^\top/\sqrt{5}, X\sim\text{Uniform}[-1,1]^5,$ and $n=100$.}
\label{fig:SimplModel}
%% label for entire figure
\end{figure}
% section section_name (end)

\section{Uniqueness of the limiting variances in Theorem~\ref{thm:Main_rate_CLSE}}\label{sec:Degenracy}

Observe that the variance of the limiting distribution (for both the heteroscedastic and homoscedastic error models) is singular. This can be attributed to the fact that $\Theta$ is a Stiefel manifold of dimension ${d-1}$ and has an empty interior in $\R^d.$

In Lemma~\ref{lem:UniqueInverse} below, we show that the limiting variances are unique, i.e., they do not depend on the particular choice of $\theta\mapsto H_\theta$.  In fact $H_{\theta_0}{I}^{-1}_{\theta_0,m_0} H_{\theta_0}^\top$ matches the lower bound obtained in~\cite{NeweyStroker93} for the single index model under only smoothness constraints. 
% Observe that when $\sigma(\cdot)\eqiv \sigma$, \sigma^4

 \begin{lemma}\label{lem:UniqueInverse}
  Suppose the assumptions of Theorem~\ref{thm:Main_rate_CLSE} hold, then the matrix $H_{\theta_0}{I}^{-1}_{\theta_0,m_0} H_{\theta_0}^\top$ is the unique Moore-Penrose inverse of
 \begin{equation}\label{eq:True_inf}
 P_{\theta_0, m_0} \big[\{(Y-m_0(\theta_0^\top X))m'_0(\theta_0^\top X)\}^2 (X-h_{\theta_0}(\theta_0^\top X))(X-h_{\theta_0}(\theta_0^\top X))^\top \big] \in \R^{d\times d}.
 \end{equation}
 \end{lemma}

\begin{proof}
Recall that 
\begin{equation}\label{eq:I_def}
{I}_{\theta,m} = H_\theta^\top \E\Big[\big(Y-m(\theta^\top X)\big) m^\prime(\theta^\top X)\big]^2   \big[ X - h_\theta(\theta^\top X)\big]  \big[ X - h_\theta(\theta^\top X)\big]^\top\Big] H_\theta.
\end{equation}
For the rest of the proof, define
\begin{equation}\label{eq:A_def}
A:=\E\Big[\big(Y-m(\theta^\top X)\big) m^\prime(\theta^\top X)\big]^2   \big[ X - h_\theta(\theta^\top X)\big]  \big[ X - h_\theta(\theta^\top X)\big]^\top\Big].
\end{equation}
In the following, we show that $G:= H_\theta(H_\theta^\top {A} H_\theta)^{-1} H_\theta^\top$ is the Moore-Penrose inverse of ${A}.$ By definition, it is equivalent to show that 
\begin{equation}\label{eq:prop}
 AGA=A,\;  GAG=G,\; (AG)^\top= AG,\text{ and } (GA)^\top= GA.
 \end{equation} 
\textbf{Proof of $\mathbf{AGA=A}$:} We will now show that $A G A= A$, an equivalent condition is that  $G A$ is idempotent and $\text{rank}(GA) = \text{rank}(A)$. Observe that $GA$ is idempotent because,
\begin{equation}\label{eq:idem_GA}
 GAGA = H_\theta(H_\theta^\top {A} H_\theta)^{-1} H_\theta^\top A H_\theta(H_\theta^\top {A} H_\theta)^{-1} H_\theta^\top A= H_\theta(H_\theta^\top {A} H_\theta)^{-1} H_\theta^\top A= GA.
\end{equation}
Note that $H_\theta^\top A G A= H_\theta^\top A.$ Thus $\text{rank}(H_\theta^\top A) \le \text{rank}(GA)$. However, 
\[ GA = H_\theta(H_\theta^\top {A} H_\theta)^{-1} H_\theta^\top A = \big[H_\theta(H_\theta^\top {A} H_\theta)^{-1}\big] H_\theta^\top A.\]
 Thus $\text{rank}(GA) = \text{rank}(H_\theta^\top A).$ Thus to prove $\text{rank}(GA) = \text{rank}(A)$ it enough to show that  $\text{rank}(H_\theta^\top A) = \text{rank}(A)$.
% We know need to show that $\text{rank}(GA) = \text{rank}(A)$. By definition,  $\text{rank}(GA) \le \text{rank}(A)$. We only need to show the reverse.
% \begin{equation}\label{eq:rank_inew}
% GA =  H_\theta(H_\theta^\top {A} H_\theta)^{-1} H_\theta^\top A
% \end{equation}
We will prove that the nullspace of $H_\theta^\top A$ is the same as that of $A$. Since $Ax = 0 $ implies that $H_\theta^\top Ax = 0$, it follows
that
\[ N(A):=\{x : Ax = 0\} \subseteq \{x : H_\theta^\top Ax = 0\}:= N(H_\theta^\top A).\]
We will now prove the reverse inclusion by contradiction.   Suppose there exists a vector $x$ such that $Ax\neq 0$ and $H_\theta^\top Ax = 0.$ Set $y = Ax$.  Then we have that $H_\theta^\top y=0$. Thus by Lemma~1 of~\cite{Patra16}, we have that $y = c\theta$ for some constant $c\neq 0$. (If $c = 0$, then $y = Ax = 0$, a contradiction). This implies that there exists $x$ such that $Ax = c\theta$ or in particular $\theta^{\top}Ax = c\neq 0$, since $\norm{\theta} = 1$. This, however, is a contradiction since $A$ is symmetric and
\begin{align}\label{eq:Atheta}
\begin{split}
 A\theta &= \E\bigg[ \big[(y-m(\theta^\top x)) m^\prime(\theta^\top x)\big]^2  H_\theta^\top \left\lbrace x - h_\theta(\theta^\top x)\right\rbrace\left\lbrace x - h_\theta(\theta^\top x)\right\rbrace^\top \theta\bigg]\\
&=E\bigg[ \big[(y-m(\theta^\top x)) m^\prime(\theta^\top x)\big]^2  H_\theta^\top \left\lbrace x - h_\theta(\theta^\top x)\right\rbrace\left\lbrace \theta^\top x - \E(\theta^\top X|\theta^\top x)\right\rbrace^\top \bigg]\\
&=\mathbf{0}_d.
\end{split}
\end{align}
\textbf{Proof of $\mathbf{GAG=G}$:} It is easy to see that 
\begin{equation}\label{eq:GAG}
 GAG = H_\theta(H_\theta^\top {A} H_\theta)^{-1} H_\theta^\top A H_\theta(H_\theta^\top {A} H_\theta)^{-1} H_\theta^\top= H_\theta(H_\theta^\top {A} H_\theta)^{-1} H_\theta^\top = G.
\end{equation}
\textbf{Proof of $\mathbf{(AG)^\top=AG}$:}
\begin{equation}\label{eq:GAtrans}
 (AG)^\top = (A H_\theta(H_\theta^\top {A} H_\theta)^{-1} H_\theta^\top)^\top = H_\theta(H_\theta^\top {A} H_\theta)^{-1} H_\theta^\top A^\top  =  H_\theta(H_\theta^\top {A} H_\theta)^{-1} H_\theta^\top A ,
\end{equation}
as $A$ is a symmetric matrix. 
Recall that $H_\theta\in \R^{d\times (d-1)}$ and the columns of $H_\theta$ are orthogonal to $\theta$. Thus let us define $\overline{H}_\theta \in \R^{d\times d}$, by adding $\theta$ as an additional column to $H_\theta$, i.e., $\overline{H}_\theta= [H_\theta ,\theta]$. Recall that by definition of $H_\theta$, $\theta^\top H_\theta= \mathbf{0}_{d-1}$ and~\eqref{eq:Atheta}, we have that $\theta^\top A = A \theta=\mathbf{0}_{d-1}$. Multiplying $(AG)^\top$ by $\overline{H}^\top_\theta$ on the left and  $\overline{H}_\theta$ on the right we have,
\begin{align}\label{eq:AGtrans}
\begin{split}
\overline{H}^\top_\theta (AG)^\top \overline{H}_\theta &=\overline{H}^\top_\theta H_\theta(H_\theta^\top {A} H_\theta)^{-1} H_\theta^\top A \overline{H}_\theta\\
 &= \begin{bmatrix}
    H^\top_\theta H_\theta(H_\theta^\top {A} H_\theta)^{-1} H_\theta^\top A H_\theta &
   H^\top_\theta H_\theta(H_\theta^\top {A} H_\theta)^{-1} H_\theta^\top A \theta \\
   \theta^\top H_\theta(H_\theta^\top {A} H_\theta)^{-1} H_\theta^\top A H_\theta & 
   \theta H_\theta(H_\theta^\top {A} H_\theta)^{-1} H_\theta^\top A \theta
   \end{bmatrix}\\
&= \begin{bmatrix}
    H^\top_\theta H_\theta &
   \mathbf{0}_{d-1} \\
   \mathbf{0}_{d-1}^\top& 
   0
   \end{bmatrix}.
\end{split}
\end{align}
Multiplying $AG$ by $\overline{H}^\top_\theta$ on the left and  $\overline{H}_\theta$ on the right we have,

\begin{align}\label{eq:AG}
\begin{split}
\overline{H}^\top_\theta AG \overline{H}_\theta &=\overline{H}^\top_\theta A H_\theta (H_\theta^\top {A} H_\theta)^{-1} H_\theta^\top \overline{H}_\theta\\
&= \begin{bmatrix}
    H^\top_\theta  A H_\theta(H_\theta^\top {A} H_\theta)^{-1} H_\theta^\top \overline{H}_\theta\\
   \mathbf{0}_d^\top
   \end{bmatrix}\\
 &= \begin{bmatrix}
    H^\top_\theta H_\theta &
   \mathbf{0}_{d-1} \\
   \mathbf{0}_{d-1}^\top& 
   0
   \end{bmatrix},
\end{split}
\end{align}
here the second equality is true, since $ \overline{H}^\top_\theta A =[H_\theta^\top A ,\theta^\top A]^\top = [H_\theta^\top A, \mathbf{0}_d].$
Thus, $\overline{H}^\top_\theta (AG)^\top \overline{H}_\theta =\overline{H}^\top_\theta AG \overline{H}_\theta$. Since $ \overline{H}_\theta $ is a nonsingular matrix, we have that $(AG)^\top =AG.$ Proof of $(GA)^\top= GA$ follows similarly.
\end{proof}

\section{Additional simulation studies}\label{app:add_simul}
% This section presents some additional simulation studies.
\subsection{A simple model}
% {\begin{remark}[Pre-binning]\label{rem:bin}
% The matrices involved in the algorithm above have entries depending on fractions such as $1/(t_{i+1} - t_i)$. Thus if there are ties in $\{t_i\}_{1\le i \le n}$, then the matrix $A$ is incomputable.  Moreover, if $t_{i+1} - t_i$ is very small, then the fractions can force the matrices involved to be ill-conditioned (for the purposes of numerical calculations). Thus  to avoid ill-conditioning of these matrices, in practice one might have to pre-bin the data which leads to a diagonal matrix $Q$ with different diagonal entries. One common method of pre-binning the data is to take the means of all data points for which the $t_i$'s are close. To be more precise, if we choose a tolerance of $\eta = 10^{-6}$ and suppose that $0 < t_2 - t_1 < t_3 - t_1 < \eta$, then we combine the data points $(t_1, y_1), (t_2, y_2), (t_3, y_3)$ by taking their mean and set $Q_{1,1} = 3$; the total number of data points is now reduced to $n-2$.
% \end{remark}}
% \todo[inline]{ Move this to supplementary?}

In this section we give a simple illustrative (finite sample) example. We observe $100$ i.i.d.~observations from the following homoscedastic model:
\begin{equation}\label{eq:SimpleMOdel}
Y=(\theta_0^\top X)^2+ N(0,1),\text{ where }\theta_0=(1,1,1,1,1)/\sqrt{5}\;\text{and}\; X\sim\text{Uniform}[-1,1]^5.
\end{equation}
In Figure~\ref{fig:SimplModel}, we have a scatter plot of $\{(\theta_0^\top X_i,Y_i)\}_{ i=1}^{ 100}$ overlaid with prediction curves $\{(\tilde{\theta}^\top X_i,\tilde{m}( \tilde{\theta}^\top X_i)\}_{ i=1}^{ 100}$ for the proposed estimators obtained from \textit{one} sample from~\eqref{eq:SimpleMOdel}. Table~\ref{tab:theta_SimplModel} displays all the corresponding estimates of $\theta_0$ obtained from the same data set. To compute the function estimates for \texttt{EFM} and \texttt{EDR} approaches we used cross-validated smoothing splines to estimate the link function using their estimates of $\theta_0$.

\begin{table}[!ht]
\caption[Estimates of $\theta_0$ for a simple model]{\label{tab:theta_SimplModel}Estimates of $\theta_0$, ``Theta Error''$:=\sum_{i=1}^5 |\tilde{\theta}_i-\theta_{0,i}|,$  ``Func Error''$:=\|\tilde{m}\circ \theta_0-m_0\circ\theta_0\|_n$, and { ``Est Error''}$:=\|\tilde{m}\circ \tilde{\theta}-m_0\circ\theta_0\|_n$ for one sample from~\eqref{eq:SimpleMOdel}.}
\centering
\begin{tabular}{lcccccccc}
  \toprule
Method & $\theta_1$ & $\theta_2$ & $\theta_3$ & $\theta_4$ & $\theta_5$ & Theta Error & Func Error & Est Error \\
  \midrule
Truth & 0.45 & 0.45 & 0.45 & 0.45 & 0.45 & --- & --- & --- \\
  \texttt{Smooth} & 0.38 & 0.49 & 0.41 & 0.50 & 0.45 & 0.21 & 0.10 & 0.10 \\
  %\texttt{CvxPen} & 0.36 & 0.50 & 0.42 & 0.47 & 0.47 & 0.21 & 0.12 & 0.13 \\
  \texttt{CvxLip} & 0.35 & 0.50 & 0.43 & 0.48 & 0.46 & 0.21 & 0.13 & 0.15 \\
  \texttt{CvxLSE} & 0.36 & 0.50 & 0.43 & 0.45 & 0.48 & 0.20 & 0.18 & 0.15 \\
  \texttt{EFM} & 0.35 & 0.49 & 0.41 & 0.49 & 0.47 & 0.24 & 0.10 & 0.11 \\
  \texttt{EDR} & 0.30 & 0.48 & 0.46 & 0.43 & 0.53 & 0.29 & 0.12 & 0.15 \\
   \bottomrule
\end{tabular}
\end{table}
\begin{figure}[!ht]
% \captionsetup[subfigure]{labelformat=empty}
\centering
\includegraphics[width=.9\textwidth]{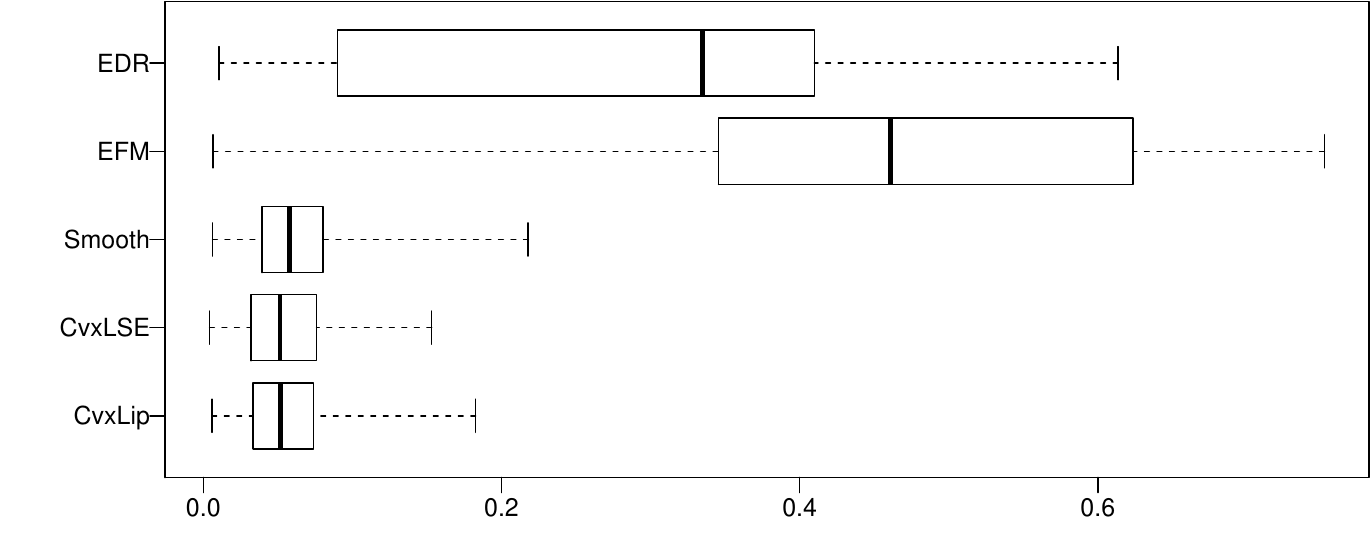}\\ %%
\caption[Boxplots of estimates when the truth is a piecewise affine convex function]{Box plots of $\sum_{i=1}^6|\tilde{\theta}_i-\theta_{0,i}|$ for the model~\eqref{eq:V_mod}. Here $d=6,$ $n=200$ and we have $500$ replications.}
\label{fig:LiPatilea15}
%% label for entire figure

\end{figure}
\subsection{Piecewise affine function and dependent covariates}
% Chapters~\ref{chap:smoothsim} and~\ref{chap:cvX_sim} study the asymptotic properties of the estimators when the true link function is smooth.
 To understand the performance of the estimators when the truth is convex but not smooth, we consider the following model:
\begin{equation}\label{eq:V_mod}
Y=|\theta_0^\top X| + N(0, .1^2),
\end{equation}
where $X\in\R^6$ is generated according to the following law: $(X_1, X_2) \sim \text{Uniform} [-1,1]^2$, $ X_3:=0.2 X_1+ 0.2 (X_2+2)^2+0.2 Z_1$, $X_4:=0.1+ 0.1(X_1+X_2)+0.3(X_1+1.5)^2+  0.2 Z_2$, $X_5 \sim \text{Ber}(\exp(X_1)/\{1+\exp(X_1)\}),$ and  $ X_6 \sim \text{Ber}(\exp(X_2)/\{1+\exp(X_2)\})$. Here $(Z_1, Z_2)\sim\text{Uniform} [-1,1]^2$ is independent of $(X_1,X_2)$ and $\theta_0$ is $(1.3, -1.3, 1, -0.5, -0.5, -0.5)/\sqrt{5.13}$.  The distribution of the covariates is similar to the one considered in Section~V.2 of~\cite{LiPatilea15}. The performances of the estimators is summarized in~Figure~\ref{fig:LiPatilea15}.  Observe that as the truth is not smooth, the convex constrained least squares estimators (\texttt{CvxLip} and \texttt{CvxLSE}) have slightly improved performance compared to the (roughness) penalized least squares estimator (\texttt{Smooth}). Also observe that both \texttt{EFM} and \texttt{EDR} fail to estimate the true parameter $\theta_0.$
\subsection{Investigation of adaptation of the CLSE} % (fold)
\label{sub:investigation_of_the_adaptation_of_the_clse}
{ In this subsection, we present a brief simulation study to illustrate the adaptive behavior of the CLSE when $m_0$ is a piecewise linear convex function.  We generate 400 replications of $n$ i.i.d.~observations from the following model:
\begin{equation}\label{eq:qqplot_sup}
Y=|\theta_0^\top X| + N(0, .1^2), \qquad \text{where}\; X \sim \text{Uniform} [-1, 1]^2 \quad \text{and}  \quad \theta_0 = (1,1)/\sqrt{2},
\end{equation}
for $n$ increasing geometrically from $100$ to $2000$. To investigate the adaptive properties of the CLSE, we compute average estimation error ($\|\check{m}(\check \theta ^\top X) -m_0(\theta_0^\top X)\|_n^2$) as sample size increases and plot $\|\check{m}(\check \theta ^\top X) -m_0(\theta_0^\top X)\|_n^2$ versus $n$ in a log-log scale; we use $L=10$. If the rate of convergence of the CLSE is $n^{-\alpha}$ then the  slope of the best fitting line should be close to $-\alpha$. In the left panel of Figure~\ref{fig:adap}, the best fitting line has a slope of $-0.95$, suggesting a near parametric rate of convergence for the CLSE; cf. the slope of $-0.8$ expected from the worst case rate in Theorem~\ref{thm:rate_m_theta_CLSE}. Additionally, the right panel shows the Q-Q plot of $\sqrt{n}(\check{\theta}-\theta_0)$. Notice that $\mathrm{Var}(\sqrt{n}(\check{\theta}-\theta_0))$ does not stabilize with the sample size, suggesting non-standard behavior for the CLSE. This kind of behavior is not well understood and can be observed in other shape constrained semiparametric models when the estimate of nonparametric component  exhibits  a near parametric rate of convergence.}
\begin{figure}[!ht]
% \captionsetup[subfigure]{labelformat=empty}
\centering
\includegraphics[width=.9\textwidth]{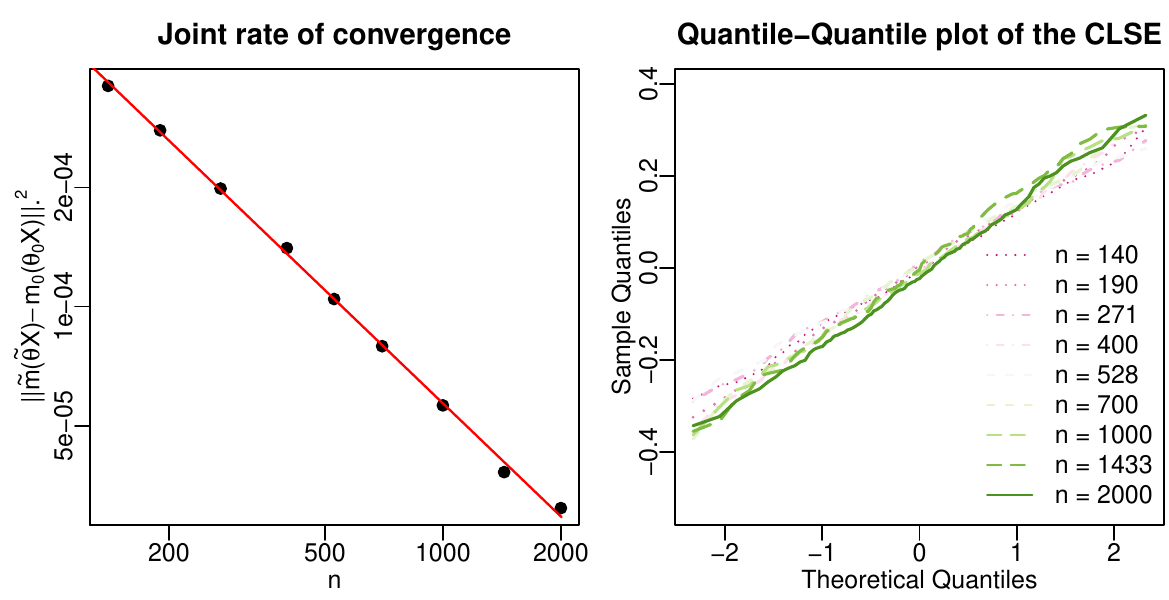}\\ %%
\caption[Boxplots of estimates when the truth is a piecewise affine convex function]{Asymptotic behavior of the CLSE when $m_0$ is a piecewise linear convex function. Left panel: plot of $\log(\|\check{m}(\check \theta ^\top X) -m_0(\theta_0^\top X)\|_n^2)$ vs $\log n$ overlaid with best fitting line (in red). The line has a slope of $-0.95$. Right panel: Q-Q plots of $\sqrt{n}(\check{\theta}-\theta_0)$ as the sample size grows from $100$ to $2000$. All simulations are based on 400  random samples.}
\label{fig:adap}
%% label for entire figure

\end{figure}
% subsection investigation_of_the_adaptation_of_the_clse (end)
% \appendix

\section{Continuing the discussion of identifiability from Section~\ref{sec:ident}}\label{sec:proof:Estims}
\subsection[Proof of (\ref{eq:true_mimina})]{Proof of \eqref{eq:true_mimina}}\label{proof:lem:Ident_cvxSim}
 In the following we show that $(m_0,\theta_0)$ is the  minimizer of $Q$ and is well-separated, with respect to the $L_2(P_X)$ norm, from $\{(m, \theta): m\circ\theta \in L_2(P_X)\}\setminus \{(m,\theta):\,\|m\circ\theta - m_0\circ\theta_0\| \le \delta\}$. Pick any $(m,\theta)$ such that $ m\circ\theta \in L_2(P_X)$ and $\|m\circ\theta - m_0\circ\theta_0\|^2 > \delta^2$. Then
\[
Q(m,\theta) = \E [Y-m_0(\theta_0^\top X)]^2 + \E[m_0(\theta_0^\top X) - m(\theta^\top X)]^2,
\]
since $\E(\epsilon|X) = 0$. Thus we have that $Q(m,\theta) > Q(m_0, \theta_0)+ \delta^2.$
%%%%%%%%%%%%%%%%%%%%%%%%%%%%%%%%%%%%%%%$
%%%%%%%%%%%%%%%%%%%%%%%%%%%%%%%%%%%%%%%$
\subsection{Discussion on the identifiability of separated parameters}\label{sec:IdentDisc}
{ The goal of the subsection is to describe various conditions on $\theta_0$, $m_0$, and the distribution of $X$ under which the model parameters can be identified separately. One of the most general sufficient conditions we could find in the literature on identifiability is from~\citet[Theorem 4.1]{ICHI93}. The author shows that $m_0$ and $\theta_0$ are separately identifiable if:

\begin{enumerate}[label=\bfseries (I)]
\item The function $m_0(\cdot)$ is non-constant, non-periodic, and a.e.~differentiable\footnote{Note that all convex functions are almost everywhere differentiable and are not periodic.} and $|\theta_0|=1$.  The components of the covariate $X =(X_{1}, \ldots, X_{d})$ do not have a perfect linear relationship.  There exists an integer $d_1\in \{1,2,\ldots,d\}$ such that  $X_{1},\ldots, X_{d_1}$ have continuous distributions  and  $X_{ d_1+1},\ldots, X_{d}$ are discrete random variables. The first non-zero coordinate of $\theta_0= (\theta_{0,1},\ldots, \theta_{0,d})$ is positive and at least one of  $\theta_{0,1}, \ldots, \theta_{0,d_1}$ is non-zero.  Furthermore, there exist an open interval $\mathcal{I}$ and non-random vectors  $c_0,c_1,\ldots, c_{d-d_1} \in \mathbb{R}^{d-d_1}$ such that \label{a0sup}
\begin{itemize}
\item $c_l-c_0$ for $l\in \{1 ,\ldots, d-d_1\}$ are linearly independent,
\item $\mathcal{I} \subset \bigcap_{l=0}^{d-d_1} \big\{ \theta_0^\top x:  x\in \rchi \text{ and } (x_{d_1+1}, \ldots x_d) =c_l \big\}.$
\end{itemize}
\end{enumerate}
 An alternative and perhaps simpler  condition for identifiability of $(m_0,\theta_0)$ is given in~\citet[Theorem 1]{LinKul07}:
\begin{enumerate}[label=\bfseries (I$'$)]
\item The support of $X$ is a bounded convex set in $\R^d$ with non-empty interior. The link function  $m_0$ is non-constant and continuous. The first non-zero coordinate of $\theta_0$ is positive and $|\theta_0|=1$.\label{a0primesup}
\end{enumerate}
Assumptions~\ref{a0sup} and~\ref{a0primesup} are necessary for identifiability in their own way; see~\cite{liracine07} and~\cite{LinKul07} for details. Also see~\cite{yuan2011identifiability}. However we prefer~\ref{a0sup} to~\ref{a0primesup}, because~\ref{a0sup} allows for discrete covariates (a common occurrence in practice).}

\section{Minimax lower bound} % (fold)
\label{sec:minimax_lower_bound}
 In the following proposition, we prove that when the single index model in~\eqref{eq:simsl} satisfies assumptions~\ref{aa1_new}--\ref{aa2} and  the errors are Gaussian random  variables  (independent of the covariates) then  $n^{-2/5}$ is a minimax lower bound on the rate of convergence for estimating $m_0\circ\theta_0$. 
Thus  $\check{m}_{L}\circ\check{\theta}_L$ is minimax rate optimal when $q\ge 5$.

\begin{prop}[Minimax lower bound]\label{lem:MinimaxLowerbound}
Suppose that $\{(X_i,Y_i)\}_{i=1}^n$ are i.i.d.~observations from~\eqref{eq:simsl} such that assumptions~\ref{aa1_new}--\ref{aa2} are satisfied and $\theta_0^\top X\sim\text{Uniform}[0,1]$. Moreover, suppose that the errors are independent of the covariates and  $\epsilon \sim N(0, \sigma^2)$ for some $\sigma>0$.  Then there exist positive constants $k_1$  and $k_2$, depending only on $\sigma$ and $L_0$, such that
\begin{equation}\label{eq:Minmaxbound}
 \inf_{\hat{f}} \sup_{(m_0,\theta_0) \in \M_{L_0}\times \Theta} \p\left(n^{2/5} \|\hat{f} -m_0\circ\theta_0\|  > k_1\right) \ge k_2 >0,
\end{equation}
% \begin{equation}\label{eq:Minmaxbound}
% \inf_{\hat{f}} \sup_{(m_0,\theta_0) \in \M_{L_0}\times \Theta} \p(\|\hat{f} -m_0\circ\theta_0\| \ge k_1 n^{-2/5}) > k_2,
% \end{equation}
where the infimum is over all estimators of $m_0\circ\theta_0$ based on $\{(X_i,Y_i)\}_{i=1}^n$. 
\end{prop}
\begin{proof}

%%%%%%%%%%%%%%%%%%%%%%%%%%%%%%%%%%%%%%%%%%%%%%%%%%%%%
Recall that for this proposition we assume that, we have i.i.d.~observations $\{(X_i,Y_i)\}_{i=1}^n$ from~\eqref{eq:simsl} such that assumptions~\ref{a0}--\ref{aa2} are satisfied and $\theta_0^\top X\sim\text{Uniform}[0,1]$. Moreover, we assume that the errors are independent of the covariates and  $\epsilon \sim N(0, \sigma^2)$, where $\sigma>0$.
Consider $\theta_0^{(2)},\ldots, \theta_0^{(d)}$ in $\R^d$ such that $\{\theta_0,\theta_0^{(2)},\ldots, \theta_0^{(d)}\}$ form an orthogonal basis of $\R^d$. We denote the matrix with $\theta_0, \theta_0^{(2)},\ldots,\theta_0^{(d-1)}$, and $\theta_0^{(d)}$ as columns by $\mathbf{O}$.  Let $Z= (Z^{(1)}, \ldots, Z^{(d)})= \mathbf{O} X$ and $Z_2^d :=(Z^{(2)}, \ldots, Z^{(d)})$. In the proof of Theorem 2 in~\cite[Page 561]{MR2369025} the authors show that if $\hat{g}$ is an estimator for $m_0\circ\theta_0$ in the model~\eqref{eq:simsl}, then 
\begin{equation}\label{eq:Equiv_1d_reg}
\|\hat{g}-m_0\circ\theta_0\|^2  \ge \int_{0}^1 \big[\hat{f}(t)-m_0(t)\big]^2 dt,
\end{equation}
where $\hat{f}:= \int \hat{g} (\mathbf{O}^{-1} Z) P_{Z_2^{d}| Z^{(1)}} (dz_2^d|z^{(1)})$. Thus we have that for any $k>0,$
\begin{align}\label{eq:minimaxbound_1}
\begin{split}
 % \qquad
  \inf_{\hat{g}}&\sup_{(m_0,\theta_0) \in \M_{L_0}\times \Theta}\mathbb{P}\left[n^{2/5}\norm{\hat{g} - m_0\circ\theta_0} \ge k \right] 
\ge{} \inf_{\hat{f}}\sup_{f_0\in \M_{L_0}}\mathbb{P}\left[n^{2/5}\|\hat{f} - f_0\|_{\Lambda}  \ge k \right],
\end{split}
\end{align}
where for any $f:[0,1] \to \R$, $\|f\|_{\Lambda} := \int_{0}^{1} f^2(t) dt$ and the infimum on the right is over all estimators of $f_0$ based on the data satisfying the assumptions with $d=1$, i.e., univariate regression. The following lemma completes the proof of Proposition~\ref{lem:MinimaxLowerbound} by  establishing an lower bound (see~\eqref{eq:1dConvexRegBound}) for the quantity on the right.
\end{proof}
% On this set of probability distributions, estimating $m\circ\theta_0$ is same as univariate convex regression (since $\theta_0$ is known). From the minimax lower bounds on univariate convex regression, it is known that {\color{blue}the convex function cannot be estimated at a rate better than $n^{2/5}$ in $\norm{\cdot}$-norm without additional assumptions on $m$} (!!! need a reference !!!). This implies that our rate is actually minimax optimal. 

\begin{lemma}\label{lem:1dConvexreg}
 Suppose we have an i.i.d.~sample from the following model:
\begin{equation}\label{eq:1dimConvex}
Z= f(W) +\xi,
\end{equation}
where $W\sim\text{Uniform}[0,1]$, $\xi\sim N(0, \sigma^2)$, and $\xi$'s are independent of the covariates. Let $f: [0,1]\to \R$ be a uniformly Lipschitz convex function with Lipschitz constant $L_0$. Then there exists a constant $k_1, k_2>0$ (depending only on $\sigma$ and $L_0$) such that  \begin{equation}\label{eq:1dConvexRegBound}
\inf_{\hat{f}} \sup_{f \in \M_{L_0}} \p_{f}\left(n^{2/5} \|\hat{f}-f\|_2 \ge k_1 \right) \ge k_2 > 0,
\end{equation} 
where the infimum is over all estimators of $f.$
 \end{lemma}
\begin{proof}
To prove the above lower bound we will follow the general reduction scheme described in Section 2.2 of~\citet[Page 79]{MR2724359}. Fix $n$ and let $m := c_0 n^{1/5}$, where $c_0$ is a constant to be chosen later (see~\eqref{eq:c_0_def}) and let  $M := 2^{m/8}$. Let us assume that there exist $f_0, \ldots,f_M \in \M_L$ such that, for all $0\le j\ne k \le M$,
\begin{equation}\label{eq:prop-1&2f_j}
\|f_j -f_k\| \ge 2 s\quad\text{where }s :=A m^{-2}\quad\text{and}\quad A:= \frac{\kappa_1}{88 c_0^2}(b-a)^{5/2}.
 \end{equation}
Let  $P_j$ denote the  joint distribution of $(Z_1,W_1),\ldots, (Z_n,W_n)$ for $f= f_{j}$ (in~\eqref{eq:1dimConvex}) and $\E_{W_1, \ldots, W_n}$ denote the expectation with respect to the joint distribution of $W_1, \ldots, W_n$. Let $\hat f$ be any estimator. Observe that  \begin{align}
  % \inf_{\hat{f}}\sup_{f \in \M_L} \E_{f} \big[n^{4/5}\|\hat{f}-f\|_2^2\big] 
  \sup_{f \in \M_L} \p_f\big(\|\hat{f}-f\|_2 \ge A c_0^{-2} n^{-2/5}\big)
% \ge& \max_{f \in \{f_{0}, \ldots, f_{M}\}} \E_{f} \big[n^{4/5}\|\hat{f}-f\|_2^2\big]\\
\ge&\max_{f \in \{f_{0}, \ldots, f_{M}\}} \p_f\big(\|\hat{f}-f\|_2 \ge A c_0^{-2} n^{-2/5}\big)\nonumber\\
\ge &\frac{1}{M+1}\sum_{j=0}^{M}\p_{j}\big(\|\hat{f}-f_j\|_2 \ge s\big)\label{eq:MultiHyp}\\
=& \frac{1}{M+1}\sum_{j=0}^{M}\E_{W_1, \ldots, W_n}\Big[\p_{j}\big(\|\hat{f}-f_j\|_2 \ge s\big\vert W_1, \ldots, W_n\big)\Big]\nonumber\\
=&  \E_{W_1, \ldots, W_n}\left[\frac{1}{M+1}\sum_{j=0}^{M}\p_{j}\big(\|\hat{f}-f_j\|_2 \ge s\big\vert W_1, \ldots, W_n\big)\right].\nonumber
\end{align}
Consider the $M+1$ hypothesis elements $f_0,\ldots, f_M$. Any test in this setup is a measurable function $\psi:\{(Z_1,W_1),\ldots, (Z_n,W_n)\} \to \{0,\ldots, M\}$.  {Let us define $\psi^* $ to be the minimum distance test, i.e., $$\psi^* := \arg \min_{0 \le k \le M} \|\hat f - f_k\|_2.$$ Then if $\psi^* \ne j$ then  $\|\hat{f}-f_j\|_2 \ge \|\hat{f}-f_\psi^*\|_2$ and \[
  2s \le \|f_j -f_\psi^*\|_2 \le \|\hat{f}-f_j\|_2 + \|\hat{f}-f_\psi^*\|_2 \le 2 \|\hat{f}-f_j\|_2.
\] Thus
\[\p_{j}\big(\|\hat{f}-f_j\|_2 \ge s\big\vert W_1, \ldots, W_n\big) \ge \p_j (\psi^* \ne j\big\vert W_1, \ldots, W_n) \qquad \text{for all } 0\le j \le M.\]}
Combining~\eqref{eq:MultiHyp} with the above display, we get
\begin{align*}
\sup_{f \in \M_L} \p_f\big(\|\hat{f}-f\|_2 \ge A c_0^{-2} n^{-2/5}\big)
& \ge  \E_{W_1, \ldots, W_n} \left[\frac{1}{M+1}\sum_{j=0}^{M} \p_j (\psi^* \ne j\big\vert W_1, \ldots, W_n) \right] \\
& \ge  \E_{W_1, \ldots, W_n}  \left[\inf_{\psi} \frac{1}{M+1}\sum_{j=0}^{M} \p_j (\psi \ne j\big\vert W_1, \ldots, W_n) \right],
\end{align*}
where the infimum is over all possible tests. Moreover, as the right side of the above display does not depend on $\hat f$, we have
\begin{equation}\label{eq:final_condi}
\inf_{\hat f} \sup_{f \in \M_L} \p_f\big(\|\hat{f}-f\|_2 \ge A c_0^{-2} n^{-2/5}\big)\ge  \E_{W_1, \ldots, W_n}  \left[\inf_{\psi} \frac{1}{M+1}\sum_{j=0}^{M} \p_j (\psi \ne j\big\vert W_1, \ldots, W_n) \right].
\end{equation}
Let $P_j^*$ denote the  joint  distribution of $Z_1,\ldots, Z_n$  (conditional on $W_1, \ldots, W_n$) for $f= f_{j}$ (in~\eqref{eq:1dimConvex}). Let us assume that there exists an $\alpha \in (0, 1/8)$ (that does not depend on $W_1, \ldots, W_n$)  such that 
\begin{equation}\label{eq:FanoCond}
\frac{1}{M+1} \sum_{j=1}^{M}K(P^*_j, P^*_0) \le\alpha \log M\quad \text{for all}  \quad W_1,\ldots, W_n,
\end{equation}
 where $K(Q^{*}, P^*)$ denotes Kullback-Leibler divergence between the conditional distributions $Q^*$ and $P^*$.
 Then by Fano's Lemma (see e.g.,~\citet[Corollary 2.6]{MR2724359}), we have
 \begin{align}\label{eq:bound_x}
 \begin{split}
 \inf_{\psi} \frac{1}{M+1}\sum_{j=0}^{M} \p_j (\psi \ne j\big\vert W_1, \ldots, W_n)  &= \inf_\psi \frac{1}{M+1}\sum_{j=0}^{M}P_j^*\big(\psi \neq j\big) \\
 &\ge \frac{\log(M+1)-\log 2}{\log M} -\alpha > 0,
 \end{split}
 \end{align}
for $M$ such that $\log M \ge (1 - \alpha)^{-1}\log 2$. Note that $A_0$ and $c_0$ are constant. Thus combining~\eqref{eq:final_condi} and~\eqref{eq:bound_x}, we have that 
\[\inf_{\hat f} \sup_{f \in \M_{L_0}} \p_{f} \big[n^{4/5}\|\hat{f}-f\|_2^2 \ge A ^2 c_0^{-4}\big] \ge \frac{\log(M+1)-\log 2}{\log M} -\alpha > 0.\]
\textbf{Construction of the $M+1$ hypotheses.} In the following, we complete the proof by constructing $f_0, \ldots, f_M \in \M_{L_0}$ that satisfy~\eqref{eq:prop-1&2f_j} and~\eqref{eq:FanoCond}.
Let $f_0$ be any function in $\M_{L_0}$ that satisfies 
\begin{equation}\label{eq:f_0def}
  0< \kappa_1 \le f_0''(t)\le \kappa_2 <\infty,\qquad\text{for all}\quad  t\in[a,b],
\end{equation}
where $0<a <b<1$ and $\kappa_1$ and $\kappa_2$ are two arbitrary constants.   Note that  $f_0(x) = L_0x^2/2$ will satisfy~\eqref{eq:f_0def} with $a = 0, b = 1$ and $\kappa_1 = \kappa_2 = L_0$. However in the following proof, we keep track of $a, b, \kappa_1,$ and $\kappa_2$.\footnote{The final result with the ``general'' constants can be easily used to establish a ``local'' minimax rate lower bound for convex functions  satisfying~\eqref{eq:f_0def}; see Section~5 and Theorem~5.1 of~\cite{MR3405621} }. Next we construct $f_1, \ldots, f_M$.
 Recall that  $m = 8 \log M/ \log 2$. For $i=0, \ldots,m$, let $t_i := a+ (b-a)i/m$. For $1\le i\le m$, let $\alpha_i: [0,1]\to \R$ define the following affine function
\begin{equation}\label{eq:alpha_def}
\alpha_i(x) := f_0(t_{i-1}) + \frac{f_0(t_i)- f_0(t_{i-1})}{t_i-t_{i-1}} (x-t_{i-1})\qquad\text{for } x \in [0,1].
\end{equation}
Note that $(\cdot,\alpha_i(\cdot))$ is straight line through $(t_{i-1},f_0(t_{i-1}))$ and $(t_{i},f_0(t_i))$. For each $\tau= (\tau_1, \ldots, \tau_m) \in \{0,1\}^m,$ let us define \[
  f_\tau(x) := \max\Big(f_0(x), \max_{i:\tau_i=1}\alpha_i(x)\Big)\qquad \text{for } x\in [0,1].\footnote{The above construction is borrowed from Section 3.2 of~\cite{MR3405621}.}
\]
 As $f_\tau$ is a pointwise maximum of $L$-Lipschitz convex functions, $f_\tau$ is itself a $L$-Lipschitz convex function. Moreover  we have 
\begin{equation}\label{eq:alpha_diff}
f_\tau(x) =\begin{cases}\alpha_i(x)& \text{ if } \tau_i=1\\
f_0(x) & \text{ if } \tau_i=0.\end{cases}  \quad \text{for}\quad  x\in[t_{i-1}, t_i].
\end{equation}
% In Section~3.2 of \cite{MR3405621}, the authors show that $\sup_{x\in[0,1]} |f_\tau(x)-f_0(x)| \le (b-a)^2 \kappa_2 / (8 m^2).$ Thus we have that $ \|f_\tau(x)-f_0(x)\|_2 \le \delta.$
 We will next show that for $\tau, \tau' \in \{0,1\}^m$, the distance between $f_\tau$ and $f_{\tau'}$ can be  bounded from below (up to constant factors) by $\rho(\tau, \tau'):= \sum_{i}\{\tau_i\neq \tau_i'\}$. 
Observe that by~\eqref{eq:alpha_diff}, we have that 
\begin{equation}\label{eq:dist_lowerbound}
\|f_\tau-f_\tau'\|^2_2 = \sum_{i:\tau_i\neq\tau_i'} \|f_0-\max(f_0, \alpha_i)\|^2_2 \ge \rho(\tau, \tau') \min_{1\le i \le m} \|f_0-\max(f_0, \alpha_i)\|^2_2. 
\end{equation}
We will now find a lower bound for $\|f_0-\max(f_0, \alpha_i)\|^2_2$.  Since $\alpha_i(x) \ge f_0(x)$ for $x\in [t_{i-1}, t_i]$ and  $\alpha_i(x) \le f_0(x)$ for $x\notin [t_{i-1}, t_i]$, we have that
\begin{align}\label{eq:bound_1}
\begin{split}
\|f_0-\max(f_0, \alpha_i)\|^2_2&=  \int_{t_{i-1}}^{t_i}(f_0(x)-\alpha_i(x))^2 dx\\
&\ge\frac{\kappa_1^2}{4}\int_{t_{i-1}}^{t_i}\big[(x-t_{i-1})(t_i- x) \big]^2dx\\
&=\frac{\kappa_1^2}{120} (t_i-t_{i-1})^5= \frac{\kappa_1^2}{120} \frac{(b-a)^5}{m^5},
\end{split}
\end{align}
where the first inequality follows from the fact that for every $x \in [t_{i-1},t_i]$, there exists $t_x\in [t_{i-1},t_i]$ such that 
\[ |f_0(x)- \alpha_i(x)| = \frac{1}{2} (x-t_{i-1}) (t_i-x) f_0''(t_x) \ge \frac{\kappa_1}{2} (x-t_{i-1}) (t_i-x).\]
Note that the bound in~\eqref{eq:bound_1} does not depend on $i$. Thus from~\eqref{eq:dist_lowerbound}, we have that 
\begin{equation}\label{eq:lowerbounddiff}
\|f_\tau-f_\tau'\|_2  \ge \frac{\kappa_1}{11} \frac{(b-a)^{5/2}}{m^{5/2}}\sqrt{\rho(\tau, \tau')}.
\end{equation}
Since $m = 8 \log M/ \log 2$, by Varshamov-Gilbert lemma (Lemma 2.9 of~\citet[Page 104]{MR2724359}) we have that there exists a set $\{\tau^{(0)}, \ldots,\tau^{(M)}\}\subset \{0,1\}^m$ such that $\tau^{(0)} =(0, \ldots, 0)$ and  $\rho(\tau^{(k)}, \tau^{(j)}) \ge m/8$ for all $0\le k< j\le M.$ Further, recall that $f_{\tau^{(0)}}$ is $f_0$ by definition.  Thus if  we define $f_j := f_{\tau^{(i)}}$  for all $1\le i \le M$, then $f_0,\ldots,f_M$ satisfy~\eqref{eq:prop-1&2f_j}. 

We will now show that $P^*_0, \ldots, P^*_M$ satisfy~\eqref{eq:FanoCond}. Let us fix $W_1, \ldots, W_n$. Let $p^*_j$ denote the joint density with respect to the Lebesgue measure on $\R^n$. Since $\xi_1,\ldots, \xi_n$ are Gaussian random variables with mean $0$ and variance $\sigma^2$, we have that
\[ p^*_j(u_1,\ldots, u_n) =\Pi_{i=1}^n \phi_\sigma(u_i -f_j(W_i))\quad \text{and}\quad  p^*_0(u_1,\ldots, u_n) =\Pi_{i=1}^n \phi_\sigma(u_i -f_0(W_i)),\] where $\phi_\sigma$ is the density (with respect to the Lebesgue measure) of a mean zero Gaussian random variable with variance $\sigma^2.$ Thus by equation (2.36) of~\citet[Page 94]{MR2724359}, we have that 
\begin{align*}
  K(P^*_j, P^*_0)  \le \frac{1}{2 \sigma^2} \sum_{i=1}^n (f_0(W_i)- f_j(W_i))^2.
  \end{align*}  Note that for any $1\le k\le M$ and $0\le i\le m$, 
  \[|f_0(x)-f_k(x)| \le |f_0(x) -\alpha_i(x)|\qquad \text{for } x\in [t_{i-1}, t_i].\]  For every  $j \in \{1,\ldots, M\} $, we have 
  \begin{align*}
  K(P^*_j, P^*_0)  &\le \frac{1}{2 \sigma^2} \sum_{i=1}^n (f_0(W_i)- f_j(W_i))^2\\
  &\le \frac{1}{2 \sigma^2} \sum_{k=1}^m \sum_{W_i \in [t_{k-1}, t_k]} (f_0(W_i)- \alpha_k(W_i))^2\\
  &\le  \frac{\kappa_2^2(b-a)^4}{128 m^4 \sigma^2}\sum_{k=1}^m \sum_{W_i \in [t_{k-1}, t_k]} 1 \\
  &=  \frac{\kappa_2^2(b-a)^4}{128 m^4 \sigma^2}\text{Card}\{i: W_i\in [a, b]\}\\
  &\le \frac{ \kappa_2^2(b-a)^4}{128 m^4 \sigma^2} n,
  \end{align*}  where the third inequality holds since for every $x \in [t_{i-1},t_i]$, there exists $t_x\in [t_{i-1},t_i]$ such that 
\[ |f_0(x)- \alpha_i(x)| = \frac{1}{2} (x-t_{i-1}) (t_i-x) f_0''(t_x) \le \frac{\kappa_2}{2} (x-t_{i-1}) (t_i-x) \le \frac{\kappa_2}{8} (t_{i}- t_{i-1})^2= \frac{\kappa_2}{8 m^2} (b-a)^2.\]
Recall that $ n= m^5 c_0^{-5}$ and $m = 8 \log M/ \log 2$, thus
\begin{align*}
\frac{1}{M+1} \sum_{j=1}^{M}K(P^*_j, P^*_0) &\le \frac{ \kappa_2^2(b-a)^4}{128 m^4 \sigma^2} n \le \frac{ \kappa_2^2(b-a)^4}{128 \sigma^2 c_0^5} m \le  \frac{ \kappa_2^2(b-a)^4}{16 \sigma^2 c_0^5 \log 2} \log M.
\end{align*}
Let us fix \begin{equation}\label{eq:c_0_def}
c_0= \left[\frac{\kappa_2^2(b-a)^4 }{\sigma^2 \log 2 }\right]^{1/5},
\end{equation}
then we have that
\[\frac{1}{M+1} \sum_{j=1}^{M}K(P^*_j, P^*_0) \le \frac{1}{16}\log M.
\]
Thus $f_0, \ldots, f_M$  satisfy~\eqref{eq:prop-1&2f_j} and~\eqref{eq:FanoCond}. \qedhere
\end{proof}

\section{Proof of existence of $\check{m}_{{L}}$ and $\check{\theta}_{{L}}$} \label{app:proof:existanceCLSE}
\begin{prop}\label{thm:existanceCLSE}
The minimizer in~\eqref{eq:CLSE} exists.
\end{prop}
\begin{proof}

We consider the estimator
\[
(\check{m}_n, \check{\theta}_n) = \argmin_{(m,\theta)\in \M_{L}\times \Theta}Q_n(m, \theta).
\]
Fix $\theta\in \Theta$ and $n\ge 1$. For $m_1, m_2\in \M_L$, let 
\[d^*_n(m_1, m_2) := \sqrt{\frac{1}{n}\sum_{i=1}^n \big(m_1(\theta^\top X_i) -m_2(\theta^\top X_i)\big)^2}.
\] Observe that  $m \in \M_L \mapsto \sqrt{Q_n(m, \theta)}$ is a coercive continuous convex function (with respect to the topology induced by $d^*_n(\cdot, \cdot)$) on a convex domain.
Thus for every $\theta\in \Theta$, the global minimizer of $m \in \M_L \mapsto Q_n(m, \theta)$ exists. Let us define

  \begin{equation}\label{eq:t_theta_def}
  m_{\theta}:= \argmin_{m\in\M_L}Q_n(m,\theta)\;\text{ and }\; T(\theta) :=Q_n(m_\theta,\theta).
  \end{equation}

Observe  that $\check{\theta}_n:= \argmin_{\theta \in \Theta} T(\theta).$ As $\Theta$ is a compact set, the existence of the minimizer $\theta\mapsto T(\theta)$ will be established if we can show that $T(\theta)$ is a continuous function on $\Theta$. We will now prove that $\theta \mapsto T(\theta)$ is a continuous function. But first we will show that  for every $\theta \in \Theta$,  $\|m_\theta\|_\infty\le C$, where the constant $C$ depends only on $\{(X_i,Y_i)\}_{i=1}^n, L,$ and $T.$  Observe that $\sum_{i=1}^n(Y_i - m_\theta(\theta^\top X_i))^2 \le \sum_{i=1}^nY_i^2$ and the constant function $0$  belongs to $\M_L$. Thus
\begin{align*}
\sum_{i=1}^n \left[m_\theta(\theta^\top X_i)\right]^2 &\le 2\sum_{i=1}^n Y_im_\theta(\theta^\top X_i) \le 2\left(\sum_{i=1}^n Y_i^2\right)^{1/2}\left(\sum_{i=1}^n \left[m_\theta(\theta^\top X_i)\right]^2\right)^{1/2}.
\end{align*}
Hence, we have $|m_\theta(\theta^\top X_1)| \le 2 \sqrt{\sum_{i=1}^n Y_i^ 2}.$ As $m_\theta$ is uniformly $L$-Lipschitz, we have that for any $t\in D$, \[
|m_{\theta}(t)| \le |m_\theta(\theta^\top X_1)| + L |t-\theta^\top X_1| \le \sqrt{4 \sum_{i=1}^n Y_i^ 2} + L T=:C. \]
 As $C$ does not depend on $\theta$, we have that $\sup_{\theta \in \Theta} \|m_{\theta}\|_{\infty} \le C.$ As a first step of proving $\theta\mapsto T(\theta)$ is continuous, we will show that the class of functions %{ to do for paper }
\[\{ \theta\mapsto Q_n(m,\theta) : m \in\M_L,\;\|m\|_{\infty}\le C\}\] is uniformly equicontinuous. Observe that for $\theta, \eta\in\Theta$, we have
\begin{align*}
 n|Q_n(m,\theta) - Q_n(m,\eta)| &= \left|\sum_{i=1}^n (Y_i - m(\theta^{\top}X_i))^2 - \sum_{i=1}^n (Y_i - m(\eta^{\top}X_i))^2\right|\\
 &= \left|\sum_{i=1}^n (m(\eta^{\top}X_i) - m(\theta^{\top}X_i))(2Y_i - m(\theta^{\top}X_i) - m(\eta^{\top}X_i))\right|\\
&\le \sum_{i=1}^n |m(\eta^{\top}X_i) - m(\theta^{\top}X_i)|\times|2Y_i - m(\theta^{\top}X_i) - m(\eta^{\top}X_i)|\\
 &\le L\sum_{i=1}^n |\eta^{\top}X_i - \theta^{\top}X_i|\times 2 \left(|Y_i|+C\right)\\
 &\le 2nLT\left(\max_i|Y_i| + C\right)|\theta - \eta|.
 \end{align*}
 Thus, we have that
\begin{equation}\label{eq:q_equi}
\sup_{\left\{m \in\M_L : \;\|m\|_{\infty}\le C\right\} }|Q_n(m,\theta) -Q_n(m,\eta)| \le C_3 |\theta-\eta|,
\end{equation}
where $C_3$ is a constant depending only on $\{Y_i\}_{i=1}^{n}$ and $C$. Next we show that $|T(\theta)-T(\eta)|\le 2 C_3 |\theta-\eta|$. Recall that $T(\theta)= Q_n(m_{\theta}, \theta)$. By~\eqref{eq:t_theta_def}, we have
\[
 Q_n(m_{\theta}, \theta)- Q_n(m_{\theta}, \eta)= T(\theta)- Q_n(m_{\theta}, \eta)  \le T(\theta)-T(\eta)
 \]
  and \[
T(\theta)-T(\eta)  \le Q_n(m_{\eta}, \theta)-T(\eta)= Q_n(m_{\eta}, \theta)-Q_n(m_{\eta},\eta).
 \]
Thus
\begin{align*}
 |T(\theta)-T(\eta)| \le |Q_n(m_{\eta}, \theta)-Q_n(m_{\eta}, \eta) |+ |Q_n(m_{\theta}, \theta)- Q_n(m_{\theta}, \eta)|\le 2C_3 |\theta-\eta|.
 \end{align*}

\end{proof}
\section{Maximal inequalities for heavy-tailed multiplier processes} % (fold)
\label{sec:maximal_inequalities_for_heavy_tailed}

% subsection maximal_inequalities_for_heavy_tailed (end)
In this section, we collect some maximal inequalities for multiplier processes with heavy-tailed heteroscedastic multipliers. These are useful for verifying some steps in the proof of semiparametric efficiency. The standard tools from empirical process theory (see e.g.,~\cite{VandeGeer00,VdVW96}) require either bounded or sub-Gaussian/sub-exponential multipliers (Lemmas 3.4.2--3.4.3 of~\cite{VdVW96}). The main ideas in the proofs of the these results are: (i) employ a truncation device on the (heavy-tailed) errors and apply the Hoffmann-J{\o}rgensen's inequality to control the remainder (see Lemma~\ref{lem:TruncatedMaximal}); (ii) use generic chaining to obtain maximal inequalities on the truncated (bounded) empirical process (see Lemma~\ref{lem:MixedTailMaximal}; also see~\cite[Theorem 3.5]{Dirksen} and~\cite[Theorem~2.2.23]{MR3184689}).

%The above results (Lemma ???) are derived based on the results of \cite{Dirksen}. %First, we reduce the maximal inequality with unbounded errors to a maximal inequality with bounded errors.
\begin{lemma}\label{lem:TruncatedMaximal}
Suppose that $\{(\eta_i,X_i)\}_{i=1}^n$ are i.i.d.~observations from $\R \times \rchi$ with $X_i \sim P_X$. Define
\[
C_{\eta} := 8\mathbb{E}\left[\max_{1\le i\le n}|\eta_i|\right],\qquad\mbox{and}\qquad \overline{\eta} := \eta\mathbbm{1}_{\{|\eta| \le C_{\eta}\}}.
\] Let $\mathcal{F}$ be a class of bounded real-valued functions on $\rchi$ such that $\sup_{f\in\mathcal{F}}\norm{f}_{\infty} \le \Phi$.  Then
\begin{equation}\label{eq:UnboundedToBounded1}
\mathbb{E}\left[\sup_{f\in\mathcal{F}}\left|\mathbb{G}_n\left[\eta f\right]\right|\right] \le \mathbb{E}\left[\sup_{f\in\mathcal{F}}\left|\mathbb{G}_n\left[\overline{\eta} f\right]\right|\right] + \frac{2\Phi C_{\eta}}{\sqrt{n}}.
\end{equation}
\end{lemma}
\begin{proof}
This lemma is similar to Lemma S.1.4 of~\cite{kuchibhotla2018moving}. As $\eta = \overline{\eta} + (\eta - \overline{\eta})$, by the triangle inequality,
\begin{equation*}
|\mathbb{G}_n[\eta f]| \le |\mathbb{G}_n[\overline{\eta} f]| + |\mathbb{G}_n[(\eta - \overline{\eta})f]|.
\end{equation*}
Thus, we have 
\begin{align}\label{eq:FirstIneqHJ}
\mathbb{E}\big[\sup_{f\in\mathcal{F}}|\mathbb{G}_n[\eta f]|\big]&\le \mathbb{E}\left[\sup_{f\in\mathcal{F}}|\mathbb{G}_n[\overline\eta f]|\right] +\mathbb{E}\left[\sup_{f\in\mathcal{F}}|\mathbb{G}_n[(\eta -\overline{\eta})f]|\right].
\end{align}
We will first simplify the second term on the right of the above inequality. Let $R_1, R_2, \ldots, R_n$ be $n$ i.i.d.~Rademacher random variables\footnote{A Rademacher random variable takes value $1$ and $-1$ with probability $1/2$ each.} independent of $\{(\eta_i, X_i), 1\le i\le n\}$. 
Using symmetrization  (Corollary 3.2.2 of \cite{Gine16}), we have that \[
  \mathbb{E}\left[\sup_{f\in\mathcal{F}}|\mathbb{G}_n[(\eta -\overline{\eta})f]|\right] \le  2 \sqrt{n} \mathbb{E}\left[\sup_{f\in\mathcal{F}}|\mathbb{P}_n[R(\eta -\overline{\eta})f]|\right].
\]
Observe that for any $f\in\mathcal{F}$,
\begin{align}\label{eq:FirstBoundHJ}
\sup_{f\in\mathcal{F}} |\mathbb{P}_n\left[R(\eta - \overline{\eta})f\right]| &= \sup_{f\in\mathcal{F}} |\mathbb{P}_n\left[R\eta\mathbbm{1}_{\{|\eta| > C_{\eta}\}}f\right]| \le \frac{\Phi}{n}\sum_{i=1}^n |\eta_i|\mathbbm{1}_{\{|\eta_i| > C_{\eta}\}}.
\end{align}
%Thus,
%\begin{equation}\label{eq:FirstBoundHJ}
%\sup_{f\in\mathcal{F}}\left|\mathbb{G}_n\left[R\eta\mathbbm{1}_{\{|\eta| > C_{\eta}\}}\right]\right| \le \frac{\Phi}{\sqrt{n}}\sum_{i=1}^n |\eta_i|\mathbbm{1}_{\{|\eta_i| > C_{\eta}\}}.
%\end{equation}
Also, note that
\begin{align*}
%\mathbb{P}\left(\sup_{f\in\mathcal{F}}\left|\mathbb{G}_n\left[R\eta\mathbbm{1}\{|\eta| > C_{\eta}\}f\right]\right| > 0\right) &\le
\mathbb{P}\left(\sum_{i=1}^n |\eta_i|\mathbbm{1}_{\{|\eta_i| > C_{\eta}\}} > 0\right)
&\le \mathbb{P}\left(\max_{1\le i\le n}|\eta_i| > C_{\eta}\right) \le \frac{\mathbb{E}\left[\max_{1\le i\le n}|\eta_i|\right]}{C_{\eta}}\le \frac{1}{8}
\end{align*}
where the last inequality follows from the definition of $C_{\eta}.$ Hence by Hoffmann-J{\o}rgensen's inequality (Proposition 6.8 of \cite{LED91} with $t_0 = 0$), we get
\begin{equation}\label{eq:SecondBoundHJ}
\mathbb{E}\left[\sum_{i=1}^n |\eta_i|\mathbbm{1}_{\{|\eta_i| > C_{\eta}\}}\right] \le 8\mathbb{E}\left[\max_{1\le i\le n}|\eta_i|\right] = C_{\eta}.
\end{equation}
Combining inequalities~\eqref{eq:FirstBoundHJ} and~\eqref{eq:SecondBoundHJ}, it follows that
\[
\mathbb{E}\left[\sup_{f\in\mathcal{F}}\left|\mathbb{P}_n\left[R\eta\mathbbm{1}_{\{|\eta| > C_{\eta}\}} f\right]\right|\right] \le \frac{\Phi C_{\eta}}{n}.
\]
Substituting this bound in~\eqref{eq:FirstIneqHJ}, we get
\[
\mathbb{E}\left[\sup_{f\in\mathcal{F}}\left|\mathbb{G}_n\left[\eta f\right]\right|\right] \le \mathbb{E}\left[\sup_{f\in\mathcal{F}}\left|\mathbb{G}_n\left[\overline{\eta} f\right]\right|\right] + \frac{2\Phi C_{\eta}}{\sqrt{n}}.\qedhere
\]
\end{proof}
%To control the empirical process with bounded errors in~\eqref{eq:UnboundedToBounded}, we will use Theorem 3.5 of~\cite{Dirksen}. A similar result can be found in~\cite[Theorem~2.2.23]{MR3184689}.
\begin{lemma}\label{lem:MixedTailMaximal}
Suppose that $\{(\eta_i,X_i)\}_{i=1}^n$ are i.i.d.~observations from $\R \times \rchi$ with $X_i \sim P_X$ such that 
\begin{equation}\label{eq:eta_cond}
 \mathbb{E}\left[\bar{\eta}^2|X\right] \le \sigma^2_{{\eta}} \qquad P_X \text{ almost every } X, \quad \text{and} \quad \mathbb{P}(|\bar{\eta}| > C_{\eta}) = 0,
\end{equation}
for some constant $C_{\eta}$. Let $\mathcal{F}$ be a class of bounded real-valued functions on $\rchi$ such that
\begin{equation}\label{eq:cond_Lema_chain_1}
\sup_{f\in\mathcal{F}}\norm{f}_{\infty} \le \Phi,\qquad\sup_{f\in\mathcal{F}}\norm{f} \le \kappa,\,\qquad\mbox{and}\qquad \log N(\nu, \mathcal{F}, \norm{\cdot}_{\infty}) \le \Delta\nu^{-\alpha},
\end{equation}
for some constant $\Delta$ and $\alpha \in (0, 1)$, where $\norm{f}^2 := \int_{\rchi} f^2(x)dP_X(x)$ and $N(\nu, \mathcal{F},\|\cdot\|_\infty)$ is the $\nu$-covering number of $\mathcal{F}$ in the $\|\cdot\|_\infty$ metric (see Section 2.1.1 of~\cite{VdVW96} for its formal definition). Then
\begin{equation}\label{eq:Maximal}
\mathbb{E}\left[\sup_{f\in\mathcal{F}}\left|\mathbb{G}_n[\bar{\eta} f]\right|\right] \le 2\sigma_{\eta}\kappa + \frac{c_2\sqrt{2\Delta}\sigma_{\eta}(2\kappa)^{1 - \alpha/2}}{1 - \alpha/2} + \frac{c_1 2\Delta C_{\eta}(2\Phi)^{1 - \alpha}}{\sqrt{n}(1 - \alpha)},
\end{equation}
where $c_1$ and  $c_2$ are universal constants.
\end{lemma}
\begin{proof} 
Define the process $\{S(f):\,f\in\mathcal{F}\}$ by
$S(f) := \mathbb{G}_n\left[\bar{\eta} f(X)\right].$
For any two functions $f_1, f_2\in\f$,
\begin{equation}\label{eq:BernsteinBound}
|\bar{\eta}(f_1 - f_2)(X)| \le C_{\eta}\norm{f_1 - f_2}_{\infty},
\end{equation}
and
\begin{equation}\label{eq:BernsteinVar}
\mbox{Var}(\bar{\eta}(f_1 - f_2)) \le \mathbb{E}\left[\bar{\eta}^2(f_1 - f_2)^2(X)\right] \le \sigma^2_{{\eta}}\norm{f_1 - f_2}^2.
\end{equation}
Since
\[
|S(f_1) - S(f_2)| = \left|\mathbb{G}_n\left[\bar{\eta}(f_1 - f_2)(X)\right]\right|,
\] and for all $m\ge2$, we have
\[ \E \left[ \big|\bar{\eta}(f_1-f_2) -\E(\bar{\eta}(f_1-f_2))\big|^m\right] \le (2C_{{\eta}} \|f_1-f_2\|_{\infty})^{m-2} \mbox{Var}(\bar{\eta}(f_1-f_2)),\]
 Bernstein's inequality (Theorem 1 of \cite{Geer13}) implies that
\begin{align*}
\mathbb{P}&\left(|S(f_1) - S(f_2)| \ge \sqrt{t}d_2(f_1, f_2) + td_1(f_1, f_2)\right)\le 2\exp(-t),
\end{align*}
for all $t\ge 0$, where
\begin{align*}
d_1(f_1, f_2) := 2C_{\eta}\norm{f_1 - f_2}_{\infty}/\sqrt{n},\qquad\mbox{and}\qquad
d_2(f_1, f_2) := \sqrt{2}\sigma_{\eta}\norm{f_1 - f_2}.
\end{align*}
Hence by Theorem 3.5 and inequality (2.3) of \cite{Dirksen}, we get
\begin{equation}\label{eq:DirksenBound}
\begin{split}
\mathbb{E}\left[\sup_{f\in\mathcal{F}}|S(f)|\right] &\le 2\sup_{f\in\mathcal{F}}\mathbb{E}\left|S(f)\right| + c_2\int_0^{2\sqrt{2}\sigma_{\eta}\kappa} \sqrt{\log N(u, \mathcal{F}, d_2)}du\\
&\quad + c_1\int_0^{4C_{\eta}\Phi/\sqrt{n}}\log N(u, \mathcal{F}, d_1)du,
\end{split}
\end{equation}
for some universal constants $c_1$ and $c_2$. It is clear that $\mathbb{E}\left[\mathbb{G}_n[\bar{\eta} f(X)]\right] = 0$ and so,
\begin{align*}
\mathbb{E}\big[\left|S(f)\right|\big] &\le \sqrt{\mbox{Var}(S(f))} = \sqrt{\mbox{Var}(\mathbb{G}_n[\bar{\eta} f(X)])}\\ &\le \sqrt{\mbox{Var}(\bar{\eta} f(X))} \le \sigma_{\eta}\norm{f} \le \sigma_{\eta}\kappa.
\end{align*}
Thus,
\begin{equation}\label{eq:Part1}
\sup_{f\in\mathcal{F}}\mathbb{E}\big[|S(f)|\big] \le \sigma_{\eta}\kappa.
\end{equation}
To bound the last two terms of~\eqref{eq:DirksenBound}, note that
\begin{align*}
N(u, \mathcal{F}, d_2) &= N\left(\frac{u}{\sqrt{2}\sigma_{\eta}}, \mathcal{F}, \norm{\cdot}\right) \le N\left(\frac{u}{\sqrt{2}\sigma_{\eta}}, \mathcal{F}, \norm{\cdot}_{\infty}\right),\\% \le \exp(\Delta\left(\frac{\sqrt{2}\sigma_{\eta}}{\nu}\right)^{\alpha}),\\
N(u, \mathcal{F}, d_1) &= N\left(\frac{u\sqrt{n}}{2C_{\eta}}, \mathcal{F}, \norm{\cdot}_{\infty}\right).% \le \exp\left(\Delta\left(\frac{2C_{\eta}}{\sqrt{n}u}\right)^{\alpha}\right).
\end{align*}
Thus by~\eqref{eq:cond_Lema_chain_1}, we get
\begin{equation}\label{eq:Part2}
\begin{split}
\int_0^{2\sqrt{2}\sigma_{\eta}\kappa} \sqrt{\log N(u, \mathcal{F}, d_2)}du &= \int_0^{2\sqrt{2}\sigma_{\eta}\kappa} \sqrt{\log N\left(\frac{u}{\sqrt{2}\sigma_{\eta}}, \mathcal{F}, \norm{\cdot}_{\infty}\right)}du\\
&= \int_0^{2\sqrt{2}\sigma_{\eta}\kappa} \sqrt{\Delta}(\sqrt{2}\sigma_{\eta})^{\alpha/2}\frac{1}{u^{-\alpha/2}}du\\
&= \sqrt{\Delta}(\sqrt{2}\sigma_{\eta})^{\alpha/2}\frac{(2\sqrt{2}\sigma_{\eta}\kappa)^{1 - \alpha/2}}{(1 - \alpha/2)}\\
&= \frac{\sqrt{2\Delta}\sigma_{\eta}(2\kappa)^{1 - \alpha/2}}{1 - \alpha/2},
\end{split}
\end{equation}
and
\begin{equation}\label{eq:Part3}
\begin{split}
\int_0^{4C_{\eta}\Phi/\sqrt{n}}\log N(u, \mathcal{F}, d_1)du &= \int_0^{4C_{\eta}\Phi/\sqrt{n}} \log N\left(\frac{u\sqrt{n}}{2C_{\eta}}, \mathcal{F}, \norm{\cdot}_{\infty}\right)du\\
&= \int_0^{4C_{\eta}\Phi/\sqrt{n}} \Delta\left(\frac{2C_{\eta}}{\sqrt{n}}\right)^{\alpha}\frac{1}{u^{\alpha}}du\\
&= \Delta\left(\frac{2C_{\eta}}{\sqrt{n}}\right)^{\alpha}\left(\frac{4C_{\eta}\Phi}{\sqrt{n}}\right)^{1 - \alpha}\frac{1}{1 - \alpha}\\
&= \frac{2\Delta C_{\eta}(2\Phi)^{1 - \alpha}}{\sqrt{n}(1 - \alpha)}.
\end{split}
\end{equation}
Substituting inequalities~\eqref{eq:Part1}, \eqref{eq:Part2} and~\eqref{eq:Part3} in the bound~\eqref{eq:DirksenBound}, we get
\[
\mathbb{E}\left[\sup_{f\in\mathcal{F}}\left|\mathbb{G}_n[\bar{\eta} f]\right|\right] \le 2\sigma_{\eta}\kappa + \frac{c_2\sqrt{2\Delta}\sigma_{\eta}(2\kappa)^{1 - \alpha/2}}{1 - \alpha/2} + \frac{c_1 2\Delta C_{\eta}(2\Phi)^{1 - \alpha}}{\sqrt{n}(1 - \alpha)}.
\]
\end{proof}

Combining Lemmas~\ref{lem:TruncatedMaximal} and~\ref{lem:MixedTailMaximal} we get the following theorem. We will use the following result in the next section to prove Theorem~\ref{thm:rate_m_theta_CLSE}.
\begin{thm}\label{thm:MaximalMoment}
Suppose that $\{(\eta_i,X_i)\}_{i=1}^n$ are i.i.d.~observations from $\R \times \rchi$ with $X_i \sim P_X$ such that 
\begin{equation}\label{eq:eta_cond_thm}
 \mathbb{E}(\eta|X)=0, \quad\text{and}\quad \mbox{Var}(\eta|X) \le \sigma^2_{\eta}, \quad P_X \text{ almost every } X.
 \end{equation}
Let $\mathcal{F}$ be a class of bounded measurable functions on $\rchi$ such that
\begin{equation}\label{eq:cond_Lema_chain}
\sup_{f\in\mathcal{F}}\norm{f}_{\infty} \le \Phi,\quad\sup_{f\in\mathcal{F}}\norm{f} \le \kappa,\,\quad\mbox{and}\quad \log N(\nu, \mathcal{F}, \norm{\cdot}_{\infty}) \le \Delta\nu^{-\alpha},
\end{equation}
for some constant $\Delta$ and $0 < \alpha < 1$, where $\norm{f}^2 := \int_{\rchi} f^2(x)dP_X(x)$.
Then
\[
\mathbb{E}\left[\sup_{f\in\mathcal{F}}|\mathbb{G}_n\left[\eta f\right]|\right] \le 2\sigma_{\eta}\kappa + \frac{k_2\sqrt{2\Delta}\sigma_{\eta}(2\kappa)^{1 - \alpha/2}}{1 - \alpha/2} + \frac{k_1 2\Delta C_{\eta}(2\Phi)^{1 - \alpha}}{\sqrt{n}(1 - \alpha)} + \frac{2\Phi C_{\eta}}{\sqrt{n}},
\]
where $k_1, k_2$ are universal constants and $C_{\eta} := 8\mathbb{E}\left[\max_{1\le i\le n}|\eta_i|\right].$ 
In particular if $\mathbb{E}\left[|\eta|^q\right] < \infty$, then $C_{\eta} \le 8n^{1/q}\norm{\eta}_q$.
\end{thm}
\begin{proof}
By Lemma~\ref{lem:TruncatedMaximal},
\[
\mathbb{E}\left[\sup_{f\in\mathcal{F}}\left|\mathbb{G}_n\left[\eta f\right]\right|\right] \le \mathbb{E}\left[\sup_{f\in\mathcal{F}}\left|\mathbb{G}_n\left[\overline{\eta} f\right]\right|\right] + \frac{2\Phi C_{\eta}}{\sqrt{n}},
\]
where $|\overline{\eta}| \le C_{\eta}$ with probability 1 and $\mathbb{E}\left[\overline{\eta}^2|X\right] \le \mathbb{E}[\eta^2|X] \le \sigma^2_{\eta}.$
Since $\overline{\eta}$ is bounded by $C_\eta$ and $\mathbb{E}[\overline{\eta}^2|X] \le \sigma_\eta^2$, the result follows by an application of Lemma~\ref{lem:MixedTailMaximal}.
\end{proof}
\subsection{Maximal inequality for heavy-tailed errors via classical tools} % (fold)
\label{sub:maximal_inequality_for_heavy_tailed_errors_via_classical_tools}

% subsection maximal_inequality_for_heavy_tailed_errors_via_classical_tools (end)
Note that the previous results require a bound on $N(\nu, \mathcal{F}, \norm{\cdot}_{\infty})$.  However, such a bound can be hard to obtain for certain function classes. The following result provides a  maximal inequality when we only have a bound on $N_{[\,]}(\nu, \mathcal{F}, \norm{\cdot})$; here $\norm{\cdot}$ denotes the $L_2$ norm.
\begin{lemma}\label{lem:Maximal342}
Suppose that $\{(\eta_i,X_i)\}_{i=1}^n$ are i.i.d.~observations from $\R \times \rchi$ with $X_i \sim P_X$ such that 
\begin{equation}
 \mathbb{E}(\eta|X)=0, \quad\text{and}\quad \mbox{Var}(\eta|X) \le \sigma^2_{\eta}, \quad P_X \text{ almost every } X.
 \end{equation}
Let $\mathcal{F}$ be a class of bounded measurable functions on $\rchi$ such that
 $\|f\| \le \delta$ and $\norm{f}_{\infty} \le \Phi$ for every $f\in\mathcal{F}$. 
Then
\begin{equation}\label{eq:OldMaximal_newform}
\mathbb{E}\left[\sup_{f\in\mathcal{F}}|\mathbb{G}_n[{\eta}f]|\right] \lesssim \sigma_{\eta}J_{[\,]}(\delta, \mathcal{F}, \norm{\cdot})\left(1 + \frac{\sigma_{\eta}J_{[\,]}(\delta, \mathcal{F}, \norm{\cdot})\Phi C_{\eta}}{\delta^2\sqrt{n}}\right)+ \frac{2\Phi C_{\eta}}{\sqrt{n}},
\end{equation}
where $C_{\eta} := 8\mathbb{E}\left[\max_{1\le i\le n}|\eta_i|\right]$ and for any class of functions $\mathcal{F}$, $J_{[\,]}$ (the entropy integral) is defined as 
\be \label{eq:def_J}
J_{[\,]}(\delta, \mathcal{F}, \norm{\cdot}) := \int_0^{\delta}\sqrt{1 + \log N_{[\,]}(\nu, \mathcal{F}, \norm{\cdot})}d\nu.
\ee
\end{lemma}
\begin{proof}
Set
$\overline{\eta} := \eta\mathbbm{1}_{\{|\eta| \le C_{\eta}\}}.$
By Lemma~\ref{lem:TruncatedMaximal}, we have
\begin{align*}
\mathbb{E}\left[\sup_{f\in\mathcal{F}}\left|\mathbb{G}_n\left[\eta f\right]\right|\right] &\le \mathbb{E}\left[\sup_{f\in\mathcal{F}}\left|\mathbb{G}_n\left[\overline{\eta} f\right]\right|\right] + \frac{2\Phi C_{\eta}}{\sqrt{n}}.
\end{align*}
Since $\norm{\overline{\eta} f}_{\infty} \le C_{\eta}\Phi$ and
\[
\mathbb{E}\left[\overline{\eta}^2f^2(X)\right] \le \mathbb{E}\left[\eta^2 f^2(X)\right] \le \mathbb{E}\left[\mbox{Var}(\eta|X)f^2(X)\right] \le \sigma_{\eta}^2\delta^2.
\]
Let $[f_1^L, f_1^U], \ldots, [f_{N_{\nu}}^L, f_{N_{\nu}}^U]$ form $\nu$-brackets of $\mathcal{F}$ with respect to the $\norm{\cdot}$-norm. Fix a function $f\in\mathcal{F}$ and let $[f_1^L, f_1^U]$ be the bracket for $f$. Then a bracket for $\overline{\eta}f$ is given by
\[
\left[f_1^L\overline{\eta}^+ - f_1^U\overline{\eta}^-, f_1^U\overline{\eta}^+ - f_1^L\overline{\eta}^-\right],
\]
and the $\norm{\cdot}$-width of this bracket is given by
\[
\norm{(f_1^U - f_1^L)|\overline{\eta}|} = \sqrt{\mathbb{E}\left[\overline{\eta}^2(f_1^U - f_1^L)^2(X)\right]} \le \sigma_{\eta}\norm{f_1^U - f_1^L} \le \sigma_{\eta}\nu.
\]
Hence
\[
N_{[\,]}(\sigma_{\eta}\nu, \overline{\eta}\mathcal{F}, \norm{\cdot}) \le N_{[\,]}(\nu, \mathcal{F}, \norm{\cdot}).
\]
Therefore, by Lemma 3.4.2 of \cite{VdVW96}, we have
\[
\mathbb{E}\left[\sup_{f\in\mathcal{F}}|\mathbb{G}_n[\overline{\eta}f]|\right] \lesssim \sigma_{\eta}J_{[\,]}(\delta, \mathcal{F}, \norm{\cdot})\left(1 + \frac{\sigma_{\eta}J_{[\,]}(\delta, \mathcal{F}, \norm{\cdot})\Phi C_{\eta}}{\delta^2\sqrt{n}}\right).\qedhere
\]
\end{proof}

\section{Proofs of results in Sections~\ref{sec:AsymRegFunEstimate} and~\ref{sec:SepPara}} % (fold)
\label{app:AsymCLSE_Proof}
To  find the rate of convergence of $\check{m}_{{L}}\circ\check{\theta}_{{L}}$, we apply Theorem~3.1 of~\cite{KuchiPatra19}. For this purpose, we need covering numbers for the class of uniformly Lipschitz convex functions. We do not know of such results without an additional uniform boundedness assumption. To accomplish this, we first prove that it is enough to consider the class of uniformly bounded, uniformly Lipschitz convex functions.
\begin{lemma}\label{lem:Upsilion_ep}
Under assumption \ref{aa2}, we have that $\|\check{m}_{{L}}\|_\infty=O_p(1).$
Moreover, for every $n \ge 1,$
\begin{equation}\label{eq:bound_m_check}
\p\Big(\check{m}_{{L}} \notin \M_{M_L', L}\mbox{ for some }L \ge L_0\Big) \le \frac{\sigma^2}{n},
\end{equation}
where
\begin{equation}\label{eq:M_ep}
M_L' :=L\diameter(D) + M_0 + 1.
\end{equation} and for any $M>0$, we define
\begin{equation}\label{eq:bounded_M}
\M_{M,L} := \{m \in \M_L :\, \|m\|_{\infty}\le M\}.
\end{equation}
\end{lemma}
\begin{proof}
Recall that
\[
(\check{m}_{{L}},\check{\theta}_{{L}}) := \argmin_{(m,\theta)\in\M_{L} \times\Theta}\,\frac{1}{n}\sum_{i=1}^n \{Y_i-m(\theta^{\top}X_i)\}^2
\]
For simplicity, we drop the subscript $n$ in the estimator $(\check{m}_{{L}},\check{\theta}_{{L}})$.
By definition, we have
\[
\sum_{i=1}^n(Y_i - \check{m}_L(\check{\theta}_L^\top X_i))^2 \le \sum_{i=1}^n(Y_i - m(\check{\theta}_L^\top X_i))^2,
\] 
for all $m\in \M_L.$ Since any constant function belongs to $\M_L$, for any fixed real $\kappa$, we have
\[
\sum_{i=1}^n(Y_i - \check{m}_L(\check{\theta}_L^\top X_i))^2 \le \sum_{i=1}^n(Y_i - \check{m}_L(\check{\theta}_L^\top X_i) + \kappa)^2.
\]
A simplification of the  above inequality gives us:
\begin{equation}\label{eq:zero_sum}
2\kappa\sum_{i=1}^n(Y_i - \check{m}_L(\check\theta_L^\top X_i)) + n\kappa^2 \ge 0,\;\mbox{for all }\kappa\quad\Rightarrow\quad \sum_{i=1}^n(Y_i - \check{m}_L(\check{\theta}_L^\top X_i)) = 0.
\end{equation}

Thus for any $t\in D$, we have
\begin{align*}
|\check{m}_L(t)|&\le \left|\check{m}_L(t) - \frac{1}{n}\sum_{j=1}^n \check{m}_L(\check{\theta}_L^\top X_j)\right| + \left|\frac{1}{n}\sum_{j=1}^n \check{m}_L(\check{\theta}_L^\top X_j)\right|\\
&\le \frac{1}{n}\sum_{j=1}^n \left|\check{m}_L(t) - \check{m}_L(\check{\theta}_L^\top X_j)\right| + \left|\frac{1}{n}\sum_{j=1}^n \{m_0({\theta}_0^\top X_j) + \epsilon_j\}\right| \quad \text{(by~\eqref{eq:zero_sum})}\\
&\le \frac{1}{n}\sum_{j=1}^n L|t - \check{\theta}_n^\top X_j| + \frac{1}{n}\sum_{j=1}^n|m_0({\theta}_0^\top X_j)| + \left|\frac{1}{n}\sum_{j=1}^n \epsilon_j\right|\\
&\le L\diameter(D) + M_0 + \left|\frac{1}{n}\sum_{j=1}^n \epsilon_j\right|,
\end{align*}
where $M_0$ is the upper bound on $m_0$; see \ref{aa1_new}. The third inequality in the above display is true because $\check{m}_L$ is $L$--Lipschitz. Therefore, 
\begin{equation}\label{eq:BoundDeterministic}
\norm{\check{m}_L}_{\infty} \le L\diameter(D) + M_0 + \left|\frac{1}{n}\sum_{i=1}^n \epsilon_i\right|,\quad\mbox{for all}\quad L \ge L_0.
\end{equation}
Now observe that
\begin{align*}
&\mathbb{P}\left(\norm{\check{m}_L}_{\infty} \ge M_0 + L\diameter(D) + 1\mbox{ for some $L \ge L_0$}\right)\\ {}&~\overset{(a)}{\le}~ \mathbb{P}\left(\left|\frac{1}{n}\sum_{i=1}^n \epsilon_i\right| \ge 1\right) ~\overset{(b)}{\le}~ \mathbb{E}\left[\left(\frac{1}{n}\sum_{i=1}^n \epsilon_i\right)^2\right] ~\overset{(c)}{\le}~ \frac{\sigma^2}{n},
\end{align*}
where inequality (a) follows from~\eqref{eq:BoundDeterministic}, (b) follows from Markov's inequality and (c) follows from~\ref{aa2}. Therefore, for all $n\ge 1$,
\[
\mathbb{P}\left(\check{m}_L \notin \mathcal{M}_{M'_L, L}\mbox{ for some }L \ge L_0\right) \le \frac{\sigma^2}{n}.\qedhere
\]
\end{proof}
The intuition for the use of Lemma~\ref{lem:Upsilion_ep} is as follows. Since $\check{m}_L$ belongs to $\mathcal{M}_{M'_L, L}$ with ``high'' probability, we get that
\[
\left(\check{m}_L, \check{\theta}_L\right) = \argmin_{(m,\theta)\in\mathcal{M}_{M'_L, L}\times\Theta}\,\frac{1}{n}\sum_{i=1}^n \left(Y_i - m\left(\theta^{\top}X_i\right)\right)^2~\qquad \text{with high probability.}
\]
This estimator can be easily studied because of the existence of covering number results for the function class $\mathcal{M}_{M, L}$. Define
\[
\h_{M, L} := \{  m\circ\theta-m_0\circ\theta_0 : (m ,\theta) \in \M_{M, L} \times \Theta\}.
\]
Then the following covering number result holds.
\begin{lemma}\label{ent10}
There exist a positive constant $c$ and $\nu_0$, such that, for every $M,L > 0$ and $\nu \le \nu_0(M + L\diameter(D))$
\begin{equation}\label{eq:entropY_H}
\log N(\nu,\mathcal{H}_{M,L},\|\cdot\|_{\infty})= \log N(\nu,\{m\circ\theta : (m, \theta) \in\M_{M,L}\times \Theta\},\|\cdot\|_{\infty}) \le \frac{\mathcal{K}_{M,L}} {\sqrt{\nu}},
\end{equation}
where
\begin{equation}\label{eq:DefMathcalK}
\mathcal{K}_{M,L} := c\left[(2M + 2L\diameter(D))^{1/2} + 2d(6LT)^{1/2}\right].
\end{equation}
\end{lemma}
\begin{proof}
To prove this lemma, we use the covering number for the class of uniformly bounded and uniformly Lipschitz convex functions obtained in \cite{ADBO}.
% \begin{lemma}[Theorem 3.2, \cite{ADBO}]\label{cent}
% Let $\mathcal{F}$ denote the class of real-valued convex functions defined on $[a,b]$ that are uniformly bounded in absolute value by $B_0$ and uniformly Lipschitz with constant $L$. Then there exist positive constants $c$ and $\nu_0$ such that for every $B_0,L>0$ and $b>a$, we have
% \[\log N(\nu,\mathcal{F},\|\cdot\|_{\infty})\le c\left(\frac{B_0+L(b-a)}{\nu}\right)^{1/2}
% \]
% for every $\nu\le\nu_0(B_0 + L(b-a))$.
% \end{lemma}
By Theorem 3.2 of \cite{ADBO} and Lemma 4.1 of \cite{Pollard90} for $\nu \in (0,1)$,  we have
\begin{align}\label{eq:entoprY_ind_lip}
\log N_{[\;]}(\nu,\M_{M,L},\|\cdot\|_{\infty})&\le c\left(\frac{M+L\diameter(D)}{\nu}\right)^{1/2},\\
\log N(\nu,\Theta,|\cdot|)&\le d \log\left(\frac{3}{\nu}\right),
\end{align}
where $c$ is a constant that depends only on $d.$
% as $|\ theta_1 - \theta_2| \le 2$ for all $\theta_1,\theta_2\in\Theta.$

Recall that $\sup_{x\in\rchi} |x| \le T$; see \ref{aa1}. Let $\{\theta_1,\theta_2,\ldots,\theta_p\}$ be a $\nu/(2LT)$-cover (with respect to the Euclidean norm) of $\Theta$ and $\{m_1,m_2,\ldots,m_q\}$ be a $\nu/2$-cover (with respect to the $\|\cdot\|_\infty$-norm) for $\M_{M,L}$. In the following we will show that the set of functions $\{m_i\circ\theta_j-m_0\circ\theta_0\}_{1\le i\le q, 1\le j\le p}$ form a $\nu$-cover for $\h_{M,L}$ with respect to the $\|\cdot\|_{\infty}$-norm. For any given $m\circ\theta-m_0\circ\theta_0\in\h_{M,L}$, we can get $m_i$ and $\theta_j$ such that $\|m-m_i\|_{\infty} \le \nu/2$ and $|\theta - \theta_j| \le \nu/(2LT).$ Therefore, for any $x\in\rchi$
\begin{align*}
|m(\theta^{\top}x) - m_i(\theta_j^{\top}x)|&\le|m(\theta^{\top}x) - m(\theta_j^{\top}x)| + |m(\theta_j^{\top}x) - m_i(\theta_j^{\top}x)|\\
&\le L|x||\theta - \theta_j| + \|m - m_i\|_{\infty} \le \frac{L|x|\nu}{2LT}+\frac{\nu}{2} \le \nu.
\end{align*}
Thus for $\nu \le \nu_0(M + L\diameter(D)),$
\[
\log N(\nu, \mathcal{H}_{M,L}\circ\Theta, \norm{\cdot}_{\infty}) \le c\left[\left(\frac{2M + 2L\diameter(D)}{\nu}\right)^{1/2} + d\log\left(\frac{6LT}{\nu}\right)\right].
\]
Hence, using $\log x\le 2\sqrt{x}$ for all $x > 0$,
\begin{align*}
\log N(\nu, \mathcal{H}_{M,L}(\delta), \norm{\cdot}_{\infty}) &\le c\left[\left(\frac{2M + 2L\diameter(D)}{\nu}\right)^{1/2} + 2d\left(\frac{6LT}{\nu}\right)^{1/2}\right]\\
&= \frac{c}{\sqrt{\nu}}\left[(2M + 2L\diameter(D))^{1/2} + 2d(6LT)^{1/2}\right],
\end{align*}
for some universal constant $c> 0$.
% The result now follows as the covering number is equal to the bracketing number for the  sup-norm.
\end{proof}
\subsection{Proof of Theorem~\ref{thm:rate_m_theta_CLSE}}\label{proof:rate_m_theta_CLSE}
In the following, we fix $n\ge 1$ and use $L$ to denote $L_n$. The proof will be an application of Theorem~3.1 of~\cite{KuchiPatra19}. However, the class of functions $\M_L\times\Theta$ is not uniformly bounded. Thus $\check{m}_{L} \circ\check{\theta}_L$ and  $\M_L\times\Theta$ do not satisfy the conditions of Theorem~3.1 of~\cite{KuchiPatra19}. To circumvent this, consider a slightly modified LSE:
\begin{equation}\label{eq:hatm}
(\hat{m}_{L},\hat{\theta}_{L}) := \argmin_{(m,\theta)\in\mathcal{F}}\,\frac{1}{n}\sum_{i=1}^n (Y_i-m(\theta^\top X_i))^2,
\end{equation}
 where $\mathcal{F}:= \M_{M'_{L}, L} \circ \Theta$ with $M'_{L}$ is defined in \eqref{eq:M_ep}. However, by Lemma~\ref{lem:Upsilion_ep}, we have that 
\[
\p\Big(\check{m}_{L} \circ\check{\theta}_L \not\equiv \hat{m}_{L} \circ\hat{\theta}_L \Big) = \p\Big(\check{m}_{L} \notin \M_{M_{L}', L}\Big) \le \frac{\sigma^2}{n},
\]
when $L \ge L_0$. Thus for any every $r_n\ge 0$ and $M\ge 0$, we have
\begin{align}\label{eq:split123}
\begin{split}
&\mathbb{P}\left(r_n\norm{\check{m}_{L}\circ\check{\theta}_{L} - m_0\circ\theta_0} \ge 2^M\right)\\ 
\le{}& \mathbb{P}\left(r_n\norm{\hat{m}_{L}\circ\hat{\theta}_{L} - m_0\circ\theta_0} \ge 2^M\right) +  \mathbb{P}\left(\hat{m}_{L}\circ\hat{\theta}_{L} \not\equiv \check{m}_{L}\circ\check{\theta}_{L}\right) \\
{}\le{}& \mathbb{P}\left(r_n\norm{\hat{m}_{L}\circ\hat{\theta}_{L} - m_0\circ\theta_0} \ge 2^M\right) + \frac{\sigma^2}{n}.
\end{split}
\end{align}

We will now apply Theorem~3.1~\cite{KuchiPatra19} $\mathcal{F} = \M_{M'_{L}, L} \circ \Theta$ and $\hat{m}_{L}\circ\hat{\theta}_{L}$.  Note that
\[
\log N(u, \mathcal{F}, \norm{\cdot}_{\infty}) \le \frac{\mathcal{K}_{M'_L,L}}{\sqrt{\nu}}, \qquad \sup_{f\in\mathcal{F}}\norm{f}_{\infty} \le M'_L, \quad \text{ and }\quad  \|f_0\| \le M_0,
\]
where $\mathcal{K}_{M'_L,L}=  c\left[(2M + 2L\diameter(D))^{1/2} + 2d(6LT)^{1/2}\right]$ for some universal constant $c>0$ (see~\eqref{eq:DefMathcalK}). Observe that by~\ref{aa2}, $\mbox{Var}(\epsilon|X) \le \sigma^2$ and $\mathbb{E}\big[|\epsilon|^q\big]$.  
% Thus $\epsilon$ and $\F$ satisfy assumptions $(\mathcal{E}_q)$ and~$(L_\infty)$ of~\cite{KuchiPatra19}. 

% { In fact, if $\rchi$ is a bounded interval in $\R$ and the functions in $\F$ are uniformly Lipschitz with Lipschitz constant $L$ then $\|F_\delta\|_{\infty} \le 2 L^{1/3} \delta^{2/3}$; see Lemma~7 of~\cite{shen1994convergence} for a proof of this.
% }
% \todo[inline]{We can't do this since $s$ depends on $P_X$. The arguments of Lemma 7 in~\cite{shen1994convergence} works only for Lebesgue like measures for finding $s$}
Thus the assumptions of~{Theorem~3.1}~\cite{KuchiPatra19} are satisfied with 
\begin{equation}\label{eq:notation_matching}
\Phi = M'_L \vee M_0 \le M'_L+M_0, \quad A =\mathcal{K}_{M'_L, L}, \quad \alpha =1/2, \quad \text{and} \quad K_q^q =\E(|\epsilon|^q). 
\end{equation}
Thus
\begin{equation}\label{eq:Tail_majorgeneral}
\mathbb{P}\left(r_n\norm{\hat{m}_{L}\circ\hat{\theta}_{L} - m_0\circ\theta_0} \ge 2^M\right)\le  \frac{C}{2^{qM}},
\end{equation}
where 
\begin{equation}\label{eq:rn_Major}
r_n:= \min \left\{  \frac{n^{2/5}}{(\mathcal{K}_{M'_L, L} (M'_L+M_0)^2)^{2/5}},  \frac{n^{1/2-1/2q}}{( M'_L+M_0)^{(3q+ 1)/(4q)}} \right\},\end{equation}
and $C$ is constant depending only on $K_q$, $\sigma,$ and $q$. Recall that $M_L'=L\diameter(D) + M_0 + 1,$ thus 
\begin{align}\label{eq:denoms}
\begin{split}
\big[\mathcal{K}_{M'_L, L} (M'_L+M_0)^2\big]^{2/5} &\asymp d^{2/5}L\quad \text{ and }\quad
( M'_L+M_0)^{(3q+ 1)/(4q)} \asymp L^{(3q+ 1)/(4q)}
\end{split}
\end{align}
where for any $a,b\in\R$, we say $a\asymp b$ if there exist constants $c_2 \ge c_1 > 0$ depending only on $\sigma, M_0, L_0,$ and $T$ such that $c_1 b \le a \le c_2b$.
Therefore by combining~\eqref{eq:split123},~\eqref{eq:rn_Major}, and~\eqref{eq:denoms}, we have that there exists a constant $\mathfrak{C}$ depending only on $\sigma, M_0, L_0, T,$ and $K_q$ and a constant $C$ depending only $K_q, \sigma,$ and $q$ such that for all $M\ge 0$ 
 \[
 \mathbb{P}\left( r_n'\norm{\check{m}_L\circ\check{\theta}_L - m_0\circ\theta_0} \ge \mathfrak{C}2^M\right) \le \frac{C}{2^{qM}}+ \frac{\sigma^2}{n}.
 \]
 where 
 \begin{equation}\label{eq:rnprime}
 r_n'= \min \left\{  \frac{n^{2/5}}{d^{2/5}L},  \frac{n^{1/2-1/2q}}{L^{(3q+ 1)/(4q)}} \right\}.
 \end{equation} Note that above finite sample bound depends on the parameters $m_0$ and $\theta_0$ and the joint distribution of $\epsilon$ and $X$ only through the constants $\sigma, M_0, L_0, T,$ and $K_q$. Thus we have that 
\[ \sup_{\theta_0, m_0,\epsilon, X}  \mathbb{P}\left( r_n'\norm{\check{m}_L\circ\check{\theta}_L - m_0\circ\theta_0} \ge \mathfrak{C}2^M\right) \le \frac{C}{2^{qM}}+ \frac{\sigma^2}{n},\]
where the supremum is taken over all joint distributions of $\epsilon$ and $X$ and parameters $m_0$ and $ \theta_0\in \Theta$ for which assumptions~\ref{aa1_new}--\ref{aa2} are satisfied with constants $\sigma, M_0, L_0, T,$ and  $K_q.$

\subsection{Proof of Theorem~\ref{thm:UniformLRate}}\label{proof:UniformLRate}
The theorem (Theorem~\ref{thm:UniformLRate_exact}) stated and proved below is a more precise version Theorem~\ref{thm:UniformLRate}. The following result provides tail bounds for the quantity of interest.  The auxiliary results used in the proof below are given in Section~\ref{sub:lemmas_used_in_the_proof_of_theorem_thm:uniformlrate}.
\begin{thm}\label{thm:UniformLRate_exact}
Under the assumptions of Theorem~\ref{thm:rate_m_theta_CLSE}, for any $M \ge 1$, and $n\ge 15$, there exists a universal constant $C > 0$ such that
\begin{align}\label{eq:uniform_exact_prob}
\begin{split}
&\mathbb{P}\left(\sup_{L_0 \le L \le nL_0}\,\varphi_n(L)\norm{\check{m}_L\circ\check{\theta}_L - m_0\circ\theta_0} \ge C2^{M+1} \sqrt{\log\log_2n}\right)\\ &\qquad\qquad\le \frac{256}{2^{2M+1} C^2\log\log_2n} + \frac{e}{2^M} + \frac{\sigma^2}{n},
\end{split}
\end{align}
where 
\begin{equation}\label{eq:DefVarphi}
\varphi_n(L) := \min\left\{\frac{n^{2/5}}{3K^{(1)}L}, \frac{n^{1/2 - 1/(2q)}}{\sqrt{2K^{(2)}L}}\right\}.
\end{equation}
Here $K^{(1)}$ and $K^{(2)}$ are constants defined as
\begin{equation}\label{eq:DefK1K2}
K^{(1)} := \max\left\{\Delta^2, \Delta^{5/4}\right\},\quad\mbox{and}\quad K^{(2)} := \norm{\epsilon}_q\max\left\{\Delta^2, \Delta^3\right\},
\end{equation}
where $\Delta$ is the following constant
\begin{equation}\label{eq:LambdaDef}
\Delta := \left(\frac{M_0 + 1}{L_0} + \diameter(D)\right)^{1/2} + d\sqrt{T} + \sqrt{\sigma/L_0}.
\end{equation}
In particular,
\begin{equation}\label{eq:UniformLRate1}
\sup_{L_0 \le L \le nL_0}\,\varphi_n(L)\norm{\check{m}_L\circ\check{\theta}_L - m_0\circ\theta_0} = O_p\left(\sqrt{\log\log n}\right).
\end{equation}

\end{thm}
\begin{proof}
By Lemma~\ref{lem:Upsilion_ep}, we know that for all $n\ge 1$,
\begin{equation}\label{eq:out_m}
\mathbb{P}\left(\check{m}_L\notin \mathcal{M}_{M'_L, L}\mbox{ for some }L \ge L_0\right) \le \frac{\sigma^2}{n},
\end{equation}
where $M'_L = M_0 + 1 + L\diameter(D)$ and $\mathcal{M}_{M'_L, L}$ denotes the set of all $L$-Lipschitz convex functions bounded by $M'_L$.
Let us first define the following class of functions, for any  $0 \le \delta_1 \le \delta_2$,
\[
\mathcal{H}_{L}(\delta_1, \delta_2) := \left\{m\circ\theta - m_0\circ\theta_0:\,(m,\theta)\in\mathcal{M}_{M'_L, L}\times\Theta,\;\delta_1 \le \norm{m\circ\theta - m_0\circ\theta_0} \le \delta_2\right\}.
\]
Also, define
\begin{equation}\label{eq:L_J_M_def}
\mathcal{L}_n := [L_0, nL_0], \; \mathcal{J}_n := \mathbb{N}\cap [1, \log_2 n], \;\text{and}\;  \mathbb{M}_n(f) := \frac{2}{n}\sum_{i=1}^n \epsilon_if(X_i) - \frac{1}{n}\sum_{i=1}^n f^2(X_i).
\end{equation}
We now bound the probability in~\eqref{eq:uniform_exact_prob}. Observe that by~\eqref{eq:out_m}, we have
\begin{align}
&\mathbb{P}\left(\sup_{L\in\mathcal{L}_n}\,\varphi_n(L)\norm{\check{m}_L\circ\check{\theta}_L - m_0\circ\theta_0} \ge \delta\right)\nonumber\\
\le{}& \mathbb{P}\left(\sup_{L\in\mathcal{L}_n}\,\varphi_n(L)\norm{\check{m}_L\circ\check{\theta}_L - m_0\circ\theta_0} \ge \delta,\,\check{m}_L\in\mathcal{M}_{M'_L, L}\mbox{ for all }L\in\mathcal{L}_n\right)\nonumber\\
&\quad+\mathbb{P}\left(\check{m}_L\notin\mathcal{M}_{M'_L, L}\mbox{ for some }L\in\mathcal{L}_n\right)\nonumber\\
\le{}& \mathbb{P}\left(\sup_{L\in\mathcal{L}_n}\,\varphi_n(L)\norm{\check{m}_L\circ\check{\theta}_L - m_0\circ\theta_0} \ge \delta,\,\check{m}_L\in\mathcal{M}_{M'_L, L}\mbox{ for all }L\in\mathcal{L}_n\right)\nonumber\\
&\qquad\quad + \frac{\sigma^2}{n}.\label{eq:first_split}
\end{align}
Recall that for any $L \ge L_0$, 
\begin{align*}
(\check{m}_L, \check{\theta}_L) &:= \argmin_{(m,\theta)\in\mathcal{M}_L\times\Theta}\,\frac{1}{n}\sum_{i=1}^n \left(Y_i - m\circ\theta(X_i)\right)^2\\
&= \argmin_{(m,\theta)\in\mathcal{M}_L\times\Theta}\,\frac{1}{n}\sum_{i=1}^n \left[\left(Y_i - m\circ\theta(X_i)\right)^2 - \left(Y_i - m_0\circ\theta_0(X_i)\right)^2\right]\\
&= \argmin_{(m, \theta)\in\mathcal{M}_L\times\Theta}\,-\frac{2}{n}\sum_{i=1}^n \epsilon_i(m\circ\theta - m_0\circ\theta_0)(X_i) + \frac{1}{n}\sum_{i=1}^n \left(m\circ\theta - m_0\circ\theta_0\right)^2(X_i).
\end{align*}
Hence, we have that 
$\mathbb{M}_n (\check{m}_L \circ \check\theta_L -m\circ\theta) \ge 0$ for all $L$; where $\mathbb{M}_n(\cdot)$ is defined in~\eqref{eq:L_J_M_def}. Thus for the first probability in~\eqref{eq:first_split}, note that
\begin{align*}
&\mathbb{P}\left(\sup_{L\in\mathcal{L}_n}\,\varphi_n(L)\norm{\check{m}_L\circ\check{\theta}_L - m_0\circ\theta_0} \ge \delta,\,\check{m}_L\in\mathcal{M}_{M'_L, L}\mbox{ for all }L\in\mathcal{L}_n\right)\\
={}& \mathbb{P}\left(\exists L\in\mathcal{L}_n: \check{m}_L\circ\check{\theta}_L - m_0\circ\theta_0 \in\mathcal{H}_L\left(\frac{\delta}{\varphi_n(L)}, \infty\right)\right)\\
={}& \mathbb{P}\left(\exists (L, f)\in \mathcal{L}_n\times\mathcal{H}_L\left(\frac{\delta}{\varphi_n(L)}, \infty\right):\,\mathbb{M}_n(f) \ge 0\right)\\
={}& \mathbb{P}\left(\exists  (j, f)\in \mathcal{J}_n \times \bigcup_{2^jL_0\le L\le 2^{j+1}L_0}\mathcal{H}_{L}\left(\frac{\delta}{\varphi_n(L)}, \infty\right):\,\mathbb{M}_n(f) \ge 0\right)\\
\overset{(a)}{\le}{}& \mathbb{P}\left(\exists  (j, f)\in \mathcal{J}_n \times \mathcal{H}_{2^{j+1}L_0}\left(\frac{\delta}{2\varphi_n(2^{j+1}L_0)}, \infty\right):\,\mathbb{M}_n(f) \ge 0\right)\\
={}& \mathbb{P}\left(\exists  (j,k, f)\in \mathcal{J}_n \times \{\mathbb{N} \cup \{0\}\}\times \mathcal{H}_{2^{j+1}L_0}\left(\frac{2^k\delta}{2\varphi_n(2^{j+1}L_0)}, \frac{2^{k+1}\delta}{2\varphi_n(2^{j+1}L_0)}\right):\,\mathbb{M}_n(f) \ge 0\right).
\end{align*}
Inequality (a) above follows from Lemma~\ref{lem:BigCupSingle}. Now define

\begin{equation}\label{eq:G_def}
\mathcal{G}_{j, k} := \mathcal{H}_{2^{j+1}L_0}\left(\frac{2^k\delta}{2\varphi_n(2^{j+1}L_0)}, \frac{2^{k+1}\delta}{2\varphi_n(2^{j+1}L_0)}\right).
\end{equation}
Then for all $f\in \mathcal{G}_{j,k}$, we have
\begin{equation}\label{eq:f_bound}
\frac{2^k\delta}{2\varphi_n(2^{j+1}L_0)} \le \|f\| \le \frac{2^{k+1}\delta}{2\varphi_n(2^{j+1}L_0)}.
\end{equation}
Thus
\begin{align}\label{eq:MtoG}
\begin{split}
\mathbb{M}_n(f) =& \frac{1}{\sqrt{n}}\left(2\mathbb{G}_n\left[\epsilon f\right]  - \mathbb{G}_n[f^2]\right) - \norm{f}^2 \\
\le& \frac{1}{\sqrt{n}}\left(2\mathbb{G}_n\left[\epsilon f\right]  - \mathbb{G}_n[f^2]\right) -\frac{2^{2k}\delta^2}{4\varphi_n^2(2^{j+1}L_0)}
\end{split}
\end{align}
and so,
\begin{align*}
&\mathbb{P}\left(\sup_{L\in\mathcal{L}_n}\,\varphi_n(L)\norm{\check{m}_L\circ\check{\theta}_L - m_0\circ\theta_0} \ge \delta,\,\check{m}_L\in\mathcal{M}_{M'_L, L}\mbox{ for all }L\in\mathcal{L}_n\right)\\
\le{}&\mathbb{P}\left(\max_{j\in \mathcal{J}_n}\max_{k \ge 0}\sup_{f\in\mathcal{G}_{j,k}} \frac{4\varphi_n^2(2^{j+1}L_0)\left(2\mathbb{G}_n[\epsilon f] - \mathbb{G}_n[f^2]\right)}{\sqrt{n}\,2^{2k}\delta^2} \ge 1\right).
\end{align*}
Since $\epsilon$ is unbounded, we will use a simple truncation method to split the above probability into two components. First define
\begin{equation}\label{eq:basic_defs}
\gamma_{j,\delta} := \frac{4\varphi_n^2(2^{j+1}L_0)}{\sqrt{n}\delta^2},\quad\bar\epsilon_{i} := \epsilon_i\mathbbm{1}_{\{|\epsilon_i|\le C_{\epsilon}\}},\quad\mbox{and}\quad \epsilon^*_{i} := \epsilon_i - \bar\epsilon_{i}, 
\end{equation}
where $C_{\epsilon} := 8\mathbb{E}\left[\max_{1\le i\le n}|\epsilon_i|\right]$. Since $\epsilon_i = \bar\epsilon_{i} + \epsilon^*_{i}$, we get
\[
\mathbb{G}_n\left[\epsilon f\right] = \mathbb{G}_n\left[\bar\epsilon  f\right] + \mathbb{G}_n\left[\epsilon^* f\right].
\]
Note that $\bar\epsilon $ is bounded while $\epsilon^*$ is unbounded. Observe that
\begin{align*}
&\mathbb{P}\left(\sup_{L\in\mathcal{L}_n}\,\varphi_n(L)\norm{\check{m}_L\circ\check{\theta}_L - m_0\circ\theta_0} \ge \delta,\,\check{m}_L\in\mathcal{M}_{M'_L, L}\mbox{ for all }L \in \mathcal{L}_n\right)\\
\le{}& \mathbb{P}\left(\max_{j\in \mathcal{J}_n}\max_{k \ge 0}\sup_{f\in\mathcal{G}_{j,k}} \frac{\gamma_{j,\delta}}{2^{2k}}\left(2\mathbb{G}_n[\bar\epsilon  f] -\mathbb{G}_n[f^2]\right) \ge \frac{1}{2}\right)+ \mathbb{P}\left(\max_{j\in \mathcal{J}_n}\max_{k \ge 0}\sup_{f\in\mathcal{G}_{j,k}} \frac{\gamma_{j,\delta}}{2^{2k}} \mathbb{G}_n[2\epsilon^* f]\ge \frac{1}{2}\right)\\
\le{}& \mathbb{P}\left(\max_{j\in \mathcal{J}_n}\max_{k \ge 0}\sup_{f\in\mathcal{G}_{j,k}} \frac{\gamma_{j,\delta}}{2^{2k}}\left(2\mathbb{G}_n[\bar\epsilon  f] - \mathbb{G}_n[f^2]\right) \ge \frac{1}{2}\right)+ 4 \mathbb{E}\left(\max_{j\in \mathcal{J}_n}\max_{k \ge 0}\sup_{f\in\mathcal{G}_{j,k}} \frac{\gamma_{j,\delta}}{2^{2k}} \mathbb{G}_n[\epsilon^* f]\right),
\end{align*}
where the last inequality above follows by Markov's inequality. Our goal is to find $\delta$ such that the above probability can be made small. To make the notation less tedious, let us define \begin{equation}\label{eq:T_def}
 T_{j,\delta} := \max_{k \ge 0}\sup_{f\in\mathcal{G}_{j,k}} \frac{\gamma_{j,\delta}}{2^{2k}}\mathbb{G}_n[2\bar\epsilon  f - f^2]. 
\end{equation}
By a simple union bound, we have 
\begin{align}
&\mathbb{P}\left(\sup_{L\in\mathcal{L}_n}\,\varphi_n(L)\norm{\check{m}_L\circ\check{\theta}_L - m_0\circ\theta_0} \ge \delta,\,\check{m}_L\in\mathcal{M}_{M'_L, L}\mbox{ for all }L \in \mathcal{L}_n\right)\nonumber\\
\le{}&\mathbb{P}\left(\max_{j\in \mathcal{J}_n} T_{j,\delta} \ge \frac{1}{2}\right)+ 2\mathbb{E}\left[\max_{j\in \mathcal{J}_n}\max_{k \ge 0}\sup_{f\in\mathcal{G}_{j,k}} \frac{\gamma_{j,\delta}}{2^{2k}}\mathbb{G}_n[\epsilon^* f]\right]\nonumber\\
\le{}&\sum_{j=1}^{\log_2 n}\mathbb{P}\left( T_{j,\delta} \ge {1}/{2}\right)+ 2\mathbb{E}\left[\max_{j\in \mathcal{J}_n}\max_{k \ge 0}\sup_{f\in\mathcal{G}_{j,k}} \frac{\gamma_{j,\delta}}{2^{2k}}\mathbb{G}_n[\epsilon^* f]\right]\label{eq:Main_proof}.
\end{align}
In Lemma~\ref{lem:Talagrand}, we provide a tail bound for $T_{j,\delta}$ (a supremum of bounded empirical process) using  Talagrand's inequality (Proposition  3.1 of \cite{Gine00}). Moreover, note that the expectation in the above display is a supremum of sum of $n$ independent unbounded stochastic process and by Hoffmann-J{\o}rgensen's inequality (Proposition 6.8 of~\cite{LED91}) we can bound the expectation by a constant multiple of the expectation of the maximum of the $n$ stochastic processes. We do this in Lemma~\ref{lem:Hoff_Jorg}.

To conclude the proof note that, if we fix $\delta= 2^{{M+1}}C \sqrt{\log \log_2 n}$ (for some $M>1$), then  by Lemmas~\ref{lem:Talagrand} and~\ref{lem:Hoff_Jorg}, we have that 
\begin{equation}\label{eq:bounds_1}
\mathbb{P}(T_{j,\delta} \ge 1/2) \le e/ (2^M \log_2 n) 
\end{equation}
and 
\begin{equation}\label{eq:bounds_2}
 4\mathbb{E}\left[\max_{j\in \mathcal{J}_n}\max_{k \ge 0}\sup_{f\in\mathcal{G}_{j,k}} \frac{\gamma_{j,\delta}}{2^{2k}}\mathbb{G}_n[\epsilon^* f]\right] \le \frac{256}{2^{2M+1} C^2 \log\log_2n},
 \end{equation} respectively.

 The proof is now complete since, by substituting the upper bounds~\eqref{eq:bounds_1} and~\eqref{eq:bounds_2} in~\eqref{eq:Main_proof} and combining the result with~\eqref{eq:first_split}, we get that 
\begin{align}
&\mathbb{P}\left(\sup_{L_0 \le L \le nL_0}\,\varphi_n(L)\norm{\check{m}_L\circ\check{\theta}_L - m_0\circ\theta_0} \ge 2^{{M+1}}C \sqrt{\log \log_2 n}\right)\\
% \le{}& \sum_{j=1}^{\log_2 n} \mathbb{P} (T_{j,\delta} \ge {1}/{2}) +\le \frac{C}{2^{2M+2} \log\log_2n}\\
  \le{}& \sum_{j=1}^{\log_2 n} \frac{e}{2^M \log_2n} + \frac{256}{2^{2M+1} C^2\log\log_2n}+\frac{\sigma^2}{n}\\
  \le{}& \frac{e}{2^M} +\frac{256}{2^{2M+1} C^2 \log\log_2n}+ \frac{\sigma^2}{n}.\qedhere
\end{align}
\end{proof}
\subsection{Lemmas used in the proof of Theorem~\ref{thm:UniformLRate}} % (fold)
\label{sub:lemmas_used_in_the_proof_of_theorem_thm:uniformlrate}

% subsection lemmas_used_in_the_proof_of_theorem_thm:uniformlrate (end)

The following two Lemmas provide basic properties about the rate $\varphi_n(L)$ and the function classes $\mathcal{H}_L(\delta_1, \delta_2)$ defined in the proof of Theorem~\ref{thm:UniformLRate_exact}.
\begin{lemma}\label{lem:PropertiesVarphi}
For any $n\ge 1$,
\begin{equation}\label{eq:qTermVarphi}
\sup_{L \ge L_0}\,\frac{L\varphi_n(L)}{n} \le \frac{1}{3n^{3/5}}\min\left\{\frac{1}{\Delta^2}, \frac{1}{\Delta^{5/4}}\right\}, 
\end{equation}
and
\begin{equation}\label{eq:ComplexityTermVarPhi}
\sup_{L \ge L_0}\,\frac{L\varphi_n^2(L)C_{\epsilon}}{n}\le 4\min\left\{\frac{1}{\Delta^2}, \frac{1}{\Delta^3}\right\}.
\end{equation}
\end{lemma}
\begin{proof}
From the definition of $\varphi_n(L)$, we get that
\begin{equation}
\varphi_n(L) \le \frac{n^{2/5}}{3K^{(1)}L}\quad\Rightarrow\quad \sup_{L}\frac{L\varphi_n(L)}{n} \le \frac{1}{3K^{(1)}n^{3/5}} \le \frac{1}{3n^{3/5}}\min\left\{\frac{1}{\Delta^2}, \frac{1}{\Delta^{5/4}}\right\},
\end{equation}
and
\[
\sup_L \frac{L\varphi_n^2(L)C_{\epsilon}}{n} \le \frac{C_{\epsilon}}{2n^{1/q}K^{(2)}} \le \frac{8\norm{\epsilon}_qn^{1/q}}{2n^{1/q}K^{(2)}} \le \frac{4\norm{\epsilon}_q}{K^{(2)}} \le 4\min\left\{\frac{1}{\Delta^2}, \frac{1}{\Delta^3}\right\}.
\]
\end{proof}
\begin{lemma}\label{lem:BigCupSingle}
For any $j \ge 0$ and any constant $C > 0$,
\[
\bigcup_{2^jL_0 \le L\le 2^{j+1}L_0}\mathcal{H}_L\left(\frac{C}{\varphi_n(L)}, \infty\right) \subseteq \mathcal{H}_{2^{j+1}L_0}\left(\frac{C}{2\varphi_n(2^{j+1}L_0)}, \infty\right).
\]
\end{lemma}
\begin{proof}
We will first prove a few inequalities of $\varphi_n(\cdot)$. Since $\varphi_n(\cdot)$ is nonincreasing and so, for all $2^jL_0 \le L \le 2^{j+1}L_0$,
\[
\varphi_n(2^jL_0) \ge \varphi_n(L) \ge \varphi_n(2^{j+1}L_0)\quad\Rightarrow\quad \frac{1}{\varphi_n(2^{j+1}L_0)} \ge \frac{1}{\varphi_n(L)} \ge \frac{1}{\varphi_n(2^jL_0)}.
\]
Also, note that
\begin{align}
\varphi_n(2^{j+1}L_0) &= \min\left\{\frac{n^{2/5}}{3K^{(1)}2^{j+1}L_0}, \frac{n^{1/2 - 1/(2q)}}{\sqrt{2K^{(2)}2^{j+1}L_0}}\right\}\nonumber\\
&\ge \frac{1}{2}\min\left\{\frac{n^{2/5}}{3K^{(1)}2^{j}L_0}, \frac{n^{1/2 - 1/(2q)}}{\sqrt{2K^{(2)}2^{j}L_0}}\right\},\nonumber\\
\Rightarrow \frac{1}{\varphi_n(2^jL_0)} &\ge \frac{1}{2\varphi_n(2^{j+1}L_0)}\quad\Rightarrow\quad \frac{1}{\varphi_n(L)} \ge \frac{1}{2\varphi_n(2^{j+1}L_0)}.\label{eq:VarPhiNIneq}
\end{align}
Also note that for $L \le 2^{j+1}L_0$,
\[
\mathcal{M}_{M'_L, L} \subseteq \mathcal{M}_{M'_{2^{j+1}L_0}, 2^{j+1}L_0}\quad\Rightarrow\quad \mathcal{H}_{L}\left(\frac{C}{\varphi_n(L)}, \infty\right)\subseteq\mathcal{H}_{2^{j+1}L_0}\left(\frac{C}{\varphi_n(L)}, \infty\right).
\]
Thus,
\begin{equation*}
\bigcup_{2^jL_0 \le L\le 2^{j+1}L_0}\mathcal{H}_L\left(\frac{C}{\varphi_n(L)}, \infty\right) \subseteq \bigcup_{2^jL_0 \le L\le 2^{j+1}L_0}\mathcal{H}_{2^{j+1}L_0}\left(\frac{C}{\varphi_n(L)}, \infty\right).
\end{equation*}
It is clear that for any $L > 0$ and for $\delta_1 \le \delta_2$, $\mathcal{H}_L\left(\delta_2, \infty\right) \subseteq \mathcal{H}_L\left(\delta_1, \infty\right),$ and combining this inequality with~\eqref{eq:VarPhiNIneq}, we get for any $L \le 2^{j+1}L_0$,
\[
\mathcal{H}_{2^{j+1}L_0}\left(\frac{C}{\varphi_n(L)}, \infty\right)\subseteq\mathcal{H}_{2^{j+1}L_0}\left(\frac{C}{2\varphi_n(2^{j+1}L_0)}, \infty\right).
\]
Therefore,
\[
\bigcup_{2^jL_0 \le L\le 2^{j+1}L_0}\mathcal{H}_L\left(\frac{C}{\varphi_n(L)}, \infty\right)\subseteq \mathcal{H}_{2^{j+1}L_0}\left(\frac{C}{2\varphi_n(2^{j+1}L_0)}, \infty\right).\qedhere
\]
\end{proof}
The following two Lemmas form an integral part in the proof of~\eqref{eq:bounds_1}.

\begin{lemma}\label{lem:Talagrand}
Recall $\gamma_{j,\delta}$ and $T_{j,\delta}$ defined in~\eqref{eq:basic_defs} and~\eqref{eq:T_def}, respectively. There exists a constant $C>1$ (depending only on $d$) such that 
\begin{align}\label{eq:ET_bound}
\begin{split}
{\delta^2\mathbb{E}\left[T_{j,\delta}\right]} &\le C\left[\frac{\Delta^2\delta}{3K^{(1)}n^{1/10}} + \frac{\Delta^{5/2}\delta^{3/4}}{(3K^{(1)})^{5/4}} + \frac{\Delta^3\norm{\epsilon}_q}{2K^{(2)}} + \frac{\Delta^{5/2}}{(3K^{(1)})^2n^{1/5}}\right]\\
&\le C \left[\delta n^{-1/10} + \delta^{3/4} +2\right],
\end{split}
\end{align}
\begin{align}\label{eq:wimpY_var}
\begin{split}
\sigma_j^2 &:= \max_{k \ge 0}\sup_{f\in\mathcal{G}_{j,k}} \mbox{Var}\left(\frac{\gamma_{j,\delta}}{2^{2k}}\mathbb{G}_n[2\bar\epsilon  f - f^2]\right) \le \frac{512n^{-1/5}}{9\delta^2},
\end{split}
\end{align}
and 
\begin{align}\label{eq:U_bound}
\begin{split}
U_{j} &:= \max_{k \ge 0}\sup_{f\in\mathcal{G}_{j,k}}\max_{1\le i\le n} \frac{1}{2^{2k}}\left|\bar\epsilon_{i}f(X_i) - f^2(X_i) - \mathbb{E}\left[\bar\epsilon_{i}f(X_i) - f^2(X_i)\right]\right|\\
&\le 2 C_{\epsilon}(2M_0 + 1 + 2^{j+1}L_0\diameter(D)) + 2(2M_0 + 1 + 2^{j+1}L_0\diameter(D))^2.
\end{split}
\end{align}
% Thus by  Talagrand's moment bounds for bounded empirical process, we have that
% \begin{align}\label{eq:Talagrand_bound}
% \begin{split}
% \norm{T_{j,\delta}}_p \le C\left[\mathbb{E}[T_{j,\delta}] + p^{1/2}\sigma_j + p\gamma_{j,\delta}U_{j}\right].
% \end{split}
% \end{align}
% Thus by~\eqref{eq:ET_bound}, \eqref{eq:wimpY_var}, and~\eqref{eq:U_bound}
Thus by  Talagrand's moment bounds for bounded empirical process
, we have
\begin{equation}\label{eq:TjTailBound_lemma}
\mathbb{P}\left(|T_{j,\delta}| \ge C\left[\frac{1}{\delta n^{1/10}} + \frac{1}{\delta^{5/4}} + \frac{\sqrt{t}}{\delta^2 n^{1/5}} + \frac{t}{\delta^2}\right]\right) \le e\exp(-t).
\end{equation}
Furthermore, choosing $\delta = 2^{{M+1}}C \sqrt{\log \log_2 n}$ and $t= \log(2^M \log_2 n)$, for $n \ge 15$ and $M \ge 1$, we have that 
\begin{equation}\label{eq:TjTailBoundFinal}
\mathbb{P}\left(|T_{j,\delta}| \ge 1/2\right) \le \mathbb{P}\left(|T_{j,\delta}| \ge C\left[\frac{1}{\delta n^{1/10}} + \frac{1}{\delta^{5/4}} + \frac{\sqrt{t}}{\delta^2 n^{1/5}} + \frac{t}{\delta^2}\right]\right) \le \frac{e}{2^M \log_2n}.
\end{equation}

\end{lemma}

\begin{proof}

 The main goal of the lemma is to prove~\eqref{eq:TjTailBound_lemma}. By Proposition 3.1 of \cite{Gine00}, we get for $p\ge 1,$
\begin{equation}\label{eq:MomentBoundEmpProc}
\big(\E|T_{j,\delta}|^p\big)^{1/p} \le K\left[\mathbb{E}[T_{j,\delta}] + p^{1/2}\sigma_j + pU_{j,p}\right],
\end{equation}
where $K$ is an absolute constant, 
\[
\sigma_j^2 = \max_{k \ge 0}\sup_{f\in\mathcal{G}_{j,k}} \mbox{Var}\left(\frac{\gamma_{j,\delta}}{2^{2k}}\mathbb{G}_n[2\bar\epsilon  f - f^2]\right),\quad \text{and}\quad
U_{j,p} := \frac{\gamma_{j,\delta}}{\sqrt{n}} \mathbb{E}\left[U_j^p\right]^{1/p}.
\]
In the following, we find upper bounds for $\E(T_{j,\delta}),$ $ \sigma_j$, and $U_{j,p}$. First up is $U_{j,p}$.  Note that~\eqref{eq:U_bound} is a simple consequence of the fact that $|\bar\epsilon_{i}| \le C_{\epsilon}$ and $\norm{f}_{\infty} \le 2M_0 + 1 + 2^{j+1}L_0\diameter(D)$ for $f\in\mathcal{G}_{j,k}$; see~\eqref{eq:basic_defs} and~\eqref{eq:G_def}. Thus for $1\le j \le \log_2 n$, we have that
\begin{align*}
U_{j,p} &\le \frac{2\gamma_{j,\delta}}{\sqrt{n}}\Big[C_{\epsilon}(2M_0 + 1 + 2^{j+1}L_0\diameter(D)) + (2M_0 + 1 + 2^{j+1}L_0\diameter(D))^2\Big]\\
&= \frac{2\gamma_{j,\delta}C_{\epsilon}(2M_0 + 1 + 2^{j+1}L_0\diameter(D))}{\sqrt{n}} + \frac{2\gamma_{j,\delta}(2M_0 + 1 + 2^{j+1}L_0\diameter(D))^2}{\sqrt{n}}\\
&\le \frac{2\gamma_{j,\delta}C_{\epsilon}2^{j+1}L_0}{\sqrt{n}}\left(\frac{2M_0 + 1}{L_0} + \diameter(D)\right) + \frac{2\gamma_{j,\delta}(2^{j+1}L_0)^2}{\sqrt{n}}\left(\frac{2M_0 + 1}{L_0} + \diameter(D)\right)^2\\
&\le \frac{2\gamma_{j,\delta}C_{\epsilon}2^{j+1}L_0}{\sqrt{n}}(2\Delta^2) + \frac{2\gamma_{j,\delta}(2^{j+1}L_0)^2}{\sqrt{n}}(4\Delta^4),
\end{align*}
where $\Delta$ is as defined in~\eqref{eq:LambdaDef}. Lemma~\ref{lem:PropertiesVarphi} and the definition of  $\gamma_{j,\delta}$, imply that
\[
\frac{\gamma_{j,\delta}C_{\epsilon}2^{j+1}L_0\Delta^2}{\sqrt{n}} = \frac{4\varphi_n^2(2^{j+1}L_0)2^{j+1}L_0C_{\epsilon}\Delta^2}{n\delta^2} \le \frac{4\Delta^2}{\delta^2}\sup_{L \ge L_0}\frac{L\varphi_n^2(L)C_{\epsilon}}{n} \le \frac{16}{\delta^2},
\]
and
\begin{align*}
\frac{\gamma_{j,\delta}(2^{j+1}L_0)^2\Delta^4}{\sqrt{n}} &\le \frac{4\varphi_n^2(2^{j+1}L_0)(2^{j+1}L_0)^2\Delta^4}{n\delta^2}\\ &\le \frac{4n\Delta^4}{\delta^2}\sup_{L\ge L_0}\frac{L^2\varphi_n^2(L)}{n^2} \le \frac{4n\Delta^4}{\delta^2}\frac{1}{9n^{6/5}\Delta^4} \le \frac{4}{9n^{1/5}\delta^2}.
\end{align*}
% the definition of $\varphi_n(L)$ that
% \begin{equation}\label{eq:ComplexityTermVarPhi}
% \varphi_n(L) \le \frac{n^{2/5}}{3K^{(1)}L}\quad\Rightarrow\quad \sup_{L}\frac{L\varphi_n(L)}{n} \le \frac{1}{3K^{(1)}n^{3/5}},
% \end{equation}
% and from~\eqref{eq:qTermVarphi},
% \[
% \sup_L \frac{L\varphi_n^2(L)C_{\epsilon}}{n} \le \frac{C_{\epsilon}}{2n^{1/q}K^{(2)}} \le \frac{8\norm{\epsilon}_qn^{1/q}}{2n^{1/q}K^{(2)}} \le \frac{4\norm{\epsilon}_q}{K^{(2)}}.
% \]
Substituting these two inequalities in the bound on $U_{j,p}$, we get
\begin{equation}\label{eq:TjBoundUjp}
U_{j,p} \le \frac{64}{\delta^2} + \frac{32n^{-1/5}}{9\delta^2} = \frac{32}{\delta^2}\left(2 + n^{-1/5}\right) \le \frac{96}{\delta^2}.
% \frac{32\norm{\epsilon}_q}{K^{(2)}\delta^2}\left(\frac{2M_0 + 1}{L_0} + \diameter(D)\right) + \frac{8}{9(K^{(1)})^2n^{6/5}\delta^2}\left(\frac{2M_0 + 1}{L_0} + \diameter(D)\right)^2%\left[\frac{64(2M_0 + 1)\norm{\epsilon}_q}{K^{(2)}L_0} + \frac{16(2M_0 + 1)^2}{9n^{6/5}(K^{(1)})^2L_0^2}\right]\frac{1}{\delta^2}.
\end{equation}
We will now prove~\eqref{eq:wimpY_var}. Recall that $\E(\epsilon^2|X)\le \sigma^2$. To bound $\sigma_j^2$, observe that for $f\in\mathcal{G}_{j,k}$,
\begin{align*}
\mbox{Var}\left(\mathbb{G}_n\left[2\bar\epsilon f - f^2\right]\right) &\le \mathbb{E}\left[\left(2\bar\epsilon_{i}f(X_i) - f^2(X_i)\right)^2\right]\\ &\le 8\mathbb{E}\left[\bar\epsilon ^2f^2(X_i)\right] + 2\mathbb{E}\left[f^4(X_i)\right]\\
&\le 8\mathbb{E}\left[\epsilon^2f^2(X_i)\right] + 2\norm{f}_{\infty}^2\mathbb{E}\left[f^2(X_i)\right]\\
&\le 8\left[\sigma^2 + (2M_0 + 1 + 2^{j+1}L_0\diameter(D))^2\right]\norm{f}^2\\ &\le 16\Delta^4\frac{2^{2k+2}(2^{j+1}L_0)^2\delta^2}{2\varphi_n^2(2^{j+1}L_0)} \le 32\Delta^4\frac{2^{2k}(2^{j+1}L_0)^2\delta^2}{\varphi_n^2(2^{j+1}L_0)}. 
\end{align*}
Substituting this in the definition of $\sigma_j^2$, we get
\begin{equation*}%\label{eq:SigmaBound}
\begin{split}
\sigma_j^2 &\le \max_{k \ge 0}\sup_{f\in\mathcal{G}_{j,k}}\,\frac{\gamma_{j,\delta}^2}{2^{4k}}\frac{2^{2k}(2^{j+1}L_0)^2\delta^2}{\varphi_n^2(2^{j+1}L_0)}(32\Delta^4)\\
&\le \max_{k\ge 0}\,\frac{16\varphi_n^4(2^{j+1}L_0)}{2^{2k}n\delta^4}\frac{(2^{j+1}L_0)^2\delta^2}{\varphi_n^2(2^{j+1}L_0)}(32\Delta^4)\\ &= \frac{512\varphi_n^2(2^{j+1}L_0)(2^{j+1}L_0)^2}{n\delta^2}\Delta^4 \le \frac{512n\Delta^4}{\delta^2}\sup_{L\ge L_0}\frac{L^2\varphi_n^2(L)}{n^2}.
\end{split}
\end{equation*}
Using~\eqref{eq:ComplexityTermVarPhi}, we get,
\begin{equation}\label{eq:TjSigmaBound}
\sigma_j^2 \le \frac{512n\Delta^4}{\delta^2}\frac{1}{9n^{6/5}\Delta^4} \le \frac{512n^{-1/5}}{9\delta^2}.%\frac{64}{9n^{1/5}(K^{(1)})^2}\left[\frac{\sigma^2}{L_0^2} + \left(\frac{2M_0 + 1}{L_0} + \diameter(D)\right)^2\right]\frac{1}{\delta^2}.
\end{equation}
To bound $\mathbb{E}\left[T_{j,\delta}\right]$, note that
\begin{equation}\label{eq:FirstBoundTj}
\begin{split}
\frac{1}{\gamma_{j,\delta}}\mathbb{E}\left[T_{j,\delta}\right] &~\le~ \sum_{k = 0}^{\infty} \frac{1}{2^{2k}}\mathbb{E}\left[\sup_{f\in\mathcal{G}_{j,k}}\,\mathbb{G}_n\left[2\bar\epsilon  f - f^2\right]\right]\\
&\,\le \sum_{k = 0}^{\infty} \frac{1}{2^{2k}}\mathbb{E}\left[\sup_{f\in\mathcal{G}_{j,k}}\left|\mathbb{G}_n\left[2\bar\epsilon f\right]\right|\right] + \sum_{k = 0}^{\infty} \frac{1}{2^{2k}}\mathbb{E}\left[\sup_{f\in\mathcal{G}_{j,k}}\left|\mathbb{G}_n\left[f^2\right]\right|\right].%,\\
% &\,\le 2\sum_{k = 0}^{\infty} \frac{1}{2^{2k}}\mathbb{E}\left[\sup_{f\in\mathcal{G}_{j,k}}\left|\mathbb{G}_n\left[\bar\epsilon f\right]\right|\right] + 2\sum_{k = 0}^{\infty}\frac{1}{2^{2k}}\mathbb{E}\left[\sup_{f\in\mathcal{G}_{j,k}}\left|\mathbb{P}_n\left[Rf^2\right]\right|\right],
\end{split}
\end{equation}
By symmetrization and contraction principles for independent Rademacher random variables $R_1, \ldots, R_n$, (see arguments  leading up to (3.175) in~\cite{Gine16}), we have 
\begin{equation}\label{eq:Contraction}
\mathbb{E}\left[\sup_{f\in\mathcal{G}_{j,k}}\left|\mathbb{G}_n\left[f^2\right]\right|\right] \le 8(2M_0 + 1 + 2^{j+1}L_0\diameter(D))\mathbb{E}\left[\sup_{f\in\mathcal{G}_{j,k}}\left|\mathbb{G}_n\left[Rf\right]\right|\right].
\end{equation}
Since for any $L \ge 0$ and $\beta > 0$,
\[
\sup_{f\in\mathcal{H}_L(0, \beta)}\norm{f}_{\infty} \le M_L' + M_0,\quad\mbox{and}\quad \sup_{f\in\mathcal{H}_L(0, \beta)}\norm{f} \le \beta,
\]
we have by Lemma~\ref{ent10}, for $\beta > 0$,
\[
\log N(\nu, \mathcal{H}_L(0, \beta), \norm{\cdot}_{\infty}) \le \frac{\mathcal{K}_{M'_L, L}}{\sqrt{\nu}},
\]
and by Lemma~\ref{lem:MixedTailMaximal},
\[
\mathbb{E}\left[\sup_{f\in\mathcal{H}_L\left(0, \beta\right)}\left|\mathbb{G}_n\left[\bar\epsilon f\right]\right|\right] \le 2\sigma\beta + \frac{c_2\sqrt{2}\mathcal{K}_{M'_L, L}^{1/2}\sigma}{3/4}(2\beta)^{3/4} + \frac{2c_1\mathcal{K}_{M'_L, L}C_{\epsilon}(2(M'_L + M_0))^{1/2}}{\sqrt{n}/2}.
\]
Here
\[
\mathcal{K}_{M'_L, L} := c\left[\left(2M'_L + 2L\diameter(D)\right)^{1/2} + 2d(6LT)^{1/2}\right],
\]
for some constant $c$ depending only on $d$. Similarly,
\[
\mathbb{E}\left[\sup_{f\in\mathcal{H}_L\left(0, \beta\right)}\left|\mathbb{G}_n\left[Rf\right]\right|\right] \le 2\beta + \frac{c_2\sqrt{2}\mathcal{K}_{M'_L, L}^{1/2}}{3/4}(2\beta)^{3/4} + \frac{2c_1\mathcal{K}_{M'_L, L}(2(M'_L + M_0))^{1/2}}{\sqrt{n}/2},
\]
Noting that for $\mathcal{G}_{j,k}\subseteq\mathcal{H}_{2^{j+1}L_0}\left(0, 2^k\delta/\varphi_n(2^{j+1}L_0)\right)$, we get that
\begin{equation}
\sum_{k = 0}^{\infty} \frac{1}{2^{2k}}\mathbb{E}\left[\sup_{f\in\mathcal{G}_{j,k}}\left|\mathbb{G}_n\left[\bar\epsilon f\right]\right|\right] \le \frac{3\sigma\delta}{\varphi_n(2^{j+1}L_0)} + \frac{{5c_2\mathcal{K}_{j}^{1/2}\sigma}}{\varphi_n^{3/4}(2^{j+1}L_0)} + \frac{16c_1\mathcal{K}_{j}C_{\epsilon}(M_j + M_0)^{1/2}}{\sqrt{n}},
\end{equation}
and
\begin{equation}
\sum_{k = 0}^{\infty} \frac{1}{2^{2k}}\mathbb{E}\left[\sup_{f\in\mathcal{G}_{j,k}}\left|\mathbb{G}_n\left[Rf\right]\right|\right] \le \frac{3\delta}{\varphi_n(2^{j+1}L_0)} + \frac{{5c_2\mathcal{K}_{j}^{1/2}}}{\varphi_n^{3/4}(2^{j+1}L_0)} + \frac{16c_1\mathcal{K}_{j}(M_j + M_0)^{1/2}}{\sqrt{n}},
\end{equation}
where $M_j := M'_{2^{j+1}L_0}$, and $\mathcal{K}_j := \mathcal{K}_{M'_{2^{j+1}L_0}, 2^{j+1}L_0}$. Substituting these inequalities in~\eqref{eq:FirstBoundTj}, we get
\begin{align*}
\frac{1}{\gamma_{j,\delta}}\mathbb{E}\left[T_{j,\delta}\right] &\le \left[2\sigma + 8(2M_0 + 1 + 2^{j+1}L_0\diameter(D))\right]\left(\frac{3\delta}{\varphi_n(2^{j+1}L_0)} + \frac{5c_2\mathcal{K}_{j}^{1/2}}{\varphi_n^{3/4}(2^{j+1}L_0)}\right)\\
&\quad+ \left[2C_{\epsilon} + 8(2M_0 + 1 + 2^{j+1}L_0\diameter(D))\right]\frac{16c_1\mathcal{K}_j(M_j + M_0)^{1/2}}{\sqrt{n}}.
\end{align*}

Now observing that
\[
\mathcal{K}_j \le c\left(2^{j+1}L_0\right)^{1/2}\left[\left(\frac{2M_0 + 2}{L_0} + 4\diameter(D)\right)^{1/2} + 2d\sqrt{6T}\right],
\]
and using Lemma~\ref{lem:PropertiesVarphi}, we get for some large constant $C>0$ that
\begin{align*}
&\frac{\sqrt{n}\delta^2\mathbb{E}[T_{j,\delta}]}{C}\\
&\,\le \Delta^2\delta\left\{\varphi_n(2^{j+1}L_0)2^{j+1}L_0\right\} + \Delta^{5/2}\delta^{3/4}\left\{\varphi_n(2^{j+1}L_0)2^{j+1}L_0\right\}^{5/4}\\ &\,\quad+ \Delta^{3}\left\{\varphi_n^2(2^{j+1}L_0)2^{j+1}L_0\right\}\norm{\epsilon}_qn^{1/q - 1/2} + \Delta^{5/2}\left\{\varphi_n(2^{j+1}L_0)2^{j+1}L_0\right\}^{2}n^{-1/2}\\
&\,\le \Delta^2\delta\frac{n^{2/5}}{3K^{(1)}} + \Delta^{5/2}\delta^{3/4}\frac{n^{1/2}}{(3K^{(1)})^{5/4}} + \Delta^3\frac{n^{1/2}\norm{\epsilon}_q}{2K^{(2)}} + \Delta^{5/2}\frac{n^{3/10}}{(3K^{(1)})^2}.
\end{align*}
% where
% \[
% \Delta := \left(\frac{M_0 + 1}{L_0} + \diameter(D)\right)^{1/2} + d\sqrt{T} + \sqrt{\sigma/L_0}.
% \]
Therefore, for $j \ge 1$,
\begin{equation}\label{eq:TjExpectationBound}
{\delta^2\mathbb{E}\left[T_{j,\delta}\right]} \le C\left[\frac{\Delta^2\delta}{3K^{(1)}n^{1/10}} + \frac{\Delta^{5/2}\delta^{3/4}}{(3K^{(1)})^{5/4}} + \frac{\Delta^3\norm{\epsilon}_q}{2K^{(2)}} + \frac{\Delta^{5/2}}{(3K^{(1)})^2n^{1/5}}\right].
\end{equation}
Using the definition of $\Delta$ and substituting inequalities~\eqref{eq:TjExpectationBound},~\eqref{eq:TjBoundUjp},~and \eqref{eq:TjSigmaBound} in~\eqref{eq:MomentBoundEmpProc}, we get for $p\ge 1$,
\begin{align*}
\frac{1}{K}\norm{\delta^2T_{j,\delta}}_p &\le C\left[\frac{\Delta^2\delta}{K^{(1)}n^{1/10}} + \frac{\Delta^{5/2}\delta^{3/4}}{(K^{(1)})^{5/4}} + \frac{\Delta^3\norm{\epsilon}_q}{K^{(2)}} + \frac{\Delta^{5/2}}{(K^{(1)})^2n^{1/5}}\right]\\
&\quad+ \frac{Cp^{1/2}}{n^{1/5}} + Cp.
\end{align*}
From the definitions~\eqref{eq:DefK1K2} of $K^{(1)}$ and $K^{(2)}$, we get for $p\ge 1$,
\begin{align*}
\norm{T_{j,\delta}}_p &\le C\left[\frac{1}{n^{1/10}\delta} + \frac{1}{\delta^{5/4}} + \frac{p^{1/2}}{n^{1/5}\delta^2} + \frac{p}{\delta^2}\right].%\\
% &\le C\left[\frac{1}{\delta} + \frac{1}{\delta^{5/4}} + \frac{p^{1/2}}{\delta^2} + \frac{p}{\delta^2}\right].
\end{align*}
Therefore, by Markov's inequality for any $t\ge 0,$
\begin{equation}\label{eq:TjTailBound}
\mathbb{P}\left(|T_{j,\delta}| \ge C\left[\frac{1}{n^{1/10}\delta} + \frac{1}{\delta^{5/4}} + \frac{t^{1/2}}{n^{1/5}\delta^2} + \frac{t}{\delta^2}\right]\right) \le e\exp(-t).
\end{equation}
Fix  $\delta = 2^{{M+1}}C \sqrt{\log \log_2 n}$ and $t= \log(2^M \log_2 n)$. Then for any $M\ge 1$ and $n \ge 15$,
\begin{align*}
C\left[\frac{1}{\delta n^{1/10}} + \frac{1}{\delta^{5/4}} + \frac{\sqrt{t}}{\delta^2 n^{1/5}} + \frac{t}{\delta^2}\right] &\le \frac{1}{2^{M+1}n^{1/10}\sqrt{\log\log_2n}}\\ &\quad+ \frac{1}{2^{5(M+1)/4}(\log\log_2n)^{5/8}}\\ &\quad+ \frac{\sqrt{M\log 2 + \log\log_2n} + (M\log 2 + \log\log_2n)}{2^{2(M+1)}\log\log_2n}\\% + \frac{}{2^{2(M+1)}C^2\log\log_2n}\right]
&\le \frac{1}{2^{M+1}n^{1/10}\sqrt{\log\log_2n}} + \frac{1}{2^{5(M+1)/4}(\log\log_2n)^{5/8}}\\ &\quad+ \frac{2(M\log 2 + \log\log_2n)}{2^{2(M+1)}\log\log_2n} \le \frac{1}{2}.
\end{align*}
Therefore, for $M \ge 1$ and $n \ge 15$,
\[
\mathbb{P}\left(|T_j| \ge \frac{1}{2}\right) \le \mathbb{P}\left(|T_{j,\delta}| \ge C\left[\frac{1}{n^{1/10}\delta} + \frac{1}{\delta^{5/4}} + \frac{t^{1/2}}{n^{1/5}\delta^2} + \frac{t}{\delta^2}\right]\right) \le \frac{e}{2^M\log_2n}.
\]
\end{proof}

\begin{lemma}\label{lem:Hoff_Jorg}
By an application of the Hoffmann-J{\o}rgensen's inequality, we get 
\begin{equation}
2\mathbb{E}\left[\max_{j\in \mathcal{J}_n}\max_{k \ge 0}\sup_{f\in\mathcal{G}_{j,k}} \frac{8\varphi_n^2(2^{j+1}L_0)\mathbb{G}_n[\epsilon^* f]}{\sqrt{n}\,2^{2k}\delta^2}\right] \le \frac{256}{\delta^2}.
\end{equation}

\end{lemma}
\begin{proof}
Note that quantity of interest is the $L_1$ norm of supremum of sum of $n$ independent stochastic process. Thus by Hoffmann-J{\o}rgensen's inequality, we can bound this expectation using the quantile of the supremum of the sum stochastic process and  the  $L_1$ norms of the maximum of the individual stochastic process. We first simplify the expectation. Note that 
\begin{align}\label{eq:hof_1}
\begin{split}
& 2\mathbb{E}\left[\max_{j\in \mathcal{J}_n}\max_{k \ge 0}\sup_{f\in\mathcal{G}_{j,k}} \frac{8\varphi_n^2(2^{j+1}L_0)\mathbb{G}_n[\epsilon^* f]}{\sqrt{n}2^{2k}\delta^2}\right]\\
\le{}& 2\mathbb{E}\left[\max_{j\in \mathcal{J}_n}\max_{k \ge 0}\sup_{f\in\mathcal{G}_{j,k}} \frac{8\varphi_n^2(2^{j+1}L_0)}{2^{2k}\delta^2}\left\vert\mathbb{P}_n[\epsilon^* f]- \mathbb{E}(\mathbb{P}_n[\epsilon^* f])\right\vert\right]\\
\stackrel{(\alpha)}{\le}{}& 4\mathbb{E}\left[\max_{j\in \mathcal{J}_n}\max_{k \ge 0}\sup_{f\in\mathcal{G}_{j,k}} \frac{8\varphi_n^2(2^{j+1}L_0)}{n2^{2k}\delta^2}\sum_{i=1}^n |\epsilon^*_{i}f(X_i)|\right]\\
{\le}{}& 4\mathbb{E}\left[\max_{j\in \mathcal{J}_n}\max_{k \ge 0}\sup_{f\in\mathcal{G}_{j,k}} \frac{8\varphi_n^2(2^{j+1}L_0)}{n2^{2k}\delta^2} \|f\|_{\infty}\sum_{i=1}^n  |\epsilon^*_{i}|\right]
\end{split}
\end{align}
where the inequality-$(\alpha)$ follows from Jensen's inequality.  Since  $\sup_{f\in\mathcal{G}_{j,k}} \norm{f}_{\infty} \le 2M_0 + 1 + 2^{j+1}L_0\diameter(D)$, we have that

\begin{align}
&4\mathbb{E}\left[\max_{j\in \mathcal{J}_n}\max_{k \ge 0}\sup_{f\in\mathcal{G}_{j,k}} \frac{8\varphi_n^2(2^{j+1}L_0)}{n2^{2k}\delta^2} \|f\|_{\infty}\sum_{i=1}^n  |\epsilon^*_{i}|\right]\nonumber\\
\le{}&4\mathbb{E}\left[\max_{j\in \mathcal{J}_n}\max_{k \ge 0} \frac{8\varphi_n^2(2^{j+1}L_0)(2M_0 + 1 + 2^{j+1}L_0\diameter(D))}{n2^{2k}\delta^2}\sum_{i=1}^n  |\epsilon^*_{i}|\right]\nonumber\\
\le{}&4\mathbb{E}\left[\max_{j\in \mathcal{J}_n} \frac{8\varphi_n^2(2^{j+1}L_0)(2M_0 + 1 + 2^{j+1}L_0\diameter(D))}{n\delta^2}\sum_{i=1}^n  |\epsilon^*_{i}|\right]\nonumber\\
\le{}&\mathbb{E}\left[\sum_{i=1}^n  |\epsilon^*_{i}|\right]\max_{j\in \mathcal{J}_n} \frac{32\varphi_n^2(2^{j+1}L_0)(2M_0 + 1 + 2^{j+1}L_0\diameter(D))}{n\delta^2}.\label{eq:hof_22}
\end{align}
 We will now bound each of terms in the product. First up is $\E(\sum_{i=1 }^{n} |\epsilon^*_{i}|)$. To apply proposition 6.8 of~\cite{LED91}, we need to find the upper $1/8$'th quantile of the sum. Note that
\begin{align}
&\mathbb{P}\left(\max_{I\le n}  \sum_{i=1}^I |\epsilon^*_{i}| \ge 0 \right) \le \mathbb{P}\left(\max_{i\le n} |\epsilon^*_{i}| \ge 0\right) \le  \mathbb{P}\left(\max_{i\le n} |\epsilon_{i}| \ge C_\epsilon\right) \le \frac{\E(\max_{i\le n} |\epsilon_{i}|)}{C_\epsilon}= \frac{1}{8}.
\end{align}
Thus by (6.8) of~\cite{LED91}, we have that 
\begin{align}\label{eq:hof_2}
\begin{split}
\mathbb{E}\left[\sum_{i=1}^n |\epsilon^*_{i}|\right]\le & 8\mathbb{E}\left[\max_{1\le i\le n}\left|\epsilon^*_{i}\right|\right] = C_{\epsilon}.
\end{split}
\end{align}
Thus combining~\eqref{eq:hof_1} and~\eqref{eq:hof_22}, we have that
\begin{align}\label{eq:const}
\begin{split}& 2\mathbb{E}\left[\max_{j\in \mathcal{J}_n}\max_{k \ge 0}\sup_{f\in\mathcal{G}_{j,k}} \frac{8\varphi_n^2(2^{j+1}L_0)\mathbb{G}_n[\epsilon^* f]}{\sqrt{n}2^{2k}\delta^2}\right]\\
\le{}& C_\epsilon\max_{j\in \mathcal{J}_n} \frac{32\varphi_n^2(2^{j+1}L_0)(2M_0 + 1 + 2^{j+1}L_0\diameter(D))}{n\delta^2}\\
\le{}& 32\left(\frac{2M_0 + 1}{L_0} + \diameter(D)\right)\max_{j\in \mathcal{J}_n} \frac{2^{j+1}L_0\varphi_n^2(2^{j+1}L_0)C_\epsilon}{n\delta^2}\\
\le{}& \frac{32}{\delta^2}\left(\frac{2M_0 + 1}{L_0} + \diameter(D)\right)\max_{L\in\mathcal{L}_n} \frac{L\varphi_n^2(L) C_\epsilon}{n}\\
\le{}& \frac{32}{\delta^2}\left(\frac{2M_0 + 1}{L_0} + \diameter(D)\right) 4\min\left\{\frac{1}{\Delta^2}, \frac{1}{\Delta^3}\right\}
\le{} \frac{256}{\delta^2},
\end{split}
\end{align}
where the last two inequalities follow from~\eqref{eq:ComplexityTermVarPhi} and~\eqref{eq:LambdaDef}, respectively.
\end{proof}
%%%%%%%%%%%%%%%%%%%%%%%%%%%%%%%%%%

%%%%%%%%%%%%%%%%%%%%%%%%%%%%%%%%%%
\subsection{Proof of Theorem \ref{thm:uucons}}\label{proof:thm:uucons}
Recall that $\mathcal{M}_L$ is a class of   equicontinuous functions defined on a closed and bounded set and $\Theta$ is a compact set.  Let $\{(m_n, \theta_n)\}$ be any sequence  in $\mathcal{M}_L\times \Theta$  such that $\{m_n\}$ is uniformly bounded.   Then, by  Ascoli-Arzel\`{a} theorem, there exists a subsequence $\{(m_{n_k}, \theta_{n_k})\}$, $\theta\in\Theta,$ and $m\in \mathcal{M}_L$ such that $|\theta_{n_k}-\theta| \rightarrow 0$ and $\|m_{n_k}-m\|_{D_0}\rightarrow 0.$
Now suppose that $\|m_n\circ\theta_n - m_0\circ\theta_0\|\rightarrow0$. This implies that $\|m\circ\theta - m_0\circ\theta_0\|=0$.  Then by assumption~\ref{a0} we have that $m\equiv m_0$ and $\theta=\theta_0.$ Now recall that in Theorem~\ref{thm:rate_m_theta_CLSE} and  Lemma~\ref{lem:Upsilion_ep}, we showed that $\|\check{m}\circ\check{\theta} - m_0\circ\theta_0\| = o_p(1)$ and  $\|\check{m}\|_\infty=O_p(1)$, respectively. Thus by taking $m_n = \check{m}_{{L}}$ and $\theta_n = \check{\theta}_{{L}}$, we have that $|\check{\theta}_{{L}} - \theta_0| = o_p(1)$ and $\|\check{m}_{{L}} - m_0\|_{D_0} = o_p(1).$ The following lemma applied to $\{\check{m}\}$ completes the proof of the theorem by showing that $\|\check{m}' - m'_0\|_{C} = o_p(1)$  for any compact subset $C$ in the interior of $D_0$.

% { Rewrite}
% As $\check{m}_n\in \mathcal{M}_L$, the sequence $\{\check{m}_n\}$
%  By Ascoli-Arzel\`{a} theorem,  the sequence $\{\check{m}_n\}$ has a uniformly converging subsequence if $\{\check{m}_n\}$ belongs to a class of closed, bounded, and  equicontinuous functions. In the proof of Theorem \ref{thm:rate_m_theta_CLSE}, we showed that $\|\check{m}_n\|_\infty=O_p(1);$ see Lemma~\ref{lem:Upsilion_ep}. Moreover,  as $\check{m}_n \in \mathcal{M}_L$  the conditions of Ascoli-Arzel\`{a} theorem are satisfied  for this sequence. Furthermore, as $\{\check{\theta}_n\}$ is a sequence in compact set. Thus $()$ we have that every subsequence $\{\check{\theta}_{n_k}\}$ has a further subsequence $\{\check{\theta}_{n_{k_l}}\}$ such that $| \check{\theta}_{n_{k_l}}-\theta_1| \rightarrow 0,$ for some $\theta_1$ in $\Theta$. Now, observe that continuity and almost everywhere differentiability of the link functions imply that $\|m_1\circ\theta_1 - m_0\circ\theta_0\|=0$ is equivalent to $m_1\equiv m_0$ and $\theta_1=\theta_0$.
%  % see Lemma \ref{lem:id}.
% % The identifiablity is used here
%   Thus we have that $\|\check{m}_n-m_0\|_{D_0}=o_p(1).$ Now the final result follows by an application of the following lemma and a standard subsequence argument.
\begin{lemma}[Lemma 3.10, \cite{EMSE}]\label{sebo}
Let $\mathcal{C}$ be an open convex subset of  $\R^d$ and $f$ a convex functions which is continuous and differentiable on $\mathcal{C}$. Consider a sequence of convex functions $\{f_n\}$ which are finite on $\mathcal{C}$ such that $f_n\to f$ pointwise on $\mathcal{C}$. Then, if $C\subset\mathcal{C}$ is any compact set,
\[
\sup_{\substack{x\in C\\\xi\in\partial f_n(x)}}|\xi - \nabla f(x)|\to 0,
\]
where $\partial f_n(x)$ represents the sub-differential set of $f_n$ at $x$. \qedhere
\end{lemma}
%%%%%%%%%%%%%%%%%%%%%%%%%%%%%%%%%%%%%%%%%
%%%%%%%%%%%%%%%%%%%%%%%%%%%%%%%%%%%%%%%%%
\subsection{Proof of Theorem \ref{thm:ratestCLSE}} \label{sec:proof_ratestPLSE}
For notational convenience and to show the dependence of $\check{m}$ and $\check{\theta}$ on $n$, we use $\check{m}_{n}$ and $\check{\theta}_n$ to denote $\check{m}_{{L}}$  (or $\check{m}$) and $\check{\theta}_{{L}}$ (or $\check{\theta}$), respectively. For the proof of Theorem \ref{thm:ratestCLSE}, we use two preliminary lemmas proved in Section~\ref{subsubsec:Theorem35}. Let us define, $A_n(x): =\check{m}_n(\check{\theta}_n^{\top}x) - m_0(\theta_0^{\top}x)$ and $B_n(x):=m_0'(\theta_0^{\top}x)x^{\top}(\check{\theta}_n - \theta_0) + (\check{m}_n - m_0)(\theta_0^{\top}x).$ Observe that
\begin{align*}
A_n(x)-B_n(x)&=\check{m}_n(\check{\theta}_n^{\top}x) -m_0'(\theta_0^{\top}x)x^{\top}(\check{\theta}_n - \theta_0) - \check{m}_n(\theta_0^{\top}x).\\
&=\check{m}_n(\check\theta_n^{\top}x)   - m_0(\theta_0^{\top}x) - \{m_0'(\theta_0^{\top}x)x^{\top}(\check\theta_n - \theta_0) + (\check{m}_n-m_0)(\theta_0^{\top}x)\}.
\end{align*}
We will now show that
\begin{equation}\label{eq:rate_prob}
D_n:=\frac{1}{|\check{\theta}_{n} -\theta_0|^2} P_X|A_{n}(X)-B_{n}(X)|^2=o_p(1).
\end{equation}
It is equivalent to show that for every subsequence $\{D_{n_{k}}\}$, there exists a further subsequence  $\{D_{n_{k_l}}\}$ that converges to 0 almost surely; see Theorem 2.3.2 of \cite{Durrett}.
We showed in Theorem \ref{thm:uucons}, that $\{\check{m}_n, \check{\theta}_n\}$ satisfies assumption \eqref{eq:rate_cond} of Lemma~\ref{lem:ratest_1CLSE} in probability. Thus by  another application of Theorem 2.3.2 of \cite{Durrett}, we have that $\{\check{m}_{n_k}, \check{\theta}_{n_k}\}$ has a further subsequence $\{\check{m}_{n_{k_l}}, \check{\theta}_{n_{k_l}}\}$ that satisfies \eqref{eq:rate_cond} almost surely. Thus by Lemma \ref{lem:ratest_1CLSE}, we have
$D_{n_{k_l}}\stackrel{a.s.}{\rightarrow} 0.$ Thus $D_n=o_p(1).$

 We will now use \eqref{eq:rate_prob} to find the rate of convergence of $\{\check{m}_n, \check{\theta}_n\}$. We first find an upper bound for $P_X |B_n(X)|^2$. By a simple application of triangle inequality and \eqref{eq:rate_prob}, we have
\begin{align*}
P_X|A_n(X)|^2\ge \frac{1}{2} P_X|B_n(X)|^2- P_X|A_n(X)-B_n(X)|^2\ge\frac{1}{2} P_X|B_n(X)|^2- o_p(|\check{\theta}_n -\theta_0|^2).
\end{align*}
As $q\ge 5$, by Theorem \ref{thm:rate_m_theta_CLSE}, we have that $P_X|A_n(X)|^2 =O_p(n^{-4/5}).$
%\begin{align}\label{eq:temp1}
%\begin{split}
%O_p(\hat{\lambda}_n^2) = P_X \big|&\check{m}_n(\check{\theta}_n^{\top}X) - m_0(\theta_0^{\top}X)\big|^2 \\
%&\ge \frac{1}{2}P_X\big|m_0'(\theta_0^{\top}X)X^{\top}(\check{\theta}_n - \theta_0) + (\check{m}_n - m_0)(\theta_0^{\top}X)\big|^2 - o_p(1)|\check{\theta}_n - \theta_0|^2.\end{split}
%\end{align}
Thus we have 
\begin{equation}\label{eq:AnMinusBnBound}
\begin{split}
P_X|B_n(X)|^2 &= P_X\big|m_0'(\theta_0^{\top}X)X^{\top}(\check{\theta}_n - \theta_0) + (\check{m}_n - m_0)(\theta_0^{\top}X)\big|^2\\  &\le O_p(n^{-4/5}) +o_p(|\check{\theta}_n - \theta_0|^2).
\end{split}
\end{equation}
Now define
\begin{equation}\label{eq:def_g_1}
\gamma_n := \frac{\check{\theta}_n - \theta_0}{|\check{\theta}_n - \theta_0|},\quad g_1(x) := m_0'(\theta_0^{\top}x)x^{\top}(\check{\theta}_n - \theta_0)  \text{ and}\quad  g_2(x) := (\check{m}_n - m_0)(\theta_0^{\top}x).
\end{equation}
Note that for all $n$,
\begin{align}\label{eq:step_th6}
\begin{split}
P_Xg_1^2 &= (\check{\theta}_n - \theta_0)^{\top}P_X[XX^{\top}|m_0'(\theta_0^{\top}X)|^2](\check{\theta}_n - \theta_0)\\
&= |\check{\theta}_n - \theta_0|^2\, \gamma_n^\top P_X[XX^{\top}|m_0'(\theta_0^{\top}X)|^2]\gamma_n\\
&\ge |\check{\theta}_n - \theta_0|^2\, \gamma_n^\top\E\big[\Var(X|\theta_0^\top X) |m_0'(\theta_0^{\top}X)|^2\big]\gamma_n.
\end{split}
\end{align}
Since $\gamma_n^{\top}\theta_0$ converges in probability to zero, we get by Lemma~14 of~\cite{Patra16} and assumption~\ref{aa5_new} that with probability converging to one,
\begin{equation}\label{eq:G1SquareBound}
\frac{P_Xg_1^2}{|\check{\theta}_n - \theta_0|^2} \ge \frac{\lambda_{\min}(H_{\theta_0}^\top\E\big[\Var(X|\theta_0^\top X) |m_0'(\theta_0^{\top}X)|^2\big] H_{\theta_0})}{2} > 0.
\end{equation}
Thus we can see that proof of this theorem will be complete if we can show that
\begin{equation}\label{eq:th6_con}
P_X g_1^2 +P_X g_2^2 \lesssim  P_X\big|m_0'(\theta_0^{\top}X)X^{\top}(\check{\theta}_n - \theta_0) + (\check{m}_n - m_0)(\theta_0^{\top}X)\big|^2.
\end{equation}
We will first prove that~\eqref{eq:th6_con} completes the proof of Theorem~\ref{thm:ratestCLSE}. Note from the combination of~\eqref{eq:step_th6}, \eqref{eq:G1SquareBound} and~\eqref{eq:th6_con} that
\[
|\check{\theta}_n - \theta_0|^2 = O_p(n^{-4/5})\quad\Rightarrow\quad |\check{\theta}_n - \theta_0| = O_p(n^{-2/5}).
\]
Substituting this in~\eqref{eq:AnMinusBnBound} and using~\eqref{eq:th6_con}, we get
\[
P_Xg_2^2 = O_p(n^{-4/5})\quad\Rightarrow\quad \norm{\check{m}_n\circ\theta_0 - m_0\circ\theta_0} = O_p(n^{-2/5}).
\]
This is same as
\begin{equation}\label{eq:IntegralFtheta_0}
\int_{D_0} \left(\check{m}_n(t) - m_0(t)\right)^2dP_{\theta_0^{\top}(X)}(t) dt = O_p(n^{-4/5}).
\end{equation}
% From assumption~\ref{aa6}, $f_{\theta_0^{\top}X}(t) \ge \underline{C_d} > 0$ and so~\eqref{eq:IntegralFtheta_0} implies that
% \[
% \int_{D_0} \left(\check{m}_n(t) - m_0(t)\right)^2dt = O_p\left(n^{-4/5}\right).
% \]
% \todo[inline]{This lemma might be there if I remove the smooth thing}
Now to prove~\eqref{eq:th6_con}. Note that by Lemma 5.7 of \cite{VANC}, a sufficient condition for~\eqref{eq:th6_con} is 
% \todo[inline]{This lemma might be there if I remove the smooth thing}
\begin{equation}\label{eq:Lemma_57}
({P}_Xg_1g_2)^2 \le c{P}_Xg_1^2{P}_Xg_2^2 \quad \text{for some constant }c<1
\end{equation}
We now show that $g_1$ and $g_2$ satisfy~\eqref{eq:Lemma_57}. By Cauchy-Schwarz inequality, we have 
\begin{align*}
\big(P_X[g_1(X) g_2(X)]\big)^2
&= \big(P_X \big[ m_0'(\theta_0^{\top}X)g_2(X)E(X^{\top}(\check{\theta} - \theta_0)|\theta_0^{\top}X)\big]\big)^2\\
&\le P_X\big[\{m_0'(\theta_0^{\top}X)\}^2E^2[X^{\top}(\check{\theta} - \theta_0)|\theta_0^{\top}X]\big]P_Xg_2^2(X)\\
&= |\check{\theta} - \theta_0|^2 \gamma_n^\top P_X\big[|m_0'(\theta_0^{\top}X)|^2E[X|\theta_0^{\top}X] E[X^{\top}|\theta_0^{\top}X]\big] \gamma_n P_Xg_2^2(X)\\
&= c_n |\check{\theta} - \theta_0|^2 \gamma_n^\top P_X\big[|m_0'(\theta_0^{\top}X)|^2XX^{\top}\big] \gamma_n P_Xg_2^2(X)\\
&= c_n P_X g_1^2 P_Xg_2^2(X),
% &< P_X\big[\{m_0'(\theta_0^{\top}X)\}^2E[\{X^{\top}(\check{\theta} - \theta_0)\}^2|\theta_0^{\top}X]\big]P_Xg_2^2(X)\\
% &= P_X\big[\mathbb{E}[\{m_0'(\theta_0^{\top}X)X^{\top}(\check{\theta} - \theta_0)\}^2|\theta_0^{\top}X)]\big]P_Xg_2^2(X)\\
% &= P_X[m_0'(\theta_0^{\top}X)X^{\top}(\check{\theta} - \theta_0)]^2P_X g_2^2(X)\\
% &= P_Xg_1^2\,P_Xg_2^2.
\end{align*}
where \[c_n:= \frac{ \gamma_n^\top P_X\big[|m_0'(\theta_0^{\top}X)|^2E[X|\theta_0^{\top}X] E[X^{\top}|\theta_0^{\top}X]\big] \gamma_n }{\gamma_n^\top P_X\big[|m_0'(\theta_0^{\top}X)|^2 X X^{\top}\big] \gamma_n}.
\]
To show that with probability converging to one, $c_n < 1$, observe that
\[
1 - c_n = \frac{ \gamma_n^\top\E\big[\Var(X|\theta_0^\top X) |m_0'(\theta_0^{\top}X)|^2\big]\gamma_n}{ \gamma_n^\top\E\big[XX^{\top} |m_0'(\theta_0^{\top}X)|^2\big]\gamma_n}
\]
and by Lemma~\ref{lem:gamma_n} along with assumption~\ref{aa5_new}, with probability converging to one,
\[
1 - c_n > \frac{ 4\lambda_{\min}\left(H_{\theta_0}^{\top}\E\big[\Var(X|\theta_0^\top X) |m_0'(\theta_0^{\top}X)|^2\big]H_{\theta_0}\right)}{\lambda_{\max}\left(H_{\theta_0}^{\top}\E\big[XX^{\top} |m_0'(\theta_0^{\top}X)|^2\big]H_{\theta_0}\right)} > 0.
\]
This implies that with probability converging to one, $c_n < 1.$
\subsection{Lemmas used in the proof of Theorem~\ref{thm:ratestCLSE}}\label{subsubsec:Theorem35}
In this section, we state and prove the two preliminary lemmas used in the proof of Theorem~\ref{thm:ratestCLSE}.
%We use Lemma 5.7 of \cite{VANC} (stated in Section \ref{app:used_lemmas}) and the following lemma to prove Theorem \ref{thm:ratestPLSE}. See Section for more details.
\begin{lemma} \label{lem:ratest_1CLSE} Let $m_0$ and $\theta_0$ satisfy the assumptions \ref{aa1_new}, \ref{aa1}. Furthermore, let  $\{\theta_n\} \in \Theta$ and $\{m_n\} \in \mathcal{M}_L$ be two non-random sequences such that
\begin{equation}\label{eq:rate_cond}
|\theta_n - \theta_0| \rightarrow 0, \qquad  \|m_n - m_0\|_{D_0}  \rightarrow 0, \quad \text{and} \quad \|{m}'_n - m'_0\|_{C}   \rightarrow 0
\end{equation}
 for any compact subset $C$ of the interior of $D_0$. Then
\begin{align*}
P_X\big|&m_n(\theta_n^{\top}X)   - m_0(\theta_0^{\top}X) - \{m_0'(\theta_0^{\top}X)X^{\top}(\theta_n - \theta_0) + (m_n-m_0)(\theta_0^{\top}X)\}  \big|^2 =o(|\theta_n-\theta_0|^2) .
\end{align*}
\end{lemma}

\begin{proof}  For any convex function $f \in \mathcal{M}_L$,  denote the right derivative of $f$ by $f'$. Note that $f'$ is  a bounded  nondecreasing function. First, observe that
\begin{align*}
m_n(\theta_n^{\top}x)   - m_0(\theta_0^{\top}x) &- \big[m_0'(\theta_0^{\top}x)x^{\top}(\theta_n - \theta_0) + (m_n-m_0)(\theta_0^{\top}x)\big]\\
&= m_n(\theta_n^{\top}x) - m_n(\theta_0^{\top}x) - m_0'(\theta_0^{\top}x)x^{\top}(\theta_n - \theta_0).
\end{align*}
% When $\theta_n=\theta_0$, then it is easy to see that \[P_X\big|m_n(\theta_n^{\top}X)   - m_0(\theta_0^{\top}X) - \{m_0'(\theta_0^{\top}X)X^{\top}(\theta_n - \theta_0) + (m_n-m_0)(\theta_0^{\top}X)\}  \big|^2=0.\]
Now,
{\small \begin{align}
&\big|m_n(\theta_n^{\top}x) - m_n(\theta_0^{\top}x) - m_0'(\theta_0^{\top}x)x^{\top}(\theta_n - \theta_0)\big|^2 \nonumber\\
\;\;  \qquad ={}& \left|\int_{\theta_n^\top x}^{\theta_0^\top x} m_n^{\prime}(t)\,dt - m_0'(\theta_0^{\top}x)x^{\top}(\theta_n - \theta_0)\right|^2 \quad \text{($m_n$ is absolutely continuous)}\nonumber\\
\;\; \qquad ={}& \left|\int_{\theta_n^\top x}^{\theta_0^\top x} m_n^{\prime}(t)\,dt - m_n^{\prime}(\theta_0^{\top}x)x^{\top}(\theta_n - \theta_0)+  m_n^{\prime}(\theta_0^{\top}x)x^{\top}(\theta_n - \theta_0)- m_0'(\theta_0^{\top}x)x^{\top}(\theta_n - \theta_0)\right|^2\nonumber\\
\;\; \qquad  ={}& \left|\int_{\theta_n^\top x}^{\theta_0^\top x} m_n^{\prime}(t)\,dt - m_n^{\prime}(\theta_0^{\top}x)x^{\top}(\theta_n - \theta_0)+  (m_n^{\prime}- m_0')(\theta_0^{\top}x)x^{\top}(\theta_n - \theta_0)\right|^2\nonumber\\
\;\; \qquad \le{} & 2 \left|\int_{\theta_n^\top x}^{\theta_0^\top x} m_n^{\prime}(t)\,dt - m_n^{\prime}(\theta_0^{\top}x)x^{\top}(\theta_n - \theta_0)\right|^2+ 2\left| (m_n^{\prime}- m_0')(\theta_0^{\top}x)x^{\top}(\theta_n - \theta_0)\right|^2\label{eq:Lemma_th3_eq1}.
\end{align}}
% \begin{align*}
% \end{align*}
We will now find an upper bound for the first term on the right hand side of the above display. Observe that $m_n^\prime$ is a nondecreasing function. When $x^\top\theta_n\neq x^\top\theta_0$, we have
\[m_n^{\prime}(\theta_n^\top x) \wedge  m_n^{\prime}(\theta_0^\top x) \le \frac{ \int_{\theta_n^\top x}^{\theta_0^\top x} m_n^{\prime}(t)\,dt}{x^{\top}(\theta_n - \theta_0)}  \le    m_n^{\prime}(\theta_n^\top x) \vee m_n^{\prime}(\theta_0^\top x).\] Thus for all  $x\in \rchi$, we have
\begin{equation}\label{eq:th3_1st_bound}
\left|\int_{\theta_n^\top x}^{\theta_0^\top x} m_n^{\prime}(t)\,dt - m_n^{\prime}(\theta_0^{\top}x)x^{\top}(\theta_n - \theta_0)\right| \le |m_n^{\prime}(\theta_n^{\top}x)- m_n^{\prime}(\theta_0^{\top}x)| |x^\top (\theta_n-\theta_0)|.
\end{equation} Note that if $x^\top\theta_n=x^\top\theta_0$, then both sides of \eqref{eq:th3_1st_bound} are $0.$
% \[m_n^{\prime}(\theta_n^\top x) \wedge  m_n^{\prime}(\theta_0^\top x)= \min_{ t\in [\theta_n^\top x,\theta_0^\top x] } m_n^{\prime}(t) \le \frac{ \int_{\theta_n^\top x}^{\theta_0^\top x} m_n^{\prime}(t)\,dt}{x^{\top}(\theta_n - \theta_0)}  \le   \max_{ t\in [\theta_n^\top x,\theta_0^\top x] } m_n^{\prime}(t) = m_n^{\prime}(\theta_n^\top x) \vee m_n^{\prime}(\theta_0^\top x) \]
%
%
% Recall that, $\{m_n^{\prime}\}$ is  a sequence of bounded increasing functions. Therefore for each $n$, $m_n$ has at most countable discontinuities. Let us denote the points of discontinuities of $m_n$ by $Dis(m_n)$. It is now easy to see that $\cup_{i\ge1} Dis(m_n)$  is  at most countable. Recall that, $\theta \in B_{\theta_0}(\delta_0)$,  $\theta^\top X $ has continuous density. As $|\theta_n-\theta_0|=o(1)$ and  we have that for almost every ($P_X$) $x$ (as for all $\theta \in B_{\theta_0}(\delta)$, $\theta^\top X$ has a continuous density)
%
% \[\big|\int_{\theta_n^\top x}^{\theta_0^\top x} m_n^{\prime}(t)\,dt- m_n^{\prime}(\theta_0^{\top}x)x^{\top}(\theta_n - \theta_0)\big|= \big|\{m_n^{\prime}(\xi^{\top}x)- m_n^{\prime}(\theta_0^{\top}x)\} x^\top (\theta_n-\theta_0)\big|\]
% where $\xi^{\top}x$ lies between $\theta_n^{\top}x$ and $\theta_0^{\top}x$. As $m_n^{\prime}$ is a monotonic function, we have that
% \[|m_n^{\prime}(\xi^{\top}x)- m_n^{\prime}(\theta_0^{\top}x)| \le |m_n^{\prime}(\theta_n^{\top}x)- m_n^{\prime}(\theta_0^{\top}x)| \]
%
%
Combine \eqref{eq:Lemma_th3_eq1} and \eqref{eq:th3_1st_bound}, to conclude that
\begin{align}\label{eq:split_rate}
 \begin{split}
&P_X\big|m_n(\theta_n^{\top}X) - m_n(\theta_0^{\top}X) - m_0'(\theta_0^{\top}X)X^{\top}(\theta_n - \theta_0)\big|^2\\
 \le{}& 2 P_X\left|(m_n^{\prime}(\theta_n^{\top}X)- m_n^{\prime}(\theta_0^{\top}X)) X^{\top}(\theta_n - \theta_0)\right|^2+ 2 P_X\left| (m_n^{\prime}- m_0')(\theta_0^{\top}X)X^{\top}(\theta_n - \theta_0)\right|^2.
 \end{split}
 \end{align} 
As $\rchi$ is bounded, the two terms on the right hand side of \eqref{eq:split_rate} can be bounded as
\begin{align*}
P_X\left|(m_n^{\prime}(\theta_n^{\top}X)- m_n^{\prime}(\theta_0^{\top}X)) X^{\top}(\theta_n - \theta_0)\right|^2
\le&T^2 |\theta_n - \theta_0|^2 P_X\left|m_n^{\prime}(\theta_n^{\top}X)- m_n^{\prime}(\theta_0^{\top}X) \right|^2,\\
 P_X\left| (m_n^{\prime}- m_0')(\theta_0^{\top}X)X^{\top}(\theta_n - \theta_0)\right|^2  \le& T^2 |\theta_n - \theta_0|^2 P_X\left| (m_n^{\prime}- m_0')(\theta_0^{\top}X)\right|^2.
\end{align*}
We will now show that both $P_X\left|m_n^{\prime}(\theta_n^{\top}X)- m_n^{\prime}(\theta_0^{\top}X) \right|^2$ and $P_X\left| (m_n^{\prime}- m_0')(\theta_0^{\top}X)\right|^2$ converge to $0$ as $n\rightarrow \infty.$ First observe that
\begin{align}
P_X\left|m_n^{\prime}(\theta_n^{\top}X)- m_n^{\prime}(\theta_0^{\top}X) \right|^2 &\lesssim  P_X\left|m_n^{\prime}(\theta_n^{\top}X)- m_0^{\prime}(\theta_n^{\top}X) \right|^2 + P_X\left|m_0^{\prime}(\theta_n^{\top}X)- m_0^{\prime}(\theta_0^{\top}X) \right|^2\nonumber\\
 & \qquad \qquad+ P_X\left|m_0^{\prime}(\theta_0^{\top}X)- m_n^{\prime}(\theta_0^{\top}X) \right|^2.\label{eq:3part}
\end{align}
Recall that $m_0^\prime$ is a continuous  and bounded function; see assumption \ref{aa1_new}. Bounded convergence theorem now implies that $P_X\left|m_0^{\prime}(\theta_n^{\top}X)- m_0^{\prime}(\theta_0^{\top}X) \right|^2 \rightarrow 0,$ as $|\theta_n-\theta_0| \rightarrow 0.$ Now consider the first term on the right hand side of \eqref{eq:3part}.  {By~\ref{aa6}, we have that $\theta_0^\top X$ has a density, for any $\varepsilon >0$, we can define a compact subset $C_\varepsilon$ in the interior of $D_0$  such that $\p(\theta_0^\top X \notin C_\varepsilon)  < \varepsilon/8L^2$}. Now note that, by Theorem \ref{thm:uucons} and the fact that $\p(\theta_n^\top X \notin C_\varepsilon) \rightarrow  \p(\theta_0^\top X \notin C_\varepsilon)$, we have
\[P_X\left|m_n^{\prime}(\theta_n^{\top}X)- m_0^{\prime}(\theta_n^{\top}X) \right|^2 \le \sup_{t\in C_\varepsilon} |m_n^\prime(t)- m_0'(t)|^2 + 4L^2  P(\theta_n^\top X \notin C_{\varepsilon}) \le \varepsilon,\] as $n\rightarrow \infty.$ Similarly, we can see that
\[P_X\left|m_0^{\prime}(\theta_0^{\top}X)- m_n^{\prime}(\theta_0^{\top}X) \right|^2 \le \sup_{t\in C_\varepsilon} |m_n^\prime(t)- m_0'(t)|^2 + 4L^2  P(\theta_0^\top X \notin C_{\varepsilon}) \le \varepsilon,\] as $n\rightarrow \infty$. Combining the results, we have shown that for every $\varepsilon>0$
\[
P_X\big|m_n(\theta_n^{\top}X) - m(\theta_0^{\top}X) - m_0'(\theta_0^{\top}X)X^{\top}(\theta_n - \theta_0)\big|^2 \le T^2 |\theta_n - \theta_0|^2 \varepsilon,\] for all sufficiently large $n.$ Thus the result follows.\qedhere

%  Now find an upper bound for the second term on the right hand of \eqref{eq:split_rate}. Observe that by Theorem \ref{thm:uucons}, we have  $\|m_n^{\prime}-m_0^\prime\|_C=o_p(1), $ for every compact subset $C$ of $D_0$. Moreover, as $\theta_0^\top X$ has continuous density and both  $m_n^{\prime}$ and $m_0^\prime$ are bounded functions, we have
% $P_X\left| (m_n^{\prime}- m_0')(\theta_0^{\top}X)x^{\top}(\theta_n - \theta_0)\right|^2 =o_p(1)|\theta_n - \theta_0|^2.$
% Thus we have
% \[P_X\big|m_n(\theta_n^{\top}X)   - m_0(\theta_0^{\top}X) - \{m_0'(\theta_0^{\top}X)X^{\top}(\theta_n - \theta_0) + (m_n-m_0)(\theta_0^{\top}X)\}  \big|^2 =o_p(1)|\theta_n - \theta_0|^2\]
\end{proof}
\begin{lemma}\label{lem:gamma_n}
Suppose $A\in\mathbb{R}^{d\times d}$ and let $\{\gamma_n\}$ be any sequence of random vectors in $S^{d - 1}$ satisfying $\theta_0^{\top}\gamma_n =o_p(1).$ Then 
\[
\p\bigg(0.5 \lambda_{\min}\left(H_{\theta_0}^{\top}AH_{\theta_0}\right) \le \gamma_n^{\top}A\gamma_n \le 2 \lambda_{\max}\left(H_{\theta_0}^{\top}AH_{\theta_0}\right)\bigg) \to 1,
\]
where for any symmetric matrix $B$, $\lambda_{\min}(B)$ and $\lambda_{\max}(B)$ denote, respectively, the minimum and the maximum eigenvalues of $B.$
\end{lemma}
\begin{proof}
Note that $\text{Col}(H_{\theta_0}) \oplus\{\theta_0\} =\R^d$, thus
\begin{equation}\label{eq:OrthogonalDecomposition}
\gamma_n = \left(\gamma_n^{\top}\theta_0\right)\theta_0 + H_{\theta_0}\left(H_{\theta_0}^{\top}\gamma_n\right).
\end{equation}
Therefore,
\begin{align*}
\gamma_n^{\top}A\gamma_n &= \left[\left(\gamma_n^{\top}\theta_0\right)\theta_0 + H_{\theta_0}\left(H_{\theta_0}^{\top}\gamma_n\right)\right]^{\top}A\left[\left(\gamma_n^{\top}\theta_0\right)\theta_0 + H_{\theta_0}\left(H_{\theta_0}^{\top}\gamma_n\right)\right]\\
&= \left(\gamma_n^{\top}\theta_0\right)^2\theta_0^{\top}A\theta_0 + \left(\gamma_n^{\top}\theta_0\right)\theta_0^{\top}AH_{\theta_0}\left(H_{\theta_0}^{\top}\gamma_n\right)\\ &\qquad+ \left(\gamma_n^{\top}\theta_0\right)\left(H_{\theta_0}^{\top}\gamma_n\right)^{\top}H_{\theta_0}^{\top}A\theta_0 + \left(H_{\theta_0}^{\top}\gamma_n\right)^{\top}H_{\theta_0}^{\top}AH_{\theta_0}\left(H_{\theta_0}^{\top}\gamma_n\right).
\end{align*}
Note that $H_{\theta_0}^{\top}\gamma_n$ is a bounded sequence of vectors. Because of $\gamma_n^{\top}\theta_0$ in the first three terms above, they converge to zero in probability and so,
\[
\left|\gamma_n^{\top}A\gamma_n - \left(\gamma_n^{\top}H_{\theta_0}\right)H_{\theta_0}^{\top}AH_{\theta_0}\left(H_{\theta_0}^{\top}\gamma_n\right)\right| = o_p(1).
\]
Also, note that from~\eqref{eq:OrthogonalDecomposition},
\[
|H_{\theta_0}^{\top}\gamma_n|^2 - 1 = |\gamma_n|^2 - \left(\gamma_n^{\top}\theta_0\right)^2 - 1 = - \left(\gamma_n^{\top}\theta_0\right)^2 = o_p(1).
\]
Therefore, as $n\to\infty$,
\begin{equation}\label{eq:InterimoPOneStatement}
\left|\gamma_n^{\top}A\gamma_n - \frac{\left(\gamma_n^{\top}H_{\theta_0}\right)H_{\theta_0}^{\top}AH_{\theta_0}\left(H_{\theta_0}^{\top}\gamma_n\right)}{|H_{\theta_0}^{\top}\gamma_n|^2}\right| = o_p(1).
\end{equation}
By the definition of the minimum and maximum eigenvalues,
\[
\lambda_{\min}\left(H_{\theta_0}^{\top}AH_{\theta_0}\right) \le \frac{\left(\gamma_n^{\top}H_{\theta_0}\right)H_{\theta_0}^{\top}AH_{\theta_0}\left(H_{\theta_0}^{\top}\gamma_n\right)}{|H_{\theta_0}^{\top}\gamma_n|^2} \le \lambda_{\max}\left(H_{\theta_0}^{\top}AH_{\theta_0}\right).
\]
Thus using~\eqref{eq:InterimoPOneStatement} the result follows.
\end{proof}
\subsection{Proof of Theorem \ref{thm:rate_derivCLSE}}\label{proof:DerivCLSE}
% For any $\delta_1 < \delta_2$ and any $t$, we have
% \[
% \frac{(\delta_2 - \delta_1)t + \delta_1(t - \delta_2)}{\delta_2} = t - \delta_1.
% \]
% Convexity of a function ${m}$ implies that
% \[
% {m}(t - \delta_1) \le \frac{\delta_2 - \delta_1}{\delta_2}m(t) + \frac{\delta_1}{\delta_2}m(t - \delta_2),
% \]
% which by rearranging implies that
% \[
% \frac{m(t - \delta_1) - m(t)}{\delta_1} \le \frac{m(t - \delta_2) - m(t)}{\delta_2}.
% \]
% Thus the function
% \[
% f_t(\delta) = \frac{m(t - \delta) - m(t)}{\delta},
% \]
% for any fixed $t$ is an increasing function of $\delta$.
\textbf{Proof of~\eqref{eq:deriv_theta_0}:}
We first show  the first part of~\eqref{eq:deriv_theta_0}. Let $\delta_n$ be a sequence of positive numbers decreasing to $0$. Let $a,b \in \R$ such that $D_0 = [a,b]$. Define $C_n := [a + 2\delta_n, b - 2\delta_n]$.  By \ref{aa6}, $f_{\theta_0^\top X},$ the density of $\theta^\top_0X$ is bounded from above. Recall that $\overline{C}_d$ denotes the maximum of $f_{\theta_0^{\top}X}(\cdot)$. Because $\check{m}$ is a convex function,  we have
\[ \frac{\check{m}(t)- \check{m}(t-\delta_n)}{ \delta_n} \le \check{m}^\prime(t-) \le \check{m}^\prime(t+)  \le \frac{\check{m}(t+\delta_n)- \check{m}(t)}{ \delta_n},\]
for all $t\in C_n$, where $\check{m}^\prime(t+)$ and $\check{m}^\prime(t-)$ denote the right and left derivatives of $\check{m}$ at $t$, respectively. Observe that
\begin{align}
\int_{t\in C_n} &\left[ \frac{\check{m}(t+\delta_n)- \check{m}(t)}{ \delta_n}- \frac{m_0(t+\delta_n)- m_0(t)}{ \delta_n}\right]^2 f_{\theta_0^\top X}(t)dt \nonumber\\
&\qquad=\frac{2}{\delta_n^2}\int_{t\in C_n}\{\check{m}(t + \delta_n) - {m}_0(t + \delta_n)\}^2 f_{\theta_0^\top X}(t)dt + \frac{2}{\delta_n^2}\int_{t\in C_n}\{\check{m}(t) - m_0(t)\}^2 f_{\theta_0^\top X}(t)dt\nonumber\\
&\qquad=\frac{2}{\delta_n^2}\int_{t\in [a + 3\delta_n, b - \delta_n]}\{\check{m}(t) - m_0(t)\}^2 f_{\theta_0^\top X}(t) dt  + \frac{2}{\delta_n^2}\int_{t\in C_n}\{\check{m}(t) - m_0(t)\}^2 f_{\theta_0^\top X}(t) dt\nonumber\\
% &\qquad\le\frac{2}{K \delta_n^2}\int_{t\in [a + 3\delta_n, b - \delta_n]}\{\check{m}(t) - m_0(t)\}^2 f_{\theta_0^{\top}X}(t) dt  + \frac{2}{K\delta_n^2}\int_{t\in C_n}\{\check{m}(t) - m_0(t)\}^2f_{\theta_0^{\top}X}(t)dt\nonumber\\
&\qquad= \frac{1}{\delta_n^2} O_p(n^{-4/5}),\label{eq:int_uppper_bound}
\end{align}
where the last equality follows from Theorem \ref{thm:ratestCLSE} (as $q\ge5$ and $L$ is fixed).
 % and assumption \ref{aa6}, we have
Similarly, it can be shown that
\begin{equation}\label{eq:int_uppper_bnd_neg}
\int_{t\in C_n} \left[ \frac{\check{m}(t)- \check{m}(t-\delta_n)}{ \delta_n}- \frac{m_0(t)- m_0(t-\delta_n)}{ \delta_n}\right]^2 f_{\theta_0^\top X}(t) dt
= \frac{1}{\delta_n^2} O_p(n^{-4/5}).
\end{equation}
Now observe that, $|m_0^{\prime} (t)- m_0^{\prime} (X_{t_n})|\le L_1 \delta_n^{1/2}$ whenever $x_{t_n} \in [t-\delta_n, t]$, we have
\begin{align}\label{eq:integrnd_bd_pos}
\begin{split}
\alpha_n^+(t):=\left[ \frac{\check{m}(t+\delta_n)- \check{m}(t)}{ \delta_n}- \frac{m_0(t+\delta_n)- m_0(t)}{ \delta_n}\right] &\ge \check{m}^\prime(t+) - m_0^\prime(x_{t_n})\\
 &\ge \check{m}^\prime(t+) - m_0^\prime(t)+  m_0^\prime(t) -m_0^\prime(x_{t_n})\\
 &\ge  \check{m}^\prime(t+) - m_0^\prime(t) -L_1 \delta_n^{1/2},
\end{split}
\end{align}
where $x_{t_n}$ lies between $t$ and $t+\delta_n$. Moreover,
\begin{align}\label{eq:integrnd_bd_neg}
\begin{split}
\alpha_n^-(t):=\left[ \frac{\check{m}(t)- \check{m}(t-\delta_n)}{ \delta_n}- \frac{m_0(t)- m_0(t-\delta_n)}{ \delta_n}\right] &\le \check{m}^\prime(t+) - m_0^\prime(x'_{t_n})\\
 &\le \check{m}^\prime(t+) - m_0^\prime(t)+  m_0^\prime(t) -m_0^\prime(x'_{t_n})\\
 &\le  \check{m}^\prime(t+) - m_0^\prime(t) +L_1 \delta_n^{1/2},
\end{split}
\end{align}
where $x'_{t_n}$ lies between $t-\delta_n$ and $t$. Combining the above two results, we have
\begin{align}\label{eq:bnd_both_side}
\begin{split}
 &\alpha_n^-(t)-L_1 \delta_n^{1/2}  \le \check{m}^\prime(t+) - m_0^\prime(t)  \le\alpha_n^+(t) + L_1 \delta_n^{1/2};
\end{split}
\end{align}
see proof of Corollary 1 of \cite{DumbgenEtAL04} for a similar inequality. Thus for every $t \in C_n$, we have   $[\check{m}^\prime(t+) - m_0^\prime(t) ]^2 \le 2 L_1^2 \delta_n  +  2\max\left\{ [\alpha_n^-(t)]^2 ,   [\alpha_n^+(t)]^2 \right\}.$ By \eqref{eq:int_uppper_bound} and \eqref{eq:int_uppper_bnd_neg}, we have
\begin{equation}\label{eq:Proof_done}
 \int_{t\in C_n} [\check{m}^\prime(t+) - m_0^\prime(t) ]^2  f_{\theta_0^\top X}(t) dt \le 2 L_1^2 \delta_n +\frac{1}{\delta_n^2} O_p(n^{-4/5}),
\end{equation}
as
\begin{align*}
 \int_{t\in C_n}   \max\left\{ [\alpha_n^-(t)]^2 ,   [\alpha_n^+(t)]^2 \right\} f_{\theta_0^\top X}(t) dt 
&\le \int_{t\in C_n}  \{\alpha_n^-(t)\}^2 f_{\theta_0^\top X}(t) dt+ \int_{t\in C_n}  \{\alpha_n^+(t)\}^2 f_{\theta_0^\top X}(t)dt\\
% &\le K'\int_{t\in C_n}  \{\alpha_n^-(t)\}^2 dt+ K' \int_{t\in C_n}  \{\alpha_n^+(t)\}^2dt\\
&=\frac{1}{\delta_n^2}O_p(n^{-4/5}).
\end{align*}
% where $K'$ is an upper bound on $f_{\theta^{\top}X}(\cdot)$ in $D_0$.
Moreover, note that $\|\check{m}^\prime\|_\infty \le L$ and  $\|m_0^\prime\|_\infty \le L_0 \le L$. Thus
\begin{align}\label{eq:rate_final_deriv}
\begin{split}
\int_{t\in D_0}\{\check{m}^\prime(t+) - m_0^\prime(t)\}^2 f_{\theta_0^\top X}(t) dt ={}&\int_{t\in C_n}\{\check{m}^\prime(t+) - m_0^\prime(t)\}^2 f_{\theta_0^\top X}(t) dt \\
&+\int_{t\in D_0\cap C_n^c}\{\check{m}^\prime(t+) - m_0^\prime(t)\}^2 f_{\theta_0^\top X}(t)dt \\
={}& 2 L_1^2 \delta_n +\frac{1}{\delta_n^2} O_p(n^{-4/5})+ 4L^2 4 \delta_n.
% \le{}& 2 \kappa^2 \delta_n^2 +\frac{1}{\delta_n^2} O_p(n^{-4/5})+ 16 L^2 \delta_n.
\end{split}
\end{align}
The tightest upper bound for the left hand side of the above display is achieved when $\delta_n =n^{-4/15}$. With this choice of $\delta_n$, we have
\begin{equation}\label{eq:for_use198}
\int\{\check{m}^\prime(t+) - m_0^\prime(t)\}^2 f_{\theta_0^\top X}(t) dt  \le 2 L_1^2 n^{-4/15} +O_p(n^{-4/15})+ 16 L^2 n^{-4/15} = O_p(n^{-4/15}).
\end{equation}
% We can find a similar upper bound for $\int_{t\in D_0}\{\check{m}^\prime(t+) - m_0^\prime(t)\}^2 dt:$
% \begin{align}
% \qquad\int_{t\in D_0}\{\check{m}^\prime(t+) - m_0^\prime(t)\}^2 dt ={}&\int_{t\in C_n}\{\check{m}^\prime(t+) - m_0^\prime(t)\}^2 dt +\int_{t\in D_0\cap C_n^c}\{\check{m}^\prime(t+) - m_0^\prime(t)\}^2 dt \label{LIP:eq:rate_final_deriv}\\
% \le{}& 2 \frac{\kappa^2 \delta_n^2}{K} +\frac{1}{K \delta_n^2} O_p(n^{-4/5})+ 16 L^2 \delta_n,\nonumber
% % \end{split}
% \end{align}
% where the first term on the right hand side of~\eqref{LIP:eq:rate_final_deriv} is bounded above via~\eqref{eq:Proof_done} and $K$ is the minimum of $f_{\theta_0^\top X}$.
% When $\delta_n =n^{-4/15}$, we have
% \begin{equation}\label{eq:bnd_111}
% \int_{t\in D_0}\{\check{m}^\prime(t+) - m_0^\prime(t)\}^2 dt \le 2\kappa^2 n^{-8/15} +O_p(n^{-4/15})+ 16 L^2 n^{-4/15} = O_p(n^{-4/15}).
% \end{equation}
{ We will now establish the second part of~\eqref{eq:deriv_theta_0}. Note that 
\begin{equation}\label{eq:part_2_3.6}
\|\check{m}^\prime \circ \check\theta- m'_0\circ\check\theta\|^2 
% = \int_{x\in \rchi}\{\check{m}^\prime(\check\theta^\top x + ) - m_0^\prime(\check\theta^\top x)\}^2 dP_X(x)
= \int\{\check{m}^\prime(t + ) - m_0^\prime(t )\}^2 f_{\check\theta^\top X }(t) dt.
\end{equation}
Note that 
\begin{align}\label{eq:11}
\begin{split}
\left|\int\{\check{m}^\prime(t + ) - m_0^\prime(t )\}^2 \big[f_{\check\theta^\top X }(t)-f_{\theta_0^\top X }(t)\big] dt\right| &\le 4L^2 \mathrm{TV}(\check\theta^\top X, \theta_0^\top X)\\ 
&\le 4L^2\overline{C}_0T |\check\theta- \theta_0|,
\end{split}
\end{align}
where $\mathrm{TV}(\check\theta^{\top}X, \theta_0^{\top}X)$ is defined as the evaluation of the total variation distance between $\theta^{\top}X$ and $\theta_0^{\top}X$ at $\theta = \check\theta$ and hence is random. The second inequality in~\eqref{eq:11} follows, if we can show that for any $\theta$,
\begin{align}\label{eq:Tv_lip}
\mathrm{TV}(\theta^\top X, \theta_0^\top X) = \sup_{t\in \R} |\p(\theta^\top X\le t)- \p( \theta_0^\top X\le t)| \le \overline{C}_0T |{\theta}-\theta_0|.
\end{align}
We will now prove~\eqref{eq:Tv_lip}. Because $\sup_{x\in \rchi} |x| \le T$, we have that $|\theta^\top x- \theta_0^\top x| \le  T|\theta-\theta_0|$ for all $x\in \rchi$. Now 
\begin{align}\label{eq:Tv_lip_proof}
\begin{split}
\p(\theta^\top X\le t ) &= \p(\theta^\top X\le t \text{ and } |\theta^\top X- \theta_0^\top X| \le  T|{\theta}-\theta_0|)\\
&  \le \p(\theta_0^\top X\le t +  T|{\theta}-\theta_0|)\\
& = \p(\theta_0^\top X\le t) + \p( t \le \theta_0^\top X\le  t+ T|{\theta}-\theta_0|)\\
& \le  \p(\theta_0^\top X\le t) + \overline{C}_0 T |{\theta}-\theta_0|
\end{split}
\end{align}
For the other side, observe
\begin{align}\label{eq:Equation Label}
\begin{split}
\p(\theta^\top X\le t ) &= \p(\theta^\top X\le t \text{ and } |\theta^\top X- \theta_0^\top X| \le  T|{\theta}-\theta_0|)\\
&\ge \p (\theta_0^\top X \le t-T|{\theta}-\theta_0|\text{ and } |\theta^\top X- \theta_0^\top X| \le  T|{\theta}-\theta_0|)\\
&= \p (\theta_0^\top X \le t-T|{\theta}-\theta_0|)\\
& = \p(\theta_0^\top X\le t) - \p( t-T|{\theta}-\theta_0| \le \theta_0^\top X\le  t)\\
& \ge  \p(\theta_0^\top X\le t) - \overline{C}_0 T |{\theta}-\theta_0|.
\end{split}
\end{align}

\textbf{Proof of~\eqref{eq:sup_rate}:}
We will use Lemma 2 of~\cite{DumbgenEtAL04} to prove both parts of~\eqref{eq:sup_rate}. We state the lemma at the end of this section for the convenience of the reader. We will now prove the first part of~\eqref{eq:sup_rate} by contradiction. Suppose 
\[\sup_{t\in C} |\check{m} (t)- m_0 (t)| > K_n n^{-8/(25 + 5\beta)},\] for some $K_n>0.$ Then by Lemma~\ref{lem:Dumbgen_lemma}, we have that there exists an interval $[c, c+\xi_n] \subset D_0$ such that 
\begin{equation}\label{eq:inf_bound}
\inf_{t\in [c, c+\xi_n]} |\check{m} (t)- m_0 (t)| > \frac{K_n}{4} n^{-8/(25 + 5\beta)},\qquad \text{for all }  n \ge [K_n/\diameter(D_0)]^{5(5 + \beta)/16},
\end{equation}
 where $\xi_n = A \sqrt{K_n n^{-8/(5(5 + \beta))}}$  and $A:= \big(64 \|m_0''\|_{D_0}\big)^{-1/2}$.
  % { Why is the next sentence here? What does it convey?}Note that by Theorem~\ref{thm:uucons} $K_n n^{-8/25}$ must converge to $0$. 
 Thus by~\eqref{eq:inf_bound}, we have
 \begin{align*}
  \int_{D_0} (\check{m}(t)- m_0 (t))^2 dP_{\theta_0^{\top}X}(t) &\ge \int_{c}^ {c+\xi_n} |\check{m} (t)- m_0 (t)|^2 dP_{\theta_0^{\top}X}(t)  \\
  &\ge \frac{K_n^2}{16} n^{-16/(5(5 + \beta))}   \int_{c}^{c+ \xi_n} dP_{\theta_0^{\top}X}(t)\\
  &\ge \frac{K_n^2}{16} n^{-16/(5(5 + \beta))} \left[\underline{C}_d\xi_n^{1 + \beta}\right]\\
  &= \frac{K_n^2\underline{C}_d}{16}n^{-4/5},
 \end{align*}
 where the last inequality above follows from assumption~\ref{bb2'}.

 However, by Theorem~\ref{thm:ratestCLSE}, we have that $\int_{D_0} (\check{m}(t)- m_0 (t))^2dP_{\theta_0^{\top}X}(t)  =O_p(n^{-4/5})$. Thus $K_n =O_p(1)$ (i.e., $K_n$ cannot diverge to infinity with $n$) and hence, $\sup_{t\in C} |\check{m} (t)- m_0 (t)| = O_p( n^{-8/(25+5\beta)}).$ Given the first part, the second part of~\eqref{eq:sup_rate} follows directly from the proof of Corollary~1 of~\cite{DumbgenEtAL04} with $\beta = 2$ (in that paper). 
}
\begin{lemma}\label{lem:Dumbgen_lemma}
Let $F$ be a twice continuously differentiable convex function on $[a,b].$ For any $\varepsilon >0$, let $\delta :=   \big(64 \|F''\|_{[a,b]}\big)^{-1/2}\min(b-a, \sqrt{\varepsilon}).$ Then  for any convex function $F_1$, we have that 
\begin{equation}\label{eq:part_1_dumb}
\sup_{t \in [a+\delta, b-\delta]} |F_1(t) -F(t)| \ge \varepsilon
\end{equation}
implies that 
\begin{equation}\label{eq:part_2_dumb}
\inf_{t \in [c, c+\delta]} |F_1(t) -F(t)| \ge \varepsilon/4,
\end{equation}
for some  $c \in [a,b-\delta].$
\end{lemma}
\begin{remark}\label{rem:DumbgenLemma}
The above statement is a slight modification of Lemma 2~\cite{DumbgenEtAL04}. However the proof remains the same as the proof does not use the fact that $\hat{F}$ (in the original statement) is a LSE.
\end{remark}

\section[Proof of  the approximate zero equation (\ref{eq:App_score_equation})]{Proof of the approximate zero equation~\eqref{eq:App_score_equation}} % (fold)
\label{sec:proof_of_eq:app_score_equation}

\begin{thm}\label{thm:projection}
Let $\gamma$ be H\"{o}lder exponent of $m_0'$. Under the assumptions of Theorem~\ref{thm:Main_rate_CLSE}, we have
 \begin{equation}\label{eq:App_score_equation1}
   \sqrt{n}\, \p_n \psi_{\check{\theta},\check{m}} =o_p(1).
   \end{equation}
\end{thm}
\begin{proof}
% \todo[inline]{Why are we not assuming the direction $\eta$ satisfies $|\eta| = 1$?}
As described, we show that
\[
\inf_{a\in\mathcal{X}_{\check{m}}}\left|\eta^{\top}\mathbb{P}_n\psi_{\check{\theta}, \check{m}} - \mathbb{P}_n\left[(y-\check{m}(\check\theta^\top x))\{ \eta^\top  \check{m}'(\check\theta ^\top x) H_{\check\theta}^\top x - a(\check\theta^\top x)\}\right]\right| = o_p(n^{-1/2}).
\]
By definition, it is enough to show that
\[
\inf_{a\in\mathcal{X}_{\check{m}}}\left|\mathbb{P}_n\left[(y - \check{m}(\check{\theta}^{\top}x))\left\{a(\check{\theta}^{\top}x) - \check{m}'(\check{\theta}^{\top}x)\eta^{\top}H_{\check{\theta}}^{\top}h_{\theta_0}(\check{\theta}^{\top}x)\right\}\right]\right| = o_p(n^{-1/2}).
\]
For every $\eta\in\mathbb{R}^{d-1}$, define 
\begin{equation}\label{eq:G_def1}
 G_{\eta}(t) := m_0'(t)\eta^{\top}H_{\theta_0}^{\top}h_{\theta_0}(t)
 \end{equation} and
\begin{equation}\label{eq:Gbar_def}
\overline{G}_{\eta}(t) := G_{\eta}(\check{t}_j) + \frac{G_{\eta}(\check{t}_{j+1}) - G_{\eta}(\check{t}_j)}{\check{t}_{j+1} - \check{t}_j}(t - \check{t}_j), \qquad \text{ when } \quad t\in[\check{t}_{j}, \check{t}_{j+1}].
\end{equation}
as a continuous piecewise affine approximation of $G_{\eta}$ with kinks at $\{\check{t}_j\}_{j=1}^\mathfrak{p}$. This implies $\overline{G}_{\eta}\in\mathcal{X}_{\check{m}}$ and hence
\begin{align*}
&\inf_{a\in\mathcal{X}_{\check{m}}}\left|\mathbb{P}_n\left[(y - \check{m}(\check{\theta}^{\top}x))\left\{a(\check{\theta}^{\top}x) - \check{m}'(\check{\theta}^{\top}x)\eta^{\top}H_{\check{\theta}}^{\top}h_{\theta_0}(\check{\theta}^{\top}x)\right\}\right]\right|\\
&\quad\le \left|\mathbb{P}_n\left[(y - \check{m}(\check{\theta}^{\top}x))\left\{\overline{G}_{\eta}(\check{\theta}^{\top}x) - \check{m}'(\check{\theta}^{\top}x)\eta^{\top}H_{\check{\theta}}^{\top}h_{\theta_0}(\check{\theta}^{\top}x)\right\}\right]\right|\\
&\quad\le \left|\mathbb{P}_n\left[(y - \check{m}\circ\check{\theta}(x))\left\{\overline{G}_{\eta}(\check{\theta}^{\top}x) - m_0'(\check{\theta}^{\top}x)\eta^{\top}H_{\theta_0}^{\top}h_{\theta_0}(\check{\theta}^{\top}x)\right\}\right]\right|\\
&\qquad+ \left|\mathbb{P}_n\left[(y - \check{m}\circ\check{\theta}(x))\left\{\check{m}'(\check{\theta}^{\top}x)\eta^{\top}H_{\check{\theta}}^{\top}h_{\theta_0}(\check{\theta}^{\top}x) - m_0'(\check{\theta}^{\top}x)\eta^{\top}H_{\theta_0}^{\top}h_{\theta_0}(\check{\theta}^{\top}x)\right\}\right]\right|\\
&\quad\le \left|\mathbb{P}_n\left[(y - \check{m}\circ\check{\theta}(x))\left\{\overline{G}_{\eta}(\check{\theta}^{\top}x) - m_0'(\check{\theta}^{\top}x)\eta^{\top}H_{\theta_0}^{\top}h_{\theta_0}(\check{\theta}^{\top}x)\right\}\right]\right|\\
&\qquad+ \left|\mathbb{P}_n\left[(y - \check{m}\circ\check{\theta}(x))(\check{m}'(\check{\theta}^{\top}x) - m_0'(\check{\theta}^{\top}x))\eta^{\top}H_{\check{\theta}}^{\top}h_{\theta_0}(\check{\theta}^{\top}x)\right]\right|\\
&\qquad+ \left|\mathbb{P}_n\left[(y - \check{m}\circ\check{\theta}(x))m_0'(\check{\theta}^{\top}x)\eta^{\top}[H_{\check{\theta}} - H_{\theta_0}]^{\top}h_{\theta_0}(\check{\theta}^{\top}x)\right]\right|\\
&\quad= \mathbf{A} + \mathbf{B} + \mathbf{C}.
\end{align*}
The terms $\mathbf{A}, \mathbf{B}$ and $\mathbf{C}$ are all of the form $(y - \check{m}\circ\check{\theta}(x))R(x)$ for a function $R(\cdot)$ that is converging to zero. We split $Y_i - \check{m}\circ\check{\theta}(X_i)$ as $\epsilon_i + (m_0\circ\theta_0 - \check{m}\circ\check{\theta})(X_i)$ and hence,
\begin{align*}
|\mathbb{P}_n[(y - \check{m}\circ\check{\theta}(x))R(x)]| &\le |\mathbb{P}_n[\epsilon R(x)]| + |\mathbb{P}_n[(\check{m}\circ\check{\theta} - m_0\circ\theta_0)(x)R(x)]|.
\end{align*}
Based on this inequality, we write $\mathbf{A} \le \mathbf{A}_1 + \mathbf{A}_2$ and similarly for $\mathbf{B}$ and $\mathbf{C}$. Now observe that
\[
|\mathbb{P}_n[(\check{m}\circ\check{\theta} - m_0\circ\theta_0)(x)R(x)]| ~\le~ \|\check{m}\circ\check{\theta} - m_0\circ\theta_0\|_{n}\|R\|_n.
\]
Using this Cauchy-Schwarz inequality, we get
\begin{align*}
% \mathbf{A}_2 &\le \|\check{m}\circ\check{\theta} - m_0\circ\theta_0\|_{n}\|\overline{G}_{\eta} - m_0'\times \eta^{\top}H_{\theta_0}^{\top}h_{\theta_0}\|_n\qquad\quad ~\overset{(a)}{=}~ O_p(n^{-(2+\gamma)/5})\\
\mathbf{A}_2 &\le \|\check{m}\circ\check{\theta} - m_0\circ\theta_0\|_{n}\|\overline{G}_{\eta} - m_0'\times \eta^{\top}H_{\theta_0}^{\top}h_{\theta_0}\|_n\qquad\quad ~\overset{(a)}{=}~ O_p\big(n^{-2/5 [1 + 2\gamma/(4+\beta)]}\big)\\
\mathbf{B}_2 &\le \|\check{m}\circ\check{\theta} - m_0\circ\theta_0\|_{n}\|\check{m}'\circ\check{\theta} - m_0'\circ\check{\theta}\|_n\|\eta^{\top}H_{\check{\theta}}^{\top}h_{\theta_0}\|_{\infty} ~\overset{(b)}{=}~ O_p(n^{-10/15})\\
\mathbf{C}_2 &\le \|\check{m}\circ\check{\theta} - m_0\circ\theta_0\|_{n}\|H_{\check\theta} - H_{\theta_0}\|_{op}\|h_{\theta_0}\|_{2,\infty}\qquad\qquad\; ~\overset{(c)}{=}~ O_p(n^{-4/5}).
\end{align*}
Equality~(a) follows from Theorem~\ref{thm:rate_m_theta_CLSE} and Lemma~\ref{lem:t_check_max} (stated and proved in the following section) under the assumption that $m_0'$ is $\gamma$-H{\"o}lder continuous. Equality~(b) follows from Theorems~\ref{thm:rate_m_theta_CLSE} and~\ref{thm:rate_derivCLSE}. { Equality~(c) follows from Theorem~\ref{thm:rate_m_theta_CLSE}, the Lipschitzness property of $\theta\mapsto H_{\theta}$, and the boundedness of the covariates (assumption~\ref{aa1}). The calculations above imply that $n^{1/2}\max\{\mathbf{A}_2, \mathbf{B}_2, \mathbf{C}_2\} = o_p(1)$ if $\beta < 8 \gamma-4$. }

We will now prove $n^{1/2}\mathbf{A}_1, n^{1/2}\mathbf{B}_1, n^{1/2}\mathbf{C}_1$ are all $o_p(1)$.  Note that for any function $R(\cdot)$, $n^{1/2}\mathbb{P}_n[\epsilon R(x)] = \mathbb{G}_n[\epsilon R(x)]$ because $\epsilon$ has a zero conditional mean. In Lemma~\ref{lem:score_aprx_1}, we prove $n^{1/2}\mathbf{A}_1 = o_p(1)$. The proof of the other two terms are similar. 

It is easy to see that $n^{1/2}\mathbf{C}_1 = o_p(1)$ because $\|H_{\theta} - H_{\theta_0}\| = O_p(n^{-2/5})$ (by~\cite[Lemma 1, part c]{Patra16} and Theorem~\ref{thm:ratestCLSE}) and $\theta\mapsto m_0'(\theta^{\top}x)\eta^{\top}[H_{\theta} - H_{\theta_0}]^{\top}h_{\theta_0}(\theta^{\top}x)$ is a $\gamma$-H{\"o}lder continuous function which implies $\{x\mapsto m_0'(\theta^{\top}x)\eta^{\top}[H_{\theta} - H_{\theta_0}]^{\top}h_{\theta_0}(\theta^{\top}x)\}$ is a Donsker class. Similarly, one can show that $n^{1/2}\mathbf{B}_1 = o_p(1)$ because $\|\check{m}'\circ\check{\theta} - m_0'\circ\check\theta\| = o_p(1)$ (by Theorem~\ref{thm:rate_derivCLSE}) and $\{x\mapsto({m}'\circ{\theta} - m'_0\circ\theta)(x)\eta^{\top}H_{\theta}^{\top}h_{\theta_0}(\theta^{\top}x):\,\theta\in\Theta\cap B_{\theta_0}(r)\mbox{ and }m'\mbox{ nondecreasing}\}$ is a Donsker class (shown in Lemma~\ref{lem:monosim}).
\end{proof}
\subsection[Lemmas used in the proof of (\ref{eq:App_score_equation})]{Lemmas used in the Proof of~\eqref{eq:App_score_equation}}
\begin{lemma}\label{lem:projection_is_identity}
Define
\[
\mathcal{X}_{\check{m}} := \{a:D\to\mathbb{R}|\, a\mbox{ is piecewise affine continuous function with kinks at }\{\check{t}_i\}_{i=1}^q\}.
\]
Then
\[
\mathcal{X}_{\check{m}} \subseteq \{a:D\to\mathbb{R}\,|\,t\mapsto\xi_t(\cdot; a, \check{m}) \text{ is differentiable at } t=0\}.
\]
\end{lemma}
\begin{proof}
For any function $f$, let $f^L_i$ and $ f^R_i$ denote the left and right derivatives (respectively) at $\check{t}_{i}$. Let $M_a := \max_{i\le q} |a^L_i- a^R_i|$. We know that $\check{m}$ is convex thus for every $i\le q$, $a^L_i <  a^R_i$; here we have the strict inequality because $\{\check{t}_i\}_{i=1}^q$ are set of kinks of $\check{m}$. Let $C_{\check{m}} := \min_{i\le q} (a^R_i- a^L_i)$. Thus for every $|t| \le C_{\check{m}}/M_a $, we have that $ \check{m}- t a $ is convex. Thus $\xi_t(\cdot; a,m)$ is the identity function for every $|t| \le C_{\check{m}}/M_a $ and differentiable at $t=0$ by definition. 
\end{proof}
\begin{lemma}\label{lem:deriv}
For every $a\in\mathcal{X}_{\check{m}}$, we have
\[
-\frac{1}{2}\frac{\partial}{\partial t}Q_n(\zeta_t(\check\theta,\eta), \xi_t(\cdot;a, \check{m}))\Big|_{t = 0} ~=~ \mathbb{P}_n\left[\big(y-\check{m}(\check\theta^\top x)\big)\Big\{ \eta^\top  \check{m}'(\check\theta ^\top x) H_{\check\theta}^\top x - a(\check\theta^\top x)\Big\}\right].
\]
\end{lemma}
\begin{proof}
If $a\in\mathcal{X}_{\check{m}}$, $\Pi_{\mathcal{M}_L}(\check{m} - t a) = \check{m} - t a$ and hence
\begin{align}\label{eq:score_eq_compu_at_a}
\begin{split}
&\left. -\frac{1}{2} \frac{\partial}{\partial t} \big[(y-\xi_t(\zeta_t(\check\theta, \eta)^\top x;  a, \check{m})^2\big]\right|_{t=0}\\
={}&(y-\xi_t(\zeta_t(\check\theta, \eta)^\top x;  a, \check{m}) \left.\frac{\partial \xi_t(\zeta_t(\check\theta, \eta)^\top x; a,\check{m})}{\partial t}\right|_{t=0}\\
={}&\big(y-\check{m}(\check\theta^\top x)\big)\Big[ \eta^\top  \check{m}'(\check\theta ^\top x) H_{\check\theta}^\top x - a(\check\theta^\top x)\Big].\qedhere
\end{split}
\end{align}
\end{proof}

{
\begin{lemma}[Property of $\{\check{t}_i\}_{i=1}^\mathfrak{p}$]\label{lem:t_check_max} If the assumptions of~Theorem~\ref{thm:ratestCLSE} hold, then
\begin{equation}\label{eq:max_t}
n^{4/5}\sum_{i=1}^\mathfrak{p} (\check{t}_{i+1}- \check{t}_{i})^{5+\beta} =O_p(1) \qquad \text{and} \qquad \max_{1\le j\le \mathfrak{p}}|\check{t}_{j+1} - \check{t}_j| = O_p(n^{-4/(25 + 5\beta)}).
\end{equation}
% where $c_{m_0}$ is the strong convexity parameter.
Furthermore, for any function $G$ that is $\gamma$-H{\"o}lder continuous, the approximating function $\bar{G}$ defined as
\[
\bar{G}(t) = G(\check{t}_i) + \frac{G(\check{t}_{i+1}) - G(\check{t}_{i})}{\check{t}_{i+1} - \check{t}_i}(t - \check{t}_i),\quad\mbox{for}\quad t\in[\check{t}_j, \check{t}_{j+1}],
\]
satisfies
{ \[
\frac{1}{n}\sum_{i=1}^n \left(G(\check{\theta}^{\top}X_i) - \bar{G}(\check{\theta}^{\top}X_i)\right)^2 = O_p(n^{-8\gamma/(20+5\beta)})\quad\mbox{for}\quad \gamma\in[0,2].
\]}
\end{lemma}
\begin{proof}
Recall the definition of $\{\check{t}_i\}_{i=1}^\mathfrak{p}$ in Page~\pageref{eq:X_m_def} of the primary document. Note that $D_0$ is an interval, $\mathfrak{p}\le n$, and $\check{t}_i\in D_{\check{\theta}}$, for all $1\le i\le \mathfrak{p}$. However, by Theorem~\ref{thm:ratestCLSE}, we have that $|\check{\theta}-\theta_0| =O_p(n^{-2/5})$. Thus $\Lambda(\text{conv}(D_{\check{\theta}})\setminus D_{\theta_0}) = O_p(n^{-2/5})$. Thus to show~\eqref{eq:max_t}, we can assume without loss of generality that for all $1\le i\le \mathfrak{p}$, we have $\check{t}_i\in D_{\theta_0}.$

% \paragraph{Short proof of the above:}
% \begin{align}\label{eq:1}
% \begin{split}
% &n^{4/5}\sum_{i=1}^\mathfrak{p} (\check{t}_{i+1}- \check{t}_{i})^{5+\beta}\\
% ={}& n^{4/5}\sum_{i=1}^\mathfrak{p} (\check{t}_{i+1}- \check{t}_{i})^{5+\beta} \mathbf{1}_{\check{t}_{i+1},\check{t}_{i} \in D_0}+ n^{4/5}\sum_{i=1}^\mathfrak{p} (\check{t}_{i+1}- \check{t}_{i})^{5+\beta} \mathbf{1}_{\check{t}_{i+1},\check{t}_{i} \in D_{\check{\theta}}\Delta D_0} \\
% &\quad + n^{4/5}\sum_{i=1}^\mathfrak{p} (\check{t}_{i+1}- \check{t}_{i})^{5+\beta} \mathbf{1}_{\check{t}_{i+1} \in D_0 \text{ and } \check{t}_{i} \in D_{\check{\theta}}\Delta D_0}
% \end{split}
% \end{align}
% Observe that 
% \[ n^{4/5}\sum_{i=1}^\mathfrak{p} (\check{t}_{i+1}- \check{t}_{i})^{5+\beta} \mathbf{1}_{\check{t}_{i+1},\check{t}_{i} \in D_{\check{\theta}}\Delta D_0} \le n^{4/5} \mathfrak{p} [n^{-2/5}]^{5+ \beta} \le n^{9/5}  [n^{(-10-\beta)/5}]= o_p(1) \]

% Let $[a_0, b_0]$ be the boundary of $D_0$, 
% \[n^{4/5} (\check{t}_{i+1}- \check{t}_{i})^{5+\beta} \lesssim n^{4/5} [(\check{t}_{i+1}- a_0)^{5+\beta} + (a_0- \check{t}_{i})^{5+\beta})]
% \le n^{4/5} [(\check{t}_{i+1}- a_0)^{5+\beta} + n^{(-2/5)(5+\beta)}] \]
Observe that  by Theorem~\ref{thm:ratestCLSE}, we have triangle inequality, 
\[\|\check{m}\circ{\theta_0}- m_0\circ{\theta_0}\| =O_p(n^{-2/5}).\]
Thus, for every $\varepsilon>0,$ we have that there exist a $K_\varepsilon$ such that 
\begin{equation}\label{eq:K_ep_def}
 \mathbb{P}\big(\|\check{m}\circ{\theta_0}- m_0\circ{\theta_0}\| \le K_\varepsilon n^{-2/5}) \ge 1-\varepsilon.
 \end{equation}
 { Thus all of the following inequalities hold with at least $1-\varepsilon$ probability:
 \begin{align}\label{eq:max_kink}
 \begin{split}
 K_\varepsilon \ge{}& n^{4/5} \sum_{i=1}^\mathfrak{p} \int_{\check{t}_i}^{\check{t}_{i+1}} (\check{m}(t) -m_0(t))^2 dP_{\theta_0^\top X}(t)\\
 \ge{}& n^{4/5}\sum_{i=1}^\mathfrak{p} \int_{\check{t}_i}^{\check{t}_{i+1}} ( a_i + b_i t  -m_0(t))^2 dP_{\theta_0^\top X}(t),
 \end{split}
 \end{align}
 where $a_i$ and $b_i$ is such that for every $1\le i\le \mathfrak{p}$, $\check{m}(t) = a_i + b_i t$ for all $t\in [\check{t}_i, \check{t}_{i+1}).$
 Further, by the $\kappa_{m_0}$-strong convexity of $t\mapsto m_0(t) - a_i - b_i(t)$, Theorem~\ref{thm:strong-convex-L2-norm} implies
 \[
 \int_{\check{t}_{i}}^{\check{t}_{i+1}}|m_0(t) - a_i - b_i t|^2 dP_{{\theta_0}^{\top}X}(t) \ge \frac{\underline{C}_d\kappa_{m_0}^2}{2^{2+\beta}3^{10 + 2\beta}} (\check{t}_{i+1}- \check{t}_{i})^{5 + \beta} =: c_{m_0}(\check{t}_{i+1}- \check{t}_{i})^{5 + \beta},
 \]
 for a constant $c_{m_0}$ depending only on $\underline{C}_d, \kappa_{m_0}$, and $\beta$.
 The proof of first part of~\eqref{eq:max_t} is now complete, because
 \begin{align}\label{eq:max_kink_1}
 \begin{split}
 K_\varepsilon \ge{}& n^{4/5}\sum_{i=1}^\mathfrak{p} \int_{\check{t}_i}^{\check{t}_{i+1}} ( a_i + b_i t  -m_0(t))^2 dP_{\theta_0^\top X}(t)
 \ge{} c_{m_0}n^{4/5}\sum_{i=1}^\mathfrak{p} (\check{t}_{i+1}- \check{t}_{i})^{5+\beta}.
 \end{split}
 \end{align}}
% Thus completing the proof of first part of~\eqref{eq:max_t}. 
To prove the second inequality in~\eqref{eq:max_t} observe that  as $\check{t}_{i}\le  \check{t}_{i+1} $ for all $1\le i < \mathfrak{p}$, we have that
\[
n^{4/5}\max_{1\le i \le \mathfrak{p}} (\check{t}_{i+1}- \check{t}_{i})^{5+\beta} \le n^{4/5}\sum_{i=1}^\mathfrak{p}  (\check{t}_{i+1}- \check{t}_{i})^{5+\beta} =O_p(1).
\]
Thus 
$\max_{1\le i \le \mathfrak{p}} |\check{t}_{i+1}- \check{t}_{i}| = O_p(n^{-4/(25 + 5\beta)}).$%\qedhere

To prove the second part of the result, define for $t\in[\check{t}_i, \check{t}_{i+1}]$,
\[
g(t) := G(t) - \bar{G}(t) = G(t) - G(\check{t}_i) - \frac{G(\check{t}_{i+1}) - G(\check{t}_i)}{\check{t}_{i+1} - \check{t}_i}(t - \check{t}_i).
\]
If $\gamma\in(0, 1]$, then there exists $C_G\in(0, \infty)$ such that for every $t\in [ \check{t}_i, \check{t}_{i+1}]$, we have 
\begin{equation}\label{eq:gamma_Holder}
|G(t) - G(\check{t}_i)| \le C_G|t - \check{t}_i|^{\gamma}\quad\Rightarrow\quad |g(t)| \le 2C_G|{t} - \check{t}_{i}|^{\gamma} \le 2C_G|\check{t}_{i+1} - \check{t}_{i}|^{\gamma}.
\end{equation}
If $\gamma\in[1, 2]$, then there exists $C_G\in(0, \infty)$ such that
\begin{equation}\label{eq:1plusgamma_Holder}
\sup_{a\neq b}\frac{|G'(b) - G'(a)|}{|b - a|^{\gamma-1}} \le C_G\quad\Rightarrow\quad |g(t)| \le 2C_G|\check{t}_{i+1} - \check{t}_i|^{\gamma},
\end{equation}
because
\begin{align*}
|g(t)| &= \left|G(t) - G(\check{t}_i) - G'(\check{t}_i)(t - \check{t}_i) + (t - t_i)\left[G'(\check{t}_i) - \frac{G(\check{t}_{i+1}) - G(\check{t}_i)}{\check{t}_{i+1} - \check{t}_i}\right]\right|\\
&\le \left|G(t) - G(\check{t}_i) - G'(\check{t}_i)(t - t_i)\right| + |t - \check{t}_i|\times\left|G'(\check{t}_i) - \frac{G(\check{t}_{i+1}) - G(\check{t}_i)}{\check{t}_{i+1} - \check{t}_i}\right|\\
&\le C_G|t - \check{t}_i|^{\gamma} + C_G|t - \check{t}_i||\check{t}_{i+1} - \check{t}_i|^{\gamma - 1} \le 2C_G|\check{t}_{i+1} - \check{t}_i|^{\gamma}.
\end{align*}
% Therefore, for $\check{\theta}^{\top}X_i\in[\check{t}_j, \check{t}_{j+1}]$,
% \begin{align*}
% |g(\check{\theta}^{\top}X_i)| &\le C_1(f)|\check{\theta}^{\top}X_i - \check{t}_j|^2 + C_1(f)|\check{\theta}^{\top}X_i - \check{t}_j||\check{t}_{j+1} - \check{t}_j|\\
% &\le C_1(f)|\theta_0^{\top}X_i - \check{t}_j|^2
% \end{align*}
This yields for any $\gamma\in(0, 2]$,
{ \begin{align}
\frac{1}{n}\sum_{i=1}^n \left(G(\check{\theta}^{\top}X_i) - \overline{G}(\check{\theta}^{\top}X_i)\right)^2 &= \frac{1}{n}\sum_{i=1}^n g^2(\check{\theta}^{\top}X_i)\nonumber\\
&\le \left(\frac{1}{n}\sum_{i=1}^n g^{(4+\beta)/\gamma}(\check{\theta}^{\top}X_i)\right)^{2\gamma/(4 + \beta)}\quad\mbox{(since $\gamma \le 2$, $4/\gamma \ge 2$)}\nonumber\\
&= \left(\sum_{j=1}^{\mathfrak{p}}\frac{1}{n}\sum_{i:\check{\theta}^{\top}X_i\in[\check{t}_j, \check{t}_{j+1}]} g^{(4+\beta)/\gamma}(\check{\theta}^{\top}X_i)\right)^{2\gamma/(4 + \beta)}\label{eq:FirstBound_f_and_fbar}\\ 
&\le \left(\sum_{j=1}^\mathfrak{p} \frac{c_j(2C_G)^{(4+\beta)/\gamma}}{n}|\check{t}_{j+1} - \check{t}_j|^{4 + \beta}\right)^{2\gamma/(4 + \beta)}\\ &= 4C_G^2\left(\sum_{j=1}^\mathfrak{p} \frac{c_j}{n}|\check{t}_{j+1} - \check{t}_j|^{4 + \beta}\right)^{2\gamma/(4 + \beta)},\nonumber
\end{align}
}
where $c_j$ denotes the number of observations $\check{\theta}^{\top}X_i$ that fall into $[\check{t}_j, \check{t}_{j+1}]$. Because $|\check{\theta} - \theta_0| = O_p(n^{-2/5})$ by Theorem~\ref{thm:ratestCLSE}, we get that with probability converging to one, $|\check{\theta} - \theta_0| \le n^{-2/5}\sqrt{\log n}$ holds true. On this event, for any $1\le j\le \mathfrak{p}$,
\begin{align*}
\frac{c_j}{n} = \frac{1}{n}\sum_{i=1}^n \mathbbm{1}\{\check{\theta}^{\top}X_i\in[\check{t}_j, \check{t}_{j+1}]\} &\le \frac{1}{n}\sum_{i=1}^n \mathbbm{1}\{\theta_0^{\top}X_i\in[\check{t}_{j} - |(\check{\theta} - \theta_0)^{\top}X_i|, \check{t}_{j+1} + |(\check{\theta} - \theta_0)^{\top}X_i|]\}\\
&\le \frac{1}{n}\sum_{i=1}^n \mathbbm{1}\{\theta_0^{\top}X_i\in[\check{t}_j - Tn^{-2/5}\sqrt{\log n}, \check{t}_{j+1} + Tn^{-2/5}\sqrt{\log n}]\}\\
&\le P\mathbbm{1}\{\theta_0^{\top}X\in [\check{t}_j - Tn^{-2/5}\sqrt{\log n}, \check{t}_{j+1} + Tn^{-2/5}\sqrt{\log n}]\}\\ &\qquad+ 2n^{-1/2}\sup_{a\in\mathbb{R}}{|\mathbb{G}_n\mathbbm{1}\{\theta_0^{\top}X \le a\}|
}.%{\sqrt{Var(\mathbbm{1}\{\theta_0^{\top}X\in[...]\})} + n^{-1/2}}.
\end{align*}
Corollary 1 of~\cite{massart1990tight} implies that with probability converging to one, $\sup_{a}|\mathbb{G}_n\mathbbm{1}\{\theta_0^{\top}X \le a\}| \le 0.5\sqrt{\log n}$. Further~\ref{aa6} yields
\[
P\mathbbm{1}\{\theta_0^{\top}X\in [\check{t}_j - Tn^{-2/5}\sqrt{\log n}, \check{t}_{j+1} + Tn^{-2/5}\sqrt{\log n}]\} ~\le~ { \overline{C}_0}\left[\check{t}_{j+1} - \check{t}_j + 2Tn^{-2/5}\sqrt{\log n}\right].
\]
Hence with probability converging to one, simultaneously for all $1\le j\le \mathfrak{p}$, we have
\[
\frac{c_j}{n} \le \overline{C}_0|\check{t}_{j+1} - \check{t}_j| + Tn^{-2/5}\sqrt{\log n} + n^{-1/2}\sqrt{\log n} \le { \overline{C}_0}|\check{t}_{j+1} - \check{t}_j| + (T + 1)n^{-2/5}\sqrt{\log n}.
\]
% \todo[inline]{need to change the notation $\overline{C}_0$ for a bound on the density.??}
{ Therefore~\eqref{eq:FirstBound_f_and_fbar} yields with probability converging to one,
\begin{align*}
\frac{1}{n}\sum_{i=1}^n \left(G(\check{\theta}^{\top}X_i) - \overline{G}(\check{\theta}^{\top}X_i)\right)^2 &\le 4C_G^2\left({ \overline{C}_0}\sum_{j=1}^\mathfrak{p} |\check{t}_{j+1} - \check{t}_j|^{5+\beta}\right)^{2\gamma/(4 + \beta)}\\ &\quad+ 4C_G^2\left(\frac{(T + 1)\sqrt{\log n}}{n^{2/5}}\sum_{j=1}^\mathfrak{p} |\check{t}_{j+1} - \check{t}_j|^{4+\beta}\right)^{2\gamma/(4 + \beta)}\\
&= O_p(n^{-{8\gamma}/{(20 + 5\beta)}}) + O_p((\log n)^{\gamma/(4+\beta)}n^{-4\gamma(11 + 3\beta)/(5(5+\beta)(4+\beta))})\\
&= O_p(n^{-8\gamma/(20+5\beta)}).
\end{align*}
%\todo[inline]{We need to use the lower  bound from the ``convex proof file'' and Assumption (D) of~\cite{MR2369025} with $\beta\ge 0$ to show that $\frac{1}{n}\sum_{i=1}^n \left(G(\check{\theta}^{\top}X_i) - \overline{G}(\check{\theta}^{\top}X_i)\right)^2 =O_p(n^{-2\gamma/5})$ as long as $m_0$ is 2-H\"{o}lder instead of 3/2-H\"{o}lder (our current assumption).}
The first equality above holds because~\eqref{eq:max_t} yields 
$$
\sum_{j=1}^\mathfrak{p} |\check{t}_{j+1} - \check{t}_j|^{4+\beta} \le \max_{1\le j\le \mathfrak{p}}|\check{t}_{j+1} - \check{t}_j|^{3+\beta}\sum_{j=1}^\mathfrak{p} |\check{t}_{j+1} - \check{t}_j| = O_p(n^{-4(3+\beta)/(25+5\beta)}).
$$
This completes the proof.
%{ Somehow we need to conclude that $c_j/n$ behaves like $|\check{t}_{j+1} - \check{t}_j|$. If $c_j$ is defined for $\theta_0^{\top}X_i$, then this is trivially true by using KS statistic bounds. 
%(As long as $|\check{t}_{j+1} - \check{t}_j|$ does not go to zero faster than $n^{-2/5}$, the dominating term on the right hand side is $|\check{t}_{j+1} - \check{t}_j|$.)
% }
}
\end{proof}

{ 
\begin{lemma}\label{lem:special-case}
Suppose $f:[a,b] \to \mathbb{R}$ is a $\lambda$-strongly convex function such that either $\inf_{x\in[a,b]}\,f(x) \ge 0$ or $\sup_{x\in[a,b]}\,f(x) \le 0$ holds true. Let $\mu$ be any probability measure such that for some $\beta \ge 0$ and all intervals $I$, $\mu(I) \ge \underline{c}|I|^{1 + \beta}$, where $|I|$ represents the Lebesgue measure of $I$. Then
\[
\int_a^b f^2(x)d\mu(x) ~\ge~ \frac{\underline{c}\lambda^2(b - a)^{5 + \beta}}{2^{2 + \beta}3^{5 + \beta}}.
\] 
\end{lemma}
\begin{proof}
Consider the case when $\inf_{x\in[a,b]}\,f(x) \ge 0$, that is, $f(x) \ge 0$ for all $x\in[a,b]$. If $f'(a) \ge 0$, then
\[
f(x) \ge f(a) + f'(a)(x - a) \ge 0\quad\mbox{for all}\quad x \in [a,b].
\] 
Note that $x\mapsto f(a) + f'(a)(x - a)$ is non-decreasing because $f'(a) \ge 0$ and is non-negative at $x = a$; this proves the second inequality above. Therefore, 
\[
f(x) - 0 \ge f(x) - \{f(a) - f'(a)(x - a)\} \ge \frac{\lambda}{2}(x - a)^2,
\] 
where the last inequality follows from $\lambda$-strong convexity of $f$. This implies that if $f'(a) \ge 0$,
\begin{align*}
\int_a^b f^2(x)d\mu(x) &\ge \frac{\lambda^2}{4}\int_a^b (x - a)^4d\mu(x)\\ 
&\ge \frac{\lambda^2(b - a)^4}{4(81)}\mu([(2a + b)/3, b])\\
&\ge \frac{\underline{c}\lambda^2(b-a)^{5 + \beta}}{4(3)^{5 + \beta}}.
\end{align*}
If, instead, $f'(b) \le 0$, then the same argument works except for the change
\[
f(x) \ge f(b) + f'(b)(x - b) \ge 0\quad\mbox{for all}\quad x\in[a, b].
\]
If, instead, $f'(a) < 0 < f'(b)$, then there exists a point $x^*\in[a, b]$ such that $f'(x^*) = 0$. Hence,
\[
f(x) \ge f(x^*) + f'(x^*)(x - x^*) = f(x^*) \ge \inf_{x\in[a,b]}f(x) \ge 0\quad\mbox{for all}\quad x\in[a, b],
\]
which implies that
\[
f(x) - 0 \ge f(x) - \{f(x^*) - f'(x^*)(x - x^*)\} \ge \frac{\lambda}{2}(x - x^*)^2.
\]
Therefore, for $I = \{x\in[a, b]:\,|x - x^*| \ge (b - a)/3\}$
\[
\int_a^b f^2(x)d\mu(x) \ge \frac{\lambda^2}{4}\int_a^b (x - x^*)^4d\mu(x) \ge \frac{\lambda^2(b - a)^4}{4(81)}\mu(I).
\]
Note that $I\subseteq[a,b]$ is a union of at most two intervals. One of which will have Lebesgue measure of at least $(b - a)/3$. Thus, $\mu(I) \ge 2^{-\beta}\underline{c}((b-a)/3)^{1 + \beta}$. Hence,
\[
\int_a^b f^2(x)d\mu(x) \ge \frac{\lambda^2}{4}\int_a^b (x - x^*)^4d\mu(x) \ge \frac{\underline{c}\lambda^2(b - a)^{5 + \beta}}{2^{2 + \beta}(3)^{5 + \beta}}.
\] 
This completes the result when $\inf_{x\in[a,b]}\,f(x) \ge 0$.

Now consider the case where $\sup_{x\in[a,b]}f(x) \le 0$. In this case,
\[
f(x) \le \ell(x) := f(a) + \frac{f(b) - f(a)}{b - a}(x - a) = f(a)\left(\frac{b-x}{b - a}\right) + f(b)\left(\frac{x - a}{b - a}\right) \le 0.
\]
Hence using the equivalent definition $f(\alpha x + (1 - \alpha)y) \le \alpha f(x) + (1 - \alpha)f(y) - \alpha(1-\alpha)\lambda(x - y)^2/2$, we conclude
\begin{align*}
\int_a^b \{0 - f(x)\}^2d\mu(x) &\le \int_a^b \{\ell(x) - f(x)\}^2d\mu(x)\\ 
&\le \int_a^b \frac{\lambda^2(b - x)^2(x - a)^2}{4(b-a)^4}(b - a)^4d\mu(x)\\
&= \frac{\lambda^2}{4}\int_a^b (b - x)^2(x - a)^2d\mu(x)\\
&\ge \frac{\lambda^2}{4}\int_{(2a + b)/3}^{(a + 2b)/3} (b - x)^2(x - a)^2d\mu(x)\\
&\ge \frac{\lambda^2(b - a)^4}{4(81)}\mu\left(\left[\frac{2a + b}{3},\,\frac{a + 2b}{3}\right]\right)\\
&\ge \frac{\underline{c}\lambda^2(b-a)^{5 + \beta}}{4(3)^{5 + \beta}}.
\end{align*}
Combining all the cases, we conclude the proof.
\end{proof}

\begin{thm}\label{thm:strong-convex-L2-norm}
Suppose $f:[a,b]\to\mathbb{R}$ is a $\lambda$-strongly convex function. Let $\mu$ be any probability measure such that for some $\beta > 0$ and all intervals $I$, $\mu(I) \ge \underline{c}|I|^{1 + \beta}$, where $|I|$ represents the Lebesgue measure of $I$. Then
\[
\int_a^b f^2(x)d\mu(x) \ge \frac{\underline{c}\lambda^2(b - a)^{5 + \beta}}{2^{2 + \beta}3^{10 + 2\beta}}.
\]
\end{thm}
\begin{proof}
If $f(x), x\in[a,b]$ is wholly above or below zero, the result follows from Lemma~\ref{lem:special-case}. Otherwise, the function $f$ on $[a, b]$ intersects the $x$-axis at no more than two points, let they be $a'$ and $b'$; if it only intersects at one point, take $a' = b'$. The function does not change its sign in the intervals $[a, a'], [a', b']$ and $[b', b]$. By virtue, at least one of $[a, a'], [a', b']$ or $[b', a]$ has to have Lebesgue measure of at least $(b - a)/3$. Therefore applying Lemma~\ref{lem:special-case} in largest of these intervals proves the result.
\end{proof}}

\begin{lemma}\label{lem:score_aprx_1} If the assumptions of~Theorem~\ref{thm:ratestCLSE} hold and $\gamma$ is the H\"{o}lder exponent of $m_0'$, then $\sqrt{n}\mathbf{A}_1= o_p(1).$
% \todo[inline]{ What can we do here? Does monotonicity work? what about $k$-monotone? What conditions do we need for this and the other term to work.}

\end{lemma}
\begin{proof}
 For any real-valued function $h: [a, b] \to \R$, let $V_\alpha(h)$ denote the $\alpha$-variation of $h$ i.e., 
\[ V_\alpha(h) := \sup\{ \sum_{i=1}^n |h(x_i)- h(x_{i-1})|^\alpha : a= x_0 < x_1<\cdots<x_n= b, n\in \mathbb{N}\}.\]
 We will now show that both $G_\eta$ and $\overline{G}_{\eta}$ (defined in~\eqref{eq:G_def} and~\eqref{eq:Gbar_def} respectively) are bounded $\alpha$-variation functions.

Recall that $G_\eta$ is $\gamma$-H\"{o}lder and is defined on a bounded interval, thus by definition it has bounded $1/\gamma$-variation; see e.g., \citet[Page 220-221]{Gine16}.  Now, observe that $\overline{G}_{\eta}$ is a piecewise linear function with kinks at $\{\check{t}_{i}\}_{i=1}^\mathfrak{p}.$ Thus we have that 
\begin{align}
\begin{split}
V_{1/\gamma}({\overline{G}_\eta}) = \sum_{j = 1}^\mathfrak{p} |\overline{G}_\eta(\check{t}_j)- \overline{G}_\eta(\check{t}_{j+1})|^{1/\gamma}\le 2C_G\sum_{j = 1}^\mathfrak{p} |\check{t}_{j+1} - \check{t}_j| &\le 2C_G \diameter(D),
\end{split}
\end{align}
where for the second inequality we use~\eqref{eq:gamma_Holder}. 
Let
\begin{equation}\label{eq:f_def}
f_\eta(t):= \overline{G}_{\eta}(t) - G_{\eta}(t),
\end{equation}
Because $ \gamma <1$, we have that $V_{1/\gamma}({\overline{G}_\eta} -{G}_\eta ) \le 2^{1/\gamma-1} (V_{1/\gamma}({\overline{G}_\eta}) + V_{1/\gamma}({{G}_\eta}))$. Thus, $f$ has bounded $1/\gamma$-variation. For any $\alpha>1$, let us now define
\begin{align}\label{eq:F_gamma_def}
\mathcal{F}_{\alpha}(K):= \{g: \rchi\to \R|\; &g(x)= f(\theta^\top x), \theta \in \Theta \cup B_{\theta_0}(r)\\ 
& \text{ and } f: D\to \R \text{ is a bounded $\alpha$-variation function with } V_\alpha(f) \le K\}.
\end{align}
 In Lemma~\ref{lem:ent_alpa_variation}, we show that 
$\log N_{[\,]}(\eta, \f_{1/\gamma}(K), \|\cdot\|) \le C \eta^{-1/\gamma}$ for some constant $C$ depending on $K$ only. Because $1/2 <\gamma<1,$ we have that $\f_{1/\gamma}(K)$ is Donsker.  Furthermore,  by~\eqref{eq:gamma_Holder}, there exists a constant $C$ such that 
\begin{align}
\begin{split}
\int f^2(t) dt &\le 2C^2 \sum_{j=1}^\mathfrak{p}\int_{\check{t}_j}^{\check{t}_{j+1}} (t - \check{t}_j)^{2 \gamma}dt\le 2 C^2 \max_{1\le j\le \mathfrak{p}}|\check{t}_{j+1} - \check{t}_j|^{2 \gamma} \sum_{j=1}^\mathfrak{p} ( \check{t}_{j+1}- \check{t}_j)dt =O_p(n^{-8 \gamma/25}),
\end{split}
\end{align} 
and by~\eqref{eq:gamma_Holder}, we have that 
\begin{equation}\label{eq:f_inf_bound}
\|f\|_{\infty} \le 2 C \max_{j\le q} |\check{t}_{j+1} - \check{t}_j|^{\gamma} =O_p\big(n^{-4\gamma  /25}\big).
\end{equation}
Because $q\ge 5$,  by Lemma~\ref{lem:Maximal342}, we have that $\sqrt{n}\mathbf{A}_1 =o_p(1)$. \qedhere

\end{proof}

\subsection{Metric entropies for  monotone and bounded $\alpha$-variation single index model}
In Lemma~\ref{lem:score_aprx_1}, we need to find the entropy of the following class of the functions:
\begin{align}\label{LIP:w_def}
\begin{split}
\mathcal{H}^*(S)&=\{q:\rchi\to\R|\, q(x)=g(\theta^\top x),  \theta\in \Theta\cap B_{\theta_0}(r)\, \text{ and }\\ &\qquad \quad g:D\to \R \text{ is a nondecreasing function and }  \|g\|_\infty \le S \}.\\
% \w^*_{M_1} &:=\left\lbrace U_{ \theta, m}: (\theta, m) \in \mathcal{C}^*_{M_1} \right\rbrace,\\%
% \w_{M_1}(n) &:=\left\lbrace U_{ \theta, m}: (\theta, m) \in \mathcal{C}_{M_1}(n) \right\rbrace,
\end{split}
\end{align}
% where $U_{ \theta, m}(\cdot)$ is defined in \eqref{eq:def_U}.  
 % In the following two lemmas, we find the bracketing numbers (with respect to the $L_2(P_{\theta_0, m_0})$ norm) and envelope functions for these three classes of functions.

\begin{lemma}\label{lem:monosim}   $\log N_{[\;]}(\varepsilon, \mathcal{H}^*(S), L_2(P_{\theta_0, m_0})) \lesssim S \varepsilon^{-1}.$ for all $\varepsilon >0.$
\end{lemma}
\begin{proof}
First recall that by assumption~\ref{aa6}, { we have that  $ \sup_{\theta\in \Theta \cap B_{\theta_0}(r)} \|f_{\theta^\top X}\|_D \le2\overline{C}_0 < \infty,$} where $f_{\theta^\top X}$ denotes the density of ${\theta^\top X}$ with respect to the Lebesgue measure.  To compute the entropy of $\mathcal{H}^*(S)$, note that by Lemma 4.1 of \cite{Pollard90} we can get $\theta_1, \theta_2, \ldots, \theta_{N_{\eta_1}},$ with $N_{\eta_1}\le 3^d T^d \eta_1^{-d}$ such that for every $\theta\in\Theta$, there exists a $j$ satisfying $|\theta - \theta_j|\le \eta_1/T$ and
\[
|\theta^{\top}x - \theta_j^{\top}x| \le |\theta - \theta_j|\cdot|x| \le \eta_1\quad \forall x\in\rchi.
\]
Thus for every $\theta\in\Theta$, we can find a $j$ such that $\theta_j^{\top}x - \eta_1 \le \theta^{\top}x \le \theta_j^{\top}x + \eta_1, \forall x\in \rchi.$  For simplicity of notation, define
% \begin{equation}\label{eq:tj_def}
$t_j^{(1)}(x) := \theta_j^{\top}x - \eta_1,$ $t_j^{(2)}(x) := \theta_j^{\top}x + \eta_1,$
% \end{equation}
and
\[\mathcal{G}^*_S:=\{g|\, g:D\to \R \text{ is a uniformly bounded nondecreasing function and }  \|g\|_\infty \le S \}.\]  Recall that $\Lambda$ denotes the Lebesgue measure on $D$. By a simple modification of Theorem 2.7.5 of~\cite{VdVW96}, we have that \[
  N_{[\;]}(\eta_2, \mathcal{G}^*_S, L_2(\Lambda)) \le \exp\left(\frac{A S \sqrt{\text{diam}(\D)}}{\eta_2^{-1}}\right):=  M_{\eta_2},
\] for some universal constant $A.$ Thus there exist $\{[l_1, u_1]\}_{i=1}^{M_{\eta_2}}$  in $ \mathcal{G}^*_S$ with $l_i\le u_i$ and $\int_D |u_i(t) - l_i(t)|^2 dt \le \eta_2^2$  such that for every $g\in\mathcal{G}^*_S$, we can find a $m\in  \{1, \ldots,M_{\eta_2}\}$ such that $l_m \le g \le u_m$. Fix any function $g\in\mathcal{G}^*_S$ and $\theta\in\Theta$. Let $|\theta_j - \theta| \le \eta_1/T$ and let $l_k \le g \le u_k$, then for every $x\in\rchi$,
\[
l_k(t_j^{(1)}(x)) \le l_k(\theta^{\top}x) \le g(\theta^{\top}x) \le u_k(\theta^{\top}x) \le u_k(t_j^{(2)}(x)),
\]
where the outer inequalities follow from the fact that both $l_k$ and $u_k$ are nondecreasing functions. Proof of Lemma~\ref{lem:monosim} will be complete if we can show that
\[
\{[l_k\circ t_j^{(1)}, u_k\circ t_j^{(2)}]: 1\le j\le N_{\eta_1}, 1\le k\le M_{\eta_2}\},
\] form a $L_2(P_{\theta_0,m_0})$ bracket for $\mathcal{H}^*(S).$  Note that by the triangle inequality, we have
\begin{equation}\label{eq:triangleineq}
\|u_k\circ t_j^{(2)} - l_k\circ t_j^{(1)}\| \le \|u_k\circ t_j^{(2)} - l_k\circ t_j^{(2)}\| + \|l_k\circ t_j^{(2)} - l_k\circ t_j^{(1)}\|.
\end{equation}
Since the density  of $X^{\top}\theta$  with respect to the Lebesgue measure is  bounded uniformly (for $\theta \in \Theta\cap B_{\theta_0}(r)$)  by $\overline{C}_0$, we get that
\[
\|u_k\circ t_j^{(2)} - l_k\circ t_j^{(2)}\|^2 = \int \left[u_k(r) - l_k(r)\right]^2 f_{\theta_j^{\top}X}(r)dr \le \overline{C}_0\int \left[u_k(r) - l_k(r)\right]^2 dr \le \overline{C}_0\eta_2^2.
\]
For the second term in \eqref{eq:triangleineq}, we first approximate the lower bracket $l_k$ by a right-continuous nondecreasing step (piecewise constant) function. Such an approximation is possible since the set of all simple functions is dense in $L_2(P_{\theta_0,m_0})$; see Lemma 4.2.1 of \cite{Bogachev}. Since $l_k$ is bounded by $S$, we can get a nondecreasing step function $A: D\to [-S, S],$ such that $\int \{l_k(r) - A(r)\}^2dr \le \eta_2^2$. Let $v_1<\cdots< v_{A_d}$ denote an points of discontinuity of $A$.
Then for every $r\in \D,$ we can write
\[
A(r) = -S+ \sum_{i =1}^{A_d} c_i\mathbbm{1}_{\{r\ge v_i\}},\;\mbox{where}\; c_i>0 \;\text{and}\; \sum_{i=1}^{A_d}c_i \le 2S.
\]
 Using triangle inequality, we get that
\begin{align*}
\|l_k\circ t_j^{(2)} - l_k\circ t_j^{(1)}\| &\le \|l_k\circ t_j^{(2)} - A\circ t_j^{(2)}\| + \|A\circ t_j^{(2)} - A\circ t_j^{(1)}\| + \|A\circ t_j^{(1)} - l_k\circ t_j^{(1)}\|\\
&\le \sqrt{\overline{C}_0}\eta_2 + \|A\circ t_j^{(2)} - A\circ t_j^{(1)}\| + \sqrt{\overline{C}_0}\eta_2.
\end{align*}
 Now observe that
% \todo[inline]{need this in the statement of theorem}
\begin{align*}
\|A\circ t_j^{(2)} - A\circ t_j^{(1)}\|^2 &= \mathbb{E}\left[\sum_{i =1}^{A_d} c_i\left(\mathbbm{1}_{\{X^{\top}\theta_j + \eta_1 \ge v_i\}} - \mathbbm{1}_{\{X^{\top}\theta_j - \eta_1 \ge v_i\}}\right)\right]^2\\
&\le 2S\mathbb{E}\left|\sum_{i =1}^{A_d}c_i\left(\mathbbm{1}_{\{X^{\top}\theta_j + \eta_1 \ge v_i\}} - \mathbbm{1}_{\{X^{\top}\theta_j - \eta_1 \ge v_i\}}\right)\right|\\
&\le 2S\sum_{i =1}^{A_d} c_i\mathbb{P}(X^{\top}\theta_j -\eta_1 < v_i \le X^{\top}\theta_j + \eta_1)\\
&\le 2S\sum_{i =1}^{A_d} c_i{\mathbb{P}(v_i -\eta_1 \le X^{\top}\theta_j < v_i + \eta_1)}\\
&\le 2S\sum_{i =1}^{A_d} c_i(2\overline{C}_0\eta_1)\le8\overline{C}_0S^2\eta_1.
\end{align*}
Therefore by choosing $\eta_2 = \varepsilon/(6\sqrt{\overline{C}_0})$ and $\eta_1 = \varepsilon^2/(32\overline{C}_0S^2)$, we have
\[
\|u_k\circ t_j^{(2)} - l_k\circ t_j^{(1)}\| \le 3\sqrt{\overline{C}_0}\eta_2 + 2\sqrt{2\overline{C}_0}S\sqrt{\eta_1} \le \varepsilon.
\]
Hence the bracketing entropy of $\mathcal{H}^*(S)$ satisfies
\[
\log N_{[\;]}(\varepsilon, \mathcal{H}^*, \|\cdot\|) \le \frac{6AS \sqrt{\overline{C}_0 \diameter(D)}}{\varepsilon} - d\log \frac{96\overline{C}_0S^2}{\varepsilon^2} \lesssim \frac{S}{\varepsilon} ,
\]
for sufficiently small $\varepsilon$. \qedhere

\end{proof}

\begin{lemma}\label{lem:ent_alpa_variation}
Let \begin{align}\label{eq:F_gamma_def_orig}
\mathcal{F}_{\alpha}(K):= \{g: \rchi\to \R|\; &g(x)= f(\theta^\top x), \theta \in \Theta \cup B_{\theta_0}(r)\\ 
& \text{ and } f: D\to \R \text{ is a bounded $\alpha$-variation function with } V_\alpha(f) \le K\}.
\end{align}
 If $\alpha > 1$, Then 
$\log N_{[\,]}(\eta, \f_{\alpha}(K), \|\cdot\|) \le C \eta^{-\alpha}$ for some constant $C$ depending on $K$ only.
\end{lemma}
\begin{proof}
By Lemma 3.6.11~\cite{Gine16}, we have  \begin{align}\label{eq:F_gamma_def_orig1}
\mathcal{F}_{\alpha}(K)= \{x\mapsto f  (h(\theta^\top x))| \theta \in \Theta &\cup B_{\theta_0}(r), h: D \to [0, K] \text{ is a nondecreasing function, and } \\&f \text{ is a $1/\alpha$-H\"{o}lder function defined on } [0, K] \text{with H\"{o}lder constant }1  \}.
\end{align}
Thus by definition~\eqref{LIP:w_def}, we have 
\[\mathcal{F}_{\alpha}(K)= \{x\mapsto f \circ k(x)| k \in \mathcal{H}^*(K) \text{ and }f \text{ is a $1/\alpha$-H\"{o}lder function defined on } [0, K] \}.\]
Let $(k_1^L, k_1^U), \ldots, (k_{N_{\delta_1}}^L ,k_{N_{\delta_1}}^U )$ be an $L_2$-bracket of $\mathcal{H}^*(K)$ of size $\delta_1$, and let $f_1, \ldots, f_{M_{\delta_2}}$ be a $\|\cdot\|_{\infty}$ cover of size $\delta_2$ for the class of bounded $1/\alpha$-H\"{o}lder functions defined on $[0, K]$. By Lemma~\ref{lem:monosim} and Example~5.11 of~\cite{wainwright2019high}, we can choose 
\[ \log N_{\delta_1} \lesssim K \delta_1^{-1} \text{   and   }  \log M_{\delta_2} \lesssim  K \delta_2^{-\alpha}.\] 
For any $f\circ k \in \mathcal{F}_{\alpha}(K)$, assume without loss of generality that $k_1^L(x)\le k(x) \le k_1^U(x)$ and $\|f-f_1\|_{\infty} \le \delta_2.$ Because $f$ is $1/\alpha$-H\"{o}lder, we have that 
\[f\circ k(x)  \le  (k(x)- k_1^L(x))^{1/\alpha}+  f_1 \circ k_1^L(x) + \delta_2 \le  (k_1^U(x)- k_1^L(x))^{1/\alpha}+  f_1 \circ k_1^L(x) + \delta_2   \]
and 
\[f\circ k(x)  \ge - (k(x)- k_1^L(x))^{1/\alpha}+  f_1 \circ k_1^L(x) - \delta_2\ge - (k_1^U(x)- k_1^L(x))^{1/\alpha}+  f_1 \circ k_1^L(x) - \delta_2 .\]
Thus $\{- (k_1^U(x)- k_1^L(x))^{1/\alpha}+  f_1 \circ k_1^L(x) - \delta_2,  (k_1^U(x)- k_1^L(x))^{1/\alpha}+  f_1 \circ k_1^L(x) + \delta_2\}$ forms a bracket for $f\circ k. $ Now the $L_2$ width of the bracket is 
\[2  \|(k_1^U(x)- k_1^L(x))^{1/\alpha}\| + 2 \delta_2 \le 2  (\|(k_1^U(x)- k_1^L(x))\|)^{1/\alpha} + 2 \delta_2\le2  \delta_1^{1/\alpha} + 2 \delta_2.\]
Thus, if  $\delta_1 = \delta_2^{\alpha}$, then we have a $4\delta_2$ bracket of cardinality $\exp(C \delta_2^{-\alpha})$. 
\end{proof}
% section proof_of_eq:app_score_equation (end)
\section{Completing the proof of Theorem~\ref{thm:Main_rate_CLSE} in Section~\ref{app:sketchCLSE}}\label{app:stepsCLSE}
In the following three sections we give a detailed discussion of \ref{item:step3}--\ref{item:step5} in the proof of Theorem~\ref{thm:Main_rate_CLSE}.  Some of the results in this section are proved in Section~\ref{sec:proof_semi}.

\subsection{Proof of \ref{item:step3} in~Section~\ref{app:sketchCLSE}}\label{proof:lip:nobias}
We start with some notation. Recall that for any (fixed or random) $(\theta, m) \in \Theta\times\mathcal{M}_L$, $P_{\theta,m}$ denotes the joint distribution of $Y$ and $X$, where $Y=m(\theta^\top X) +\epsilon$ and $P_X$ denotes the distribution of $X$. Now, let $ P_{\theta,m}^{(Y,X)|\theta^\top X}$ denote the joint distribution of $(Y,X)$ given $\theta^\top X$. For any $(\theta,m) \in \Theta\times \mathcal{M}_L$ and $f \in L_2(P_{\theta,m})$, we have
$P_{\theta,m}[f(X)]= P_X(f(X))$ and
\begin{align}\label{eq:def_EX_X}
\begin{split}
P_{\theta,m} \big[\big(Y-m_0(\theta^\top X)\big) f(X)\big]=& P_X\big[   P_{\theta,m}^{(Y,X)|\theta^\top X} \big[f(X)\big(Y-m_0(\theta^\top X)\big)\big] \big]\\=&P_X\big[ \E( f(X)|\theta^\top X) \big(m(\theta^\top X)-m_0(\theta^\top X)\big) \big].
\end{split}
\end{align}
\begin{thm}[\ref{item:step3}] \label{LIP:thm:nobiasCLSE}
Under assumptions of Theorem~\ref{thm:Main_rate_CLSE}, $\sqrt{n} P_{\check{\theta}, m_0} \psi_{\check{\theta},\check{m}} =o_p(1).$
\end{thm}

\begin{proof}

By the above display, we have that
\begin{align}
\begin{split} \label{eq:nobias_thm}
P_{\check{\theta}, m_0} \psi_{\check{\theta},\check{m}}={}& H_{\check\theta}^\top P_{\check{\theta},m_0}\bigg[(Y-\check{m}(\check{\theta}^\top X)) \big[\check{m}^\prime(\check{\theta}^\top X) X- (\check{m}^\prime \, h_{\theta_0})(\check{\theta} ^\top X)\big] \bigg]\\
={}&H_{\check\theta}^\top P_X\bigg[(m_0-\check{m})(\check{\theta}^\top X)\check{m}^\prime(\check{\theta}^\top X) \big[ E(X|\check{\theta}^\top X)- h_{\theta_0}(\check{\theta}^\top X)\big] \bigg]
\end{split}
\end{align}
Now we will show right \eqref{eq:nobias_thm} is $o_p(n^{-1/2}).$ By \ref{aa1} and the Cauchy-Schwarz inequality, we have
\begin{align}\label{eq:nobias_23}
\begin{split}
&\big|P_X[(m_0-\check{m})(\check{\theta}^\top X) \check{m}^\prime(\check{\theta}^\top X) (E(X|\check{\theta}^\top X)-h_{\theta_0}(\check{\theta}^\top X)) ]\big|\\
\le{}& \|\check{m}^\prime\|_\infty \sqrt{P_X\big[(m_0-\check{m})^2(\check{\theta}^\top X)\big] P_X\big[ |h_{\check{\theta}}(\check{\theta}^\top X)-h_{\theta_0}(\check{\theta}^\top X)|^2\big] }\\
={}& \|\check{m}^\prime\|_\infty \|m_0\circ\check{\theta}-\check{m}\circ\check{\theta}\| \; \| h_{\check{\theta}}\circ\check{\theta}-h_{\theta_0}\circ\check{\theta}\|_{2,P_{\theta_0,m_0}}.
\end{split}
\end{align}
Combining~\eqref{eq:nobias_thm} and \eqref{eq:nobias_23},  we have that
\begin{align}\label{eq:Lip:nobias_1}
|P_{\check{\theta}, m_0} \psi_{\check{\theta},\check{m}}|&\le \|\check{m}^\prime\|_\infty \|m_0\circ\check{\theta}-\check{m}\circ\check{\theta}\| \; \| h_{\check{\theta}}\circ\check{\theta}-h_{\theta_0}\circ\check{\theta}\|_{2,P_{\theta_0, m_0}}.
\end{align}
 % By Theorem \ref{thm:rate_derivCLSE}, we have $\|\check{m}^\prime\circ\check{\theta}-m_0^\prime\circ\check{\theta}\|=O_p(n^{-2/15})$.
  Furthermore, by Theorems \ref{thm:rate_m_theta_CLSE} and  \ref{thm:ratestCLSE}  we have
 \begin{align}
\begin{split}
\|m_0\circ\check{\theta}-\check{m}\circ\check{\theta}\| \le{}&\|\check{m}\circ\check{\theta} - m_0 \circ\theta_0 \| + \|m_0 \circ\theta_0 -m_0\circ\check{\theta} \| \\
\le{}&\|\check{m}\circ\check{\theta} - m_0 \circ\theta_0 \| + L_0 T^2 |\theta_0- \check{\theta}| \\
={}& O_p(n^{-2/5}).
\end{split}
\end{align}
{ To bound the last factor on the right hand side of~\eqref{eq:Lip:nobias_1}, note that
\begin{align*}
\mathrm{TV}(\theta_0^{\top}X, \check{\theta}^{\top}X) &= \sup_{\ell:\|\ell\|_{\infty} \le 1}\left|P_X[\ell(\theta_0^{\top}X) - \ell(\check{\theta}^{\top}X)]\right|\\ 
&\ge \frac{1}{2T^2}\left|P_X\left[|h_{\check{\theta}}(\check{\theta}^\top X)-h_{\theta_0}(\check{\theta}^\top X)|^2 - |h_{\check{\theta}}(\theta_0^\top X)-h_{\theta_0}(\theta_0^\top X)|^2\right]\right|.
\end{align*}
The inequality here follows because $\ell(u) := |h_{\check{\theta}}(u)-h_{\theta_0}(u)|^2$ is upper bounded by $2T^2$ for all $u$. Therefore,
\begin{align*}
P_X|h_{\check{\theta}}(\check{\theta}^\top X)-h_{\theta_0}(\check{\theta}^\top X)|^2 &\le 2T^2\mathrm{TV}(\theta_0^{\top}X, \check{\theta}^{\top}X) + P_X|h_{\check{\theta}}(\theta_0^\top X)-h_{\theta_0}(\theta_0^\top X)|^2\\
&\le 2T^3\overline{C}_0|\check{\theta} - \theta_0| + \bar{M}|\check{\theta} - \theta_0|\\
&= O_p(n^{-1/5}).
\end{align*}
The second inequality here follows from~\eqref{eq:Tv_lip} and assumption~\ref{bb2}. Thus the right hand side of \eqref{eq:Lip:nobias_1} is $O_p(n^{-3/5})$. Thus $|P_{\check{\theta}, m_0} \psi_{\check{\theta},\check{m}}| =o_p(n^{-1/2}).$ \qedhere}
% \todo[inline]{Need B2 here}
\end{proof}

\subsection{Proof of~\ref{item:step4} in~Section~\ref{app:sketchCLSE}} % (fold)
\label{sub:proof_of_item:step4}

% subsection proof_of_item:step4 (end)
In Lemma \ref{thm:nobiasCLSE_part2}, stated and proved in Section~\ref{app:nobias_part2CLSE}, we prove that $\psi_{\check{\theta},\check{m}}$ is a consistent estimator of $\psi_{\theta_0,m_0}$ under $L_2(P_{\theta_0,m_0})$ norm. The following theorem (proved in Section~\ref{LIP:app:ConsisG_n}) completes the proof of Theorem \ref{thm:Main_rate_CLSE}.
\begin{thm}[\ref{item:step4}] \label{LIP:thm:ConsistencyofG_n}
Under assumptions of Theorem~\ref{thm:Main_rate_CLSE}, we have
\begin{equation}\label{LIP:eq:donskerity}
\g_n (\psi_{\check{\theta},\check{m}} -\psi_{\theta_0,m_0})=o_p(1).
 \end{equation}
\end{thm}
We first find an upper bound for the left side of \eqref{LIP:eq:donskerity} and then show that each of the terms converge to zero; see Lemmas~\ref{LIP:lem:tau_donskerity}  and \ref{LIP:lem:upsion_donskerity} in Section~\ref{LIP:app:ConsisG_n}.

\begin{proof}

Recall the definition~\eqref{eq:App_score}. Under model \eqref{eq:simsl},

\begin{align}
\psi_{\check{\theta},\check{m}} -\psi_{\theta_0,m_0}={}& [\epsilon+ m_0(\theta_0^\top x) -\check{m}(\check{\theta}^\top x)] H_{\check\theta}^\top [ \check{m}^\prime (\check{\theta}^\top x) \big(x- h_{\theta_0} (\check{\theta}^\top x)\big)] \\
&\quad- \epsilon H_{\theta_0}^\top m'_0(\theta_0^\top x)\big[ x- h_{\theta_0}(\theta_0^\top x)\big]\\
={}&\epsilon \Big[H_{\check\theta}^\top  \check{m}^\prime (\check{\theta}^\top x) \big[x- h_{\theta_0} (\check{\theta}^\top x)\big] - H_{\theta_0}^\top m'_0(\theta_0^\top x)\big[ x- h_{\theta_0}(\theta_0^\top x)\big]\Big] \\
&\quad +H_{\check\theta}^\top\Big[ \big[m_0(\theta_0^\top x) -\check{m}(\check{\theta}^\top x)\big] \big[ \check{m}^\prime (\check{\theta}^\top x) \big(x- h_{\theta_0} (\check{\theta}^\top x)\big)\big] \Big].\\
&{}=  \epsilon \tau_{\check{\theta},\check{m}} + \upsilon_{\check{\theta},\check{m}}, 
\label{eq:33_11}
\end{align}
where for every $(\theta,m)\in \Theta\times \M_L$, the functions $\upsilon_{\theta, m}:\rchi \rightarrow \R^{d-1}$ and  $ \tau_{\theta, m}:\rchi \rightarrow \R^{d-1}$ are defined as:
\begin{align}\label{eq:G_n_cons_part1}
\begin{split}
 \tau_{\theta, m}(x):={}& H_{\check\theta}^\top  \check{m}^\prime (\check{\theta}^\top x) \big[x- h_{\theta_0} (\check{\theta}^\top x)\big] - H_{\theta_0}^\top m'_0(\theta_0^\top x)\big[ x- h_{\theta_0}(\theta_0^\top x)\big]\\
 \upsilon_{\theta, m}(x):={}& H_\theta^\top  [m_0(\theta_0^\top x) -m(\theta^\top x)]m^\prime (\theta^\top x)  [  x- \;h_{\theta_0} (\theta ^\top x)].
\end{split}
\end{align}
We begin with some definitions.  Let $b_n$ be a sequence of real numbers such that $b_n\rightarrow \infty$ as $n\rightarrow \infty$, $b_n=o(n^{1/2})$,  and $b_n\|\check{m}-m_0\|_{D_{0}} =o_p(1).$ Note that we can always find such a sequence $b_n,$ as by Theorem \ref{thm:uucons} we have $\|\check{m}-m_0\|_{D_{0}} =o_p(1).$ For all $n \in \mathbb{N},$ define\footnote{The notations with $*$ denote the classes that do not depend on $n$ while the ones with $n$ denote shrinking neighborhoods around the truth.}
\begin{align}
% \mathcal{C}^{m*}_{M_1} &:=\Big\{ m\in\M_L:\,\|m\|_\infty \le M_1\Big\}, \\
% \mathcal{C}^m_{M_1}(n) &:=\Big\{ m \in \mathcal{C}^{m*}_{M_1} :  n^{1/5}  \int_{D_0} (m'(t)- m'_0(t))^2 dt \le 1, b_n \|m - m_0\|_{D_0} \le 1\Big\},\\
%
\mathcal{C}^*_{M_1}&:= \Big\{ (\theta,m):  m\in\M_L,\,\|m\|_\infty \le M_1, \text{ and } \theta \in \Theta\cap B_{\theta_0}(r)\Big\},\label{LIP:eq:c_n}\\
%
% \mathcal{C}^\theta(n)&:=\Big\{\theta \in \Theta\cap B_{\theta_0}(r):  n^{1/10} |\theta-\theta_0| \le1 \Big\},  \\
%
% \mathcal{C}_{M_1}(n)&:= \Big\{ (\theta,m):\theta \in \mathcal{C}^\theta(n) \text{ and } m\in \mathcal{C}^m_{M_1}(n)\Big\},
\mathcal{C}_{M_1}(n)&:= \Big\{ (\theta,m) \in \mathcal{C}^*_{M_1}: n^{1/10} |\theta-\theta_0| \le 1,\, n^{1/10}  \|m'\circ\theta- m'_0\circ\theta\|\le 1, \text{ and } b_n \|m - m_0\|_{D_0} \le 1\Big\}\nonumber\end{align}
where $r$ is defined in~\ref{aa6}.
Thus, for every fixed $M_1$, we have
\begin{align}
&\p( |\g_n ( \psi_{\check{\theta},\check{m}} - \psi_{\theta_0,m_0})|> \delta)\\
 \leq{}& \p(|\g_n ( \epsilon \tau_{\check{\theta},\check{m}} + \upsilon_{\check{\theta},\check{m}})|> \delta, (\check{\theta},\check{m}) \in \mathcal{C}_{M_1}(n) )+\p((\check{\theta},\check{m}) \notin \mathcal{C}_{M_1}(n) ) \\
\leq{}& \p\left(|\g_n ( \epsilon \tau_{\check{\theta},\check{m}} )|> \frac{\delta}{2}, (\check{\theta},\check{m}) \in \mathcal{C}_{M_1}(n) \right)\\
& \quad + \p\left(|\g_n \upsilon_{\check{\theta},\check{m}} |> \frac{\delta}{2}, (\check{\theta},\check{m}) \in \mathcal{C}_{M_1}(n) \right)+\p\left((\check{\theta},\check{m}) \notin \mathcal{C}_{M_1}(n) \right) \\
\leq{}& \p\Big(\sup_{(\theta, m) \in \mathcal{C}_{M_1}(n) } |\g_n \epsilon \tau_{\theta,m}|> \frac{\delta}{2}\Big)+ \p\Big(\sup_{(\theta, m) \in \mathcal{C}_{M_1}(n) } |\g_n \upsilon_{\theta,m} |> \frac{\delta}{2}\Big)\\
& \quad +\p\big((\check{\theta},\check{m}) \notin \mathcal{C}_{M_1}(n) \big). \label{38_1}
%& \leq \frac{1}{\delta} E\left(\sup_{f \in \Xi_{M_1}(n) } |\g_n f|\right)+ \frac{5\nu }{8} \\
%& \leq \frac{1}{\delta} J^*_{} \Big(1, \Xi_{M_1}(n) \bigm) [P^* (F_{M_1}^2)]^{1/2}+ \frac{5\nu }{8},
\end{align}
%where $F_{M_1}$ is the envelope of the class of functions $ \mathcal{C}_{M_1}.$
Recall that by Theorems~\ref{thm:rate_m_theta_CLSE}--\ref{thm:rate_derivCLSE}, we have $\p\big((\check{\theta},\check{m}) \notin \mathcal{C}_{M_1}(n) \big) =o(1)$. Thus the proof of Theorem~\ref{LIP:thm:ConsistencyofG_n} will be complete if we show that the first two terms in~\eqref{38_1} are $o(1).$ Lemmas~\ref{LIP:lem:tau_donskerity} and \ref{LIP:lem:upsion_donskerity} do this.

\end{proof}

\section{Proof of results in Section~\ref{app:stepsCLSE}} % (fold)
\label{sec:proof_semi}
\subsection{Lemma used in the proof of Theorem~\ref{LIP:thm:ConsistencyofG_n}}\label{LIP:app:ConsisG_n}

\begin{lemma}\label{LIP:lem:tau_donskerity}
Fix $M_1$ and $ \delta>0.$ Under assumptions \ref{aa1_new}--\ref{aa2},   we have
\[\p\Big(\sup_{(\theta,m ) \in \mathcal{C}_{M_1}(n)  } |\g_n \epsilon \tau_{\theta,m}|> \frac{\delta}{2}\Big) = o(1).\]
\end{lemma}
\begin{proof}
Recall that
\[\tau_{\theta, m}(x):= H_{\check\theta}^\top  \check{m}^\prime (\check{\theta}^\top x) \big[x- h_{\theta_0} (\check{\theta}^\top x)\big] - H_{\theta_0}^\top m'_0(\theta_0^\top x)\big[ x- h_{\theta_0}(\theta_0^\top x)\big].\]
% {\clg \begin{align*}
% \tau_{\theta, m}(x):= H_\theta^\top&\big\{[ m^\prime (\theta^\top x) -m_0^\prime ({\theta_0}^\top x)] x+ [ (m_0^\prime\;h_{\theta_0}) ({\theta_0}^\top x) - (m_0^\prime \;h_{\theta_0}) (\theta ^\top x) ]\big\}\\
% &\quad + (H_{\theta}^\top -H_{\theta_0}^\top) \big[m'_0(\theta_0^\top x) x- (m_0'h_{\theta_0})(\theta_0^\top x)\big],
% \end{align*}}
Let us define,
\[ \Xi_{M_1}(n):=\big\{ \tau_{\theta, m} \big| (\theta, m)\in \mathcal{C}_{M_1}(n) \big\} \quad \text{and}\quad \Xi^*_{M_1}:=\big\{ \tau_{\theta, m} \big| (\theta, m)\in \mathcal{C}^*_{M_1} \big\}.\]
We will prove Lemma~\ref{LIP:lem:tau_donskerity} by applying Lemma~\ref{lem:Maximal342} with $\mathcal{F}=\Xi_{M_1}(n)$ and $\epsilon$. Recall that as $q\ge5$, by~\ref{aa2} we have 
 \[
\mathbb{E}[\epsilon|X] = 0,\quad\mbox{Var}(\epsilon|X) \le \sigma^2,\quad\mbox{and}\quad C_{\epsilon} := 8\mathbb{E}\left[\max_{1\le i\le n}|\epsilon_i|\right] \le n^{1/5}.
\]
We will show that
\be \label{LIP:eq:entropY_xi_star}
 N_{[\,]}(\varepsilon, \Xi_{M_1}(n), \|\cdot\|_{2,P_{\theta_0,m_0}}) \le N(\varepsilon, \Xi^*_{M_1}, \|\cdot\|_{2,\infty}) \leq c \exp (c/\varepsilon) \varepsilon^{-10d},
 \ee 
 and
\begin{equation}\label{eq:bound_last_thm}
 \sup_{f \in \Xi_{M_1}(n)} \|f\|_{2,P_{\theta_0,m_0}} \le  C n^{-1/10} \qquad \text{and} \qquad \sup_{f\in\Xi_{M_1}(n)}\|f\|_{2, \infty}  \le 4LT
\end{equation} where $c$ depends only on $M_1$ and $d$ and $C$ depends only on $L,L_0, T, m_0,$ and $h_{\theta_0}.$ The second inequality in~\eqref{eq:bound_last_thm} follows trivially from the definitions. 

The first inequality of~\eqref{LIP:eq:entropY_xi_star} is trivially true. 
 To prove the second inequality, we will now construct a bracket for $\Xi_{M_1}^*$. Recall that by Lemma~\ref{lem:monosim}, we have
\begin{equation}\label{eq:mono_ent_61}
\log N_{[\,]}(\varepsilon,\{m'(\theta^\top \cdot)| (\theta,m) \in  \mathcal{C}^*_{M_1}\}, L_2(P_{\theta_0, m_0}) ) \lesssim L/\varepsilon.
\end{equation}
Moreover, by Lemma~15 of~\cite{Patra16}, we can find a $\theta_1, \theta_2, \ldots, \theta_{N_{\varepsilon}}$ with $N_{\varepsilon} \lesssim \varepsilon^{-2d}$ such that for every $\theta\in\Theta \cap B_{\theta_0}(1/2)$, there exists a $\theta_j$ such that
 \[|\theta - \theta_j| \le \varepsilon/T, \;\|H_\theta-H_{\theta_j}\|_2\le \varepsilon/T, \text{  and  } |\theta^{\top}x - \theta_j^{\top}x| \le \varepsilon, \;\forall  x\in\chi.\]
Observe that for all $x\in \rchi$, we have $H_{\theta_j}^\top x -\varepsilon \preceq H_\theta^\top x \preceq H_{\theta_j}^\top x +\varepsilon$. Thus
\begin{equation}\label{eq:H_theta_X_brac}
N_{[\,]}(\varepsilon, \{f:\rchi\rightarrow \R^d| f(x)=H_\theta^\top x, \forall x\in \rchi, \theta \in \Theta\cap B_{\theta_0}(1/2) \}, \|\cdot\|_{2,\infty}) \lesssim \varepsilon^{-2d}
\end{equation}
% {\clg Similarly as $|m_0^\prime(\theta^\top x) -m_0^\prime(\theta_j^\top x)| \le L_0 \varepsilon,$ we have
% \begin{equation}\label{eq:mprime0_x}
% N_{[\,]}(\varepsilon, \{m_0'\circ\theta: \theta \in  \Theta\cap B_{\theta_0}(1/2) \}, L_2(P_{\theta_0, m_0})) \lesssim \varepsilon^{-2d}
% \end{equation}}
Finally observe that
\begin{align}\label{eq:H_h_brac}
\begin{split}
&|H_\theta^\top h_{\theta_0}(\theta^\top x)-H_{\theta_j}^\top h_{\theta_0}(\theta_j^\top x)| \\
\le{}& |H_\theta^\top h_{\theta_0}(\theta^\top x)-H_{\theta}^\top h_{\theta_0}(\theta_j^\top x)|+ |H_{\theta}^\top h_{\theta_0}(\theta_j^\top x)-H_{\theta_j}^\top h_{\theta_0}(\theta_j^\top x)|\\
\le{}& | h_{\theta_0}(\theta^\top x)- h_{\theta_0}(\theta_j^\top x)| +\|H_{\theta}^\top -H_{\theta_j}^\top \|_2 \|h_{\theta_0}\|_{2,\infty}\\
\le{}& \| h'_{\theta_0}\|_{2, \infty}|\theta- \theta_j| T +\|H_{\theta}^\top -H_{\theta_j}^\top \|_2 \|h_{\theta_0}\|_{2,\infty} \le\varepsilon (\| h'_{\theta_0}\|_{2, \infty}|+ \|h_{\theta_0}\|_{2,\infty}/T) \lesssim \varepsilon
\end{split}
\end{align}
and
\begin{equation}\label{eq:H_h_0}
|H_\theta^\top h_{\theta_0}(\theta_0^\top x)-H_{\theta_j}^\top h_{\theta_0}(\theta_0^\top x)| \le  \|h_{\theta_0}(\theta_0^\top \cdot)\|_{2, \infty}\varepsilon/ T.
\end{equation}
Thus we have
{\small \begin{align}
 N_{[\,]}(\varepsilon, \{f:\rchi\rightarrow \R^d| f(x)=H_\theta^\top h_{\theta_0}(\theta^\top x),  \theta \in \Theta\cap B_{\theta_0}(1/2) \}, \|\cdot\|_{2,\infty}) &\lesssim \varepsilon^{-2d}, \qquad\label{eq:H_h_ent}
 % \vsap
 % N_{[\,]}(\varepsilon, \{f:\rchi\rightarrow \R^d| f(x)=H_\theta^\top h_{\theta_0}(\theta_0^\top x),  \theta \in \Theta\cap B_{\theta_0}(1/2) \}, \|\cdot\|_{2,\infty}) &\lesssim \varepsilon^{-2d}.\qquad\label{eq:H_h_ent1}
\end{align}}
Thus by applying Lemma 9.25 of \cite{Kosorok08} to sums and product of classes of functions in \eqref{eq:mono_ent_61},\eqref{eq:H_theta_X_brac}, and \eqref{eq:H_h_ent}, we have~\eqref{LIP:eq:entropY_xi_star}. Now, we will find an upper bound for $\sup_{f \in \Xi_{M_1}(n)} \|f\|_{2,P_{\theta_0, m_0}}$.  For every $(\theta, m) \in \mathcal{C}_{M_1}(n)$  and $x \in \rchi$ note that 
\begin{align}\label{eq:L_2_bound}
\begin{split}
 \|\tau_{\theta, m}(X)\|_{2, P_{\theta_0, m_0}} ={}& \Big\|H_{\theta}^\top m^\prime (\theta^\top X) \big[X- h_{\theta_0} (\theta^\top X)\big] - H_{\theta_0}^\top m'_0(\theta_0^\top X)\big[ X- h_{\theta_0}(\theta_0^\top X)\big]\Big\|_{2, P_{\theta_0, m_0}}\\
  \le{}& \Big\|(H_{\theta}^\top - H_{\theta_0}^\top)m^\prime (\theta^\top X) \big[X- h_{\theta_0} (\theta^\top X)\big] \Big\|_{2, P_{\theta_0, m_0}}\\
  &\quad + \Big\|H_{\theta_0}^\top m^\prime (\theta^\top X) \big[X- h_{\theta_0} (\theta^\top X)\big] - H_{\theta_0}^\top m'_0(\theta_0^\top X)\big[ X- h_{\theta_0}(\theta_0^\top X)\big]\Big\|_{2, P_{\theta_0, m_0}}\\
    \le{}& |\theta -\theta_0| 2 L T+ \Big\|H_{\theta_0}^\top \big[m^\prime (\theta^\top X) -m_0^\prime (\theta^\top X)\big] \big[X- h_{\theta_0} (\theta^\top X)\big]\Big\|_{2, P_{\theta_0, m_0}}\\
  &\quad + \Big\|H_{\theta_0}^\top m_0^\prime (\theta^\top X) \big[X- h_{\theta_0} (\theta^\top X)\big] - H_{\theta_0}^\top m'_0(\theta_0^\top X)\big[ X- h_{\theta_0}(\theta_0^\top X)\big]\Big\|_{2, P_{\theta_0, m_0}}\\
   \le{}& |\theta -\theta_0| 2 L T+ 2 T \big\| m^\prime (\theta^\top X ) -m_0^\prime (\theta^\top X)\big] \big\|\\
  &\quad + \Big\|H_{\theta_0}^\top \big[m_0^\prime (\theta^\top X)-m_0^\prime (\theta_0^\top X)\big] \big[X- h_{\theta_0} (\theta^\top X)\big]\Big\|_{2, P_{\theta_0, m_0}}\\
  &\quad + \Big\|H_{\theta_0}^\top m_0^\prime (\theta_0^\top X)\big[X- h_{\theta_0} (\theta^\top X)\big] - H_{\theta_0}^\top m'_0(\theta_0^\top X)\big[ X- h_{\theta_0}(\theta_0^\top X)\big]\Big\|_{2, P_{\theta_0, m_0}}\\
\le{}& |\theta -\theta_0| 2 L T+ 2 T \big\| m^\prime (\theta^\top X ) -m_0^\prime (\theta^\top X)\big] \big\| +  2 T \|m_0''\|_{\infty} |\theta -\theta_0| \\
  &\quad + \Big\|H_{\theta_0}^\top m_0^\prime (\theta_0^\top X)\big[h_{\theta_0}(\theta_0^\top X)- h_{\theta_0} (\theta^\top X)\big]\Big\|_{2, P_{\theta_0, m_0}}\\
  \le{}& { |\theta -\theta_0| 2 L T+ 2 T  \big\| m^\prime (\theta^\top X ) -m_0^\prime (\theta^\top X)\big] \big\| +  2 T \|m_0''\|_{\infty} |\theta -\theta_0|} \\
  &\quad { + L L_{h_0}|\theta_0- \theta|^{1/2}}\\
  \le{}& C_{11} n^{-1/10}
 \end{split}
\end{align}
{ where  the penultimate inequality holds, as $L_{h_0} := \sup_{u_1\neq u_2}| h_{\theta_0}(u_1)- h_{\theta_0}(u_2)|/|u_1-u_2|^{1/2}$ is finite (by~\ref{bb2}) and the last inequality follows from~\eqref{LIP:eq:c_n} and $C_{11}$ is constant depending only on $L,L_0, T, m_0,$ and $h_{\theta_0}$.}
For any $f: \rchi \to \R^{d-1}$, let $f_1,\ldots, f_{d-1}$ denote its real-valued components. For any $k\in \{1, \ldots, d-1\}$, let 
\[ \Xi_{M_1}^{(k)}(n) := \{f_k : f\in \Xi_{M_1}(n)\}.\]  By Markov's inequality, we have 
\begin{align}\label{eq:cvxdim_split}
\begin{split}
&\p\Big(\sup_{f \in \Xi_{M_1}(n) } |\g_n \epsilon f|> \frac{\delta}{2}\Big)
 \le{} 2  \delta^{-1} \sqrt{d-1}\sum_{i=1}^{d-1} \mathbb{E}  \Big(\sup_{g \in \Xi_{M_1}^{(i)}(n) } |\g_n \epsilon g|\Big).
\end{split}
\end{align}
We can bound each term in the summation of the above display by Lemma~\ref{lem:Maximal342},  since by~\eqref{LIP:eq:entropY_xi_star} and~\eqref{eq:bound_last_thm}, we have 
{ \begin{equation}\label{eq:f_i_req1}
 J_{[\,]}(\varepsilon, \Xi^{(i)}_{M_1}(n), \|\cdot\|_{P_{\theta_0,m_0}}) \lesssim \varepsilon^{1/2},\;  \sup_{f \in \Xi^{(i)}_{M_1}(n)} \|f\|_{2,P_{\theta_0,m_0}}   \le C_{11} n^{-1/10}, \quad\text{and}\quad \sup_{f\in\Xi^{(i)}_{M_1}(n)}\|f\|_{2, \infty} \le 4LT.
 \end{equation}
} 
By Lemma~\ref{lem:Maximal342}, we have
\begin{equation}\label{eq:first_step_f_i1}
\mathbb{E}  \Big[\sup_{f \in \Xi^{(i)}_{M_1}(n) } |\g_n \epsilon f|\Big] \lesssim \sigma  \sqrt{C_{11}} n^{-1/20}\left(1 + \sigma \frac{\sqrt{C_{11}} n^{-1/20} 4LT n^{1/5}}{C_{11}^2n^{-1/5} \sqrt{n}}\right)+ \frac{ 8LT n^{1/5}}{\sqrt{n}} =o(1)
\end{equation}
for all $i\in \{1,\ldots, d-1\}$. Thus we have that $\p\Big(\sup_{f \in \Xi_{M_1}(n) } |\g_n \epsilon f|> \frac{\delta}{2}\Big) =o(1).$ \qedhere

\end{proof}

\begin{lemma}\label{LIP:lem:upsion_donskerity}
Fix $M_1$ and $ \delta>0.$ For $n\in \mathbb{N},$ we have
\[\p\left(\sup_{(\theta,m ) \in \mathcal{C}_{M_1}(n) } |\g_n \upsilon_{\theta, m}|> \frac{\delta}{2}\right) =o_p(1). \]
\end{lemma}
\begin{proof}
Recall that
\[\upsilon_{\theta, m}(x):={}    H_\theta^\top  [m_0(\theta_0^\top x) -m(\theta^\top x)]m^\prime (\theta^\top x)  [  x- \;h_{\theta_0} (\theta ^\top x)].\]
We will first show that
\begin{equation}\label{eq:J_up_class}
J_{[\,]}(\nu, \{  \upsilon_{\theta, m}:(\theta,m ) \in \mathcal{C}_{M_1}(n) \}, \|\cdot\|_{2,P_{\theta_0,m_0}})
\lesssim \nu^{1/2}
\end{equation}
By Lemmas~\ref{ent10} and \ref{lem:monosim} and \eqref{eq:H_theta_X_brac} and \eqref{eq:H_h_ent}, we have
\begin{align}
N_{[\,]}(\varepsilon,\{m_0(\theta_0^\top \cdot) -m(\theta^\top \cdot)| (\theta,m) \in  \mathcal{C}^*_{M_1}\}, \|\cdot\|_\infty ) &\lesssim \exp(1/\sqrt{\varepsilon}),\nonumber\\
N_{[\,]}(\varepsilon,\{m'(\theta^\top \cdot)| (\theta,m) \in  \mathcal{C}^*_{M_1}\}, \|\cdot\| ) &\lesssim \exp(1/\varepsilon),\label{eq:62_ents}\\
N_{[\,]}(\varepsilon, \{f:\rchi\rightarrow \R^d| f(x)=H_\theta^\top x, \forall x\in \rchi, \theta \in \Theta\cap B_{\theta_0}(1/2) \}, \|\cdot\|_{2,\infty}) &\lesssim \varepsilon^{-2d}\nonumber\\
N_{[\,]}(\varepsilon, \{f:\rchi\rightarrow \R^d| f(x)=H_\theta^\top h_{\theta_0}(\theta^\top x),  \theta \in \Theta\cap B_{\theta_0}(1/2) \}, \|\cdot\|_{2,\infty}) &\lesssim \varepsilon^{-2d}.\nonumber
\end{align}
Thus by applying Lemma 9.25 of \cite{Kosorok08} to sums and product of classes of functions in \eqref{eq:62_ents},  we have
\begin{equation}\label{eq:uppsilon_Lip_ent}
N_{[\,]}(\varepsilon, \{  \upsilon_{\theta, m}:(\theta,m ) \in \mathcal{C}^*_{M_1} \}, \|\cdot\|_{2,P_{\theta_0,m_0}}) \lesssim \exp{\left( \frac{1}{\varepsilon}+  \frac{1}{\sqrt\varepsilon}\right)} \varepsilon^{-6d}.
\end{equation}
Now~\eqref{eq:J_up_class} follows  from the definition of $J_{[\,]}$ by observing that \[
J_{[\,]}(\nu, \{  \upsilon_{\theta, m}:(\theta,m ) \in \mathcal{C}_{M_1}(n) \}, \|\cdot\|_{2,P_{\theta_0,m_0}}) \le J_{[\,]}(\nu, \{  \upsilon_{\theta, m}:(\theta,m ) \in \mathcal{C}^*_{M_1} \}, \|\cdot\|_{2,P_{\theta_0,m_0}}).\]
Next find $\sup_{(\theta,m) \in \mathcal{C}_{M_1}(n) }\|\upsilon_{\theta,m}\|_{2,\infty}.$ For every $x \in \rchi$ observe that,
\begin{align*}
|\upsilon_{\theta, m}(x)|
\leq{}& \big[|m_0(\theta_0^\top x) -m({\theta_0}^\top x)| + |m(\theta_0^\top x) -m(\theta^\top x)|\big] |m^\prime (\theta^\top x)|    |x- \;h_{\theta_0} (\theta ^\top x)| \\
\leq{}& \big[\|m_0 -m\|_{D_0} + L|\theta_0^\top x -\theta^\top x|\big] |m^\prime (\theta^\top x)|    |x- \;h_{\theta_0} (\theta ^\top x)| \\
 % &\;+\cdot \cdot |m^\prime (\theta^\top x)|    |x- \;h_{\theta_0} (\theta ^\top x)|\\
% \leq{}& \|m_0 -m\|_{D_0} |m^\prime (\theta^\top x) x- m_0^\prime (\theta^\top x) h_{\theta_0} (\theta^\top x)| \\
%  &\;+ L|\theta_0^\top x -\theta^\top x| |m^\prime (\theta^\top x) x- m_0^\prime (\theta^\top x) h_{\theta_0} (\theta^\top x)| \\
 \leq{}&\big[b_n^{-1} + 2LT |\theta-\theta_0|] 2L T\\
 \le{}& C [b_n^{-1} + n^{-1/10}],
\end{align*}
where $C$ is a constant depending only on $T,L,$ and $ M_1$. Thus \[
\sup_{(\theta,m) \in \mathcal{C}_{M_1}(n) }\|\upsilon_{\theta,m}\|_{2,P_{\theta_0,m_0}}\le \sup_{(\theta,m) \in \mathcal{C}_{M_1}(n) }\|\upsilon_{\theta,m}\|_{2,\infty} \le C [b_n^{-1} + n^{-1/10}].\]
Thus using arguments similar to~\eqref{eq:cvxdim_split} and the  maximal inequality in Lemma 3.4.2 of \cite{VdVW96} (for uniformly bounded function classes), we have

\begin{align*}
&\p\left(\sup_{(\theta,m ) \in \mathcal{C}_{M_1}(n) } |\g_n \upsilon_{\theta, m}|> \frac{\delta}{2}\right) \\
 \lesssim{}&  2 \delta^{-1} \sqrt{d-1}\sum_{i=1}^{d-1} \mathbb{E} \Big(\sup_{(\theta,m ) \in \mathcal{C}_{M_1}(n) } |\g_n \upsilon_{\theta, m,i}|\Big)\\
\lesssim{}& J_{[\,]}([b_n^{-1} + n^{-1/10}],\w_{M_1}(n), \|\cdot\|_{2,P_{\theta_0,m_0}})+ \frac{J^2_{[\;]}([b_n^{-1} + n^{-1/10}],\w_{M_1}(n), \|\cdot\|_{2,P_{\theta_0,m_0}})}{[b_n^{-1} + n^{-1/10}]^2\sqrt{n}} \\
\lesssim{}&  [b_n^{-1} + n^{-1/10}]^{1/2} +\frac{ [b_n^{-1} + n^{-1/10}]}{[b_n^{-1} + n^{-1/10}]^2\sqrt{n}} \\
\lesssim{}& [b_n^{-1} + n^{-1/10}]^{1/2} +\frac{ 1}{b_n^{-1} \sqrt{n} + n^{4/10}} =o(1),
\end{align*}
 as $b_n=o(n^{1/2})$, here in the first inequality $\upsilon_{\theta, m,i}$ denotes the $i$th component of $\upsilon_{\theta, m}.$
\end{proof}

\subsection{Lemma used in the proof of~\ref{item:step5}}\label{app:nobias_part2CLSE}

The following lemma is used in the proof of~\ref{item:step5} in Theorem~\ref{thm:Main_rate_CLSE}; also see~\citet[Section 10.4]{Patra16}.

\begin{lemma} \label{thm:nobiasCLSE_part2}
If the conditions in Theorem~\ref{thm:Main_rate_CLSE} hold, then \begin{align}
P_{\theta_0, m_0} |\psi_{\check{\theta},\check{m}}- \psi_{\theta_0,m_0}|^2&=o_p(1),\label{eq:LIPL_2conv}\\
P_{\check{\theta}, m_0} |\psi_{\check{\theta},\check{m}}|^2&=O_p(1). \label{eq:LIPL_2bound}
\end{align}
\end{lemma}
\begin{proof}
We first prove \eqref{eq:LIPL_2conv}. 
By the smoothness properties of $\theta \mapsto H_\theta$; see Lemma 1 of~\cite{Patra16}, we have
% \todo[inline]{I do not see where (B2) is used. Only use $h_{\theta_0}$'s Lipschitness property.}
\begin{align}
&\; P_{\theta_0, m_0} | \psi_{\check{\theta},\check{m}}- \psi_{\theta_0,m_0}|^2\\
\qquad&= P_{\theta_0, m_0} \Big|( y-\check{m}(\check{\theta} ^\top X)) H_{\check{\theta}}^\top \big[ \check{m}^\prime(\check{\theta} ^\top X) \big(X -h_{\theta_0} (\check{\theta}^\top X)\big)\big]\\
&\quad\; - ( y-m_0(\theta_0^\top X)) H_{\theta_0}^\top\big[ m_0^\prime(\theta_0^\top X) \big(X - h_{\theta_0} (\theta_0^\top X)\big)\big]\Big|^2\\
\qquad&= P_X \Big| \big[( m_0(\theta_0 ^\top X)-\check{m}(\check{\theta}^\top X))+\epsilon\big] H_{\check{\theta}}^\top \big[\check{m}^\prime(\check{\theta} ^\top X) \big(X -  h_{\theta_0}(\check{\theta} ^\top X)\big)\big] \\
&\quad\; -\epsilon H_{\theta_0}^\top\big[m_0^\prime(\theta_0^\top X) \big(X - h_{\theta_0} (\theta_0^\top X)\big)\big] \Big|^2\\
\qquad&= P_X \Big| \big[ m_0(\theta_0 ^\top X)-\check{m}(\check{\theta}^\top X)\big] H_{\check{\theta}}^\top \big[\check{m}^\prime(\check{\theta} ^\top X) \big(X -  h_{\theta_0}(\check{\theta} ^\top X)\big)\big]\Big|^2 \\
&\quad\;  + P_{\theta_0, m_0} \bigg| \epsilon \Big[ H_{\check{\theta}}^\top \big[\check{m}^\prime(\check{\theta} ^\top X) \big(X -  h_{\theta_0}(\check{\theta} ^\top X)\big)\big] -H_{\theta_0}^\top \big[m_0^\prime(\theta_0^\top X) \big(X - h_{\theta_0} (\theta_0^\top X)\big)\big]\Big] \bigg|^2\\
\qquad&\le P_X \Big| \big[ m_0(\theta_0 ^\top X)-\check{m}(\check{\theta}^\top X)\big] \big[\check{m}^\prime(\check{\theta} ^\top X) \big(X -  h_{\theta_0}(\check{\theta} ^\top X)\big)\big]\Big|^2 \\
&\quad\; + P_{\theta_0, m_0} \bigg| \epsilon  H_{\check{\theta}}^\top  \Big[\check{m}^\prime(\check{\theta} ^\top X) \big(X -  h_{\theta_0}(\check{\theta} ^\top X)\big) -m_0^\prime(\theta_0^\top X) \big(X - h_{\theta_0} (\theta_0^\top X)\big)\Big] \bigg|^2\\
&\quad\;+ P_{\theta_0, m_0} \bigg| \epsilon \Big[ H_{\check{\theta}}^\top  -H_{\theta_0}^\top\Big] \big[m_0^\prime(\theta_0^\top X) \big(X - h_{\theta_0} (\theta_0^\top X)\big)\big] \bigg|^2\\
% \qquad&\le P_{\theta_0, m_0} \Big| \big[ m_0(\theta_0 ^\top X)-\check{m}(\check{\theta}^\top X)\big] \big[\check{m}^\prime(\check{\theta} ^\top X) \big(X -  h_{\theta_0}(\check{\theta} ^\top X)\big)\big]\Big|^2 \\
% &\quad\; +  \|\sigma^2(\cdot)\|_\infty P_X \Big|  \check{m}^\prime(\check{\theta} ^\top X) X -m_0^\prime(\theta_0^\top X) X + (m_0^\prime \; h_{\theta_0}) (\theta_0^\top X) - (\check{m}^\prime \; h_{\theta_0})(\check{\theta} ^\top X) \Big|^2\\
% &\quad\;+ 4 M_1^2 T^2 \|\sigma^2(\cdot)\|_\infty \|H_{\check\theta}-H_{\theta_0}\|_2^2 \\
\qquad&\le P_X \Big| \big[ m_0(\theta_0 ^\top X)-\check{m}(\check{\theta}^\top X)\big] \big[\check{m}^\prime(\check{\theta} ^\top X) \big(X -  h_{\theta_0}(\check{\theta} ^\top X)\big)\big]\Big|^2 \\
&\quad\; +  \|\sigma^2(\cdot)\|_\infty P_X \Big|  \check{m}^\prime(\check{\theta} ^\top X) \big(X -  h_{\theta_0}(\check{\theta} ^\top X)\big) -m_0^\prime(\theta_0^\top X) \big(X - h_{\theta_0} (\theta_0^\top X)\big) \Big|^2\\
&\quad\;+ 4 M_1^2 T^2 \|\sigma^2(\cdot)\|_\infty \|H_{\check\theta}-H_{\theta_0}\|_2^2 \\
% 
% \qquad &\leq{}   P_X \Big| \big[ m_0(\theta_0 ^\top X)-\check{m}(\check{\theta}^\top X)\big] \big[\check{m}^\prime(\check{\theta} ^\top X) \big(X -  h_{\theta_0}(\check{\theta} ^\top X)\big)\big]\Big|^2 \\
% &\quad\; + 2\sigma^2   P_X \Big|  \check{m}^\prime(\check{\theta} ^\top X) X -m_0^\prime(\theta_0^\top X) X  \Big|^2\\
\qquad&\le P_X \Big| \big[ m_0(\theta_0 ^\top X)-\check{m}(\check{\theta}^\top X)\big] \big[\check{m}^\prime(\check{\theta} ^\top X) \big(X -  h_{\theta_0}(\check{\theta} ^\top X)\big)\big]\Big|^2 \\
&\quad\; +  \|\sigma^2(\cdot)\|_\infty P_X \Big|  \check{m}^\prime(\check{\theta} ^\top X) \big(X -  h_{\theta_0}(\check{\theta} ^\top X)\big) -m_0^\prime(\theta_0^\top X) \big(X - h_{\theta_0} (\theta_0^\top X)\big) \Big|^2\\
&\quad\; + 4 M_1^2 T^2 |\check\theta-\theta_0|^2\sigma^2\\
\qquad&=\textbf{\Rome{1}} + \sigma^2 \textbf{ \Rome{2}}+  4 M_1^2 T^2  \sigma^2 |\check\theta-\theta_0|^2, \label{eq:consis_11}
\end{align}
where
\begin{align}\label{eq:roja_jane}
\begin{split}
\textbf{\Rome{1}} &:= P_X \Big| \big[ m_0(\theta_0 ^\top X)-\check{m}(\check{\theta}^\top X)\big] \big[\check{m}^\prime(\check{\theta} ^\top X) \big(X -  h_{\theta_0}(\check{\theta} ^\top X)\big)\big]\Big|^2, \\
  \textbf{ \Rome{2}} &:=   P_X \Big|  \check{m}^\prime(\check{\theta} ^\top X) \big(X -  h_{\theta_0}(\check{\theta} ^\top X)\big) -m_0^\prime(\theta_0^\top X) \big(X - h_{\theta_0} (\theta_0^\top X)\big) \Big|^2.
\end{split}
\end{align}
We will now show that both $\textbf{\Rome{1}} $ and $\textbf{ \Rome{2}},$ are $o_p(1).$  By Theorems~\ref{thm:ratestCLSE} and \ref{thm:rate_derivCLSE}, we have
{\small \begin{align*}
 \textbf{ \Rome{2}} \le{}& P_X \Big|  \check{m}^\prime(\check{\theta} ^\top X) \big(X -  h_{\theta_0}(\check{\theta} ^\top X)\big) -m_0^\prime(\theta_0^\top X) \big(X - h_{\theta_0} (\theta_0^\top X)\big) \Big|^2\\
 \le{}& P_X \Big|  \check{m}^\prime(\check{\theta} ^\top X) \big(X -  h_{\theta_0}(\check{\theta} ^\top X)\big)- \check{m}^\prime(\check{\theta} ^\top X) \big(X -  h_{\theta_0}(\theta_0 ^\top X)\big)+ \big(\check{m}^\prime(\check{\theta} ^\top X)  -m_0^\prime(\theta_0^\top X)\big) \big(X - h_{\theta_0} (\theta_0^\top X)\big) \Big|^2\\
  \le{}&  2 P_X \Big|  \check{m}^\prime(\check{\theta} ^\top X) \big(h_{\theta_0}(\theta_0 ^\top X) -  h_{\theta_0}(\check{\theta} ^\top X)\big)\Big|^2+ 2P_X \Big|\big(\check{m}^\prime(\check{\theta} ^\top X)  -m_0^\prime(\theta_0^\top X)\big) \big(X - h_{\theta_0} (\theta_0^\top X)\big) \Big|^2\\
   \le{}&  2 L^2 P_X \Big| h_{\theta_0}(\theta_0 ^\top X) -  h_{\theta_0}(\check{\theta} ^\top X)\Big|^2+ 4 T^2P_X \Big|\check{m}^\prime(\check{\theta} ^\top X)  -m_0^\prime(\theta_0^\top X)\Big|^2\\
\le{}&  2 L^2 T^2 L_{h_0} |\theta_0 -\check{\theta}|+ 4 T^2P_X \Big|\check{m}^\prime(\check{\theta} ^\top X)  -m_0^\prime(\theta_0^\top X)\Big|^2\\
\le{}&  2 L^2 T^2 L_{h_0} |\theta_0 -\check{\theta}|+ 8 T^2 \|\check{m}^\prime(\check{\theta} ^\top X)  -m_0^\prime(\check{\theta}^\top X)\|^2 +  8 T^2 \|m_0^\prime(\check{\theta}^\top X)  -m_0^\prime(\theta_0^\top X)\|^2\\
\le{}&  { 2 L^2 T^2 L_{h_0} |\theta_0 -\check{\theta}|+ 8 T^2 \|\check{m}^\prime(\check{\theta} ^\top X)  -m_0^\prime(\check{\theta}^\top X)\|^2 +  8 T^2 \|m_0''\|_{\infty} T^2 |\theta_0-\check{\theta}|^2
= o_p(1),}
\end{align*}}
{ as $L_{h_0} := \sup_{u_1\neq u_2}| h_{\theta_0}(u_1)- h_{\theta_0}(u_2)|/|u_1-u_2|^{1/2}$ is finite by~\ref{bb2}.} For $\textbf{\Rome{1}}$, observe that
\begin{equation}\label{eq:bound_1_second}
|\check{m}'(\check{\theta}^{\top}x)\big(x - h_{\theta_0}(\check{\theta}^{\top}x)\big)| \le |\check{m}'(\check{\theta}^{\top}x)x| + |m_0'(\check{\theta}^{\top}x)h_{\theta_0}(\check{\theta}^{\top}x)| \le 2LT
\end{equation}
Moreover, by Theorem \ref{thm:rate_m_theta_CLSE}, we have $\|\check{m}\circ\check{\theta} - m_0\circ\theta_0\| \stackrel {P}{\rightarrow} 0$. Thus,
\begin{align*}
\textbf{\Rome{1}}={}& P_X \big|( m_0(\theta_0^\top X)-\check{m}(\check{\theta}^\top X)) (\check{m}^\prime(\check{\theta} ^\top X) \big(X -  h_{\theta_0}(\check{\theta} ^\top X)\big))\big|^2\\
\leq{}&2 LT\|m_0\circ\theta_0 - \check{m}\circ\check{\theta}\|^2=o_p(1).
\end{align*}
% Finally, we have
% { \begin{align*}
% \textbf{\Rome{3}}={}& P_X \Big| (m_0^\prime \; h_{\theta_0}) (\theta_0^\top X) - (\check{m}^\prime \; h_{\theta_0})(\check{\theta} ^\top X) \Big|^2\\
% % \le{}&P_{\theta_0, m_0} \Big[ \|m_0^{\prime\prime} \; h_{\theta_0} + m_0^\prime \; h^\prime_{\theta_0}\|_{2,\infty} |(\theta_0 -\check{\theta}) ^\top X|\Big]^2\\
% % \le{}& \|m_0^{\prime\prime} \; h_{\theta_0} + m_0^\prime \; h^\prime_{\theta_0}\|^2_{2,\infty}  T^2 |\theta_0 -\check{\theta}| ^2=o_p(1). 
% \end{align*}}
% \todo[inline]{Need to show that this is $o_p(1)$.}
Thus proof of~\eqref{eq:LIPL_2conv} is complete. We now prove \eqref{eq:LIPL_2bound}. Note that
\begin{align}\label{eq:psi_bound}
\begin{split}
P_{\check{\theta}, m_0} |\psi_{\check{\theta},\check{m}}|^2
\le{}& P_{\check{\theta}, m_0} \Big| (Y-\check{m}(\check{\theta} ^\top X))^2 \big[ \check{m}^\prime(\check{\theta}^\top X) \big(X - h_{\theta_0}(\check{\theta} ^\top X)\big)\big] \Big|^2\\
={}& P_{\check{\theta}, m_0} \Big| \big[(m_0(\check{\theta} ^\top X)-\check{m}(\check{\theta} ^\top X))+\epsilon\big] \;\big[ \check{m}^\prime(\check{\theta}^\top X) \big(X -h_{\theta_0}(\check{\theta} ^\top X)\big)\big]\Big|^2\\
\le{}& P_{\check{\theta}, m_0} \Big| \big[(m_0(\check{\theta} ^\top X)-\check{m}(\check{\theta} ^\top X))\big] \;\big[ \check{m}^\prime(\check{\theta}^\top X) \big(X -h_{\theta_0}(\check{\theta} ^\top X)\big)\big]\Big|^2\\
& + \sigma^2 P_{\check{\theta}, m_0} \big| \check{m}^\prime(\check{\theta}^\top X) \big(X -h_{\theta_0}(\check{\theta} ^\top X)\big)\big|^2\\
\le{}&(\|m_0\|_\infty ^2 + \|\check{m}\|_\infty^2) P_{\check{\theta}, m_0} \big|  \check{m}^\prime(\check{\theta}^\top X) \big(X -h_{\theta_0}(\check{\theta} ^\top X)\big)\big|^2\\
& + P_{\check{\theta}, m_0} | \check{m}^\prime(\check{\theta}^\top X) \big(X -h_{\theta_0}(\check{\theta} ^\top X)\big)|^2\\
\le{}&(\|m_0\|_\infty ^2 + \|\check{m}\|_\infty^2 + 1) P_{\check{\theta}, m_0} |  \check{m}^\prime(\check{\theta}^\top X) \big(X -h_{\theta_0}(\check{\theta} ^\top X)\big)|^2. \qedhere
\end{split}
\end{align}
% All these facts combined prove $P_{\theta_0, m_0} | \psi_{\check{\theta},\check{m}}- \psi_{\theta_0,m_0}|^2=o_p(1).$
\end{proof}

\section{Remark on pre-binning}\label{rem:bin}
% { Need to add a line about ties..; see first paragraph of Section 5.1.}
{ The matrices involved in the optimization problem~\eqref{lipcons} and~\eqref{p1} in Section~\ref{sec:compute} have entries depending on fractions $1/(t_{i+1} - t_i)$. Thus if there are ties in $\{t_i\}_{1\le i \le n}$, then the matrix $A$ is incomputable. Moreover, if $t_{i+1} - t_i$ is very small, then the fractions can force the matrices involved to be ill-conditioned (for the purposes of numerical calculations). Thus  to avoid ill-conditioning of these matrices, in practice one might have to pre-bin the data which leads to a diagonal matrix $Q$ with different diagonal entries. One common method of pre-binning the data is to take the means of all data points for which the $t_i$'s are close. To be more precise, if tolerance  $\eta = 10^{-6}$ and $0 < t_2 - t_1 < t_3 - t_1 < \eta$, then we will combine the data points $(t_1, y_1), (t_2, y_2), (t_3, y_3)$ by taking their mean and set $Q_{1,1} = 3$. Note that the total number of data points is now reduced to $n-2$. The above pre-binning step is implemented in the accompanying package. 
}

\section{Discussion on the theoretical analysis of the {CvxLSE}} % (fold)
\label{sec:discussion_on_the_theoretical_analysis_of_the_texttt_CvxLSE}

   The \texttt{CvxLSE} defined in~\eqref{eq:ultimate_est} is a natural estimator for the convex single index model~\eqref{eq:simsl}.  We have investigated its performance in our simulation studies in Section~\ref{sec:Simul_Cvx} and \ref{app:add_simul}. However, a thorough study of the theoretical properties of the \texttt{CvxLSE} is an open research problem. The difficulties are multifaceted. A result like Theorem~\ref{thm:rate_m_theta_CLSE} (which is used throughout the paper) for the \texttt{CvxLSE} is not known. 
   The recent advancements of~\cite{2018arXiv180502542H} in the analysis of the \texttt{CvxLSE} in the one-dimensional regression  problem is encouraging. However, these techniques cannot be directly extended to our framework as the index parameter is unknown. Even if we have a result like Theorem~\ref{thm:rate_m_theta_CLSE}, deriving Theorem~\ref{thm:uucons} for the \texttt{CvxLSE} brings further challenges. In particular the standard technique (see discussion in page~\pageref{thm:uucons}) used to prove consistency of $\{{m}^\dagger_n\}_{n\ge 1}$ would require control  on $m^\dagger_n$ and its right-derivative near the boundary of its domain. Another bottleneck is deriving a result similar to Theorem~\ref{thm:rate_derivCLSE} for the \texttt{CvxLSE}.  Even in the case of 1-dimensional convex LSE, there are no results  that study the $L_2$-loss for the derivative of the LSE. Note that the derivative is an important quantity in the case of the single index model as the efficient score has $m_0'$ in its formulation; see~\cite{groeneboom2016current,2016arXiv161006026B,2017arXiv171205593B} for similar difficulties that arise in related models. However, if one can prove results similar to Theorems~\ref{thm:rate_m_theta_CLSE}--\ref{thm:rate_derivCLSE} for the convex LSE, then the techniques used in Section~\ref{sec:SemiInf} can be readily applied to prove asymptotic normality of $\theta^\dagger$. These challenges make the study of the \texttt{CvxLSE} a very interesting problem for future research.

% section discussion_on_the_theoretical_analysis_of_the_texttt_ (end)

% section discussion_on_the (end)
% \noeqref{eq:rate_final_deriv}
% \bibliographystyle{chicago}
% \bibliography{SigNoise}
% \end{document}

\end{document}